\documentclass[a4paper]{amsart}

\usepackage{fix-cm}
\usepackage[final]{microtype}
\usepackage[dvipsnames,svgnames,x11names,hyperref]{xcolor}
\usepackage{dsfont,url,graphicx,verbatim,amssymb,enumerate,stmaryrd,booktabs,lmodern,mathtools,mathabx,nicefrac}
\SetSymbolFont{stmry}{bold}{U}{stmry}{m}{n}
\usepackage[pagebackref,colorlinks,citecolor=Mahogany,linkcolor=Mahogany,urlcolor=Mahogany,filecolor=Mahogany]{hyperref}
\usepackage[capitalize]{cleveref}
\usepackage[mathscr]{euscript}
\usepackage[margin=1.33in]{geometry}
\usepackage{tikz,tikz-cd}
\usetikzlibrary{matrix,calc,positioning,arrows,decorations.pathreplacing,patterns,arrows,patterns.meta}
\tikzset{
  snake left/.style={
    rounded corners,
    to path={
      let \p1 = (\tikztostart.east),
          \p2 = (\tikztotarget.west),
          \p3 = ($(\p1)!0.5!(\p2)$),
          \n1 = {8pt} 
      in
      (\p1)
      -- (\x1 + \n1, \y1)
      -- (\x1 + \n1, \y3)
      -- (\x2 - \n1, \y3) \tikztonodes
      -- (\x2 - \n1, \y2)
      -- (\p2)
    }
  }
}

\newtheorem{theorem}{Theorem}[section]
\newtheorem*{theorem*}{Theorem}
\newtheorem{lemma}[theorem]{Lemma}
\newtheorem{proposition}[theorem]{Proposition}
\newtheorem{corollary}[theorem]{Corollary}
\newtheorem*{corollary*}{Corollary}

\newtheorem{atheorem}{Theorem}

\newtheorem{innercustomgeneric}{\customgenericname}
\providecommand{\customgenericname}{}
\newcommand{\newcustomtheorem}[2]{%
  \newenvironment{#1}[1]
  {%
   \ifdefined\crefalias\crefalias{innercustomgeneric}{#2}\fi
   \renewcommand\customgenericname{#2}%
   \renewcommand\theinnercustomgeneric{##1}%
   \innercustomgeneric
  }
  {\endinnercustomgeneric}%
  \ifdefined\crefname\crefname{#2}{#2}{#2s}\fi
}

\newcustomtheorem{customthm}{Theorem}

\newcustomtheorem{customconj}{Conjecture}
\theoremstyle{definition}
\newtheorem{definition}[theorem]{Definition}
\newtheorem*{definition*}{Definition}

\newtheorem{notation}[theorem]{Notation}
\newtheorem{convention}[theorem]{Convention}

\theoremstyle{remark}
\newtheorem{example}[theorem]{Example}
\newtheorem*{example*}{Example}
\newtheorem*{remark*}{Remark}

\newtheorem{remark}[theorem]{Remark}

\renewcommand{\sf}[1]{{\mathsf{#1}}}
\newcommand{\scr}[1]{{\mathscr{#1}}}
\renewcommand{\rm}[1]{{\mathrm{#1}}}
\newcommand{\lra}{\longrightarrow}

\newcommand{\ul}[1]{{\underline{#1}}}
\newcommand{\ol}[1]{{\overline{#1}}}
\newcommand{\op}{\rm{op}}
\DeclareMathOperator*{\colim}{colim}
\newcommand{\ds}[1]{{\mathds{#1}}}
\newcommand{\fr}[1]{{\mathfrak{#1}}}
\newcommand{\bb}[1]{{\mathds{#1}}}
\newcommand{\cal}[1]{{\mathcal{#1}}}
\renewcommand{\bf}[1]{{\mathbf{#1}}}

\makeatletter
\newcommand{\superimpose}[2]{{%
  \ooalign{%
    \hfil$\m@th#1\@firstoftwo#2$\hfil\cr
    \hfil$\m@th#1\@secondoftwo#2$\hfil\cr
  }%
}}
\makeatother

\newcommand{\Fun}{\rm{Fun}}
\newcommand{\Map}{\rm{Map}}

\newcommand{\Sp}{\scr{S}\rm{p}}
\newcommand{\aug}{\rm{aug}}

\newcommand{\St}{\rm{St}}
\newcommand{\SSt}{{\scr{S}\rm{t}}}
\newcommand{\SStH}{{\SSt^2}}
\newcommand{\SStL}{{\SSt^\infty}}
\newcommand{\ISSt}{{\scr{S}\rm{t}_{>0}}}
\newcommand{\ISStH}{{\SSt^2_{>0}}}

\newcommand{\PGL}{\rm{PGL}}
\newcommand{\GL}{\rm{GL}}
\newcommand{\alt}{\rm{alt}}
\newcommand{\Cat}{{\scr{C}\rm{at}}}
\newcommand{\Spc}{{\scr{S}\rm{pc}}}

\newcommand{\alg}{\rm{alg}}
\newcommand{\Alg}{\rm{Alg}}
\newcommand{\Op}{{\scr{O}\rm{p}}}
\newcommand{\Mon}{\rm{Mon}}
\newcommand{\coAlg}{\rm{coAlg}}
\newcommand{\Q}{\ds{Q}}
\newcommand{\N}{\ds{N}}
\newcommand{\Z}{\ds{Z}}
\newcommand{\C}{\ds{C}}

\renewcommand{\cot}{\rm{cot}}
\newcommand{\triv}{\rm{triv}}
\newcommand{\prim}{\rm{prim}}
\newcommand{\fgt}{\rm{fgt}}
\newcommand{\free}{\rm{free}}
\newcommand{\cofree}{\rm{cofree}}
\newcommand{\fil}{\rm{fil}}
\newcommand{\gc}{\rm{gc}}

\newcommand{\indec}{\rm{indec}}

\newcommand{\gr}{\rm{gr}}
\newcommand{\Fin}{\rm{Fin}}
\newcommand{\coLie}{\rm{\rm{coLie}}}

\newcommand{\BGL}{\rm{BGL}}

\newcommand{\BGLb}{\bf{BGL}}
\newcommand{\topro}{\mathrel{\mathpalette\superimpose{{\to}{\shortmid}}}}
\newcommand{\lrapro}{\mathrel{\mathpalette\superimpose{{\lra}{\shortmid}}}}
\newcommand{\levi}{\boxplus}
\newcommand{\para}{\boxbackslash}
\renewcommand{\Bar}{\rm{Bar}}
\newcommand{\Vect}{\rm{Vect}}
\newcommand{\Cobar}{\rm{Cobar}}

\newcommand{\id}{\rm{id}}
\newcommand{\inc}{\rm{inc}}
\newcommand{\pr}{\rm{pr}}

\newcommand{\StH}{\rm{St}^{2}}
\newcommand{\StL}{\rm{St}^{\infty}}
\newcommand{\FC}{\mathrm{FC}} 
\newcommand{\FI}{\mathrm{FI}} 
\newcommand{\FCR}{\mathrm{FCR}} 

\newcommand{\compactldots}{\mathinner{\ldotp\mkern-2mu\ldotp\mkern-2mu\ldotp}}
\newcommand{\compactcdots}{\mathinner{\cdotp\mkern-2mu\cdotp\mkern-2mu\cdotp}}

\newcommand{\PolyL}{\mathscr{G}} 
\newcommand{\PolyLF}{\scr{L}^\rm{f}} 
\newcommand{\MotCoLie}{\scr{L}^\rm{MTM}} 

\newcommand{\QMHTS}{\rm{MHTS}_\bb{Q}} 

\newcommand{\SGroup}{\mathfrak{S}} 
\newcommand{\BSt}{\mathrm{BSt}} 

\newcommand{\mahoganydashed}[1]{\smash{\textcolor{Mahogany}{\protect\begin{tikz}[baseline=(a.base)] \protect\node[rectangle,draw,dashed] (a) at (0,0) {#1};\end{tikz}}}}

\newcommand{\Cor}{\mathrm{Cor}} 
\newcommand{\CorG}{\mathrm{Cor}^\PolyL} 
 
\newcommand{\CorF}{\mathrm{Cor}^\mathrm{f}} 
\newcommand{\ItG}{\mathrm{I}^\PolyL}

\newcommand{\LiG}{\mathrm{Li}^\PolyL}

\newcommand{\LiH}{\mathrm{Li}^{\rm{Hod}}}
\newcommand{\dpw}{\rm{dpw}}
\newcommand{\nil}{\rm{nil}}
\newcommand{\red}{\rm{red}}
\newcommand{\DQ}{{\scr{D}_\bb{Q}}}

\newcommand{\Dec}{{\rm{Dec}}}
\newcommand{\sym}{{\rm{sym}}}

\newcommand{\smid}{\mkern1mu|\mkern1mu}
\DeclareMathSymbol{\shortminus}{\mathbin}{AMSa}{"39}

\newcommand{\circled}[1]{\raisebox{.5pt}{\textcircled{\raisebox{-.9pt} {#1}}}}

\DeclareFontFamily{U}{min}{}
\DeclareFontShape{U}{min}{m}{n}{<-> udmj30}{}

\AtBeginDocument{%
	\def\MR#1{}
}

\title{The Goncharov Lie coalgebra of a field}

\author{Alexander Kupers}
\address{Department of Computer and Mathematical Sciences, University of Toronto Scarborough, 1265 Military Trail, Toronto, ON M1C 1A4, Canada}
\email{a.kupers@utoronto.ca}

\author{Daniil Rudenko}
\address{Department of Mathematics,
University of Chicago,
5734 S. University Avenue, 
Chicago, IL 60637, USA}
\email{rudenkodaniil@uchicago.edu}

\author{Ismael Sierra}
\address{Department of Mathematics, University of Toronto. Bahen Centre 40 St. George Street, Room 6290. Toronto, Ontario Canada M5S 2E4}
\email{ismael.sierra@utoronto.ca}

\date{\today}

\begin{document}
	
\begin{abstract}This paper relates algebraic $K$-theory of fields to polylogarithms via general linear groups. We introduce the Goncharov Lie coalgebra, defined in terms of the $E_\infty$-homology of general linear groups. Using Steinberg modules, we find a presentation, compute its Lie cobracket, and construct motivic and Hodge realisations. Combining these results with the Rognes rank spectral sequence, we give symbolic descriptions of the rationalisation of the algebraic $K$-theory of fields beyond the cases studied by Matsumoto--Milnor and Bloch--Suslin: we express $\smash{K^{(3)}_4(F)}$ and the indecomposable part of $\smash{K^{(3)}_5(F)}$ in terms of Goncharov's polylogarithmic complex of weight 3.
\end{abstract}
	
\maketitle

\vspace{-.75cm}

\tableofcontents

\vspace{-.75cm} 

\section{Introduction}

Let $F$ be a field. The classifying spaces of the general linear groups $\GL_n(F)$ assemble to an $E_{\infty}$-algebra with multiplication given by block sum
\[
\BGLb(F)\simeq \bigsqcup_{n > 0} \BGL_n(F).
\]
The algebraic $K$-theory spectrum $K(F)$ is obtained from $\BGLb(F)$ by group completion, relating the groups
\[
K_n(F)\coloneq \pi_n(\BGLb(F)^\gc) \quad \text{for} \quad  n\geq 0
\]
to the stable homology of general linear groups. Its relationship to the unstable homology of general linear groups is mediated by the $E_\infty$-homology groups $\smash{H_{n,d}^{E_{\infty}}(\BGLb(F))}$, defined as the homology groups of its derived indecomposables (see \cref{sec:intro-einfty}). Working rationally, we combine vanishing results for these groups with Koszul duality to construct the following:

\begin{definition*}The \emph{Goncharov Lie coalgebra} of the field $F$ is given by
\[
\PolyL(F) \coloneq \bigoplus_{n\geq 1} \PolyL_n(F) \qquad \text{with} \qquad \PolyL_n(F) \coloneq H_{n,2n-1}^{E_{\infty}}(\BGLb(F)_\Q).
\]
\end{definition*}

Here the cobracket is a map taking $\scr{G}_n(F)$ to $\bigoplus_{n=n'+n''} \scr{G}_{n'}(F) \otimes \scr{G}_{n''}(F)$, satisfying antisymmetry and co-Jacobi properties (without any additional Koszul signs). This Lie coalgebra is the central object of our series of papers.

\subsection{The Goncharov Lie coalgebra and mixed Tate motives} Before giving a more detailed introduction to the results that appear in this paper, we outline the conjecture that informs them: through mixed Tate motives, the Goncharov Lie coalgebra should be connected to multiple polylogarithms and algebraic $K$-theory.

\subsubsection{Multiple polylogarithms} Classically, \emph{multiple polylogarithms} are multivalued functions depending on integers $n_1,\dots,n_k \in \N$ defined as
\begin{equation}\label{FormulaPolylogarithm}
\rm{Li}_{n_1,\compactldots, n_k}(a_1,\compactldots,a_k)= \sum_{0<m_1<m_2<\dots<m_k}\frac{a_1^{m_1} a_2^{m_2}\compactcdots a_k^{m_k}}{m_1^{n_1}m_2^{n_2}\compactcdots m_k^{n_k}}\,
\end{equation}
for complex arguments $|a_i|<1$ and then analytically continued. These have been studied for centuries, see \cite{HainClassical,Zagier, Gon01} and their references, for their special values, functional equations, and applications in other parts of mathematics. For example, Zagier's conjecture uses them to connect special values of $L$-functions to algebraic $K$-theory, and the works of Goncharov on this conjecture inform our papers \cite{Cathelineau,Gon95, Dup20}.

We give a polylogarithmic interpretation of the Goncharov Lie coalgebra by constructing explicit elements
\[
\LiG_{n_1,\compactldots, n_k}(a_1,\compactldots,a_k) \in \PolyL_n(F), \quad \text{for $n_i\in \N$ so that $n_1+\dots+n_k=n$ and $a_i \in F$,}
\]
proving that these generate $\scr{G}(F)$ as a $\bb{Q}$-vector space and giving a complete set of relations corresponding to a certain family of functional equations satisfied by multivalued functions \eqref{FormulaPolylogarithm} (see \cref{sec:functional equations}). A precise statement appears in \cref{thm:polyl-presentation-additive}, which uses a more convenient alternative generating set of \emph{correlators} (see the discussion following \cref{thm:polyl-presentation-cobracket}). For $n\leq 3$, we identify in \cref{thm:polyl-identification} the spaces $\PolyL_n(F)$ with objects that have appeared previously in works of Bloch, Suslin, and Goncharov and admit a similar polylogarithmic interpretation.

\subsubsection{Mixed Tate motives over a number field} \label{sec:intro-mtm-number-field}

For a number field $F$, these two instances of multiple polylogarithms can be related through the theory of \emph{mixed Tate motives}, see \cite{BD94,Levine,DG05,Dup20}. The central object of this theory is a category $\rm{MTM}_\bb{Q}(F)$ of rational mixed Tate motives, which the reader can informally think as follows: consider the universal cohomology theory for varieties and take in its target category those objects that are iterated extensions of simple pieces obtained from the cohomology projective spaces. This is a Tannakian tensor category and can be shown to be equivalent to the category $\smash{\rm{Comod}^\rm{fd}_{\scr{L}^\rm{MTM}(F)}(\rm{GrMod}_\bb{Q})}$ of finite-dimensional graded comodules over a Lie coalgebra $\scr{L}^\rm{MTM}(F)$. Goncharov constructed a family of elements \cite[(10.20)]{Gon19}
\[\rm{Li}^\rm{MTM}_{n_1,\compactldots,n_k}(a_1,\compactldots,a_k) \in \scr{L}^\rm{MTM}(F)\]
called \emph{motivic multiple polylogarithms}. Their periods are closely related to the functions \eqref{FormulaPolylogarithm}: for example, multiple polylogarithms appear in the monodromy matrix comparing de Rham and Betti realisations \cite{Ramakrishnan}, and real-valued single-valued variants of multiple polylogarithms appear when applying Goncharov's real period construction (combine \cite[Corollary 1.15]{Gon19} with \S 11.1.3 loc.cit.). In this paper we show these can be obtained from the corresponding elements in $\scr{G}_n(F)$: in \cref{thm:motivic-realisation} we will construct a functor of Tannakian categories
\[R^\rm{MTM} \colon \rm{Comod}^\rm{fd}_{\PolyL(F)}(\rm{GrMod}_\bb{Q}) \lra \rm{MTM}_\bb{Q}(F)\]
such that the induced morphism $r^\rm{MTM}$ of graded Lie coalgebras satisfies
\[r^\rm{MTM}\left(\LiG_{n_1,\compactldots,n_k}(a_1,\compactldots,a_k)\right)=\rm{Li}^\rm{MTM}_{n_1,\compactldots,n_k}(a_1,\compactldots,a_k).\]
In the sequel \cite{KRS2} we will prove that $R^\rm{MTM}$ is an equivalence, which yields an isomorphism $\scr{G}(F) \cong \scr{L}^\rm{MTM}(F)$. This has many concrete consequences: it implies all periods of mixed Tate motives over the number field $F$ can be expressed as rational linear combinations of periods of motivic multiple polylogarithms, implying the universality conjecture of Goncharov \cite[Conjecture~17a]{Gon95}, cf.~\cite{Bro12} for a similar result for mixed Tate motives over the integers $\bb{Z}$.

\subsubsection{Mixed Tate motives over a general field}\label{sec:motives-general}
In view of this, we propose the category $\smash{\rm{Comod}^\rm{fd}_{\scr{G}(F)}(\rm{GrMod}_\bb{Q})}$ of finite-dimensional graded comodules over $\PolyL(F)$ as a candidate for the abelian category of mixed Tate motives over a general field $F$, whose existence was conjectured by Beilinson and Deligne \cite{BD94} and is closely related to the Beilinson--Soul\'e vanishing conjecture \cite{Levine}:

\begin{customconj}{A.a} \label{conjecture main motivic} Let $F$ be an arbitrary field. The category of finite-dimensional graded comodules over the Goncharov Lie coalgebra $\PolyL(F)$ is equivalent to the conjectural category of mixed Tate motives over $F$.
\end{customconj}

\begin{remark*}There are other candidates for the Lie coalgebra of mixed Tate motives (or, equivalently, for its universal coenveloping Hopf algebra) \cite[\S 1.16]{Gon95b}. The first is the \emph{Hopf algebra $\mathcal{A}(F)$ of Aomoto polylogarithms}, which was constructed by Beilinson, Goncharov, Shekhtman, and Varchenko in \cite{BGSV90}. The second, in the special case $F=\C$, is the \emph{scissors congruence Hopf algebra} $\mathbb{S}_\bullet$ constructed by Goncharov in \cite{Gon99}. The third is the \emph{1-minimal model} $\cal{M}_{\cal{N}}$ due to Bloch--Kritz \cite{BlochKriz}. We intend to investigate the relationship between $\PolyL(F)$ and these candidates in future work.\end{remark*}

This conjecture predicts a relationship between $\PolyL(F)$ and the algebraic $K$-groups of $F$ that can be stated---and possibly (dis)proven---without referring to the (conjectural) motivic formalism. To state it, we use the Adams eigenspaces $\smash{K^{(r)}_i(F)_\bb{Q} \subseteq K_i(F)_\bb{Q}}$, given by those elements on which the Adams operation $\psi^p$ acts by $p^r$ \cite[IV.5]{Weibel}; we refer to $r$ as \emph{weight}. These agree with the associated graded $\rm{gr}^r_\gamma\,K_i(F)_\bb{Q}$ of the $\gamma$-filtration \cite[p.~500]{Soule}, the motivic cohomology group $H^{2r-i,r}(\rm{Spec}(F);\bb{Q})$ \cite[Lecture 19]{MVW}, and the higher Chow group $\rm{CH}^r(\rm{Spec}(F),i)_\bb{Q}$ \cite[Theorem 3.1]{LevineBloch}. Beilinson gave a conjectural formula for the algebraic $K$-theory groups in terms of Ext-groups in the category of mixed Tate motives \cite{BD94}:
\[K^{(n)}_{2n-i}(F)_{\bb{Q}} \cong \rm{Ext}^i_{\rm{MTM}_\bb{Q}(F)}(\bb{Q}(-n),\bb{Q}(0)).\]
The Tannakian formalism allows one to compute these Ext-groups in terms of the Chevalley--Eilenberg complex of the corresponding Lie coalgebra. We let $H^i(\PolyL(F))_n$ denote the $i$th homology group of the weight $n$ part of this complex, explicitly given by
\begin{equation}\label{eqn:polyl-ce}
0 \lra \PolyL_n(F) \longrightarrow (\Lambda^2 \PolyL(F))_n \longrightarrow \cdots \longrightarrow (\Lambda^{n-1} \PolyL(F))_n \longrightarrow \Lambda^n F^{\times}_\Q \lra 0.
\end{equation}
where by convention $\scr{G}_n(F)$ lies in degree $i=1$, $\smash{\Lambda^n F^\times_\bb{Q}}$ lies in degree $n$, and the differential is induced by the cobracket.

\begin{customconj}{A.b}\label{conjecture main gamma} Let $F$ be an arbitrary field. For $1\leq i \leq n$, there is an isomorphism
\[
K^{(n)}_{2n-i}(F)_{\Q} \overset{\cong}\lra H^i(\PolyL(F))_n.
\]
\end{customconj}

\stepcounter{atheorem}

If true then the right side provides a symbolic description for rationalised algebraic $K$-theory of fields, generalising those for Milnor $K$-theory due to Matsumoto and Milnor and for the Bloch group due to Bloch and Suslin. In the sequel \cite{KRS2} we provide a candidate for such an isomorphism that does not rely on the motivic formalism and prove that it is an isomorphism for $F$ a number field. More precisely, recall the rank conjecture says that the associated graded $\rm{gr}^\rm{prim}_r\,K_d(F)_\bb{Q}$ of the primitive rank filtration is isomorphic to $\smash{K_d^{(r)}(F)_\bb{Q}}$ \cite[\S 2.3]{Cathelineau} \cite[Conjecture 1.17]{Gon95b} \cite{deJeu}. We will construct maps
\[\rm{gr}^\rm{prim}_n\,K_{2n-i}(F)_\bb{Q} \lra H^i(\PolyL(F))_n\]
and show these induce isomorphisms in cases where we understand the left side well.

\subsubsection{The Goncharov Lie coalgebra and algebraic $K$-theory in low weights}
\cref{conjecture main gamma} is known to be true in several cases: (i) the case $i=n$ is a consequence of Suslin's identification of $\smash{K_n^{(n)}(F)_{\Q}}$ with the rationalised Milnor $K$-group $K_n^M(F)_\bb{Q}$ \cite[\S 2.7]{Sus84}, (ii) for $i=1, n=2$ it is equivalent to the exactness of the complex
\[
0\lra K^{(2)}_{3}(F)_{\Q}\lra B_2(F)_\Q\lra \Lambda^2 F^{\times}_\Q\lra K^{(2)}_{2}(F)_{\Q}\lra 0
\]
proved by Suslin \cite[\S 5]{Sus90}, and (iii) for any $i=n-1,n \geq 1$ a proof was given in \cite{Bol24}. In further support of \cref{conjecture main gamma}, we study the case $n=3$. We do so by combining our description of $\scr{G}(F)$ with a \emph{rank spectral sequence} related to that of Rognes \cite{Rognes}, resulting in a symbolic description of $\smash{K_4^{(3)}(F)_{\bb{Q}}}$ (note \cref{theorem weight 3} \eqref{enum:FormulaGoncharovConjecture2} is an instance of \cite{Bol24} but his proof is quite different) and the \emph{indecomposable part} of $\smash{K_5^{(3)}(F)_\bb{Q}}$:

\begin{atheorem}\label{theorem weight 3} Let $F$ be an arbitrary field. 
\begin{enumerate}
    \item \label{enum:FormulaGoncharovConjecture2}
There is an isomorphism
\[
K^{(3)}_{4}(F)_{\Q} \cong H^{2}(\PolyL(F))_3.
\]
    \item \label{enum:FormulaGoncharovConjecture3} There is an exact sequence
\[
K^{(2)}_{4}(F)_{\Q}\otimes F^{\times}_{\Q}
\lra
K^{(3)}_{5}(F)_{\Q}
\lra 
H^{1}(\PolyL(F))_3\lra 0.
\]
\end{enumerate}
\end{atheorem}

In the case $n=3$, \cref{conjecture main gamma} is equivalent to a conjecture of Goncharov \cite[Conjecture 1.15]{Gon95b} concerning his weight $3$ polylogarithmic complex $\Gamma_3(F;\bb{Q})$, given by
\[
B_3(F)_{\Q}\lra B_2(F)_{\Q}\otimes F^{\times}_{\Q} \lra \Lambda^3 F^{\times}_{\Q}.
\]
It says that homology groups of this complex are given by $K_5^{(3)}(F)_\bb{Q}$, $K_4^{(3)}(F)_\bb{Q}$, and $K^M_3(F)_\bb{Q}$ respectively. \cref{theorem weight 3} affirms a variant of this conjecture, using that under the identifications of \cref{thm:polyl-identification}, this agrees with \eqref{eqn:polyl-ce} in weight $n=3$, given by
\[\scr{G}_3(F) \lra \scr{G}_2(F) \otimes \scr{G}_1(F) \lra \Lambda^3 \scr{G}_1(F).\]

\smallskip

We now turn to a more detailed introduction to the results contained in this paper, in particular \cref{theorem weight 3}.

\subsection{$E_\infty$-homology of general linear groups} \label{sec:intro-einfty} We start by recalling the approach to studying the homology of general linear groups introduced by Galatius, Kupers, and Randal-Williams \cite{GKRW18,GKRW19,GKRW20} and how it leads to the Goncharov Lie coalgebra. This takes place in a homotopy-theoretic context, so we work in $\infty$-categories or, if the reader prefers, suitable model categories.

\subsubsection{The $E_\infty$-algebra $\BGLb(F)_\bb{Q}$}
The block sum of matrices induces a product on chains
\[C_*(\GL_n(F);\bb{Q})\otimes C_*(\GL_m(F);\bb{Q}) \lra C_*(\GL_{n+m}(F);\bb{Q})\]
which is associative and commutative up to chain homotopy. It is convenient to keep track of the dimension (or \emph{rank}) $n$ by defining a functor
\begin{align*}\BGLb(F)_\bb{Q} \colon \bb{N} &\lra \DQ \\
n &\longmapsto C_*(\BGL_n(F);\bb{Q}),\end{align*}
where $\DQ$ is the derived category of $\bb{Q}$ obtained by inverting the quasi-isomorphisms in the category of rational chain complexes. That the above product admits coherent associativity or commutativity chain homotopies means that this functor lifts to a nonunital $E_\infty$-algebra $\BGLb(F)_\bb{Q} \in \Alg_{E_\infty^\rm{nu}}(\Fun(\bb{N},\DQ))$, where $\Fun(\bb{N},\DQ)$ is endowed with the Day convolution symmetric monoidal structure. 

There is a well-developed theory of derived indecomposables of nonunital $E_\infty$-algebras such as $\BGLb(F)_\bb{Q}$, e.g.~\cite{GKRW18}. Informally this recovers its generators and relations, as homotopy-theoretically there is no distinction between these. More precisely, the functor making an object into a nonunital $E_\infty$-algebra with trivial multiplication has a left adjoint
\[\cot_{E_\infty^\rm{nu}} \colon \Alg_{E_\infty^\rm{nu}}(\Fun(\bb{N},\DQ)) \lra \Fun(\bb{N},\DQ)\]
known as the \emph{cotangent complex} or \emph{$E_\infty$-indecomposables}; we will use the latter terminology. In our rational setting, it may be computed by strictifying to a commutative dg-algebra and taking Harrison homology; for example, for a minimal commutative dg-algebra this yields its indecomposables in an underived sense, which are canonically isomorphic to its generators. 

\begin{figure}[h]
	\begin{tikzpicture}
	\begin{scope}
	\clip (-2,-1) rectangle ({2.5*4+2},6.5);
	\draw (0,0)--(10.5,0);
	\draw (0,0) -- (0,6.5);
	\foreach \s in {0,...,6}
	{
		\draw [dotted] (-1.5,\s)--(10.5,\s);
		\node [fill=white] at (-1.5,\s) [left] {\tiny $\s$};
	}
	\foreach \s in {0,...,4}
	{
		\draw [dotted] ({2.5*\s},-0.5)--({2.5*\s},6.5);
		\node [fill=white] at ({2.5*\s},-.5) {\tiny $\s$};
	}
	\draw [very thick,Mahogany,densely dotted] (2.5,1) -- (10,7);

	\node [fill=white] at (2.5,0) {$\ds{Q}$};
	\node [fill=white] at (2.5,1) {$F^\times_\ds{Q} = \scr{G}_1(F)$};
	\node [fill=white] at (2.5,2) {$\Lambda^2 F^\times_\ds{Q}$};
	\node [fill=white] at (2.5,3) {$\Lambda^3 F^\times_\ds{Q}$};
	\node [fill=white] at (2.5,4) {$\Lambda^4 F^\times_\ds{Q}$};
	\node [fill=white] at (2.5,5) {$\Lambda^5 F^\times_\ds{Q}$};
	\node [fill=white] at (2.5,6) {$\Lambda^6 F^\times_\ds{Q}$};

	\node [fill=white] at (5,3) {$\scr{G}_2(F)$};
	\node [fill=white] at (5,4) {?};
	\node [fill=white] at (5,5) {?};
	\node [fill=white] at (5,6) {?};

	\node [fill=white] at (7.5,5) {$\scr{G}_3(F)$};
	\node [fill=white] at (7.5,6) {?};
	
	\node at (-.5,-.5) {$\nicefrac{d}{n}$};
	\end{scope}
	\end{tikzpicture}
	\caption{The $E_\infty$-homology of $\BGLb(F)_\ds{Q}$, which vanishes for $d \leq 2n-2$ as long as $n \geq 2$. In this grading convention the Lie cobracket has degree $-1$, so the entries on the dashed line assemble to a Lie coalgebra.}
	\label{fig:gltable1} 
\end{figure}
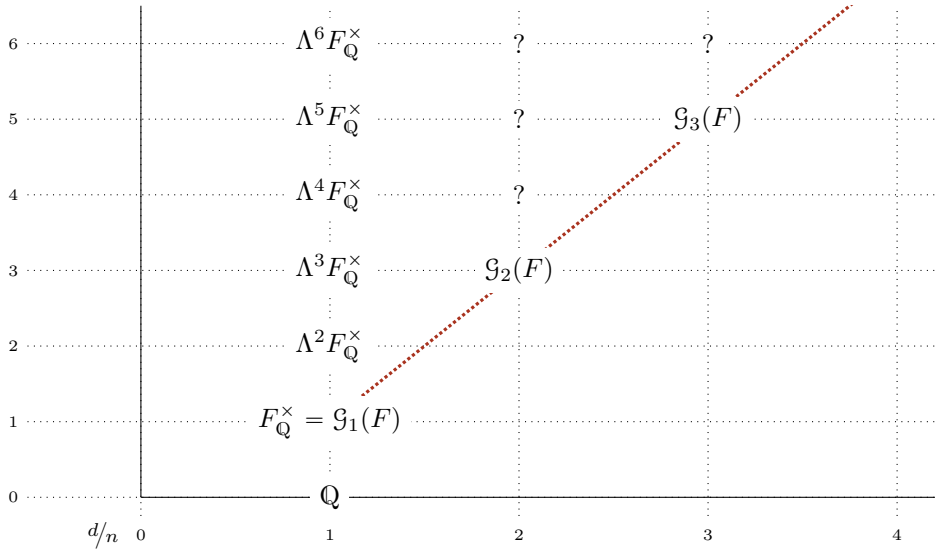

The \emph{$E_\infty$-homology groups} of $\BGLb(F)_\bb{Q}$ are then defined as
\[H^{E_\infty}_{n,d}(\BGLb(F)_\bb{Q}) \coloneq H_d(\cot_{E_\infty^\rm{nu}}(\BGLb(F)_\bb{Q})(n))\]
and one of the main results of \cite{GKRW20} is a determination of some of these groups, see \cref{fig:gltable1}: in particular, they vanish for $d \leq 2n-2$ with the exception of $(n,d) = (1,0)$. The Goncharov Lie coalgebra is defined as those entries on the critical line above where this vanishing result applies:
\[
\PolyL(F) = \bigoplus_{n\geq 1} H_{n,2n-1}^{E_{\infty}}(\BGLb(F)_\Q).
\]
Its name is justified by Koszul duality between the nonunital commutative operad and the suspended Lie cooperad, which implies that $E_\infty$-indecomposables admit the structure of a shifted Lie coalgebra: for a minimal commutative dg-algebra this is the quadratic part of the differential on generators. This endows the $E_\infty$-homology groups $\smash{H^{E_\infty}_{n,d}(\BGLb(F)_\bb{Q})}$ with a cobracket of degree $-1$, restricting to a Lie coalgebra structure on $\PolyL(F)$; the gradings work out to make so that Lie cobracket does not involve additional Koszul signs. It is a construction for $E_\infty$-algebras with a slope 2 vanishing line that is analogous to the ``stability Hopf algebra'' of \cite{RWchromatic} for $E_2$-algebras with a slope 1 vanishing line, and encodes the attaching maps for $E_\infty$-cells in the sense of \cite{GKRW18}.

\subsubsection{The presentation of the Goncharov Lie coalgebra}
Once one defines the Goncharov Lie coalgebra $\scr{G}(F)$, one should ask what it is. The first main result of this paper attempts to answer this question by giving a presentation of it.

\begin{customthm}{C.a}\label{thm:polyl-presentation-additive}
    The Lie coalgebra $\PolyL(F)$ is generated as a $\bb{Q}$-vector space by \emph{correlators}
    \[\Cor^\scr{G}(x_0,x_1,\compactldots,x_n) \in \scr{G}_n(F) \qquad \text{for $x_0,\ldots,x_n \in F$ not all equal }\]  
    subject to the following relations:
    \begin{enumerate}[\noindent (1)]
\item \label{enum:rel-goncharov-1} Homogeneity: $\Cor^\scr{G}(x_0,x_1,\compactldots,x_n)=\Cor^\scr{G}(x_0+b,x_1+b,\compactldots,x_n+b)$ for $b\in F$.
\item \label{enum:rel-goncharov-2} Cyclic symmetry: $\Cor^\scr{G}(x_0,x_1,\compactldots,x_n)=\Cor^\scr{G}(x_1,x_2,\compactldots,x_0)$.
\item \label{enum:rel-goncharov-3} Shuffle relations: 
\[\sum_{\sigma \in \rm{Sh}(n_1,n_2)} \Cor^\scr{G}(x_0,x_{\sigma(1)},\compactldots,x_{\sigma(n_1+n_2)}) =0 \quad 
\text{for $n=n_1+n_2$, $n_1,n_2>0$.}\]
\item \label{enum:rel-goncharov-4} Decomposition relations:
\begin{align*}
    &\CorG(x_0,\compactldots,x_n)-\CorG(y_0,\compactldots,y_n)\\
     &=\sum_{\iota=((i_1,j_1),\compactldots,(i_{n},j_n))\in T(n)}\rm{sign}(\iota)\,
     \CorG\left(0,\frac{x_{i_1}-x_{j_1}}{y_{i_1}-y_{j_1}},\compactldots,\frac{x_{i_n}-x_{j_n}}{y_{i_n}-y_{j_n}}\right),
\end{align*}
where we omit terms with $y_{i_k}=y_{j_k}$ for some k, and the set $T(n)$ of pairs of indices as well as the sign $\rm{sign}(\iota)$ are given by \cref{prop:universal-symbol-combinatorics}.
\end{enumerate}
\end{customthm}

\begin{remark*}The expression in the decomposition relation \eqref{enum:rel-goncharov-4} can be inductively determined, takes the same form as the symbol for formal correlators, and admits a formula as a sum over trees. In practice, it is sufficient to merely know its form, rather than the details of $T(n)$ or $\rm{sign}(\iota)$.\end{remark*}

The decomposition relations imply that $\CorG(x_0,\compactldots,x_n)$ with all $x_0,\ldots,x_n$ distinct generate $\PolyL_n(F)$ (\cref{proposition: generic correlators}), so the following completely determines the Lie cobracket:

\begin{customthm}{C.b}\label{thm:polyl-presentation-cobracket} With respect to the presentation of \cref{thm:polyl-presentation-additive},
the Lie cobracket is given by
\[\delta(\CorG(x_0,x_1,\compactldots,x_n)) = \sum_{j=0}^n \sum_{i=1}^{n-1} \CorG(x_j,x_{j+1},\compactldots,x_{j+i}) \wedge \CorG(x_j,x_{j+i+1},\compactldots,x_{j+n})\]
as long as $x_0,\ldots,x_n$ are distinct.
\end{customthm}

There are other generating sets, e.g.~the analogues $\LiG_{n_1,\compactldots,n_k}(a_1,\compactldots,a_k)$ of multiple polylogarithms. One has for example
\[\LiG_n(a)=-\CorG(1,\underbrace{0,\compactldots,0}_{n-1},a),\]
and a general formula expressing multiple polylogarithms in terms of correlators (and vice versa) can be deduced from \cref{sec:iterated-integrals-and-multiple-polylogarithms}.

It can be difficult to make use of a presentation, but \cref{thm:polyl-presentation-additive} is such that one can perform calculations in $\scr{G}(F)$: through it one can understand the action of the duality involution (it is by $(-1)^n$, see \cref{theorem: duality}) and, in characteristic $p$, the action of Frobenius endomorphism (see \cite{KRS2}). In fact, for $n \leq 3$ we can find identifications with more classical groups:

\begin{customthm}{C.c}
\label{thm:polyl-identification} There are isomorphisms
\[\begin{aligned}F^\times_\bb{Q} &\overset{\cong}\lra \scr{G}_1(F) \\
a &\longmapsto \CorG(0,a),\end{aligned} \qquad \begin{aligned}B_2(F)_\bb{Q} &\overset{\cong}\lra \scr{G}_2(F) \\
\{a\}_2 &\longmapsto -\CorG(1,0,a),\end{aligned} \qquad\begin{aligned}B_3(F)_\bb{Q} &\overset{\cong}\lra \scr{G}_3(F) \\
\{a\}_3 &\longmapsto -\CorG(1,0,0,a),\end{aligned}\]
where $B_2(F)$ is the Bloch group studied by Bloch and Suslin \cite{Bloch,Sus90} (who call it the pre-Bloch group), and $B_3(F)$ was introduced by Goncharov \cite{Gon95b}.
\end{customthm}

\stepcounter{atheorem}

\subsection{Steinberg modules} The presentation of \cref{thm:polyl-presentation-additive} is obtained by relating the $E_\infty$-homology groups of $\BGLb(F)_\bb{Q}$ to Steinberg modules. We will first recall the Steinberg modules and their double and infinite variants, giving a particularly nice generating set for the latter, and explain how these are related to the $E_\infty$-homology groups discussed in the previous subsection.

\subsubsection{Steinberg modules, double Steinberg modules, and infinite Steinberg modules} 

\begin{figure}[h]
\begin{tikzpicture}
\draw (0,0) -- (-1,4) -- (2,2) -- cycle;
\node at (0,0) [left] {$v_3$};
\node at (-1,4) [left] {$v_2$};
\node at (2,2) [right] {$v_1$};
\end{tikzpicture}
\qquad \begin{tikzpicture}
\draw (0,0) -- (-1,4) -- (3,2) -- cycle;
\node at (0,0) [left] {$v_3$};
\node at (-1,4) [left] {$v_2$};
\node at (3,2) [right] {$v_1$};
\node at (1,0) [right] {$w_1$};
\node at (-1.5,2) [left] {$w_2$};
\node at (1.5,3.5) [right] {$w_3$};
\draw (1,0) -- (-1.5,2) -- (1.5,3.5) --cycle;
\end{tikzpicture}
\caption{A generator of $\St_3(F)$ (on the left) and a generator of $\StH_3(F)$ (on the right).}
\label{fig:generators-st-sth} 
\end{figure}
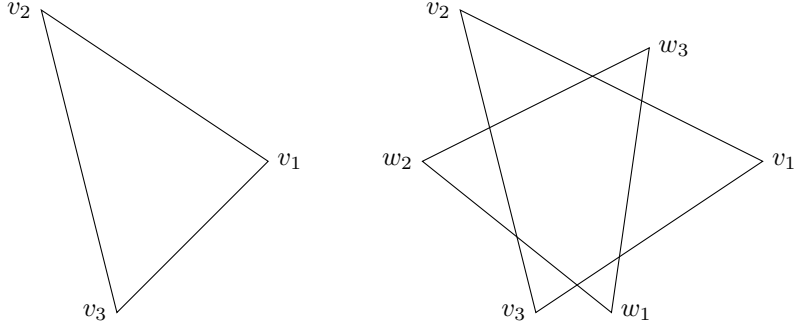

The \emph{Tits building} $T(F^n)$ is defined as the geometric realisation of the nerve of the poset of proper nonzero subspaces of $F^n$. By the Solomon--Tits theorem it is equivalent to a wedge of $(n-2)$-spheres, and its top reduced rational homology is a $\GL_n(F)$-module known as the \emph{Steinberg module}
\[\St_n(F) \coloneq \widetilde{H}_{n-2}(T(F^n);\bb{Q}).\]
This has a well-known presentation with generators given by apartment classes $[v_1,\compactldots,v_n]$ indexed by ordered bases $v_1,\compactldots,v_n$; we suggest representing these by simplices in projective space $\bb{P}^{n-1}(F)$ as on the left of \cref{fig:generators-st-sth}.

To explain further structure, we organise all Steinberg modules into a single object. Let $\Vect_F$ denote the groupoid of finite-dimensional vector spaces over $F$ and linear isomorphisms, then we can think of the Steinberg modules as a functor (with $\St(V)$ in degree $\dim(V)$)
\begin{align*}\SSt \colon \Vect_F &\lra \rm{GrMod}_\bb{Q} \\
V &\longmapsto \St(V).\end{align*}
The category $\Vect_F$ admits not only a symmetric monoidal structure given by direct sum but also a promonoidal structure given by ``flag sum'', which assemble to a produoidal category. By Day convolution, these induce on functors $\Vect_F \to \rm{GrMod}_\bb{Q}$ a symmetric monoidal structure whose tensor product $\levi$ we call the \emph{Levi tensor  product} and a monoidal structure whose tensor product $\para$ we call the \emph{parabolic tensor product}, which assemble to a duoidal category. The Steinberg modules $\SSt$ form a bialgebra, with product and coproduct
\[\parbox{4cm}{\centering $(\SSt \levi \SSt)(V) \cong$ \\
$\bigoplus_{V_1 \oplus V_2 = V} \St(V_1) \otimes \St(V_2)$} \lra \St(V) \quad \text{and} \quad \St(V) \lra \parbox{4cm}{\centering $(\SSt \para \SSt)(V) \cong $ \\ 
$\bigoplus_{U \subseteq V} \St(U) \otimes \St(V/U)$}\]
given by concatenation and splitting of apartments (see \cref{prop:st-explicit-pres} for explicit formulas).

To see this, one considers the functor $\ul{\bb{Q}}{}_{>0} \colon \Vect \to \DQ$, which takes the value $\bb{Q}$ on nonzero vector spaces and $0$ otherwise, which admits the structure of a commutative algebra with respect to $\levi$ and an associative algebra with respect to $\para$. Then computing its $E_1$-indecomposables with respect to $\para$ via a bar construction, one finds these are given by the (rationalised) Steinberg modules, but now it has a coassociative coproduct with respect to $\para$ and a remaining commutative product with respect to $\levi$. Similar algebraic structures were used in \cite{AMP,BCGP}.

A crucial property of $\SSt$ as a commutative algebra is that it is Koszul \cite{MNP,MPW23,CharltonRadchenkoRudenko}. Its Koszul dual as an associative algebra is given by the \emph{double Steinberg modules} $\SStH$ 
\[\StH_n(F) \cong \St_n(F) \otimes \St_n(F),\] 
which merely means that the associative bar complex
\[0 \to \StH_n(F) \to (\SSt^{\levi n})(F^n) \to (\SSt^{\levi n-1})(F^n) \to \cdots \to (\SSt^{\levi 2})(F^n) \to \St_n(F) \to 0\]
is exact. We suggest representing elements of $\StH_n(F)$ by pairs of simplices in projective space $\bb{P}^{n-1}(F)$ as on the right of \cref{fig:generators-st-sth}. It has a product by concatenating apartments termwise and the Koszul dual of $\SSt$ as a commutative algebra is given by the \emph{infinite Steinberg modules} $\SStL$, which can be defined as the indecomposables with respect to the algebra structure on the double Steinberg modules:
\[\StL(V) \cong \rm{coker}\left[\bigoplus_{\substack{V_1 \oplus V_2 = V \\ V_1,V_2 \neq 0}} \StH(V_1)\otimes \StH(V_2) \lra \StH(V)\right].\]
We believe that the infinite Steinberg module is a remarkable object worth studying independently. It can also be obtained from the common basis complex \cite{Rognes} or the partial decomposition poset \cite{HHS,BPW}.

By Koszul duality, $\SStL$ admits the structure of a Lie coalgebra with respect to $\levi$. Charlton, Radchenko, and Rudenko gave a presentation of the infinite Steinberg modules in which the cobracket takes a simple form \cite[Section 2]{CharltonRadchenkoRudenko}: $\StL_n(F)$ is generated by Steinberg correlators $\rm{C}[u_0: \compactcdots u_n]$ indexed by affine bases $u_0,\ldots,u_n$ (this means that $u_1-u_0,\ldots,u_n-u_0$ are a basis) defined as (see \cref{fig:generator-stl})
\[\rm{C}[u_0: \compactcdots: u_n] \coloneq \pi\big[[u_n-u_0,\compactldots,u_1-u_0] \otimes [u_n-u_0,u_{n-1}-u_n,\compactldots,u_1-u_2]\big]\]
with $\pi \colon \StH_n(F) \to \StL_n(F)$ the projection onto indecomposables. These satisfy homogeneity, cyclic symmetry, and shuffle relations. In terms of these the cobracket is given by
\[\delta\big(\rm{C}[u_0:\cdots:u_n]\big) = \sum_{j=0}^n \sum_{i=1}^{n-1} \rm{C}[u_j:\compactcdots:u_{j+i}] \wedge \rm{C}[u_j:u_{j+i+1}:u_{j+i+2}:\compactcdots:u_{j+n}].\]

\begin{figure}\begin{tikzpicture}
\draw (0,0) -- (-1,4) -- (3,2) -- cycle;
\node at (0,0) [left] {$u_3-u_0$};
\node at (-1,4) [left] {$u_2-u_0$};
\node at (3,2) [right] {$u_1-u_0$};
\node at (-.5,2) [left] {$u_2-u_3$};
\node at (1,3) [right] {$\,u_1-u_2$};
\draw (0,0) -- (-.5,2) -- (1,3) --cycle;
\end{tikzpicture}
\caption{A generator of $\StL_3(F)$.}
\label{fig:generator-stl} \end{figure}

\begin{remark*}In \cref{sec:higher-apts} we will explain how to obtain both the correlators and their cobracket by a formal procedure given only the description of Steinberg modules in terms of apartments, as what we call \emph{higher apartments}.\end{remark*}

\subsubsection{Relationship to $E_\infty$-homology}
We can relate $\BGLb(F)_\bb{Q}$ to Steinberg modules by constructing this nonunital $E_\infty$-algebra in a different manner. Taking the dimension gives a symmetric monoidal functor $\dim \colon \Vect_F \to \bb{N}$. The nonunital $E_\infty$-algebra $\BGLb(F)_\bb{Q}$ can be obtained as $\dim_!(\ul{\bb{Q}}{}_{>0})$, the left Kan extension along $\dim$ of the nonunital $E_\infty$-algebra $\bb{Q}_{>0}$. The advantage of this perspective is that constructions such as indecomposables commute with $\dim_!$. This can be used to prove that there are isomorphisms \cite[Section 6]{GKRW20}
\begin{equation}\label{eqn:he1-st-heinfty-stl} \begin{aligned} H^{E_1}_{n,d}(\BGLb(F)_\bb{Q}) & \overset{\cong}\lra H_{d-n+1}(\GL_n(F);\St_n(F)), \\
H^{E_2}_{n,d}(\BGLb(F)_\bb{Q}) & \overset{\cong}\lra H_{d-2n+2}(\GL_n(F);\StH_n(F)), \\
H^{E_\infty}_{n,d}(\BGLb(F)_\bb{Q}) & \overset{\cong}\lra H_{d-2n+2}(\GL_n(F);\StL_n(F)),\end{aligned}\end{equation}
for $n \geq 1$, where $H^{E_1}_{*,*}(-)$ and $H^{E_2}_{*,*}(-)$ are obtained as the indecomposables with respect to the nonunital $E_1$- or $E_2$-algebra structures obtained by forgetting commutativity. In particular, there is an isomorphism
\[\PolyL_n(F) = H^{E_\infty}_{n,2n-1}(\BGLb(F)_\bb{Q}) \cong H_1(\GL_n(F);\StL_n(F)).\]
The top isomorphism of \eqref{eqn:he1-st-heinfty-stl} is induced by taking $E_1$-indecomposables of $\ul{\bb{Q}}{}_{>0}$ with respect to $\para$ and uses that applying $\dim_!$ to the comparison map from $\levi$ to $\para$ yields an isomorphism, an observation due to Nesterenko and Suslin \cite[\S 1]{NesterenkoSuslin}. The bottom isomorphism of \eqref{eqn:he1-st-heinfty-stl} can similarly be obtained by first taking $E_1$-indecomposables of $\ul{\bb{Q}}{}_{>0}$ with respect to $\para$ and then $E_\infty$-indecomposables with respect to $\levi$.

This reveals that $\SStL$ has additional structure: in addition to the cobracket with respect to $\levi$ it has a compatible coproduct with respect to $\para$, though by an Eckmann--Hilton argument the cobracket determines the coproduct. The cobracket on $E_\infty$-homology is \emph{not} induced by the cobracket on $\SStL$ (it would have degree $0$ rather than degree $-1$) but rather is a secondary cobracket arising from the compatibility of cobracket and coproduct with respect to tensor products that are made equal upon applying $\dim_!$. This should call to mind the Dunn--Lurie additivity theorem \cite[5.1.2.2]{LurieHA}.

\subsection{From infinite Steinberg modules to a presentation for $\PolyL(F)$.} \label{sec:intro-presentation} 
We next explain how we obtain a presentation of $\PolyL_n(F)$ from that for $\StL_n(F)$, leading to \cref{thm:polyl-presentation-additive,thm:polyl-presentation-cobracket}. This is an outline of the arguments that comprise the first half of this paper.

\subsubsection{A projective resolution of the infinite Steinberg module} \label{sec:intro-presentation-resolution} The presentation of $\PolyL_n(F)$ is obtained by constructing a novel projective resolution of $\StL_n(F)$ and using the identification of $\scr{G}_n(F)$ as $H_1(\GL_n(F);\StL_n(F))$. 

We start by constructing a generating set. The $\GL(V)$-module $\FC(V)$ is the $\bb{Q}$-vector space generated by symbols $\FC[u_0:\compactcdots:u_n]$ for affine bases $u_0,\ldots, u_n$ of $V$, which we call \emph{formal correlators}, subject to the following relations: 
\begin{itemize}
\item Homogeneity: $\rm{FC}[u_0:\compactcdots:u_n] = \rm{FC}[u_0-u:\compactcdots:u_n-u]$ for any $u \in V$.
\item Cyclic symmetry: $\rm{FC}[u_0:u_1:\compactcdots:u_n] = \rm{FC}[u_1:u_2:\cdots:u_0]$.
\item Shuffle relations: 
\[\sum_{\sigma \in \rm{Sh}(n_1,n_2)} \rm{FC}[u_0:u_{\sigma(1)}:\compactcdots:u_{\sigma(n_1+n_2)}] =0 \quad \text{for $n=n_1+n_2$ with $n_1,n_2>0$.}\]
\end{itemize}
As a $\GL(V)$-module $\FC(V)$ is projective, and there is an evident surjection
\begin{align*}
\pr^\FC \colon \FC(V)&\lra \StL(V)\\
\FC[u_0:\compactcdots:u_n] &\longmapsto \rm{C}[u_0:\compactcdots:u_n]
\end{align*}
which is the start of a projective resolution
\begin{equation}\label{eqn:dec-resolution}\cdots \overset{d}\lra \bb{Q}[\Dec_{V}]^{\otimes 2} \otimes \FC(V) \overset{d}\lra \bb{Q}[\Dec_{V}] \otimes \FC(V) \overset{d}\lra \FC(V) \overset{\pr^\FC}\lra \StL(V)\lra 0.\end{equation}
The further terms use the key observation, based on \cite[Section 3.9]{CharltonRadchenkoRudenko}, that the vector space $\StL(V)$ has a collection of ``almost bases'' labelled by nonzero linear functionals $h\in V^{\vee}$. Namely, for every such functional, the Steinberg correlators $\rm{C}[0:v_1:\compactcdots:v_n]$ with $h(v_i) = 1$ span $\StL(V)$ as a $\Q$-vector space; they are almost a basis in the sense that all the relations between such elements for the same $h$ follow from the shuffle relations. This allows us to introduce \emph{decomposition operators}
\[
D^\FC_{h} \colon \FC(V) \lra \FC(V)
\]
which are obtained by projecting $\FC(V)$ to $\StL(V)$, expanding the obtained element as a linear combination of Steinberg correlators $\rm{C}[0:v_1:\compactcdots:v_n]$ with $h(v_i) = 1$, and then lifting these to $\FC(V)$. This is well-defined because, for fixed $h$, the only relations among the normalized correlators are the shuffle relations, which are imposed in $\FC(V)$.

\begin{example*}For $V = F^2$ the decomposition operator is given by the following formula:
\begin{align*}
D^\FC_h(\FC[u_0:u_1:u_2])
&=\,\FC\left[0:\frac{u_{1}-u_{0}}{h(u_{1}-u_0)}:\frac{u_{2}-u_{0}}{h(u_{2}-u_0)}\right]\\
&-\FC\left[0:\frac{u_{1}-u_{0}}{h(u_{1}-u_0)}:\frac{u_{2}-u_{1}}{h(u_{2}-u_1)}\right]\\
&+\FC\left[0:\frac{u_{2}-u_{0}}{h(u_{2}-u_0)}:\frac{u_{2}-u_{1}}{h(u_{2}-u_1)}\right],
\end{align*}
where if a denominator is zero that term must be dropped.
\end{example*}

In general, the formula for the decomposition operator is equally explicit but more intricate, see \cref{sec:symbol-maps}. It has the same form as the decomposition relation of \cref{thm:polyl-presentation-additive} \eqref{enum:rel-goncharov-4}:
\begin{align*}
    &D^\FC_h(\FC[u_0:\compactldots:u_n])\\
     &=\sum_{\iota=((i_1,j_1),\ldots,(i_{n},j_n))\in T(n)}\rm{sign}(\iota)\,
     \FC\left[0:\frac{u_{i_1}-u_{j_1}}{h(u_{i_1}-u_{j_1})},\compactldots,\frac{u_{i_n}-u_{j_n}}{h(u_{i_n}-u_{j_n})}\right],
\end{align*}
for the set $T(n)$ of pairs of indices and signs $\rm{sign}(\iota)=\pm 1$ of \cref{prop:universal-symbol-combinatorics}.

Let $\Dec_{V}$ denote the set of nonzero linear functionals on $V$, and define \eqref{eqn:dec-resolution} as
\[\cdots \overset{d}\lra \bb{Q}[\Dec_{V}]^{\otimes 2} \otimes \FC(V) \overset{d}\lra \bb{Q}[\Dec_{V}] \otimes \FC(V) \overset{d}\lra \FC(V) \overset{\pr^\FC}\lra \StL(V)\lra 0\]
with differential is akin to that in the bar construction, given by an alternating sum of forgetting functionals and applying a decomposition operator: for example, $d([h_1|h_2] \otimes \FC[u_0 : \compactcdots : u_n])$ is given by $[h_2] \otimes \FC[u_0 : \compactcdots : u_n]-[h_1] \otimes \FC[u_0 : \compactcdots : u_n]+[h_1] \otimes D_{h_2}^\FC(\FC[u_0 : \compactcdots : u_n])$.

\subsubsection{Presentation of $\PolyL_n(F)$ as a $\bb{Q}$-vector space} We now use this projective resolution to give a presentation for $\PolyL_n(F)$. The natural map $H_0(\GL(V);\FC(V)) \to H_0(\GL(V);\StL(V))$ is an isomorphism: both sides vanish for $\dim V>1$, and for $\dim V=1$ both are canonically isomorphic to $\bb{Q}$. By taking coinvariants we obtain from \eqref{eqn:dec-resolution} an exact sequence
\[
    \bigl(\bb{Q}[\Dec_{V}]^{\otimes 2}\otimes \FC(V)\bigr)_{\GL(V)} \lra \bigl(\bb{Q}[\Dec_{V}]\otimes \FC(V)\bigr)_{\GL(V)}  \lra \PolyL_n(F)\lra 0.
\]
The projection of the element $[h]\otimes \FC[u_0:\compactcdots:u_n]$ to the coinvariants $\PolyL_n(F)$ depends only on the elements $h(u_0), \ldots, h(u_n)\in F$ and its image is by definition the correlator
\[\CorG(h(u_0),\compactldots,h(u_n)) \in \PolyL_n(F)\] 
that appears in \cref{thm:polyl-presentation-additive}. The remainder of that theorem is now a consequence of the above exact sequence and the formula for the decomposition operator.

\begin{example*} \label{ex:decomposition weigh 2} For $V =F^2$, the element $[h_1|h_2]\otimes \FC[u_0:u_1:u_2]$ yields the relation
\begin{align*}
    &\CorG(x_0,x_1,x_2)-\CorG(y_0, y_1,y_2)\\
     &=\CorG\Bigl(0,\frac{x_1-x_0}{y_1-y_0},\frac{x_2-x_0}{y_2-y_0}\Bigr)-\CorG\Bigl(0,\frac{x_1-x_0}{y_1-y_0},\frac{x_2-x_1}{y_2-y_1}\Bigr)+\CorG\Bigl(0,\frac{x_2-x_0}{y_2-y_0},\frac{x_2-x_1}{y_2-y_1}\Bigr),
\end{align*}
where $x_i=h_1(u_i)$ and  $y_i=h_2(u_i)$. Under the isomorphism of \cref{thm:polyl-identification} between $\scr{G}_2(F)$ and the classical Bloch group $B_2(F)_\bb{Q}$, this is equivalent to the 5-term relation. 
\end{example*}

\subsubsection{Presentation of $\scr{G}(F)$ as a Lie coalgebra} Now that we have the presentation of \cref{thm:polyl-presentation-additive}, we give a formula for the aforementioned secondary cobracket in terms of correlators. The result is \cref{thm:polyl-presentation-cobracket}, and it is obtained by similarly reducing the computation to one on coinvariants. 

The starting point is the bigraded Hopf algebra $\rm{H}^+ \coloneq \bigoplus_{n \geq 0} H_*(\GL_n(F);\St_n(F))$, where we place $\St_n(F)$ in homological degree $n$. This vanishes for $d<2n-2$ (\cref{thm:steinberg-homology-improved-vanishing}) and one can obtain the Lie coalgebra $\bigoplus_{n \geq 0} H_*(\GL_n(F);\StL_n(F))$, where we place $\StL_n(F)$ in homological degree $2n$, as its indecomposables $\rm{H}/\rm{H}^2$ and the cobracket as induced by the antisymmetrisation of its reduced coproduct. We then model the zigzag
\[\rm{H}/\rm{H}^2 \overset{\pi}\longleftarrow \rm{H}/\rm{H}^3 \xrightarrow{\ol{\Delta}-\sigma \circ \ol{\Delta}} \rm{H}/\rm{H}^2 \otimes \rm{H}/\rm{H}^2\]
with leftwards map surjective, by applying $\dim_!$ to a zigzag of chain complexes
\[[\St^\infty \to 0] \longleftarrow [\St^\infty \to \Lambda^2 \St^\infty] \lra [0 \to \St^\infty \para \St^\infty].\]
This zigzag is then resolved in a similar manner as above by formal correlators and decomposition operators. When computing the cobracket in the generic case with these resolutions, one gets to make choices of lifts and the subtlety lies in making a good choice there. We defer to \cref{sec:cobracket-outline} for a more detailed outline.

\subsection{Realisations} \label{sec:intro-realisations} In his works on conjectures relating multiple polylogarithms to mixed Tate motives and algebraic $K$-theory, Goncharov defines several $\bb{Q}$-vector spaces of formal polylogarithms subject to an inductively defined collections of relations. A closely related and detailed construction appears in work of Charlton, Matveiakin, Radchenko, and Rudenko \cite{CMRR24}, who construct a Lie coalgebra $\scr{L}^\rm{f}(F)$ of \emph{formal multiple polylogarithms}. This is generated by formal correlators $\rm{Cor}^\rm{f}(x_0,\compactldots,x_n)$ and there is a well-defined map of Lie coalgebras
\begin{align*} r^\rm{f} \colon \scr{G}(F) &\lra \scr{L}^\rm{f}(F) \\
\CorG(x_0,\compactldots,x_n) &\longmapsto \rm{Cor}^\rm{f}(x_0,\compactldots,x_n)\end{align*}
because the formal correlators satisfy the relations in \cref{thm:polyl-presentation-additive}, though possibly many more, and cobrackets are given by the same formulas. This is the \emph{formal realisation}.

Depending on the field $F$, there are realisation maps from $\scr{L}^f(F)$ to other Lie coalgebras. Firstly, for a number field $F$ there is a \emph{formal-to-motivic realisation} $\scr{L}^\rm{f}(F) \to \scr{L}^\rm{MTM}(F)$ with target the motivic Lie coalgebra as in \cref{sec:intro-mtm-number-field}. It is uniquely determined by sending the formal correlators $\CorF(x_0,\compactldots,x_n)$ to the \emph{motivic correlators} $\rm{Cor}^\rm{MTM}(x_0,\compactldots,x_n) \in \scr{L}^\rm{MTM}_n(F)$ \cite[\S 10.4.1]{Gon19}. Precomposing with formal realisation we obtain a \emph{motivic} realisation.

\begin{customthm}{D.a}\label{thm:motivic-realisation} For a number field $F$, there exists a unique functor of Tannakian categories
\[R^\rm{MTM} \colon \rm{Comod}^\rm{fd}_{\PolyL(F)}(\rm{GrMod}_\bb{Q}) \lra \rm{MTM}_\bb{Q}(F)\]
such that the induced morphism $r^\rm{MTM}$ of graded Lie coalgebras satisfies
\[r^\rm{MTM}(\CorG(x_0,\compactldots,x_n))=\Cor^\rm{MTM}(x_0,\compactldots,x_n).\]
\end{customthm}

Since multiple polylogarithms admit the same expression in terms of correlators in the Goncharov Lie coalgebra and formal Lie coalgebra, it follows that 
\[r^\rm{MTM}\left(\LiG_{n_1,\compactldots,n_k}(a_1,\compactldots,a_k)\right)=\rm{Li}^\rm{MTM}_{n_1,\compactldots,n_k}(a_1,\compactldots,a_k).\]

Secondly, given an embedding $\sigma \colon F \to \bb{C}$ there is a \emph{formal-to-Hodge realisation} $\scr{L}^{f}(F) \to \scr{L}^\rm{Hod}$. Consider the Tannakian category of rational mixed Hodge structures, constructed by Deligne \cite{Del71b} and let $\QMHTS$ be its Tannakian subcategory of mixed Hodge--Tate structures; this is equivalent to the category of finite-dimensional graded comodules over the Lie coalgebra $\scr{L}^\rm{Hod}$ of framed Hodge-Tate structures. The formal-to-Hodge realisation is uniquely determined by sending the formal correlators $\CorF(x_0,\compactldots,x_n)$ to the \emph{Hodge correlators} $\rm{Cor}^\rm{Hod}(\sigma(x_0),\compactldots,\sigma(x_n)) \in \scr{L}_n^\rm{Hod}$ \cite[\S 3.2]{Gon19}. Precomposing with formal realisation we obtain a Hodge realisation:

\begin{customthm}{D.b}\label{thm:hodge-realisation} For an embedding $\sigma\colon F\hookrightarrow \C$, there exists a unique functor of Tannakian categories
\[
R^{\rm{Hod}}_{\sigma}\colon \rm{Comod}^\rm{fd}_{\scr{G}(F)}(\rm{GrMod}_\bb{Q}) \lra \QMHTS
\]
such that the induced morphism $r^{\rm{Hod}}_{\sigma}$ of graded Lie coalgebras satisfies
\[
r^{\rm{Hod}}_{\sigma}(\CorG(x_0,\compactldots,x_n))=\Cor^\rm{Hod}(\sigma(x_0),\compactldots,\sigma(x_n)).
\]
\end{customthm}

\stepcounter{atheorem}

As for motivic realisation, it follows that
\[
r^{\rm{Hod}}_{\sigma}\left(\LiG_{n_1,\compactldots,n_k}(a_1,\compactldots,a_k)\right)=\LiH_{n_1,\compactldots,n_k}(\sigma(a_1),\compactldots,\sigma(a_k))
\]
with right-hand side the \emph{Hodge multiple polylogarithms} \cite[\S 6]{Gon01} \cite[\S 11.1.4]{Gon19}.

More concretely, the real period construction of Goncharov yields maps \cite[\S 1.11]{Gon19}
\[p_\bb{R} \colon \scr{L}^\rm{Hod}_n \lra \bb{R} \qquad \text{for $n \geq 1$},\]
that extract real numbers from Hodge correlators or Hodge multiple polylogarithms; in the latter case, it yields the values of generalisations of the Bloch--Wigner dilogarithm, a single-valued real-valued variant of the dilogarithm \cite[\S 11.1.3]{Gon19} \cite[Appendix]{Mal20}. One can use these real period maps to show, for example, that $\scr{G}_n(\bb{C})$ is uncountable for $n \geq 1$. In the sequel \cite{KRS2} we will use this to give a novel cocycle representing the Borel regulator class.

\subsection{The Rognes rank spectral sequence} Rognes constructed a filtration of the algebraic $K$-theory spectrum $K(F)$ by rank (that is, dimension) \cite{Rognes}. An equivalent spectral sequence was constructed by Galatius, Kupers, and Randal-Williams by filtering the group completion of $\BGLb(F)$ \cite[Section 13.8]{GKRW18}. Rationally it takes the form
\[E^1_{n,d} \cong H^{E_\infty}_{n,d}(\BGLb(F)_\bb{Q}) \Longrightarrow K_d(F)_\bb{Q}\]
with $d^r$-differential of bidegree $(-r,-1)$; we call this the \emph{Rognes rank spectral sequence}. The $E^1$-page looks as in \cref{fig:gltable1} but to use it for computations we use more information (see \cref{thm:rank-ss-omnibus}) obtained from our understanding of the Goncharov Lie coalgebra:
\begin{enumerate}[(i)]
    \item It is compatible with the splitting induced by a scaling action.
    \item It is compatible with the duality involution (which gives the Adams operation $\psi^{-1}$ on the target $K_*(F)_\Q$).
    \item The $d^1$-differential is determined by the part of the Lie cobracket on $E_\infty$-homology given by the ``$\sigma$-component''
    \[\delta_\sigma \colon H^{E_\infty}_{n,d}(\BGLb(F)_\bb{Q}) \lra H^{E_\infty}_{n-1,d-1}(\BGLb(F)_\bb{Q}) \otimes \bb{Q}\{\sigma\},\]
    where $\sigma$ denotes the generator of $H^{E_\infty}_{1,0}(\BGLb(F)_\bb{Q})$.
    \item For elements of the $E_\infty$-homology lying on the critical line $\PolyL(F)$, (part of) the $\sigma$-component 
    \[\delta_\sigma \colon \PolyL(F) \lra H_2(\GL_{n-1}(F),\StL_{n-1}(F)) \otimes \Q\{\sigma\}\]
    may be expressed in terms of the cobracket $\delta: \PolyL(F) \to \Lambda^2 \PolyL(F)$, see \cref{sec:sigma-component}, and hence can be computed using our presentation. 
\end{enumerate}
The resulting $E^1$- and $E^2$-pages are displayed in \cref{fig:e1page} and \cref{fig:e2page} (to explain the notation in those figures, in the displayed range there is an isomorphism $H^*(\scr{G}(F))_n \cong H^*(\Gamma_n(F))$ with the right-hand side the Goncharov's polylogarithmic complex of weight $n$). Before discussing this spectral sequence more generally, let us first deduce some concrete consequences. These are limited to low weights, as this spectral sequence is inconclusive in higher weights because we lack a good understanding of $H_*(\GL_n(F);\StL_n(F))$ for $* > 1$.

Combining the $E^2$-page and known bounds on the weights, one may prove \cref{theorem weight 3} as follows. Looking at the $(-1)$-eigenspaces of the duality involution on the $E^2$-page and discarding the Milnor $K$-theory contribution, the row $d=4$ yields an isomorphism
\[K_4^{(3)}(F)_{\bb{Q}} \cong H^2(\PolyL(F))_3,\]
which is precisely \cref{theorem weight 3} \eqref{enum:FormulaGoncharovConjecture2}. Similarly, the row $d=5$ yields an exact sequence
\[K^{(2)}_{4}(F)_{\Q}\otimes F^{\times}_{\Q}
\lra
K^{(3)}_{5}(F)_{\Q}
\lra 
H^{1}(\PolyL(F))_3\lra 0,\]
which is precisely \cref{theorem weight 3} \eqref{enum:FormulaGoncharovConjecture3}.
As we observed in \cref{sec:motives-general}, the weight $3$ part of the Chevalley--Eilenberg complex for $\scr{G}(F)$ agrees with the weight 3 polylogarithmic complex under the isomorphisms of \cref{thm:polyl-identification}, so one may rephrase these results in terms of the latter. 

Observe that if $F$ satisfies $\smash{K^{(2)}_4(F)_\Q}=0$ as predicted by the Beilinson--Soul\'e vanishing conjecture, it follows from \cref{theorem weight 3} \eqref{enum:FormulaGoncharovConjecture3} that there is an isomorphism $\smash{K^{(3)}_{5}(F)_{\Q} \cong H^{1}(\PolyL(F))_3}$ as predicted by \cref{conjecture main gamma}. Using localisation sequences, we prove in \cref{cor: weight 3 goncharov iff beilinson soule} that the converse holds as well: if $\smash{K^{(3)}_{5}(F)_{\Q} \cong
H^{1}(\PolyL(F))_3}$ holds for all fields $F$ then the instance $\smash{K^{(2)}_4(E)_\Q}=0$ of the Beilinson--Soul{\'e} vanishing conjecture holds for all fields $E$ too (in fact, this statement holds characteristic-wise).  We further prove that these statements are equivalent to certain conjectures regarding ``homotopy-invariance'' for $B_3(F)_\Q$. 

We return to rank spectral sequences in general. Their use goes back to Quillen's work on the algebraic $K$-theory groups of rings of integers \cite{QuillenFiniteGeneration}, but it was Rognes who recognised their potential for studying algebraic $K$-theory more generally \cite{Rognes,RognesMotivic,RognesWeight}. Rognes made several conjectures about his rank spectral sequence for a field. The first is that there is a relationship between the entry $E^1_{3,5}$ and trilogarithms, which is provided by our identification of $E^1_{3,5}$ as $\scr{G}_3(F) \cong B_3(F)_\bb{Q}$ in \cref{thm:polyl-identification}. The second is that there is a slope 2 vanishing line, which was proven by Galatius, Kupers, and Randal-Williams \cite{GKRW20} and is a crucial input to the results stated above. The third is that it collapses at the $E^2$-page---or may even coincide with the motivic spectral sequence---but in future work we will use our understanding of this spectral sequence to disprove this.

\subsection{Functional equations for multiple polylogarithms} \label{sec:functional equations} We next explain how our work sheds light on functional equations for multiple polylogarithms and their variants.

\subsubsection{The decomposition relation as universal functional equation}
\cref{thm:polyl-presentation-additive} suggests that the decomposition relation  
\begin{align*}
    &\CorG(x_0,\compactldots,x_n)-\CorG(y_0,\compactldots,y_n)\\
     &=\sum_{\iota=((i_1,j_1),\dots,(i_{n},j_n))\in T(n)}\rm{sign}(\iota)\,
     \CorG\left(0,\frac{x_{i_1}-x_{j_1}}{y_{i_1}-y_{j_1}},\compactldots,\frac{x_{i_n}-x_{j_n}}{y_{i_n}-y_{j_n}}\right),
\end{align*}
together with its degenerations is a natural candidate for the \emph{universal} functional equation for correlators, and hence for multiple polylogarithms. As discussed in an example above, the case $n=2$ yields the 5-term relation. For the convenience of the reader, we spell the decomposition relation in the case $n=3$: 

\begin{align*}
\CorG&\bigl(x_0,x_1,x_2,x_3\bigr)
-\CorG\bigl(y_0,y_1,y_2,y_3\bigr)\\
&=\CorG\Bigl(
0,
\frac{x_1-x_0}{y_1-y_0},
\frac{x_2-x_0}{y_2-y_0},
\frac{x_3-x_0}{y_3-y_0}
\Bigr)
-\CorG\Bigl(
0,
\frac{x_1-x_0}{y_1-y_0},
\frac{x_2-x_1}{y_2-y_1},
\frac{x_3-x_0}{y_3-y_0}
\Bigr)\\
&+\CorG\Bigl(
0,
\frac{x_2-x_0}{y_2-y_0},
\frac{x_2-x_1}{y_2-y_1},
\frac{x_3-x_0}{y_3-y_0}
\Bigr)
-\CorG\Bigl(
0,
\frac{x_1-x_0}{y_1-y_0},
\frac{x_3-x_0}{y_3-y_0},
\frac{x_2-x_1}{y_2-y_1}
\Bigr)\\
&
+\CorG\Bigl(
0,
\frac{x_1-x_0}{y_1-y_0},
\frac{x_3-x_1}{y_3-y_1},
\frac{x_2-x_1}{y_2-y_1}
\Bigr)
-\CorG\Bigl(
0,
\frac{x_3-x_0}{y_3-y_0},
\frac{x_3-x_1}{y_3-y_1},
\frac{x_2-x_1}{y_2-y_1}
\Bigr)\\
&
-\CorG\Bigl(
0,
\frac{x_1-x_0}{y_1-y_0},
\frac{x_2-x_0}{y_2-y_0},
\frac{x_3-x_2}{y_3-y_2}
\Bigr)
+\CorG\Bigl(
0,
\frac{x_1-x_0}{y_1-y_0},
\frac{x_2-x_1}{y_2-y_1},
\frac{x_3-x_2}{y_3-y_2}
\Bigr)\\
&-\CorG\Bigl(0,\frac{x_2-x_0}{y_2-y_0},
\frac{x_2-x_1}{y_2-y_1},
\frac{x_3-x_2}{y_3-y_2}
\Bigr)
+\CorG\Bigl(0,
\frac{x_2-x_0}{y_2-y_0},
\frac{x_3-x_0}{y_3-y_0},
\frac{x_2-x_1}{y_2-y_1}
\Bigr)\\
&-\CorG\Bigl(0,\frac{x_2-x_0}{y_2-y_0},\frac{x_3-x_2}{y_3-y_2},
\frac{x_2-x_1}{y_2-y_1}\Bigr)
+\CorG\Bigl(0,\frac{x_3-x_0}{y_3-y_0},
\frac{x_3-x_2}{y_3-y_2},
\frac{x_2-x_1}{y_2-y_1}\Bigr)\\
&+\CorG\Bigl(0,
\frac{x_1-x_0}{y_1-y_0},
\frac{x_3-x_0}{y_3-y_0},
\frac{x_3-x_2}{y_3-y_2}
\Bigr)
-\CorG\Bigl(0,
\frac{x_1-x_0}{y_1-y_0},
\frac{x_3-x_1}{y_3-y_1},
\frac{x_3-x_2}{y_3-y_2}
\Bigr)\\
&+\CorG\Bigl(0,
\frac{x_3-x_0}{y_3-y_0},
\frac{x_3-x_1}{y_3-y_1},
\frac{x_3-x_2}{y_3-y_2}
\Bigr).
\end{align*}

The decomposition relation is not entirely new: from the viewpoint of
iterated integrals, it is a natural instance of the classical change-of-variable formalism for iterated integrals. In some sense, it is the simplest general functional equation for multiple polylogarithms/iterated integrals/correlators one can write down. To explain the idea behind it, we work with iterated integrals as multivalued functions, defined by 
\[
I_\gamma(a_0;a_1,\compactldots,a_n;a_{n+1})
=
\int_{0\le t_1\le \cdots \le t_n\le 1}
\frac{\gamma'(t_1)\,dt_1}{\gamma(t_1)-a_1}
\cdots
\frac{\gamma'(t_n)\,dt_n}{\gamma(t_n)-a_n}
\]
for a path $\gamma \colon [0,1]\to \mathbb C\setminus\{a_1,\ldots,a_n\}$ with start point $\gamma(0)=a_0$, end point $\gamma(1)=a_{n+1}$. It is well-known, e.g.~\cite[Theorem 2.1]{Gon01}, that the iterated integral
\[
 I(f_0(t);f_1(t),\compactldots,f_n(t);f_{n+1}(t)) \qquad \text{for rational functions $f_i(t)\in \C(t)$}
\] 
can be rewritten, after choosing branches and possibly enlarging the set of singularities, as a linear combination of hyperlogarithms in the variable $t$, i.e.~of iterated integrals of the form
\[
I(a_0;a_1,\compactldots,a_n;t) \qquad \text{for constants $a_0, \ldots,a_{n}$.}
\] 
This idea has appeared in the literature before, for instance, in the work of Wojtkowiak \cite{Wojtkowiak}. The same heuristic applies to analytic, Hodge, motivic, or formal iterated integrals and correlators. In $\scr{G}(F)$, the resulting expansion for
\[
\CorG(x_0+y_0t,\compactldots,x_n+y_nt)
\]
is precisely the decomposition relation upon evaluation at $t=0,1$. For the formal, Hodge, and motivic realisations, this argument is made precise in \cref{prop: formal realization}.

\subsubsection{Functional equations for polylogarithms from relations in
infinite Steinberg modules}

The identification
\[
   \PolyL_n(F)\cong H_1(\GL_n(F);\StL_n(F))
\]
gives a new source of functional equations for polylogarithms: every
linear relation among Steinberg correlators gives, after evaluation by a linear functional, an identity among polylogarithmic correlators. More precisely, one has the following statement:

\begin{atheorem}\label{thm:polyl-relations-from-stl-relations}
Let $V$ be an $n$-dimensional vector space over $F$ and suppose that
\[
   \sum_i a_i\, \mathrm C
   \Bigl[u^{(i)}_0:\compactcdots:u^{(i)}_n\Bigr]=0
   \qquad\text{in } \StL(V),
\]
where $a_i\in \bb{Q}$, and where each
$u^{(i)}_0,\ldots,u^{(i)}_n$ is an affine basis of $V$. Then the element
\[
   \sum_i a_i\,
   \CorG\bigl(
      h\bigl(u^{(i)}_0\bigr),\compactcdots,h\bigl(u^{(i)}_n\bigr)
   \bigr)
   \in \PolyL_n(F)
\]
is independent of the choice of nonzero linear functional
$h\in V^\vee$. Equivalently, for any two nonzero linear functionals $h_1,h_2\in V^\vee$, one has the functional equation
\[
   \sum_i a_i\,
   \CorG\Bigl(
      h_1\bigl(u^{(i)}_0\bigr),\compactcdots,h_1\bigl(u^{(i)}_n\bigr)
   \Bigr)
   =
   \sum_i a_i\,
   \CorG\Bigl(
      h_2\bigl(u^{(i)}_0\bigr),\compactcdots,h_2\bigl(u^{(i)}_n\bigr)
   \Bigr)
   \qquad\text{in } \PolyL_n(F).
\]
\end{atheorem}
The decomposition relation is a special case of \cref{thm:polyl-relations-from-stl-relations}.

\subsection*{Acknowledgments} We would like to thank Steven Charlton for helping us to establish the results of \cref{section: trilog} with computer-assisted methods.  AK would like to thank S\o ren Galatius and Oscar Randal-Williams for many helpful discussions, Max Blans, Thomas Blom, and Gijs Heuts for answering questions about Koszul duality, Peter Scholze for some conversations regarding Goncharov's programme, and Elden Elmanto and Nick Rozenblyum for answering questions about algebraic $K$-theory and higher category theory. DR would like to thank Danylo Radchenko for numerous helpful discussions; in particular, regarding the results in \cref{section: duality computation}. DR would also like to thank  Alexander Beilinson, Cl{\'e}ment Dupont, S\o ren Galatius, Alexander Goncharov, Richard Hain, David Kazhdan, Jeremy Miller, and Peter Patzt for many helpful discussions. IS would like to thank Oscar Randal-Williams for many helpful conversations.

AK acknowledges the support of the Natural Sciences and Engineering Research Council of Canada (NSERC) [funding reference number 512156 and 512250]. DR was supported by NSF grant DMS-2502729.

\section{Steinberg modules} In this section we will define the Steinberg modules as well as the double and infinite variants that appeared in \cite{GKRW20,CharltonRadchenkoRudenko}. Our definitions will be in terms of Koszul duality, clarifying the algebraic structures present on these objects, followed by descriptions in terms of presentations. We next discuss the ``almost bases'' for $\StL(V)$ mentioned in the introduction and the corresponding decomposition operators. Finally, with an eye towards the computation of the cobracket on the Goncharov Lie coalgebra, we give (a) consequences of a duoidal Eckmann--Hilton argument for the coproduct on double Steinberg modules (\cref{lem:sth-coproduct-symmetry}) and the cobracket on infinite Steinberg modules (\cref{lem:stl-cobracket-vanishing}), and (b) a lift of the coproduct on Steinberg modules to the coLie cobar complex of $\SStL$ (\cref{sec:lift-coproduct}).

\begin{convention}\label{conv:shorter-notation} We fix a field $F$ and suppress it from the notation unless there is a risk of confusion, e.g.~write $\GL_n$ for $\GL_n(F)$. We work with rational coefficients and suppress this from the notation unless there is a risk of confusion, e.g.~write $H_*(-)$ for $H_*(-;\bb{Q})$.\end{convention}

\subsection{The Levi and parabolic tensor products} Before we can state the definitions of the Steinberg modules and their variants, we need to define several categories with additional structure.

\smallskip

Let $\Vect$ be the 1-category whose objects are finite-dimensional vector spaces over $F$ and whose morphisms are isomorphisms. It admits a symmetric monoidal structure by direct sum. It also admits a second ``flag sum'' promonoidal structure (a notion going back to \cite[Section 3]{Day}, where it was called a ``premonoidal'' structure) given as follows: the $k$-fold iterated tensor product for $k \geq 0$ is given by a profunctor (recall a profunctor $A \topro B$ is a functor $A \times B^\op \to \rm{Set}$)
\begin{align*} \obackslash_k \colon \Vect^k  &\lrapro \Vect \\
(V_1,\ldots,V_k,W) &\longmapsto \left\{\text{\parbox[c]{8cm}{\centering flags of subspaces $0 = W_0 \subseteq W_1 \subseteq \cdots \subseteq W_k = W$ with identifications $W_i/W_{i-1} \cong V_i$}}\right\},\end{align*}
where by convention $V_0 = 0$. Note for $k=2$ this is the same as the set of short exact sequences $0 \to V_1 \to W \to V_2 \to 0$ and for $k=0$ this assigns to $W \in \Vect$ the empty set unless $W \cong 0$ in which case it assigns a singleton.

For comparison, ``direct sum'' considered as a promonoidal structure has $k$-fold iterated tensor product for $k \geq 0$ given by a profunctor
\begin{align*} \oplus_k \colon \Vect^k &\lrapro \Vect \\
(V_1,\ldots,V_k,W) &\longmapsto \left\{\text{\parbox[c]{6cm}{\centering injections $V_i \to W$ so that the map $V_1 \oplus \cdots \oplus V_k \to W$ is an isomorphism}}\right\}.\end{align*}
Note that this has the same monoidal unit as ``flag sum'': for $k=0$ this assigns $W \in \Vect$ the empty set unless $W \cong 0$ in which case it assigns a singleton.

The identity lifts to a unital lax promonoidal functor from direct sum to flag sum. The natural transformations of this lax promonoidality are given by sending an ordered collection of summands $V_1,\ldots,V_k \subseteq W$ to the flag $0\subseteq \rm{im}(V_1) \subseteq \cdots \subset \rm{im}(V_1 \oplus \cdots \oplus V_{k-1}) \subseteq \rm{im}(V_1 \oplus \cdots \oplus V_k) = W$
with standard identifications. 

Moreover, we can assemble both to the structure of a normal $(E_\infty,E_1)$-produoidal category $(\Vect,\oplus,\obackslash)$ (see \cref{def:duoidal}); it has a promonoidal structure $\obackslash$ and a symmetric promonoidal structure $\oplus$ sharing the same units, related by a pronatural interchange natural transformation (see \eqref{eqn:duoidal-zeta}). The latter has components
\[(U_1 \obackslash U_2) \oplus (V_1 \obackslash V_2) \lra (U_1 \oplus V_1) \obackslash (U_2 \oplus V_2)\]
that are given, upon mapping into $W$, by the map sending a sum of short exact sequences $(0 \to U_1 \to W_1 \to U_2 \to 0) \oplus (0 \to V_1 \to W_2 \to V_2 \to 0)$ with $W_1 \oplus W_2 = W$ to the short exact sequence $0 \to (U_1 \oplus V_1) \to W \to (U_2 \oplus V_2) \to 0$. By taking $U_2 = 0 = V_1$ one recovers the aforementioned binary natural transformation of lax promonoidality on the identity functor: the identity is a lax functor of produoidal $(E_\infty,E_1)$-categories $(\Vect,\oplus,\oplus) \to (\Vect,\oplus,\obackslash)$ in the sense of \cite[Definition 2.10]{BataninMarkl}.

The category of functors from a symmetric (pro)monoidal category to a presentable symmetric monoidal category admits a Day convolution tensor product and this generalises to the produoidal setting (see \cref{sec:day}). Thus we can make the following definition:

\begin{definition}Let $\scr{C}$ be a presentable symmetric monoidal category.
\begin{enumerate}
    \item Direct sum induces through Day convolution a \emph{Levi tensor product}
    \[\levi \colon \Fun(\Vect,\scr{C}) \times \Fun(\Vect,\scr{C}) \lra \Fun(\Vect,\scr{C}).\]
    \item Flag sum induces through Day convolution a \emph{parabolic tensor product}
    \[\para \colon \Fun(\Vect,\scr{C}) \times \Fun(\Vect,\scr{C}) \lra \Fun(\Vect,\scr{C}).\]
\end{enumerate}
\end{definition}

These assemble to an $(E_\infty,E_1)$-duoidal category $(\Fun(\Vect,\scr{C}),\levi,\para)$ and the identity lifts to a lax monoidal functor $\id \colon (\Fun(\Vect,\scr{C}),\levi) \to (\Fun(\Vect,\scr{C}),\para)$ with associated natural transformation $\ol{\zeta} \colon X \levi Y \to X \para Y$: here we have used that in this case the map $\ol{\zeta}$ (see \eqref{eqn:duoidal-ol-ul-zeta}) obtained from the interchange transformation is part of a lax monoidality on the identity.

\begin{example}\label{exam:levi-para-explicit-formulas} Explicit formulas are given by
\begin{align*}(X \levi Y)(F^n) &\simeq \bigsqcup_{k=0}^n \rm{Ind}^{\GL_n(F)}_{\GL_k(F) \times \GL_{n-k}(F)} \left(X(F^k) \otimes_\scr{C} Y(F^{n-k})\right) \\
(X \para Y)(F^n) &\simeq \bigsqcup_{k=0}^n \rm{Ind}^{\GL_n(F)}_{\rm{P}_{n,k}(F)} \left(X(F^k) \otimes_\scr{C} Y(F^{n-k})\right)\end{align*}
where $\rm{P}_{n,k}(F) \subset \GL_n(F)$ is the parabolic subgroup 
that preserves the flag $F^k \subset F^n$ and $\GL_k(F) \times \GL_{n-k}(F) \subset \GL_n(F)$ is the Levi subgroup that preserves the splitting $F^k \oplus F^{n-k} = F^n$. The map $\ol{\zeta}$ is induced in terms of the above formulas by the inclusion $\GL_k(F) \times \GL_{n-k}(F) \subset \rm{P}_{n,k}(F)$.\end{example}

Recall that $E_\infty^\rm{u}$ denotes the unital $E_\infty$-operad (see \cref{def:einfty}), whose algebras are unital $E_\infty$-algebras. In the duoidal category $(\Fun(\Vect,\Spc),\levi,\para)$, one may use the interchange natural transformation $\zeta$ to endow the category of $E^\rm{u}_\infty$-algebras with respect to the Levi tensor product $\levi$ with a parabolic tensor product that we will also denote $\para$. Informally, for a presentable category $\scr{C}$ and $\bf{A},\bf{B} \in \Alg_{E^\rm{u}_\infty,\levi}(\Fun(\Vect),\scr{C})$, we have that $\bf{A} \para \bf{B}$ is again an $E^\rm{u}_\infty$-algebra whose product is given by
\[(\bf{A} \para \bf{B}) \levi (\bf{A} \para \bf{B}) \overset{\zeta}\lra (\bf{A} \levi \bf{A}) \para (\bf{B} \levi \bf{B}) \lra \bf{A} \para \bf{B}.\]
Precisely, there is a monoidal structure on $\Alg_{E^\rm{u}_\infty,\levi}(\Fun(\Vect,\scr{C}))$ so that the forgetful functor is symmetric monoidal \cite[Definition 3.7]{ToriiMult}. We can thus define the category of \emph{$(E_\infty,E_1)$-algebras} as $\Alg_{E_1^u,\para}(\Alg_{E_\infty^\rm{u},\levi}(\Fun(\Vect,\Spc)))$; this is a mild generalisation of \cite[Section 6.3]{ToriiDuoidal} in the setting of \cite{ToriiHigher} (see \cite[Section 3]{ToriiMult} for a summary).

\subsection{Steinberg modules and higher variants}\label{sec:steinberg-modules-defs}
We now define the Steinberg modules and their variants, through Koszul duality. For our purposes it is sufficient to specialise the previous discussion to rational chain complexes by taking $\scr{C} = \DQ$, the derived category of $\bb{Q}$ obtained by inverting the quasi-isomorphisms on the 1-category $\rm{Ch}_\bb{Q}$ of chain complexes over $\bb{Q}$.  Note that it contains the category $\rm{GrMod}_\bb{Q}$ of graded vector spaces over $\bb{Q}$ as a full subcategory. The category $\DQ$ has a symmetric monoidal structure given by tensor product, and using the Day convolution construction we obtain a $(E_\infty,E_1)$-duoidal category $\Fun(\Vect,\DQ)$.

The constant functor $\ul{\bb{Q}}$ admits the structure of an $(E^\rm{u}_\infty,E^\rm{u}_1)$-algebra in $\Fun(\Vect,\DQ)$. To see this, recall that $C_*(-;\bb{Q}) \colon \Spc \to \DQ$ is symmetric monoidal, and construct $\bb{Q}$ as the image under the induced functor
\begin{align*} C_*(-;\bb{Q})^\Alg \colon \Alg_{E^\rm{u}_1,\para}(\Alg_{E^\rm{u}_\infty,\levi}(\Fun(\Vect,\Spc))) &\lra  \Alg_{E^\rm{u}_1,\para}(\Alg_{E^\rm{u}_\infty,\levi}(\Fun(\Vect,\DQ))) \\
\bf{A} &\longmapsto C_*(\bf{A};\bb{Q})\end{align*}
of $\ul{\ast}$, which as a terminal object has a unique lift to an $(E^\rm{u}_\infty,E^\rm{u}_1)$-algebra. As $\ul{\ast}$ is given by the terminal object on the zero vector space, it admits a unique augmentation, inducing an augmentation on $\bb{Q}$.

We will first phrase Koszulity properties in terms of iterated bar constructions (see \cref{sec:iterated-bar-constructions}), though at the end of this subsection we will rephrase them in terms of indecomposables (see \cref{sec:bar-cobar-operads-algebras}). Our first Koszulity property is the following: as an augmented $E^\rm{u}_1$-algebra with respect to $\para$, $\ul{\bb{Q}}$ is Koszul in the sense that
\begin{equation}\label{hyp:ast-koszul}H_*(\Bar_{\para}(\ul{\bb{Q}})(V)) = 0 \qquad \text{unless $* = \dim(V)$.}\end{equation}
This is equivalent to the Solomon--Tits theorem: the bar construction in \eqref{hyp:ast-koszul} is \emph{isomorphic} to the reduced simplicial chains of the double simplicial suspension of the Tits building $T(F^n)$, see \cref{cor:koszul-st}. Recalling that $\ul{\bb{Q}}$ is an augmented $(E^\rm{u}_\infty,E^\rm{u}_1)$-algebra, by \cref{lem:bar-underlying} there is a preferred lift of its bar construction to
\[\Bar_{\para}(\ul{\bb{Q}}) \in \coAlg^\aug_{E_1^\rm{u},\para}(\Alg^\aug_{E_\infty^\rm{u},\levi}(\Fun(\Vect,\DQ))).\]
For now, we forget the $E^\rm{u}_1$-coalgebra structure, retaining only the $E^\rm{u}_\infty$-algebra structure. Our second Koszulity property is the following: as an augmented $E_1^\rm{u}$-algebra with respect to $\levi$, $\Bar_{\para}(\ul{\bb{Q}})$ is Koszul as an $E_1^\rm{u}$-algebra:
\begin{equation}\label{hyp:st-koszul} H_*(\Bar_\levi(\Bar_{\para}(\ul{\bb{Q}}))(V)) = 0 \text{ unless $* = 2\,\dim(V)$}.\end{equation}
This was proven in \cite{MNP}, \cite{MPW23}, and \cite{CharltonRadchenkoRudenko}, see \cref{cor:koszul-st} for the middle proof though we will use the latter in \cite{KRS2}: by the proof of \cite[Theorem 18]{CharltonRadchenkoRudenko} there is an isomorphism of chain complexes $B^\rm{As}(\SSt) \cong C_*(T(F^n);\bb{Q}) \otimes \St_n(F)$ up to a grading shift, with left side as in \cref{def:bar-as}, and the result follows from the Solomon--Tits theorem.

Over $\bb{Q}$, the commutative bar complex $B^\rm{Com}(\SSt)$ is a summand of the associative bar complex $B^\rm{As}(\SSt)$ by \cref{prop:barr-splitting} and it follows that $\Bar_{\para}(\ul{\bb{Q}})$ is also Koszul as an $E_\infty^\rm{u}$-algebra: 
\[\colim_{k \to \infty} H_*(\Sigma^{-k+1} \Bar^k_\levi(\SSt)(V)) = 0 \text{ unless $* = 2\dim(V)$}.\] 
In fact, the left side is given by the indecomposables of $H_*(\Bar_\levi(\Bar_{\para}(\ul{\bb{Q}}))(V)) $ with respect to the remaining product structure.

\begin{definition}The \emph{Steinberg modules} are defined as
\begin{align*}\SSt &\coloneq H_*(\Bar_{\para}(\ul{\bb{Q}})) \in \Fun(\Vect,\rm{GrMod}_\bb{Q}).\\
\intertext{The \emph{double Steinberg modules} are defined as}
        \SStH &\coloneq H_*(\Bar_\levi(\SSt)) \in \Fun(\Vect,\rm{GrMod}_\bb{Q}).\\
\intertext{The \emph{infinite Steinberg modules} are defined as}
        \SStL &\coloneq \colim_{k \to \infty} H_{*-k+1}(\Bar^k_\levi(\SSt))\in \Fun(\Vect,\rm{GrMod}_\bb{Q}).\end{align*}
\end{definition}

\begin{remark}It would have been equivalent to define $\SSt$, $\SStH$, and $\SStL$ as objects in the $\Fun(\Vect,\DQ)$ instead, as lying in $\Fun(\Vect,\rm{GrMod}_\bb{Q})$ is a property. We have opted for the latter to stress that we are in a case where the Koszulity hypotheses \eqref{hyp:ast-koszul} and \eqref{hyp:st-koszul} hold.\end{remark}

\begin{notation}For $V \in \Vect$ we define  $\GL(V)$-representations over $\bb{Q}$
\begin{align*}\St(V) &\coloneq H_{\dim(V)}(\Bar_{\para}(\ul{\ds{Q}})(V))\\
\StH(V) &\coloneq H_{2\dim(V)}(\Bar_\levi(\SSt)(V)) \\
\StL(V) &\coloneq \colim_{k \to \infty} H_{2\dim(V)-k+1}(\Bar^k_\levi(\SSt)(V)).
\end{align*}
Note that these are \emph{not} graded, and if we do want to consider them as graded, $\St(V)$ must be placed in grading $\dim(V)$, and $\StH(V)$ and $\StL(V)$ in grading $2\,\dim(V)$. We will also use the abbreviation
\[\St_n(F) \coloneq \St(F^n), \quad \StH_n(F) \coloneq \StH(F^n), \quad \text{and} \quad \StL_n(F) \coloneq \StL(F^n).\]
\end{notation}

\begin{remark}Let us spell out the case $V = 0$. There are preferred identifications $\St(0) \cong \bb{Q}$ and $\StH(0) \cong \bb{Q}$, and $\StL(0) = 0$. Occasionally we will need to pass to the subobjects 
\[\SSt_{>0} \subset \SSt \quad \text{and} \quad \SSt^2_{>0}  \subset \SStH\] 
which differ only in that their value on $V=0$ is now equal to $0$; these are the augmentation ideals for augmentations to be discussed later.
\end{remark}

Using \cref{thm:indec-is-bar} these can also be described in terms of indecomposables
\begin{align*}\SSt &\simeq (\Sigma\, \indec_{E_1^\rm{nu}}(\ul{\bb{Q}}_{>0}))^+ \\
\SStH &\simeq (\Sigma\, \indec_{E_1^\rm{nu}}(\SSt_{>0}))^+ \\
\SStL &\simeq (\Sigma\, \indec_{E_\infty^\rm{nu}}(\SSt_{>0})),\end{align*}
where $(-)^+$ denotes unitalisation (see \cref{sec:augmentation-ideals}). To see this, recall that a $t$-structure allows one to talk about the connectivity of objects and that functor categories inherit $t$-structures from the standard $t$-structure on $\DQ$, once we fix a function specifying how connected the value of a functor should be at each object of $\Vect$; in \cite{GKRW18} such a function was called an \emph{abstract connectivity}. In our case it will be convenient to consider the $t$-structure specified by the function $c_0 \colon \rm{ob}(\Vect) \to [-\infty,\infty]$ given by $c_0(0) = 0$ and $c_0(V) = -\infty$ for $V \neq 0$. This is left-compatible with respect to either (symmetric) monoidal structure, i.e.~a Levi or parabolic tensor product of $c_0$-connective functors is $c_0$-connective; consequently, $c_0$-connected objects have a well-behaved theory of Koszul duality. This allows us to apply Koszul duality results \cref{thm:kd-connected} and \cref{thm:indec-is-bar} to identify, up to suspensions, bar constructions with the unitalisations of indecomposables, in case of the $c_0$-connected nonunital algebras $\ul{\bb{Q}}_{>0}$ and $\SSt_{>0}$.

\subsection{Steinberg modules} We now discuss in more detail the Steinberg modules, namely algebraic structures present on them and explicit presentations.

\subsubsection{Algebraic structures on Steinberg modules} By definition we have 
\[\SSt = H_*(\Bar_{\para}(\ul{\ds{Q}})) \qquad \text{with} \qquad \Bar_{\para}(\ul{\ds{Q}}) \in \coAlg_{E_1^\rm{u}}(\Alg_{E_\infty^\rm{u}}(\Fun(\Vect,\DQ))).\]
Noting that because we work with rational coefficients there are no $\rm{Tor}$-terms when we apply the K\"unneth theorem, we get from a mild generalization of \cref{lem:bialgebra-structure} to the duoidal setting:

\begin{proposition}The Steinberg modules $\SSt$ come equipped with the structure of a commutative bialgebra in the normal $(E_\infty,E_1)$-duoidal category $(\Fun(\Vect,\rm{GrMod}_\ds{Q}),\levi,\para)$.\end{proposition}

Explicitly, this means it has a commutative (in the graded sense, i.e.~using the symmetry in $\rm{GrMod}_\bb{Q}$ with Koszul sign) unital associative product and counital coassociative coproduct
\[\mu \colon \SSt \levi \SSt \lra \SSt \qquad \text{and} \qquad \Delta \colon \SSt \lra \SSt \para \SSt\]
so that $\mu$ is a map of coalgebras. This means that the following diagram commutes
    \[\begin{tikzcd} \SSt \levi \SSt \rar{\mu} \dar[swap]{\Delta \levi \Delta} & \SSt \arrow{dd}{\Delta} \\[-3pt] 
    (\SSt \para \SSt) \levi (\SSt \para \SSt) \dar[swap]{\zeta} &  \\[-3pt]
    (\SSt \levi \SSt) \para (\SSt \levi \SSt) \rar{\mu \para \mu} & \SSt \para \SSt.\end{tikzcd}\]
It is worth pointing out a different perspective on this. Recall that we have an equivalence 
\[\Bar_{\para}(\ul{\ds{Q}}) \simeq \big(\Sigma\,\indec_{E_1^\rm{nu}}(\ul{\ds{Q}}_{>0})\big)^+ \quad \text{with} \quad \indec_{E_1^\rm{nu}}(\ul{\ds{Q}}) \in \coAlg^{\nil}_{BE_1^\rm{nu}}(\Fun(\Vect,\DQ)))\] 
and the latter admits, after suspension, a conilpotent $E^\rm{nu}_1$-coalgebra structure as $BE_1^\rm{nu}$ is equivalent to the operadic suspension of the linear dual cooperad $DE_1^\rm{nu}$.

\subsubsection{Presentations and formulas for Steinberg modules} We now give a well-known presentation of the Steinberg modules, see e.g.~\cite[Section 1]{KahnSun} \cite[Section 3.1]{CharltonRadchenkoRudenko}. Recall that an ordered basis $v_1,\ldots,v_n$ of a vector space $V$ of dimension $n$ gives an \emph{apartment class}
\[[v_1,\ldots,v_n] \in \St(V),\]
though it is more convenient to allow \emph{any} ordered collection $v_1,\ldots,v_n$ of vectors in $V$ and then set $[v_1,\ldots,v_n]$ to be zero if $v_1,\ldots,v_n$ do not span. The elements $[v_1,\ldots,v_n] \in \St(V)$ satisfy the following relations
\begin{enumerate}[(1)] \setcounter{enumi}{-1}
    \item $[v_1,\compactldots,v_n] = 0$ if $v_1,\ldots,v_n$ are linearly dependent,
    \item $[v_{\sigma(1)},\compactldots,v_{\sigma(n)}] = (-1)^\sigma[v_1,\compactldots,v_n]$ for $\sigma \in \fr{S}_n$,
    \item $[\lambda v_1,\compactldots,v_n] = [v_1,\compactldots,v_n]$ for $\lambda \in F^\times$,
    \item $\sum_{i=0}^n (-1)^i [v_0,\compactldots,\widehat{v}_i,\compactldots,v_n] = 0$ for an ordered collection of vectors $v_0,\ldots,v_n$ in $V$.
\end{enumerate}
These give a presentation of the Steinberg modules:

\begin{proposition}\label{prop:st-explicit-pres} The following map of $\bb{Q}[\GL(V)]$-modules is an isomorphism
\[\frac{\bb{Q}[[v_1,\compactldots,v_n] \text{ for ordered collections $v_1,\ldots,v_n$}]}{\text{(0)--(3)}}\overset{\cong}\lra \St(V).\]
\end{proposition} 

We can now ask for a description of the product (with respect to $\levi$) and coproduct (with respect to $\para$) in terms of this presentation: this was done in \cite{AMP,CharltonRadchenkoRudenko}. To do so, it is convenient to include $V$ in the notation of the apartment and write $[V|v_1,\ldots,v_n] \coloneq [v_1,\ldots,v_n] \in \St(V)$.

\begin{proposition}\label{prop:st-explicit-prod-coprod} With respect to the presentation of \cref{prop:st-explicit-pres}, the product and coproduct on $\SSt$ are given by
\begin{equation}\label{eqn:st-product}\begin{aligned} \mu \colon \St(V) \otimes \St(W) &\lra \St(V \oplus W) \\
[V|v_1,\compactldots,v_n] \otimes [W|w_1,\compactldots,w_{n'}] &\longmapsto [V\oplus W|v_1,\compactldots,v_n,w_1,\compactldots,w_{n'}]\end{aligned}\end{equation}
and \begin{equation}\label{eqn:st-coproduct}\begin{aligned} \Delta \colon \St(V) &\lra \bigoplus_{W \subseteq V} \St(W) \otimes \St(V/W) \\
[V|v_1,\compactldots,v_n] &\longmapsto \sum_{I \subseteq \ul{n}} (-1)^{\sigma_I} [V_I|v_{i_1},\compactldots,v_{i_k}] \otimes [V/V_I|\overline{v}_{j_1},\compactldots,\overline{v}_{j_{n-k}}]\end{aligned}\end{equation}
\end{proposition}

Here in the formula \eqref{eqn:st-product} for the product, the right vectors in $V$ or $W$ are to be considered as vectors in $V \oplus W$ through the inclusion. It follows from relation (1) that \eqref{eqn:st-product} is graded-commutative, as we put $\St(V)$ in grading $\dim(V)$. The formula \eqref{eqn:st-coproduct} for the coproduct should be interpreted as follows: there are only nonzero terms for $W$ of the form $V_I \coloneq \rm{span}(v_{i_1},\ldots,v_{i_k})$ for a subset $I = \{i_1,\ldots,i_k\}$ (with induced order) and $\overline{v}_j$ for $j \in I^c$ in its complement (with induced order) denote the projection of $v_j$ to $V/V_I$. Finally, $\sigma_I$ is the shuffle permutation that shuffles $I$ to the front of $I^c$ preserving the induced ordering. 

\begin{proof}[Proof of \cref{prop:st-explicit-prod-coprod}] We will give two proofs, as we believe this is clarifying. In the first proof, for the product we cite \cite[Lemma 6.8]{GKRW20}. For the coproduct, we use that the induced coproduct and the formula \eqref{eqn:st-coproduct} both have the property that they make $\St$ into a bialgebra in the normal duoidal category $(\Fun(\Vect,\rm{GrMod}_\bb{Q}),\levi,\para)$. As $\SSt$ is generated under products by $\St(F) \cong \bb{Q}$ placed in degree 1, and the product is a homomorphism for the coproduct, it suffices to verify that the coproducts agree on $\St(F)$, which follows from counitality.

In the second proof, we pass to the 1-category of rational chain complexes and use \cref{sec:ass-algebras} to obtain an identification
\[\SSt \cong H_*(\Bar_{\para}(\ul{\bb{Q}})) \cong H_*(B^\rm{As}(\ul{\bb{Q}}))\]
as a commutative bialgebra in $\Fun(\Vect,\rm{GrMod}_\bb{Q})$. Unwinding the construction, we recognise that the underlying object of $B^\rm{As}(\ul{\bb{Q}})$ is given by $V \mapsto \smash{\widetilde{C}_*(D^1(V))}$. By the proof of \cite[Proposition 25]{CharltonRadchenkoRudenko} the shuffle product on this yields \eqref{eqn:st-product} in terms of apartments, and similarly by the discussion following \cite[Proposition 21]{CharltonRadchenkoRudenko} the deconcatenation coproduct on this yields \eqref{eqn:st-coproduct} in terms of apartments. 
\end{proof}

\begin{remark} Koszul duality yields several useful resolutions. Firstly, the equivalence $\ul{\bb{Q}} \simeq \Omega^\rm{As}_{\para}(\SSt)$, with right side the cobar construction of \cref{def:cobar-as}, yields upon evaluation at $V$ the following resolution of $\bb{Q}$ by Steinberg modules: there is a chain complex
\[\St(V) \to \bigoplus_{0 \subsetneq V_1 \subsetneq V} \St(V_1) \otimes \St(V/V_1) \to \cdots \to \bigoplus_{0 \subsetneq V_1 \subsetneq \cdots \subseteq V_{n-1} \subsetneq V} \St(V_1) \otimes \cdots \otimes \St(V/V_{n-1})\]
with differentials induced by the coproduct, and augmentation to $\bb{Q}$ that is a quasi-isomorphism. (The dual equivalence $\SSt \simeq B^\rm{As}_{\para}(\bb{Q})$ is merely the definition of the Steinberg module.) Secondly, one can use the acyclicity of the twisted tensor products as in \cite[Proposition 2.2.13]{LodayVallette}.
\end{remark}

\subsection{Double Steinberg modules} We now do the same for the double Steinberg modules.

\subsubsection{Algebraic structure on double Steinberg modules} By definition we have
\[\SStH = H_*(\Bar_\levi(\SSt)) \qquad \text{and} \qquad \Bar_{\levi}(\SSt) \in \Alg^\rm{aug}_{E_\infty^\rm{u}}(\Fun(\Vect,\DQ)).\]

\begin{proposition}The double Steinberg modules $\SStH$ comes equipped with the structure of a commutative bialgebra in the symmetric monoidal category $(\Fun(\Vect,\rm{GrMod}_\ds{Q}),\levi)$. 
\end{proposition}

Explicitly, this means it has a graded-commutative unital associative product and counital coassociative coproduct
\[\mu^{(2)} \colon \SStH \levi \SStH \lra \SStH \qquad \text{and} \qquad \Delta^{(2)} \colon \SStH \lra \SStH \levi \SStH\]
so that $\mu^{(2)}$ is a map of coalgebras. This means that the following diagram commutes
     \[\begin{tikzcd} \SStH \levi \SStH \rar{\mu^{(2)}} \dar{\Delta^{(2)} \levi \Delta^{(2)}} &[10pt] \SStH \arrow{dd}{\Delta^{(2)}} \\[-3pt] 
    (\SStH \levi \SStH) \levi (\SStH \levi \SStH) \dar{\id \levi \sigma \levi \id} &  \\[-3pt]
    (\SStH \levi \SStH) \levi (\SStH \levi \SStH) \rar{\mu^{(2)} \levi \mu^{(2)}} & \SStH \levi \SStH.\end{tikzcd}\]
As in the previous case, the coalgebra structure also arises through the equivalence
\[\Bar_\levi(\SSt) \simeq \big(\Sigma\,\indec_{E_1^\rm{nu}}(\SSt_{>0})\big)^+\]
where $\indec_{E_1^\rm{nu}}(\SSt_{>0})$ admits, after a suspension, a conilpotent $E^\rm{nu}_1$-coalgebra structure.

\medskip

In fact, $\SStH$ admits what at first seems like additional structure, though we will see momentarily this is determined by the previous data: a coproduct with respect to $\para$, which is compatible with $\Delta^{(2)}$. As $\SSt$ is an augmented bialgebra, informally its coproduct $\Delta \colon \SSt \to \SSt \para \SSt$ is a map of augmented commutative algebras and hence induces the left map in map
\[\Bar_\levi(\SSt) \lra \Bar_\levi(\SSt \para \SSt) \lra \Bar_\levi(\SSt) \para \Bar_\levi(\SSt)\]
while the right map is induced by $\zeta$ using that $\para$ preserves sifted colimits in each entry.

To make this precise, we use a generalisation of \cref{lem:bar-underlying} to the duoidal setting: the category $\coAlg^\rm{aug}_{E_1^\rm{u},\para}(\Fun(\Vect,\DQ))$ admits a symmetric monoidal structure so that the forgetful functor $\fgt_{E_1} \colon \coAlg_{E_1^\rm{u}}(\Fun(\Vect,\DQ)) \to \Fun(\Vect,\DQ)$ is symmetric monoidal \cite[Proposition 3.10]{ToriiMult}. As $\fgt_{E_1}$ admits a right adjoint given by $\cofree_{E_1}$, it preserves colimits and $\Bar_\levi(\SSt)$ has a preferred lift to an object
\[\Bar_\levi(\SSt) \in \coAlg^\rm{aug}_{E_1^\rm{u},\levi}(\Alg^\rm{aug}_{E_\infty^\rm{u}}(\coAlg^\rm{aug}_{E_1^\rm{u},\para}(\Fun(\Vect,\DQ)))).\]
Forgetting the algebra structure and passing to homology, we get on $\SStH$ a counital coassociative coproduct, which we will also denote $\Delta$, compatible with $\Delta^{(2)}$ in the sense that the following diagram commutes:
\[\begin{tikzcd} \SStH \rar{\Delta^{(2)}} \arrow{dd}{\Delta} &[20pt] \SStH \levi \SStH \dar{\Delta \levi \Delta} \\[-3pt]
& (\SStH \para \SStH) \levi (\SStH \para \SStH) \dar{\zeta} \\[-3pt]
\SStH \para \SStH \rar{\Delta^{(2)} \para \Delta^{(2)}} & (\SStH \levi \SStH) \para (\SStH \levi \SStH).\end{tikzcd}\]
This provides the input for the duoidal Eckmann--Hilton argument of \cref{cor:duoidal-eh-sym} and we conclude that this additional structure is determined by $\Delta^{(2)}$ as $\overline{\zeta} \circ \Delta^{(2)} = \Delta = \overline{\zeta} \circ \sigma \circ \Delta^{(2)}$. Let us record the consequence for $\Delta^{(2)}$:

\begin{lemma}\label{lem:sth-coproduct-symmetry} The coproduct \eqref{eqn:sth-coproduct} on $\SStH$ has the following symmetry property: 
\[\overline{\zeta} \circ \Delta^{(2)} = \overline{\zeta} \circ \sigma \circ \Delta^{(2)}.\]\end{lemma}

\subsubsection{Presentations and formulas for double Steinberg modules} 
We will next give a presentation of $\SStH$, and describe the product and coproduct in terms of this presentation. The outlined proof of Koszulity for Steinberg modules induces an isomorphism
\[\SStH \overset{\cong}\lra \SSt \odot \SSt\]
where $\odot$ denotes the pointwise tensor product (also known as the Hadamard tensor product). We hence obtain from the presentation for $\SSt$ in \cref{prop:st-explicit-pres} a presentation for $\SStH$, simply by taking a tensor product. For example, for $n=3$, it is generated by pairs of triangles in the projective space $\mathbb{P}^2(F)$:
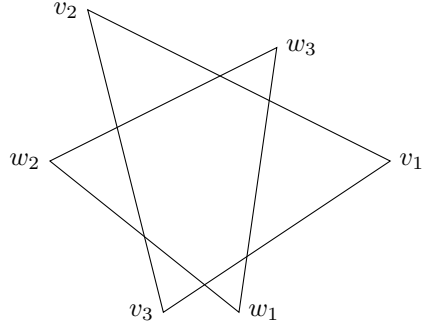
\[\begin{tikzpicture}
\draw (0,0) -- (-1,4) -- (3,2) -- cycle;
\node at (0,0) [left] {$v_3$};
\node at (-1,4) [left] {$v_2$};
\node at (3,2) [right] {$v_1$};
\node at (1,0) [right] {$w_1$};
\node at (-1.5,2) [left] {$w_2$};
\node at (1.5,3.5) [right] {$w_3$};
\draw (1,0) -- (-1.5,2) -- (1.5,3.5) --cycle;
\end{tikzpicture}\]
With respect to this, the explicit formulas for the product and coproduct were determined in \cite[Section 3.3]{CharltonRadchenkoRudenko} in terms of those in \cref{prop:st-explicit-prod-coprod}:

\begin{proposition}\label{prop:sth-explicit-prod-coprod} With respect to the presentation induced by \cref{prop:st-explicit-pres}, the product and coproduct on $\SStH$ are given by
\begin{equation}\label{eqn:sth-product}\begin{aligned} \mu^{(2)} \colon \StH(V) \otimes \StH(W) &\lra \StH(V \oplus W) \\
(a \otimes a') \otimes (b \otimes b') &\longmapsto \mu(a\otimes b) \otimes \mu(a' \otimes b') \end{aligned}\end{equation}
and
 \begin{equation}\label{eqn:sth-coproduct}\begin{aligned} \Delta^{(2)} \colon \StH(V) &\lra \bigoplus_{V = V_1 \oplus V_2} \StH(V_1) \otimes \StH(V_2) \\
a \otimes a' &\longmapsto \Psi(\Delta(a) \otimes \tau\Delta(a')).\end{aligned}\end{equation}
\end{proposition}

\begin{remark} Observe that \eqref{eqn:sth-product} from \cite[Proposition 20]{CharltonRadchenkoRudenko} involves no Koszul sign, in contrast with the proof of \cite[Theorem 6.9]{GKRW20}: this is due to how exactly $\SStH$ is identified with $\SSt \odot \SSt$, and these references differ by multiplication with $(-1)^{n(n-1)/2}$ in degree $2n$.\end{remark}

Here, first, $\Psi$ is the transformation 
\[(A \para B) \odot (C \para D) \lra (A \odot C) \levi (B \odot D)\]
given on subspaces $V_1,V_2 \subseteq V$ by the zero map unless $V = V_1 \oplus V_2$, where it is given by
\begin{align*}A(V_1) \otimes B(V/V_1) \otimes C(V_2) \otimes D(V/V_2) &\lra A(V_1) \otimes D(V_1) \otimes B(V_2) \otimes C(V_2) \\
a \otimes b \otimes c \otimes d &\longmapsto a \otimes \pi^{-1}_2(d) \otimes \pi^{-1}_1(c) \otimes d,
\end{align*}
where $\pi_1 \colon V_1 \to V \to V/V_2$ and $\pi_2 \colon V_2 \to V \to V/V_1$ are the induced isomorphisms. Second, $\tau \colon \SSt \levi \SSt \to \SSt \levi \SSt$ is a sign, given by $[V|a] \otimes [W|b] \mapsto (-1)^{\dim(V)\dim(W)} [V|a] \otimes [W|b]$.

\begin{remark}This may be part of a ``trioidal'' structure on $\Fun(\Vect,\rm{GrVect}_\bb{Q})$ but we will not pursue this further in this paper. \end{remark}

\begin{proof}[Proof of \cref{prop:sth-explicit-prod-coprod}] For the coproduct, we observe that we can compute the bar construction as a \emph{coalgebra} using the bar complex $B^\rm{As}(\SSt)$ (see \cref{def:bar-as}), which means that the injection of chain complexes known as the associative symbol map (see \cref{sec:symbol-maps})
\[s^\rm{As} \colon \SStH \lra  B^\rm{As}(\SSt)\]
is one of coalgebras, where the target has the deconcatenation product. This criterion uniquely determines the formula in \cite[Proposition 21]{CharltonRadchenkoRudenko}.

We will give three proofs for the product. For the first proof, we cite \cite[Section 6.3]{GKRW20}. For the second proof, we use that the product on $B^\rm{As}(\SSt)$ is given by the shuffle product and this criterion determines the product of \cite[Proposition 20]{CharltonRadchenkoRudenko}. For the third proof, to determine the product we use that the induced product and the formula both have the property that, together with the above coproduct, they make $\SStH$ into a bialgebra. As the $n$-fold iterated reduced coproduct $\StH(F^n) \to \StH(F)^{\levi n}$ is injective and the product is a homomorphism for the coproduct, it suffices to verify that the products agree on $\StH(0) \otimes \StH(F)$, which follows from unitarity. 
\end{proof}

 Note that \eqref{eqn:sth-coproduct} does not obviously have the symmetry property of \cref{lem:sth-coproduct-symmetry}, and this in fact involves a surprising amount of cancellation. We illustrate that with two examples. 

\begin{example}\label{ex:symmetry-para-decon_1}Consider a generic pair of simplexes, i.e., an element
\[x= [v_1, \compactldots, v_n] \otimes [w_1, \compactldots, w_n] \in \StH(V) \cong \St(V) \otimes \St(V)\]
such that subspaces $V_I = \rm{span}(v_i \mid i \in I)$ and $W_J = \rm{span}(w_j \mid j \in J)$ for $I, J \subseteq \ul{n}$ intersect transversally. In this case, $\overline{\zeta} \circ \Delta^{(2)}(x)$ is supported in summands $\StH(V_I) \otimes \StH(V/V_I)$ of $(\SStH \para \SStH)(V)$ while $\overline{\zeta} \circ \sigma \circ \Delta^{(2)}(x)$ is supported in summands $\StH(W_J) \otimes \StH(V/W_J)$, and these summands cannot agree by the genericity hypothesis. Nonetheless \cref{lem:sth-coproduct-symmetry} is true: both terms $\overline{\zeta} \circ \Delta^{(2)}(x)$ and $\overline{\zeta} \circ \sigma \circ \Delta^{(2)}(x)$ vanish identically simply because $\Delta([v_1,\compactldots,v_n])$ will be supported in direct summands of the form $\St(V_I) \otimes \St(V_{I^c})$ while $\Delta([w_1,\compactldots,w_n])$ will be supported in direct summands of the form $\St(W_J) \otimes \St(W_{J^c})$, and these cannot agree by the genericity hypothesis so $\Psi$ takes $\Delta([v_1,\compactldots,v_n]) \otimes \tau \Delta([w_1,\compactldots,w_n])$ to zero.
\end{example}

\begin{example}\label{ex:symmetry-para-decon_2}
Consider next a pair of simplexes
\[x= [v_1, v_2, v_3] \otimes [w_1, w_2, w_3] \in \StH(V)\]
which is not generic because they have a common vertex $v_1=w_1$ (but otherwise are generic). The projection of the term $\overline{\zeta} \circ \Delta^{(2)}(x)$ to the summand $\StH(V_1)\otimes \StH(V/V_1) \subseteq (\SStH \para \SStH)(V)$ is equal to 
\begin{equation}\label{eqn: component v1 part 1}
\Bigl([v_1]\otimes [v_1]\Bigr) \otimes \Bigl([\overline{V_{12}\cap W_{23}},  \overline{V_{13}\cap W_{23}}]\otimes [\overline{w_2},\overline{w_3}] \Bigr).
\end{equation}
The projection of the term $\overline{\zeta} \circ \sigma \circ \Delta^{(2)}(x)$ to the same summand is equal to 
\begin{equation}\label{eqn: component v1 part 2}
\Bigl([w_1]\otimes [w_1]\Bigr) \otimes \Bigl([\overline{v_2}, \overline{v_3}] \otimes [ \overline{W_{12}\cap V_{23}}, \overline{W_{13}\cap V_{23}}] \Bigr).
\end{equation}
Terms $\overline{\zeta} \circ \Delta^{(2)}(x)$ and $\overline{\zeta} \circ \sigma \circ \Delta^{(2)}(x)$ give the same contribution to the aforementioned component as \eqref{eqn: component v1 part 1} and \eqref{eqn: component v1 part 2} coincide.
\end{example}

\subsection{Infinite Steinberg modules} Finally, we do the same for infinite Steinberg modules.

\subsubsection{Algebraic structure on infinite Steinberg modules}\label{sec:alg-str-sstl} As the nonunital commutative algebra $\St_{>0} \in \Fun(\Vect,\rm{GrMod}_{\ds{Q}})$ is Koszul as a nonunital associative algebra with Koszul dual given by $\SSt^2_{>0}$, it is also Koszul as a nonunital commutative algebra with Koszul dual given by the quotient $\SStL$ of $\SSt^2_{>0}$ given by the indecomposables (cf.~the notation of \cite[Section 3.1]{CharltonRadchenkoRudenko}) with cobracket given as follows
\[\begin{tikzcd} \SSt^2_{>0} \rar{\Delta^{(2)}-\sigma \circ \Delta^{(2)}} \dar[two heads] &[30pt] \SSt^2_{>0} \levi \SSt^2_{>0} \dar[two heads] \\[-5pt]
\SStL \rar[dashed]{\delta} & \SStL \levi \SStL,\end{tikzcd}\]
where by construction $\delta$ takes values in the summand $\Lambda^2 \SStL$, including $a \wedge b$ as $\frac{1}{2}(a \otimes b - b \otimes a)$. The latter is by definition $\bb{Q}_\rm{sign} \otimes_{\fr{S}_2} (\SStL \levi \SStL)$ and we may suppress $\levi$ from the notation since it is not possible to perform this construction with respect to the monoidal structure $\para$ as it lacks a symmetry.

\begin{proposition}The infinite Steinberg modules $\SStL$ come equipped with the structure of a graded Lie coalgebra in the symmetric monoidal category $(\Fun(\Vect,\rm{GrMod}_\ds{Q}),\levi)$.
\end{proposition}

In this case the perspective using indecomposables is arguably better: iterating the bar construction on a $c_0$-connected augmented $E^\rm{u}_\infty$-algebra is up to a suspension the same as taking $E^\rm{nu}_\infty$-indecomposables of its augmentation ideal by \cref{thm:indec-is-bar}, so we have
\[\SStL = H_*(\Sigma \, \indec_{E^\rm{nu}_\infty}(\SSt_{>0})) \quad \text{where} \quad \indec_{E^\rm{nu}_\infty}(\SSt_{>0}) \in \coAlg^\nil_{BE_\infty^\rm{nu}}(\Fun(\Vect,\DQ)))\]
with $BE_\infty^\rm{nu} \simeq s \,\rm{coLie}$, so $\indec_{E^\rm{nu}_\infty}(\SSt_{>0})$ admits up to a suspension a conilpotent Lie coalgebra structure by Koszul duality.

\medskip

The following is a direct consequence of \cref{lem:sth-coproduct-symmetry}, as it tells us that after composing with $\overline{\zeta}$ the cobracket is not just skew-commutative but also commutative. We give a different proof using an explicit formula for the cobracket in \cref{lem:fc-cobracket-symmetry}. There is a map
\begin{equation}\label{eqn:zeta-alt}\begin{aligned}\zeta^\alt = \ol{\zeta}-\ul{\zeta} \colon \Lambda^2 X &\lra X \para X \\
x \wedge y &\longmapsto \frac{1}{2}(x \para \overline{y} - y \para \overline{x}),\end{aligned}\end{equation}
where the overlines indicate that we use the naturality of $X$ with respect to the isomorphism $V_2 \overset{\cong}\lra V/V_1$ arising from a splitting $V_1 \oplus V_2 \overset{\cong}\lra V$. It can also be thought of as the inclusion $\inc \colon \Lambda^2 X \to X \levi X$ followed by $\ol{\zeta} \colon X \levi X \to X \para X$.

\begin{lemma}\label{lem:stl-cobracket-vanishing} The cobracket on $\SStL$ has the following vanishing property: 
\[\zeta^\alt \circ \delta = 0.\]
\end{lemma}

\subsubsection{Presentations and formulas for infinite Steinberg modules}\label{sec:pres-form-stl} We continue with a presentation for the infinite Steinberg modules and a formula for the cobracket, which by \cref{sec:alg-str-sstl} are obtained by passing to the indecomposables in $\SStH$ and antisymmetrising the (reduced) cobracket. If the dimension of $V$ is $n$ then a pair of apartments $[v_1,\ldots,v_n],[w_1,\ldots,w_n] \in \St(V)$ gives rise to an element
\[[v_1,\compactldots,v_n] \otimes [w_1,\compactldots,w_n] \in \StH(V) \cong \St(V) \otimes \St(V).\]
The following elements will play an important role \cite[Definition 28]{CharltonRadchenkoRudenko}:

\begin{definition}\label{def:steinberg-iterated-integral} Let $v_1,\ldots,v_n$ be a basis of $V$, then the \emph{Steinberg iterated integral} is the element
\[\rm{I}[v_1,\compactldots,v_n] \coloneq (-1)^n [v_n,v_{n-1},\compactldots,v_1] \otimes [v_n,v_{n-1}-v_n,\compactldots,v_1-v_2] \in \StH(V).\]
\end{definition}

It is convenient to think of these geometrically as a pair of simplexes in projective spaces that are \emph{not} in general position. For example, for $n=3$, we have two triangles in the projective space $\mathbb{P}^2(F)$ intersecting in the following pattern:
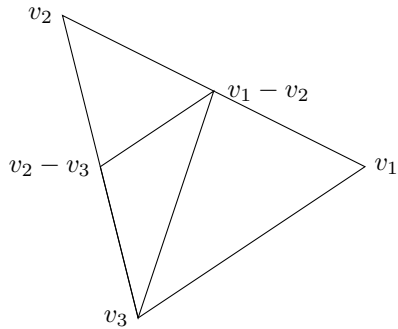
\[\begin{tikzpicture}
\draw (0,0) -- (-1,4) -- (3,2) -- cycle;
\node at (0,0) [left] {$v_3$};
\node at (-1,4) [left] {$v_2$};
\node at (3,2) [right] {$v_1$};
\node at (-.5,2) [left] {$v_2-v_3$};
\node at (1,3) [right] {$\,v_1-v_2$};
\draw (0,0) -- (-.5,2) -- (1,3) --cycle;
\end{tikzpicture}\]

Taking the quotient by the decomposables yields a projection map $\pi \colon \SStH \to \SStL$ and we will combine this with the Steinberg iterated integrals to construct elements of $\StL(V)$ that depend on a so-called \emph{affine basis}: if $V$ is $n$-dimensional then this is a collection of $n+1$ vectors $u_0,\ldots,u_n$ in $V$ so that $u_1-u_0,\ldots,u_n-u_0$ are linearly independent. The following are defined in \cite[Section 3.8]{CharltonRadchenkoRudenko}:

\begin{definition}\label{def:steinberg-correlator} Let $u_0,\ldots,u_n$ be an affine basis, then the \emph{Steinberg correlator} is the element 
\[\rm{C}[u_0:\compactcdots:u_n] \coloneq \pi\big[(-1)^n \rm{I}[u_1-u_0,\compactldots,u_n-u_0]\big] \in \StL(V).\]
\end{definition}

In \cite[Section 3.8]{CharltonRadchenkoRudenko} it is shown these have the following properties:
\begin{enumerate}[(1)]
\item \label{enum:stl-relations-i}  They are homogeneous $\rm{C}[u_0:\cdots:u_n] = \rm{C}[u_0-u:\compactcdots:u_n-u]$ for any $u \in V$.
\item \label{enum:stl-relations-ii} They are cyclically symmetric: $\rm{C}[u_0:u_1:\compactcdots:u_n] = \rm{C}[u_1:u_2:\compactcdots:u_0]$.
\item \label{enum:stl-relations-iii} They satisfy the shuffle relations: 
\[\sum_{\sigma \in \rm{Sh}(n_1,n_2)} \rm{C}[u_0:u_{\sigma(1)}:\compactcdots:u_{\sigma(n_1+n_2)}] =0 \quad \text{for $n=n_1+n_2$.}\]
\end{enumerate}
We will see momentarily, in \eqref{eqn:decomposition-operator}, that there is a linear map $D^\St_h \colon \StL(V) \to \StL(V)$ for each nonzero functional $h \colon V \to F$, and the only remaining required relations can be written in terms of these:
\begin{enumerate}[(1)] \setcounter{enumi}{3}
\item \label{enum:stl-relations-iv} The Steinberg correlators satisfy 
\[\rm{C}[u_0:\compactcdots:u_n] = D^\St_h\big(\rm{C}[u_0:\compactcdots:u_n]) \quad \text{for $h \in V^\vee \setminus \{0\}$.}\] 
\end{enumerate}

\begin{proposition}\label{prop:stl-explicit-pres}The following map of $\bb{Q}[\GL(V)]$-modules is an isomorphism
\[\frac{\bb{Q}[\rm{C}[u_0:\compactcdots:u_n] \text{ for affine bases $u_0,\ldots,u_n$}]}{\text{\eqref{enum:stl-relations-i}-\eqref{enum:stl-relations-iv}}} \overset{\cong}\lra \StL(V).\]
\end{proposition}

\begin{proof}Referring forward to \cref{prop:stl-resolution} for details, there is an exact sequence
\[\bb{Q}[\Dec_V] \otimes \FC(V) \overset{d}{\lra} \FC(V) \overset{\pr^\FC} \lra \StL(V) \lra 0\]
where $\FC(V)$ is generated by ``formal correlators'' satisfying \eqref{enum:stl-relations-i}--\eqref{enum:stl-relations-iii} and \eqref{enum:stl-relations-iv} is imposed by the differential $d$.
\end{proof}

We can now ask for a description of the Lie coalgebra structure with respect to this presentation: this was done in \cite[Section 3.8]{CharltonRadchenkoRudenko}.

\begin{proposition}\label{prop:stl-explicit-cobracket} With respect to this presentation, the cobracket on $\SStL$ is given by
\begin{equation}\label{eqn:stinfty-cobracket}
\delta\big(\rm{C}[u_0:\cdots:u_n]\big) = \sum_{j=0}^n \sum_{i=1}^{n-1} \rm{C}[u_j:\compactcdots:u_{j+i}] \wedge \rm{C}[u_j:u_{j+i+1}:u_{j+i+2}:\compactcdots:u_{j+n}]\end{equation}
where indices are to be interpreted cyclically.
\end{proposition}

For $I \subseteq \{0,\ldots,n\}$ we write $V_I \coloneq \rm{span}(u_i-u_j \mid i,j \in I)$, and see that the cobracket of $\rm{C}[u_0:\compactcdots:u_n]$ is concentrated on summands $\StL(V_I) \otimes \StL(V_J) \subseteq (\SStL \levi \SStL)(V)$ for direct sum decomposition $V_I \oplus V_J \smash{\overset{\cong}\lra} V$ where both subspaces are of this form.

It is convenient to represent these formulas pictorially. A Steinberg correlator $\rm{C}[u_0:\cdots:u_n]$ can be represented by a $(n+1)$-sided polygon whose vertices are decorated cyclically in clockwise order by $u_0,\ldots,u_n$. For each choice of vertex $j$ and nonadjacent edge $(j+i,j+i+1)$ (that is, $j \neq j+i,j+i+1$), we draw a cut from the vertex to the edge. We then interpret both sides as $(i+1)$- and $(n-i+1)$-sided polygons with vertices labelled cyclically by a subset of $u_0,\ldots,u_n$ (the vertex through which we cut will be duplicated, appearing once in both polygons) and take the corresponding Steinberg correlator, where the one clockwise from the cut vertex appears first in the wedge product: for example,
\[\begin{tikzpicture}[baseline={([yshift=-.5ex]current bounding box.center)}]
   \draw (0:1) \foreach \x in {60,120,...,360} {  -- (\x:1) };
   \foreach \x/\l/\p in
     { 60/{$u_2$}/above,
      120/{$u_1$}/above,
      180/{$u_0$}/left,
      240/{$u_5$}/below,
      300/{$u_4$}/below,
      360/{$u_3$}/right
     }
     \node[inner sep=1pt,circle,draw,fill,label={\p:\l}] at (\x:1) {};
     \draw (360:1) -- (90:1);
\end{tikzpicture} \quad \leadsto \quad \text{$\rm{C}[u_3,u_4,u_5,u_0,u_1] \wedge \rm{C}[u_3,u_2]$}\]
corresponds to the term $j=3$ and $i=4$ in the formula for $\delta(\rm{C}[u_0:\cdots:u_5])$.

\medskip

We can apply the map $\zeta^\alt$ from \eqref{eqn:zeta-alt} to $\delta$: then, in $\zeta^\alt \circ \delta$ we will have twice the number of terms:
\begin{align*}(\zeta^\alt \circ \delta)\big(\rm{C}[u_0:\cdots:u_n]\big) &= \sum_{j=0}^n \sum_{i=1}^{n-1} \tfrac{1}{2} \rm{C}[u_j:\cdots:u_{j+i}] \para \overline{\rm{C}[u_j:u_{j+i+1}:v_{j+i+2}:\cdots:u_{j+n}]} \\
&\quad - \sum_{j=0}^n \sum_{i=1}^{n-1} \tfrac{1}{2} \rm{C}[u_j:u_{j+i+1}:v_{j+i+2}:\cdots:u_{j+n}] \para \overline{\rm{C}[u_j:\cdots:u_{j+i}]}.\end{align*}
In terms of the pictorial description, we mark one of the two polygons to be ``projected'', and add a required negative sign when this marked polygon appears clockwise from the cut vertex. We proved in \cref{lem:stl-cobracket-vanishing} using an abstract argument that $\zeta^\alt \circ \delta$ vanishes and in \cref{lem:fc-cobracket-symmetry} will give a computational proof.

\subsection{Symbol maps and decomposition operators}\label{sec:symbol-maps} We next recall from \cite[Section 3.9]{CharltonRadchenkoRudenko} how to give bases for $\StH(V)$ and $\StL(V)$ and how to pass between different bases. This uses the symbol and to explain its combinatorial structure it is useful to first define a ``universal symbol''.

\subsubsection{Universal symbol}
Consider a set $S$ and letters $\omega_{ij}$ for $i\neq j \in S$ subject to the relations $\omega_{ij}=\omega_{ji}$. 

\begin{definition}\emph{Universal symbols} are the elements
\[
\mathbf{S}(s_0,\compactldots,s_n)\in \free_{\rm{As}^\rm{u}}(\Q\{\omega_{ij} \mid i,j\in S\}) \qquad  \text{for $n\geq 1$ and $s_0,\compactldots,s_n \in S$,}
\]
defined inductively by the formulas $S(s_0,s_1) = \omega_{s_0s_1}$, and 
\begin{equation}\label{eqn:symbol-combinatorics}
\begin{aligned}\mathbf{S}(s_0,\compactldots,s_n) &=\mathbf{S}(s_0,\compactldots,s_{n-1})\otimes \omega_{s_0s_n} \\
&\quad -\sum_{i=1}^{n-1} \mathbf{S}(s_0,\compactldots,\widehat{s}_{i+1},\compactldots,s_n)\otimes \omega_{s_is_{i+1}} \\
&\quad+\sum_{i=1}^{n-1} \mathbf{S}(s_0,\compactldots,\widehat{s}_i,\compactldots,s_n)\otimes \omega_{s_is_{i+1}}.\end{aligned}
\end{equation}
\end{definition}

It is sometimes convenient to have an expanded formula for the universal symbol, rather than an inductive one. The inductive definition of the universal symbol implies that $\bf{S}(s_0,\ldots,s_n)$ is a sum of symbols of the form
\[
\pm \omega_{i_1j_1}\otimes \dots \otimes \omega_{i_{n},j_n}
\]
where pairs $\{i_1j_1\} ,\dots \{i_{n},j_n\}$ can be identified with the set of edges of a spanning tree of the complete graph with vertices $s_0,\dots,s_n$. It is easy to see from \eqref{eqn:symbol-combinatorics} that each set of edges appears at most once. We conclude that:

\begin{proposition}\label{prop:universal-symbol-combinatorics} There exists a unique subset $T(n)$ of the set of ordered tuples of edges of the complete graph with vertex set $\ul{n} = \{0,1,\ldots,n\}$ and a unique sign function $\rm{sign}\colon T(n) \to \{\pm 1\}$ such that 
\[
\mathbf{S}(s_0,\compactldots,s_n)=\sum_{\iota=(\{i_1,j_1\},\dots,\{i_{n},j_n\})\in T(n)}\rm{sign}(\iota)\ \omega_{i_1j_1}\otimes \dots \otimes  \omega_{i_nj_n}.
\]
\end{proposition} 

\begin{example}\label{exam:universal-symbol} We have 
\[ \mathbf{S}(s_0,s_1,s_2) = \omega_{s_0s_1}\otimes \omega_{s_0s_2} - \omega_{s_0s_1}\otimes \omega_{s_1s_2} + \omega_{s_0s_2}\otimes \omega_{s_1s_2}.\]
That is, $T(2)$ consists of tuples $(\{0,1\},\{0,2\}), (\{0,1\},\{1,2\}),$ and $(\{0,2\},\{1,2\})$ with the signs
\[
\rm{sign}\bigl((\{0,1\},\{0,2\})\bigr)=1, \quad \rm{sign}\bigl((\{0,1\},\{1,2\})\bigr)=-1, \quad \rm{sign}\bigl((\{0,2\},\{1,2\})\bigr)=1.
\]\end{example}

\subsubsection{A property of the universal symbol related to partitions} We continue with a discussion that will allow us to prove cancellation of certain terms in the symbol; it suffices to do this for the universal symbol.

Consider an equivalence relation $\sim$ on the set $S$ and a collection of letters $a_{ij}$ for $i, j\in S$ together with an additional letter $a$. Define a map of associative algebras
\begin{align*}
\rho_{\sim}\colon \free_{\rm{As}^\rm{u}}(\Q\{\omega_{ij}| i,j\in S\}) &\lra \free_{\rm{As}^\rm{u}}(\Q\{a_{ij} \mid i,j\in S\}\oplus \Q\{a\}) \\
\omega_{ij} &\longmapsto
\begin{cases}
 a_{ij} & \text{ if } i\sim j, \\
 a  & \text{ if } i \not \sim j. 
\end{cases}\end{align*}

\begin{lemma}\label{lem:symbol and equivalence relation} Assume that $s_0 \not \sim s_n$. We have 
\[
\rho_{\sim}\bigl(\bf{S}(s_0,\dots,s_n)\bigr)=a^{\otimes n}.
\]
\end{lemma}
\begin{proof}
We argue by induction on $n$. The base case $n=1$ is trivial as we assume $s_0 \not \sim s_1$. By the induction assumption,  we have
\begin{align*}
\rho_{\sim}\left (\sum_{i=1}^{n-1} \bf{S}(s_0,\compactcdots,\widehat{s}_i,\compactldots,s_n)\otimes \omega_{s_is_{i+1}}\right) &= a^{\otimes n-1}\otimes \sum_{i=1}^{n-1} \rho_{\sim}(\omega_{s_is_{i+1}}) \\
\rho_{\sim}\left(\sum_{i=1}^{n-2} \bf{S}(s_0,\compactldots,\widehat{s}_{i+1},\compactldots,s_n)\otimes \omega_{s_is_{i+1}} \right) &=a^{\otimes n-1}\otimes \sum_{i=1}^{n-2} \rho_{\sim}(\omega_{s_is_{i+1}}),\end{align*}
so using \eqref{eqn:symbol-combinatorics} we have 
\begin{align*}
\rho_{\sim}\bigl(\bf{S}(s_0,\dots,s_n)\bigr)&=
\rho_{\sim}\bigl(\bf{S}(s_0,\compactldots,s_{n-1})\bigr)\otimes a
\\
&-\rho_{\sim}\bigl(\bf{S}(s_0,\compactldots,s_{n-1})\big)\otimes \rho_{\sim}(\omega_{s_{n-1}s_n})\\
&+a^{\otimes (n-1)}\otimes \rho_{\sim}(\omega_{s_{n-1}s_n}).
\end{align*}
If $s_{n-1} \sim s_{n}$, then $s_{n-1} \not \sim s_0$ and we have 
\[
\rho_{\sim}(\bf{S}(s_0,\compactldots,s_{n-1}))=a^{\otimes (n-1)}
\]
by the induction hypothesis. In this case, the above formula simplifies to 
\[
\rho_{\sim}\bigl(\bf{S}(s_0,\dots,s_n)\bigr)=
a^{\otimes (n-1)}\otimes a
-a^{\otimes (n-1)}\otimes a_{s_{n-1}s_n}
+a^{\otimes (n-1)}\otimes a_{s_{n-1}s_n}=a^{\otimes n}.
\]
If on the other hand $s_{n-1} \not \sim s_{n}$, the above formula simplifies to
\[
\rho_{\sim}\bigl(\bf{S}(s_0,\dots,s_n)\bigr)=
\rho_{\sim}(\bf{S}(s_0,\compactldots,s_{n-1}))\otimes a
-\rho_{\sim}(\bf{S}(s_0,\compactldots,s_{n-1}))\otimes a
+a^{\otimes (n-1)}\otimes a= a^{\otimes n}.
\]
This completes the proof of the induction step.
\end{proof}

\begin{example}For the discrete equivalence relation this says that the number of terms in the universal symbol counted with sign is 1. For other equivalence relations it gives similar but more refined information: for example, if $1 \sim 2$ but no other elements are identified it tells us that the number of terms in the universal symbol containing $\omega_{12}$ counted with sign is $0$.\end{example}

\subsubsection{Symbol maps}
As a consequence of the Koszulity of Steinberg modules from \eqref{hyp:st-koszul}, the description $\SStH$ and $\SStL$ as the associative and commutative Koszul duals of $\SSt$ yields injections in the bar constructions of \cref{def:bar-as,def:bar-comm}
\begin{align*}s^\rm{As} \colon \StH(V) &\lra (\rm{B}^\rm{As}\SSt)_n(V) \\
s^\rm{Com} \colon \StL(V) &\lra (\rm{B}^\rm{Com}\SSt)_n(V)\end{align*}
for $V$ of dimension $n$, that we will refer to as the \emph{associative symbol map} and \emph{commutative symbol map} respectively. The former induces the latter, in that the following diagram commutes
\[\begin{tikzcd} \StH(V) \dar[two heads] \rar[hook]{s^\rm{As}} &[20pt]  (\rm{B}^\rm{As}\SSt)_n(V) \dar[two heads]\\[-5pt]
\StL(V)  \rar[hook]{s^\rm{Com}}  & (\rm{B}^\rm{Com}\SSt)_n(V)\end{tikzcd}\]
where the vertical maps are the natural maps that take the quotient by decomposables.

\begin{notation}Unless there is a risk of confusion, we abbreviate both symbol maps to $s$.
\end{notation}

To perform computations, we use the inductive formula for the associative symbol of a Steinberg iterated integral \cite[Lemma 30]{CharltonRadchenkoRudenko}: we have 
\[s(\rm{I}[v]) = -[v]\] by definition, and then set
\begin{equation}\label{eqn:symbol-iterated-integral}\begin{aligned}s(\rm{I}[v_1,\compactldots,v_n]) &=- s(\rm{I}[v_1,\compactldots,v_{n-1}])\otimes [v_n] \\
&\quad +\sum_{i=1}^{n-1} \big(s(\rm{I}[v_1,\compactldots,\widehat{v}_{i+1},\compactldots,v_n]-s(\rm{I}[v_1,\compactldots,\widehat{v}_i,\compactldots,v_n]\big)\otimes [v_{i+1}-v_i].\end{aligned}\end{equation}
Since the Steinberg correlators are up to sign given by projecting Steinberg iterated integrals, the formula for their symbol is essentially the same: $s(\rm{C}[u_0:u_1]) = [u_1-u_0]$, and 
\begin{equation}\label{eqn:symbol-correlator}\begin{aligned}s(\rm{C}[u_0:\compactcdots:u_n]) &=s(\rm{C}[u_0:\compactcdots:u_{n-1}])\otimes [u_n-u_0] \\
&\quad -\sum_{i=1}^{n-1} s(\rm{C}[u_0:\compactcdots:\widehat{u}_{i+1}:\compactcdots:u_n])\otimes [u_{i+1}-u_i] \\
&\quad+\sum_{i=1}^{n-1} s(\rm{C}[u_0:\compactcdots:\widehat{u}_i,\compactcdots:u_n])\otimes [u_{i+1}-u_i]\end{aligned}\end{equation}
where we implicitly project the right side to $(B^\rm{Com}\SSt)_n(V)$. Comparing the inductive formula \eqref{eqn:symbol-correlator} to the inductive formula \eqref{eqn:symbol-combinatorics} for the universal symbol, we see that
\[s(\rm{C}[u_0:\compactcdots: u_n]) = \pi\left(\sum_{\iota=(\{i_1,j_1\},\dots,\{i_{n},j_n\})\in T(n)}\rm{sign}(\iota)\ [u_{j_1}-u_{i_1}] \otimes \dots \otimes  [u_{j_n}-u_{i_n}]\right)\]
with $T(n)$ and $\rm{sign}(\iota)$ as in \cref{prop:universal-symbol-combinatorics}.

\begin{example}For $n=2$, we have \cite[Example 31]{CharltonRadchenkoRudenko} (compare to \cref{exam:universal-symbol}) 
\begin{align*}s(\rm{I}[v_1,v_2]) &= [v_1|v_2]-[v_1|v_2-v_1]+[v_2|v_2-v_1],\\
s(\rm{C}[u_0:u_1:u_2]) &= \pi\big([u_1-u_0|u_2-u_0]-[u_1-u_0|u_2-u_1]+[u_2-u_0 \mid u_2-u_1]\big).\end{align*}
For $n=3$, see \cref{sec:decomp-rels-wt-3}.
\end{example}

\subsubsection{Decomposition operators}
Unwinding the definitions, the target of the commutative symbol map is given by
\begin{align*}(\rm{B}^\rm{Com}\SSt)_n(V) &\cong (\rm{B}^\rm{Ass}\SSt)_n(V) \otimes_{\fr{S}_n} \rm{coLie}_n \\
&\cong \left(\bigoplus_{V = P_1 \oplus \cdots \oplus P_n} \St(P_1) \otimes \cdots \otimes \St(P_n)\right) \otimes_{\fr{S}_n} \rm{coLie}_n.\end{align*}
Given a hyperplane $H \subset V$ we can project onto those summands where none of the lines $P_i$ is contained in $H$, yielding a map
\[\pi_H \colon (\rm{B}^\rm{Com}\SSt)_n(V) \lra (\rm{B}^\rm{Com}\SSt)^H_n(V) \coloneq \left(\bigoplus_{\substack{V = P_1 \oplus \cdots \oplus P_n \\ P_1,\ldots,P_n \not \subseteq H}} \St(P_1) \otimes \cdots \otimes \St(P_n)\right) \otimes_{\fr{S}_n} \rm{coLie}_n.\]
This can be used to construct a basis for $\StL(V)$ \cite[Proposition 44]{CharltonRadchenkoRudenko}:
\begin{proposition}For every hyperplane $H \subset V$ the following is an isomorphism
\[s_H \coloneq (\pi_H \circ s) \colon \StL(V) \lra (\rm{B}^\rm{Com}\SSt)^H_n(V).\]
\end{proposition}

If we pick a (necessarily nonzero) linear functional $h \colon V \to F$ so that $H = \ker(h)$, we can construct an inverse by
\begin{align*}\rm{C}^\St_h \colon (\rm{B}^\rm{Com}\SSt)^H_n(V) &\lra \StL(V)\\
[P_1|\compactcdots|P_n] &\longmapsto \rm{C}[0:v_1:\compactcdots:v_n]\end{align*}
where $v_i$ is the unique vector in $P_i$ so that $h(v_i) = 1$. This is well-defined because the Steinberg correlators satisfy the shuffle relations. Given that $s_H$ is an isomorphism, to see that its inverse is given by $\smash{\rm{C}^\St_h}$ it suffices to verify:

\begin{lemma}\label{lem:ch-inverse-of-symbol} We have $s_H \circ \rm{C}^\St_h = \rm{id}$.
\end{lemma}

\begin{proof}Given a Steinberg correlator $\rm{C}[0:v_1:\compactcdots:v_n]$ in the image of $\rm{C}^\St_h$, in the inductive formula \eqref{eqn:symbol-iterated-integral} each term in the second sum has $\rm{span}(v_{i+1}-v_i) \subseteq H$. Thus only the term $[v_1|\compactcdots|v_n]$ of the symbol survives.
\end{proof}

Consequently if we define a \emph{decomposition operator} as
\begin{equation}\label{eqn:decomposition-operator} D^\St_h \coloneq (\rm{C}^\St_h \circ s_H) \colon \StL(V) \lra \StL(V)\end{equation}
it satisfies $D^\St_h = \id$. This explains and justifies relation \eqref{enum:stl-relations-iv}. Note that replacing $h$ with $\lambda h$ for $\lambda \neq 0$ replaces $v_i$ by $\lambda^{-1} v_i$, so $D_h^\St$ only depends on $H$; we will not use this.

\subsection{A lift of the coproduct to cobar complexes}\label{sec:lift-coproduct} In this subsection we will work in the 1-category $\rm{Ch}_\ds{Q}$ of chain complexes over $\ds{Q}$; this yields the category $\DQ$ upon inverting the quasi-isomorphisms, see \cref{sec:rect-dg}. Koszul duality yields equivalences in terms of the cobar constructions of \cref{def:cobar-as,def:cobar-colie}
\[\Sigma^{-1}\Omega^\rm{coAs} \SStH \simeq \SSt \qquad \text{and} \qquad \Sigma^{-1}\Omega^\rm{coLie} \SStL \simeq \SSt,\]
where the cobar constructions are constructed using the tensor product $\levi$, and the gradings are so that $\SStH(V)$ and $\SStL(V)$ are in degree $2\,\dim(V)$ and $\SSt(V)$ is in degree $\dim(V)$. In this section we will \emph{explicitly} lift to these resolutions the coproduct
\[\Delta \colon \SSt \lra \SSt \para \SSt.\]
This will eventually be used to give a formula for the cobracket on $H_1(\GL;\SStL)$. 

\begin{proposition}\label{prop:coproduct-on-cobar} There exist coassociative counital coproducts
\begin{align*}\Delta &\colon \Sigma^{-1}\Omega^\rm{coAs}(\SStH) \lra \Sigma^{-1}\Omega^\rm{coAs}(\SStH) \para \Sigma^{-1}\Omega^\rm{coAs}(\SStH) \\
\Delta &\colon \Sigma^{-1}\Omega^\rm{coLie}(\SStL) \lra \Sigma^{-1}\Omega^\rm{coLie}(\SStL) \para \Sigma^{-1}\Omega^\rm{coLie}(\SStL)\end{align*}
so that the augmentation maps to $\SSt$ are maps of coassociative counital coalgebras.
\end{proposition}

For $X \in \rm{Fun}(\Vect,\rm{Ch}_{\ds{Q}})$ we define the free unital associative and commutative algebras
\[T^\bullet X \coloneq \rm{As}^\rm{u} \circ X \qquad \text{and} \qquad S^\bullet X \coloneq \rm{Com}^\rm{u} \circ X,\]
whose underlying objects decompose as sums
\begin{align*}T^\bullet X &= \bigoplus_{p \geq 0} T^p X \qquad \text{where } T^p X = X^{\levi p},\\
S^\bullet X &= \bigoplus_{p \geq 0} S^p X \qquad \text{where } S^p X = \ds{Q}_\rm{triv} \otimes_{\fr{S}_p} X^{\levi p}.\end{align*} 
Recall that now $\para$ induces a monoidal structure on the categories of unital associative and unital commutative algebras with respect to $\levi$, so we get unital associative and unital commutative algebras $T^\bullet X \para T^\bullet X$ and $S^\bullet X \para S^\bullet X$, with product in both cases given by $(x \para y)(x' \para y') = (-1)^{|x'||y|} (x \levi x') \para (y \levi y')$. We can then extend the assignment $x \mapsto x \para 1 + 1 \para x$ uniquely to homomorphisms
\[\xi^\rm{As} \colon T^\bullet X \lra T^\bullet X \para T^\bullet X \qquad \text{and} \qquad \xi^\rm{Com} \colon S^\bullet X \lra S^\bullet X \para S^\bullet X.\]
For example, we have 
\begin{equation}\label{eqn:xi-computation}\xi^\rm{As}(x\levi y) = \xi^\rm{Com}(x \levi y) = (x \levi y) \para 1 + x \para y + (-1)^{|x||y|} y \para x + 1 \para (x \levi y).\end{equation}
Note that these maps are counital and coassociative, because the identities that one needs to verify are between maps of unital associative or unital commutative algebras whose domain is free and they are easily verified on generators.

\begin{lemma}\,
\begin{enumerate}[(i)]
    \item Let $\bf{C}$ be a nonunital dg-coassociative coalgebra in $\Fun(\Vect,\rm{GrMod}_{\ds{Q}})$ with respect to $\levi$ concentrated in even degrees. Then
        \[\xi^\rm{As} \colon \Sigma^{-1}\Omega^\rm{coAs} \bf{C} \lra \Sigma^{-1}\Omega^\rm{coAs} \bf{C} \para \Sigma^{-1}\Omega^\rm{coAs} \bf{C}\]
    is a map of chain complexes if and only if $\ol{\zeta} \circ \Delta_\bf{C}  = \ol{\zeta} \circ \sigma \circ \Delta_\bf{C}$.
    \item Let $\bf{L}$ be a dg-coLie coalgebra in $\Fun(\Vect,\rm{GrMod}_{\ds{Q}})$ with respect to $\levi$ concentrated in even degrees. Then
        \[\xi^\rm{Com} \colon \Sigma^{-1}\Omega^\rm{coLie} \bf{L} \lra \Sigma^{-1}\Omega^\rm{coLie} \bf{L} \para \Sigma^{-1}\Omega^\rm{coLie} \bf{L}\]
    is a map of chain complexes if and only if $\zeta^\alt \circ \delta_\bf{L} =0$.
\end{enumerate}
\end{lemma}

\begin{proof}We give the details in the second case, as it is most relevant to this paper and the first case is entirely analogous. Recall from \cref{def:cobar-colie} that the Lie coalgebra cobar complex $\Sigma^{-1}\Omega^\rm{coLie} \bf{L}$ is given by $S^\bullet \Sigma^{-1} \bf{L}$ with differential $d = d_\Omega$ is determined as the unique derivation extending a map on $\Sigma^{-1} \bf{L}$ (there is no nonzero internal differential $d_\bf{L}$ since $\bf{L}$ has none). We need to verify that 
\[\xi^\rm{Com} \circ d = (d \para \id+\epsilon\, \id \para d) \circ \xi^\rm{Com}\] where $\epsilon$ is an appropriate sign coming from the tensor product of chain complexes, acting in bidegree $(p,q)$ by $(-1)^p$. Note that both sides are derivations, so it suffices to verify they agree on the generators $\Sigma^{-1}\bf{L}$. Fixing $\ul{x} \in \Sigma^{-1} \bf{L}$ the desuspension of $x \in \bf{L}$, if $\delta_\bf{L}(x) = \sum_i x_i \levi x'_i$ we have $d\ul{x} = \sum_i \ul{x}_i \levi \ul{x}'_i$.

Writing the middle two terms in \eqref{eqn:xi-computation} as $(\ol{\zeta}+\ol{\zeta}\circ \sigma)(x \levi y)$, we compute 
\begin{align*} \xi^\rm{Com} \circ d(\ul{x}) &= d(\ul{x}) \para 1+(\ol{\zeta}+\ol{\zeta}\circ \sigma) \circ \delta_\bf{L}(\ul{x})+1 \para d(\ul{x}), \\
(d \para \id+\epsilon\, \id \para d) \circ \xi^\rm{Com}(\ul{x}) &= d(\ul{x}) \para 1+1 \para d(\ul{x}),\end{align*}
where there is no sign $\epsilon$ because $1$ is in even degree. This establishes the result after incorporating a minus sign from the desuspension; $\sigma \circ \delta_\bf{L}(\ul{x}) = - \ul{\sigma \circ \delta_\bf{L}(x)}$ when $\bf{L}$ is concentrated in even degree so $(\ol{\zeta}+\ol{\zeta} \circ \delta) \circ \delta_\bf{L}(\ul{x})$ vanishes if and only if $\zeta^\alt \circ \delta(x)$ does.
\end{proof}

\begin{proof}[Proof of \cref{prop:coproduct-on-cobar}] Using \cref{lem:sth-coproduct-symmetry} for $\SStH$ with coproduct $\Delta^{(2)}$ or \cref{lem:stl-cobracket-vanishing} for $\SStL$ with cobracket $\delta$, the previous lemma provides the maps $\Delta$ as soon as we verify these induce the correct coproduct on $\SSt$. We give the details in the second case, as it is most relevant to this paper and the first case is entirely analogous. We will use the description from \cref{prop:st-explicit-prod-coprod} of the coproduct on $\SSt$ in terms of a shuffle coproduct on apartment classes, whose proof does not rely on this proposition.

In these terms, the augmentation map for the Lie coalgebra cobar complex is the unique multiplicative extension
\[S^\bullet 1_! \ds{k}[1] \lra \SSt\]
of the identification $\ds{k}[1] \cong \StL_1(F)$. When evaluated on $F^n$, the left side is spanned by symbols $[L_1,\compactldots,L_n]$ indexed by decompositions $L_1 \oplus \cdots \oplus L_n = F^n$ into lines, up to the relation $[L_{\sigma(1)},\compactldots,L_{\sigma(n)}] = (-1)^\sigma [L_1,\compactldots,L_n]$, and the map is given by sending this to the corresponding apartment class. We now observe that the shuffle coproduct of \eqref{eqn:st-coproduct} is the unique map $\SSt \to \SSt \para \SSt$ of commutative algebras extending $[v] \mapsto [v] \para 1+ 1\para [v]$ and hence is compatible with the maps $\xi^\rm{Com}$.
\end{proof}

\section{Higher apartments} The purpose of this section is to give a conceptual explanation of the origin of the Steinberg correlators that generate $\StL_n$. It may be skipped on a first reading; it is used later only to justify some formulas, but that can also be done by hand. Given an affine basis of $F^n$, we construct an ``apartment class'' map $\rm{apt} \colon  \rm{Lie}^\vee_{n+1} \otimes \bb{Q}_\rm{sign} \to \St_n$ and then extract from this by Koszul duality a ``higher apartment class'' map  $\rm{apt}_{\coLie} \colon \rm{cycLie}^\vee_n \otimes \bb{Q} \to \StL_n$
whose image contains Steinberg correlator $\rm{C}[u_0:\compactcdots:u_n]$.
This construction may be used to justify the formula for the cobracket of Steinberg correlators of \cref{prop:stl-explicit-cobracket}, as presaged in \cite[\S 5.1]{Goncharov01} (it is related to the genus zero case of \cite[\S 6, \S 8]{Gon19}). We follow \cref{conv:shorter-notation}.

\subsection{Restriction to subspaces of a fixed module}For a vector space $V$, we let $\rm{Sub}(V)$ be the discrete category of subspaces of $V$. This admits a symmetric promonoidal structure with $k$-fold iterated tensor products, for $k \geq 0$, given by the ``internal sum if disjoint'' functor
\begin{align*}\oplus_k \colon \rm{Sub}(V)^k &\lrapro \rm{Sub}(V) \\
(V_1,\ldots,V_k,W) &\longmapsto \begin{cases} \ast & \text{if $V_1 \oplus \cdots \oplus V_k \overset{\cong}\lra W$,} \\
\varnothing & \text{otherwise.}\end{cases}\end{align*}
There is a lax symmetric promonoidality on the functor 
\begin{align*}j_V \colon \rm{Sub}(V) & \lra \Vect \\
(W \subseteq V) &\longmapsto W \end{align*}
which induces a lax symmetric monoidality, with respect to the Day convolution symmetric monoidal structures, on the restriction functor
\[j_V^* \colon \Fun(\Vect,\scr{C}) \lra \Fun(\rm{Sub}(V),\scr{C}).\]

\begin{lemma}$j_V^*$ is symmetric monoidal.\end{lemma}

\begin{proof}It is clear that $j^*_V$ preserves the unit, and unwinding the definitions the lax monoidality is given on the subspace $W \subseteq V$ by the map
\[\bigsqcup_{\substack{V_1,V_2 \subseteq W,\\V_1 \oplus V_2 \xrightarrow{\cong} W}} F(V_1) \otimes G(V_2) \lra \colim_{\substack{V_1,V_2 \to W,\\V_1 \oplus V_2 \xrightarrow{\cong} W}} F(V_1) \otimes G(V_2)\]
which is an isomorphism by a cofinality argument.
\end{proof}

As a symmetric monoidal left adjoint, $j^*_V$ induces an equivalence
\[j^*_V(\SStL) \simeq \indec^\nil_{E^\rm{nu}_\infty}(j^*_V(\SSt))\]
in $\coAlg^{\nil}_\coLie(\Fun(\rm{Sub}(V),\DQ))$. It is thus possible to understand the cobracket on $\SStL$ in terms of that on $j^*_V(\SStL)$. The advantage of passing to this restricted setting is that we may choose $V$ to come equipped with additional structure, in this case an affine basis.

\subsection{Lie algebras of trees and derivations} \label{sec:lie-alg-trees-derivations} In this section we describe some algebraic constructions that are used below to understand better the Drinfeld--Kohno Lie algebra.

\begin{remark}Some of the references used below take completions of the free associative and free Lie algebras, but the algebraic structures preserve the free associative and free Lie algebras sitting inside their completions, restricting to these.
\end{remark}

\subsubsection{Associative variant} \label{sec:trees-as} For the associative variant, we mostly follow \cite{AKKN}. We fix a finite set $S$ and consider the free nonunital associative algebra on the finite-dimensional vector space $V_S \coloneq \bb{Q}\{X_s \mid s \in S\}$ of symbols $X_s$ for $s \in S$: 
\[\sf{as}_S \coloneqq \free^{\rm{As}^\rm{nu}}(V_S).\]
This admits an equivalent description in terms of trees: let $\sf{rtree}_S$ be the vector space spanned by isomorphism classes of rooted finite trivalent planar trees whose leaves are labelled by elements from the set of symbols $\{X_s \mid s \in S\}$, modulo the $IH$ (associativity) relation: for a subtree ($\circled{0},\circled{1},\circled{2},\circled{3}$ denote the remainder of the rooted tree, with root in $\circled{0}$) we have
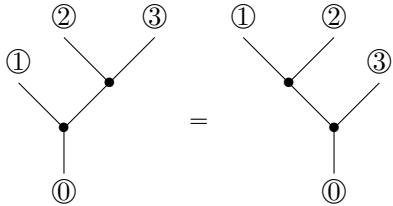
\[\begin{tikzpicture}[scale=.6,baseline]
    \node at (0,0) {$\bullet$};
    \node at (-1,1) [above] {$\circled{1}$};
    \node at (1,1) {$\bullet$};
    \node at (0,2) [above] {$\circled{2}$};
    \node at (2,2) [above] {$\circled{3}$};
    \node at (0,-1) [below] {$\circled{0}$};
    \draw (0,-1) --(0,0) -- (-1,1);
    \draw (0,0) -- (1,1) -- (0,2);
    \draw (1,1) -- (2,2);
\end{tikzpicture} \, = \, \begin{tikzpicture}[scale=.6,baseline]
    \node at (0,0) {$\bullet$};
    \node at (1,1) [above] {$\circled{3}$};
    \node at (-1,1) {$\bullet$};
    \node at (-2,2) [above] {$\circled{1}$};
    \node at (0,2) [above] {$\circled{2}$};
    \node at (0,-1) [below] {$\circled{0}$};
    \draw (0,-1) --(0,0) -- (1,1);
    \draw (0,0) -- (-1,1) -- (0,2);
    \draw (-1,1) -- (-2,2);
\end{tikzpicture}\]
Equivalence classes of such trees are rooted planar corollas and these correspond to words in the alphabet $\{X_s \mid s \in S\}$, in the sense that multiplication in the order described by the tree induces an isomorphism of vector spaces
\[\sf{rtree}_S \overset{\cong}\lra \sf{as}_S.\]
Under this isomorphism, the multiplication is given on the left side by taking the disjoint union of trees and joining their roots into a trivalent vertex:
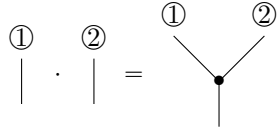
\[\begin{tikzpicture}[scale=.6, baseline]
    \node at (0,.5) [above] {$\circled{1}$};
    \draw (0,-.5) -- (0,.5);
\end{tikzpicture} \, \cdot \, \begin{tikzpicture}[scale=.6,baseline]
    \node at (0,.5) [above] {$\circled{2}$};
    \draw (0,-.5) -- (0,.5);
\end{tikzpicture} =  \begin{tikzpicture}[scale=.6,baseline]
    \node at (0,0) {$\bullet$};
    \node at (-1,1) [above] {$\circled{1}$};
    \node at (1,1)  [above] {$\circled{2}$};
    \draw (0,-1) --(0,0) -- (-1,1);
    \draw (0,0) -- (1,1);
\end{tikzpicture}\]
There is a Lie algebra of derivations $\sf{der}_S$ of $\sf{as}_S$ with two distinguished Lie subalgebras. Firstly, the \emph{tangential derivations} are given by
\[\sf{tder}_S \coloneq \{D \mid \text{for all $s \in S$ we have $D(X_s) = [a_s,X_s]$ for some $a_s \in \sf{as}_S$}\} \subseteq \sf{der}_S,\]
and secondly the \emph{special derivations} are given by
\[\sf{sder}_S \coloneq \{D \mid D({\textstyle \sum_{s \in S} X_s}) = 0\} \subseteq \sf{tder}_S.\]
A tangential derivation $D$ is uniquely determined by the elements $a_s$ such that $D(X_s) = [a_s,X_s]$, which yields a split short exact sequence \cite[Remark 2.1]{AKKN}
\[0 \lra {\textstyle \prod_{s \in S} \sf{as}_s} \lra {\textstyle \prod_{s \in S} \sf{as}_S} \lra \sf{tder}_S \lra 0\]
with left map induced by the inclusions $\sf{as}_s \to \sf{as}_S$ for $s \in S$, inducing an isomorphism $\rho$ between $\sf{tder}_S$ and tuples $(a_s)_{s \in S}$ so that $a_s$ is the kernel of the map sending all $X_{s'}$ for $s' \neq s$ to zero. Doing so, for $D = \rho((a_s)_{s \in S})$ and $D' = \rho((a'_s)_{s \in S})$ we have $[D,D'] = \rho((D(a'_s)-D'(a_s)-[a_s,a'_s])_{s \in S})$; the difference in sign with the reference is due to them rather defining $a_s$ by $D(X_S) = [X_S,a_s]$, cf.~\cite[p.~7]{AKKN}.

The Lie subalgebra $\sf{sder}_S$ of special derivations can be identified with the summand of the space of cyclic words in the alphabet $\{X_s \mid s \in S\}$, given by $|\sf{as}_S| \coloneq \sf{as}_S/[\sf{as}_S,\sf{as}_S]$, of those words in which at least two distinct letters occur. Additively, $|\sf{as}_S|$ can be described using that the nonunital associative operad $\rm{As}^\rm{nu}$ extends to a cyclic operad $\rm{cycAs}^\rm{nu}$ and taking 
\[|\sf{as}_S| \coloneq \bigoplus_{n \geq 1} \rm{cycAs}^\rm{nu}_n \otimes_{\fr{S}_{n+1}} V_S^{\otimes (n+1)}.\]
That is, there is an isomorphism \cite[Lemma 8.3]{AKKN} (using Remark 2.1 loc.cit.)
\begin{align*} |\sf{as}_S|/{\textstyle \prod_{s \in S} |\sf{as}_s|} &\overset{\cong}\lra \sf{sder}_S \\
|a| & \longmapsto (a_s)_{s \in S} \end{align*}
where the tuple $(a_s)_{s \in S}$ is determined by $N(|a|) = \sum_{s \in S} X_s a_s$, where $N \colon |\sf{as}_S| \to \sf{as}_S$ is the symmetrisation map sending a cyclic word $|z| = z_1\cdots z_k$ to $\sum_{1 \leq j \leq k} z_j \cdots z_k z_1\cdots z_{j-1}$. Under this isomorphism, the bracket of derivations is given by the (negative of the) Kirillov--Kostant--Souriau Lie bracket with explicit formula as in \cite[p.~25]{AKKN}; if $|z| = z_0\cdots z_k$ and $|w| = w_0\cdots w_\ell$ then
\begin{equation}\label{eqn:kks-bracket} [|z|,|w|] = -\sum_{i=0}^k\sum_{j=0}^\ell \delta_{z_i,w_j} \left(\parbox{6.5cm}{\centering $w_1\cdots w_{j-1} z_{i+1}\cdots z_k z_1 \cdots z_i w_{j+1} \cdots w_\ell$ \\
$- w_1\cdots w_{j-1} z_{i}\cdots z_k z_1 \cdots z_{i-1} w_{j+1} \cdots w_\ell$}\right).\end{equation}
This in turn admits an equivalent description in terms of cyclic trees: let $\sf{tree}_S$ be the vector space spanned by isomorphism classes of finite trivalent planar trees whose leaves are labelled by elements from the set of symbols $\{X_s \mid s \in S\}$, modulo the $IH$ relation (now without root). Equivalence classes of such cyclic trees are in bijection with planar corollas corresponding to cyclic words in the alphabet $\{X_s \mid s \in S\}$, and multiplication in cyclic order induces an isomorphism
\[\sf{tree}_S \overset{\cong}\lra |\sf{as}_S| \qquad \qquad \begin{tikzpicture}[scale=.6,baseline]
    \node at (0,0) {$\bullet$};
    \foreach \i in {0,...,5} 
    {
    \draw (0,0) -- (360/6*\i:1);
    \node at (-360/6*\i:1.5) {$X_{z_\i}$};
    }
\end{tikzpicture}  \longmapsto X_{z_0}X_{z_1}\cdots X_{z_5}.\]
The Lie bracket on cyclic words is given in terms of cyclic trees as follows: for two cyclic trees $T_1,T_2$, $[T_1,T_2] = T_1 \circ T_2 - T_2 \circ T_1$ where $T_1 \circ T_2$ is given by the sum over $s \in S$, all leaves of $T_1$ labelled by $X_s$, and all leaves of $T_2$ labelled by $X_s$, of $T_1$ and $T_2$ joined at these leaves to a trivalent vertex labelled by $X_s$. The following is an example:
\[\begin{tikzpicture}[scale=.6,baseline]
    \node at (0,0) {$\bullet$};
    \foreach \i in {0,...,5} 
    {
    \draw (0,0) -- (360/6*\i:1);
    \node at (-360/6*\i:1.5) {$X_{z_\i}$};
    }
\end{tikzpicture} \circ \begin{tikzpicture}[scale=.6,baseline]
    \node at (0,0) {$\bullet$};
    \foreach \i in {0,...,3} 
    {
    \draw (0,0) -- (360/4*\i:1);
    \node at (-360/4*\i:1.5) {$X_{w_\i}$};
    }
\end{tikzpicture} = \sum_{i=0}^5 \sum_{j=0}^3 \delta_{z_i,w_j}  
\begin{tikzpicture}[scale=.6,baseline]
    \node at (0,0) {$\bullet$};
    \foreach \i in {0,...,5} 
    {
    \draw (0,0) -- (360/6*\i:1);
    }
    \foreach \i in {1,...,5} 
    {
    \node at (-360/6*\i:1.5) {$X_{z_{i+\i}}$};
    }
    \begin{scope}[xshift=4cm]
    \node at (0,0) {$\bullet$};
    \foreach \i in {0,...,3} 
    {
    \draw (0,0) -- (360/4*\i:1);
    }
    \node at (-360/4*-1:1.5) {$X_{w_{j-1}}$};
    \node at (-360/4*0:1.5) {$X_{w_{j}}$};
    \node at (-360/4:1.5) {$X_{w_{j+1}}$};
    \end{scope}
    \draw (0,0) -- (4,0);
    \node at (2,0) {$\bullet$};
    \draw (2,0) -- (2,-2);
    \node at (2,-2) [below] {$X_{z_i} = X_{w_j}$};
\end{tikzpicture}
\]
where we remark that the right picture can be collapsed to a 9-valent corolla.

Moreover, $\sf{sder}_S$ acts on $\sf{as}_S$ and the description in terms of trees is as follows: for a cyclic tree $T$ and rooted tree $R$, $D_T(R)$ is given by the sum over $s \in S$, all leaves of $T$ labelled by $X_s$, and all leaves $R$ labelled by $X_s$, of these leaves joined into a single edge.

\subsubsection{The Lie representation and its dual} \label{sec:lie-rep} Before moving on the Lie variant, we recall some representations that will appear in that setting. First, let us construct a pair of $\bb{Z}[\fr{S}_n]$-modules:
\begin{definition}\,
\begin{enumerate}[(i)]
    \item Let $\rm{Lie}_n$ be the span of those Lie words containing each generator exactly once in the free Lie algebra on generators $X_1,\ldots,X_n$.
    \item Let $\rm{Lie}_n^\vee$ denote its linear dual $\rm{Hom}_\bb{Z}(\rm{Lie}_n,\bb{Z})$.
\end{enumerate}
\end{definition}

The former can be identified with the $\bb{Z}[\fr{S}_n]$-module of rooted trivalent trees with a cyclic order at each vertex whose leaves are in bijection with $\ul{n}$, modulo the AS (anti-symmetry) and IHX (Jacobi) relations; we call these \emph{Lie trees}. These relations apply to a subtree ($\circled{0},\circled{1},\circled{2},\circled{3}$ denote the remainder of the rooted tree, with root in $\circled{0}$) and are given respectively by:
\[\begin{tikzpicture}[scale=.6,baseline]
    \node at (0,0) {\rotatebox{-45}{$\circlearrowright$}};
    \node at (-1,1) [above] {$\circled{1}$};
    \node at (1,1) [above] {$\circled{2}$};
    \node at (0,-1) [below] {$\circled{0}$};
    \draw (0,-1) --(0,0) -- (-1,1);
    \draw (0,0) -- (1,1);
\end{tikzpicture} \, = -\, \begin{tikzpicture}[scale=.6,baseline]
    \node at (-1,1) [above] {$\circled{2}$};
    \node at (1,1) [above] {$\circled{1}$};
    \node at (0,-1) [below] {$\circled{0}$};
    \draw (0,-1) --(0,0) -- (-1,1);
    \draw (0,0) -- (1,1);
    \node at (0,0) {\rotatebox{45}{$\circlearrowleft$}};
\end{tikzpicture}\]
\[\begin{tikzpicture}[scale=.6,baseline]
    \draw (-1,-1) -- (0,-.4) -- (1,-1);
    \draw (1,1) -- (0,.4) -- (-1,1);
    \draw (0,-.4) -- (0,.4);
    \node at (0,-.4) {\rotatebox{45}{$\circlearrowleft$}};
    \node at (0,.4) {\rotatebox{45}{$\circlearrowleft$}};
    \node at (1,1) [above right] {$\circled{2}$};
    \node at (-1,1) [above left] {$\circled{1}$}; 
    \node at (1,-1) [below right] {$\circled{3}$};
    \node at (-1,-1) [below left] {$\circled{0}$};
\end{tikzpicture}=
\begin{tikzpicture}[scale=.6,baseline]
    \draw (-1,-1) -- (-.4,0) -- (-1,1);
    \draw (1,-1) -- (.4,0) -- (1,1);
    \draw (-.4,0) -- (.4,0);
    \node at (-.4,0) {\rotatebox{45}{$\circlearrowleft$}};
    \node at (.4,0) {\rotatebox{45}{$\circlearrowleft$}};
    \node at (1,1) [above right] {$\circled{2}$};
    \node at (-1,1) [above left] {$\circled{1}$}; 
    \node at (1,-1) [below right] {$\circled{3}$};
    \node at (-1,-1) [below left] {$\circled{0}$};
\end{tikzpicture}
\,-\,\begin{tikzpicture}[scale=.6,baseline]
    \draw (-1,-1) -- (1,1);
    \draw [line width=1ex,white] (1,-1) -- (-1,1);
    \draw (1,-1) -- (-1,1);
    \draw (-.5,-.5) -- (.5,-.5);
    \node at (-.5,-.5) {\rotatebox{45}{$\circlearrowleft$}};
    \node at (.5,-.5) {\rotatebox{45}{$\circlearrowleft$}};
    \node at (1,1) [above right] {$\circled{2}$};
    \node at (-1,1) [above left] {$\circled{1}$}; 
    \node at (1,-1) [below right] {$\circled{3}$};
    \node at (-1,-1) [below left] {$\circled{0}$};
\end{tikzpicture}
\]
If $\circled{0}$ is the root, $\circled{1} = X,\circled{2} = Y,\circled{3}=Z$, these say $[X,Y] = -[Y,X]$ and $[[X,Y],Z]= [X,[Y,Z]] -[Y,[X,Z]]$. Its linear dual can be identified with the quotient of the free $\bb{Z}[\fr{S}_n]$-module on a single generator by $(i,n-i)$-shuffles for $1 \leq i \leq n-1$ \cite[Theorem 2.6]{Whitehouse}. (Over a field $\ds{k}$ of characteristic zero, there is a noncanonical isomorphism $\rm{Lie}_n \cong \rm{Lie}_n^\vee$ but this is not the case over the integers, and we find it illuminating to distinguish them.)

By construction, we may identify $\rm{Lie}_n$ with the space of $n$-ary operations of the operad $\rm{Lie}$. Since the operad $\rm{Lie}$ extends to a cyclic operad $\rm{cycLie}$ \cite[3.9(c)]{GetzlerKapranov}, $\rm{Lie}_n$ must be the restriction of a $\bb{Z}[\fr{S}_{n+1}]$-module $\rm{cycLie}_n$ and dually $\rm{Lie}^\vee_n$ must be the restriction of a $\bb{Z}[\fr{S}_{n+1}]$-module $\rm{cycLie}^\vee_n$. The former can be described as $\bb{Z}[\fr{S}_{n+1}]$-module of trivalent trees with a cyclic order at each vertex whose leaves are in bijection with the set $[n] = \{0,1,\ldots,n\}$, modulo AS and IHX relations; we call these \emph{cyclic Lie trees}. The latter was described by Whitehouse: $\rm{cycLie}^\vee_n$ is the quotient of the free $\bb{Z}[\fr{S}_{n+1}]$-module on a single generator by $(i,n-i)$-shuffles for $1 \leq i \leq n-1$ and cyclic symmetry \cite[p.~317, 319]{Whitehouse}.

\subsubsection{Lie variant} \label{sec:trees-lie} There is a closely related Lie variant of the associative construction in \cref{sec:trees-as}, and in explaining it we mostly follow \cite{AlekseevTorossian}. We again fix a finite set $S$ and consider the free Lie algebra on the vector space $V_S = \bb{Q}\{X_s \mid s \in S\}$ of symbols $X_s$ for $s \in S$:
\[\fr{lie}_S \coloneq \free_{\rm{Lie}}(V_S).\]
This admits a description in terms of Lie trees: let $\fr{rtree}_S$ be the vector space spanned by isomorphism classes of rooted finite trivalent trees with a cyclic order at each vertex whose leaves are labelled by elements from the set of symbols $\{X_s \mid s \in S\}$, modulo the AS and IHX relations. Equivalence classes of Lie trees correspond to Lie words in the alphabet $\{X_s \mid s \in S\}$, in the sense that sending a trivalent vertex with positive cyclic order to the Lie bracket induces an isomorphism of vector spaces
\[\fr{rtree}_S \overset{\cong}\lra \fr{lie}_S.\]
Under this isomorphism, the Lie bracket is given on the left side by taking the disjoint union of trees and joining their roots into a trivalent vertex with positive cyclic order:
\[\left[\begin{tikzpicture}[scale=.6, baseline]
    \node at (0,.5) [above] {$\circled{1}$};
    \draw (0,-.5) -- (0,.5);
\end{tikzpicture} \, , \, \begin{tikzpicture}[scale=.6,baseline]
    \node at (0,.5) [above] {$\circled{2}$};
    \draw (0,-.5) -- (0,.5);
\end{tikzpicture}\right] =  \begin{tikzpicture}[scale=.6,baseline]
    \node at (0,0) {\rotatebox{45}{$\circlearrowleft$}};
    \node at (-1,1) [above] {$\circled{1}$};
    \node at (1,1)  [above] {$\circled{2}$};
    \draw (0,-1) --(0,0) -- (-1,1);
    \draw (0,0) -- (1,1);
\end{tikzpicture}\]

There is a Lie algebra $\fr{der}_S$ of derivations of $\fr{lie}_S$ with two distinguished Lie subalgebras. Firstly, the \emph{tangential derivations} are given by
\[\fr{tder}_S \coloneq \{D \mid \text{for all $s \in S$ we have $D(X_s) = [a_s,X_s]$ for some $a_s \in \fr{lie}_S$}\}\subseteq \fr{der}_S,\]
and secondly the \emph{special derivations} are given by
\[\fr{sder}_S \coloneq \{D \mid D(\textstyle{\sum_{s \in S} X_s})=0\} \subseteq \fr{tder}_S.\]
A tangential derivation $D$ is uniquely determined by the elements $a_s$ such that $D(X_s) = [a_s,X_s]$, which yields a split short exact sequence \cite[Remark 3.3]{AlekseevTorossian}
\[0 \lra {\textstyle \prod_{s \in S} \fr{lie}_s} \lra {\textstyle \prod_{s \in S} \fr{lie}_S} \lra \fr{tder}_S \lra 0\]
inducing an isomorphism $\rho$ between $\fr{tder}_S$ and tuples $(a_s)_{s \in S}$ so that $a_s$ is in the kernel of the map sending all $X_{s'}$ for $s' \neq s$ to zero. Doing so, for $D= \rho((a_s)_{s \in S})$ and $D' = \rho((a'_s)_{s \in S})$ we have $[D,D'] = \rho((D(a'_s)-D'(a_s)-[a_s,a_{s'}])_{s \in S})$ \cite[p.~423]{AlekseevTorossian}; the sign difference has the same origin as in the associative variant.

The Lie subalgebra $\fr{sder}_S$ of special derivations can be identified with the summand of the space of cyclic Lie words in the alphabet $\{X_s \mid s \in S\}$ of those words in which at least two distinct letters occur (this only rules out the case where no bracketings are present). Additively, this space is given by first recalling that the Lie operad $\rm{Lie}$ extends to a cyclic operad $\rm{cycLie}$ and then taking 
\[|\fr{lie}_S| \coloneq \bigoplus_{n \geq 1} \rm{cycLie}_n \otimes_{\fr{S}_{n+1}} V_S^{\otimes (n+1)}.\]
That is, there is an isomorphism
\begin{align*} |\fr{lie}_S|/{\textstyle \prod_{s \in S} |\fr{lie}_s|} &\overset{\cong}\lra \fr{sder}_S \\
    |a| &\longmapsto (a_s)_{s \in S}.\end{align*}
This admits a description in terms of cyclic Lie trees: let $\fr{tree}_S$ be the vector space spanned by isomorphism classes of finite trivalent trees with a cyclic order at each vertex whose leaves are labelled by elements from $\{X_s \mid s \in S\}$, modulo the AS and IHX relations. Equivalence classes of cyclic Lie trees correspond to cyclic Lie words in the alphabet $\{X_s \mid s \in S\}$, in the sense that there is an isomorphism of vector spaces
\[\fr{tree}_S \overset{\cong}\lra |\fr{lie}_S|.\]
There is an isomorphism 
\begin{align*}\fr{tree}_S/{\textstyle \prod_{s \in S} \fr{tree}_s} & \overset{\cong}\lra \fr{sder}_S \\
\Gamma &\longmapsto (a_s)_{s \in S} \end{align*}
given on a tree $\Gamma$ by letting $a_s$ be given by a sum over all leaves of $\Gamma$ with label $s$, of the Lie tree obtained by interpreting that leaf as a root \cite[Section 2.2.1]{AlekseevTorossian} (see also \cite[Appendix A]{Felder}). Under this isomorphism, the Lie bracket on $\fr{sder}_S$ is given as follows: for two cyclic Lie trees $T_1,T_2$, the bracket $[T_1,T_2]$ is given by the sum over all $s \in S$, all leaves of $T_1$ labelled by $X_s$, and all leaves of $T_2$ labelled by $X_s$, of $T_1$ and $T_2$ joined at these leaves to a trivalent vertex with positive orientation and remaining leaf labelled by $X_s$. Similarly, the action of $\fr{sder}_S$ on $\fr{lie}_S$ is given as follows: for a cyclic Lie tree $T$ and a rooted Lie tree $R$, $D_T(R)$ is given by the sum over $s \in S$, all leaves of $T$ labelled by $X_s$, and all leaves of $R$ labelled by $X_s$, of those leaves joined into a single edge.

\subsubsection{Comparison of the associative and Lie variant}\label{sec:trees-comp} As the similarity of the descriptions of the associative and Lie variants suggests, there are inclusions
\[\fr{lie}_S \lra \sf{as}_S \qquad \text{and} \qquad \fr{sder}_S \lra \sf{sder}_S\]
given by including a free Lie algebra into its universal enveloping algebra and uniquely extending a special derivation of a free Lie algebra to the one of its universal enveloping algebra. It is clear from the constructions that these maps are compatible with the Lie brackets on the free algebras and special derivations, and the actions of special derivations on the corresponding free algebras \cite[p.~7]{AKKN}.

\subsubsection{Relationship to universal traces} The above discussion is related to universal invariant bilinear forms, due to Drinfeld \cite[p.~856]{Drinfeld}: for a Lie algebra $\fr{g}$ define a vector space
\[\mathscr{F}(\fr{g}) \coloneqq \frac{\fr{g} \otimes \fr{g}}{[X,Y] \otimes Z - X \otimes [Y,Z] \text{ and } X\otimes Y-Y \otimes X}.\]
In \cite[Section 2.2]{AlekseevTorossianNote} it is proven that there is an isomorphism of vector spaces 
\[\mathscr{F}(\fr{lie}_S) \overset{\cong}\lra \fr{tree}_S\]
given on a pair $X \otimes Y$ of Lie words in the $X_s$ by interpreting them as rooted trees and connecting their roots. In \cite[p.~857]{Drinfeld}, Drinfeld described the composite homomorphism
\[\mathscr{F}(\fr{lie}_S) \overset{\cong}\lra \fr{tree}_S \lra \fr{sder}_S\]
which is a surjection of Lie algebras (the only elements in the kernel are of the form $X_s \otimes X_s$) if the domain is given the Kirillov bracket, see also \cite[Section 6]{KontsevichFormal}. This gives another construction of the Lie bracket of special derivations.

\begin{remark}Additively, the construction of Drinfeld is a special instance of one for cyclic operads due to Getzler--Kapranov: inspecting \cite[Definition 4.7]{GetzlerKapranov} we note there is an equality $\mathscr{F}(\fr{g}) = \lambda(\rm{cycLie},\fr{g})$. The left side is defined for any cyclic operad $\scr{P}$ and $\scr{P}$-algebra $A$. If the latter is a free $\scr{P}$-algebra one gets an isomorphism \cite[Proposition 4.9]{GetzlerKapranov}, originally outlined in the cases $\mathscr{P} = \rm{Lie},\rm{Comm},\rm{Ass}$ by Kontsevich \cite[Section 4, 5]{KontsevichFormal}:
\[\bigoplus_{n \geq 0} \scr{P}(n) \otimes_{\fr{S}_{n+1}} V^{\otimes (n+1)} \overset{\cong}\lra \lambda(\scr{P},\free_\scr{P}(V))\]
induced by thinking of $\scr{P}(n) \otimes V^{\otimes (n+1)}$ as $V \otimes (\scr{P}(n) \otimes V^{\otimes n})$ and mapping this to $\free_\scr{P}(V) \otimes \free_\scr{P}(V)$ in the evident manner.\end{remark}

\subsubsection{Dual Lie coalgebras} \label{sec:sder-dual} 
Linearly dualising the discussions in \cref{sec:trees-as,sec:trees-lie}, degreewise by putting the generators in degree $1$, we obtain from \cref{sec:trees-comp} a pair of surjective maps of Lie coalgebras
\[\sf{as}_S^\vee \lra \fr{lie}_S^\vee \qquad \text{and} \qquad \sf{sder}_S^\vee \lra \fr{sder}_S^\vee\]
which are compatible with coactions of the duals of special derivations of the duals of the corresponding free algebras. Let us use this to deduce formulas for the cobracket and coaction.

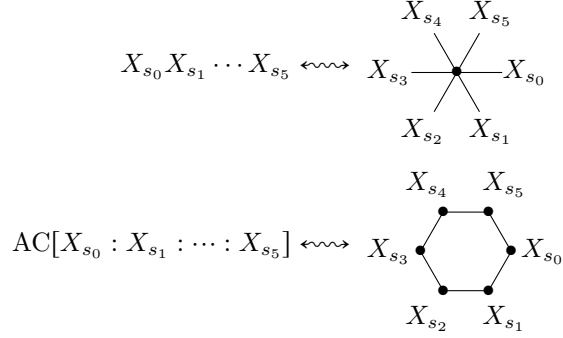
\begin{figure}
\begin{align*}X_{s_0}X_{s_1}\cdots X_{s_5} &\leftrightsquigarrow \begin{tikzpicture}[scale=.6,baseline]
    \node at (0,0) {$\bullet$};
    \foreach \i in {0,...,5} 
    {
    \draw (0,0) -- (360/6*\i:1);
    \node at (-360/6*\i:1.5) {$X_{s_\i}$};
    }
\end{tikzpicture} \\ 
\rm{AC}[X_{s_0}:X_{s_1}:\compactcdots:X_{s_5}] &\leftrightsquigarrow \begin{tikzpicture}[scale=.6,baseline]
    \foreach \i in {0,...,5} 
    {
    \node at (-360/6*\i:1) {$\bullet$};
    \draw (360/6*\i+360/6:1) -- (360/6*\i:1);
    \node at (-360/6*\i:1.7) {$X_{s_\i}$};
    }
\end{tikzpicture}\end{align*}
\caption{A pictorial interpretation of duals to cyclic words, compatible with our discussion of infinite Steinberg modules, is as polygons.}
\end{figure}

Firstly, $\sf{as}_S^\vee$ is additively generated by symbols dual to nonempty words in the $\{X_s \mid s \in S\}$ or equivalently planar rooted corollas with leaves labelled by these symbols. That is, a nonempty word $X_{s_1}\compactcdots X_{s_n}$ for $n \geq 1$ has a dual symbol $\rm{AI}[X_{s_1},\compactldots,X_{s_n}]$ and in terms of these the dual cobracket, which we will denote  $\delta_\rm{dec}$, is simply deconcatenation:
\[\delta_\rm{dec}(\rm{AI}[X_{s_1},\compactldots,X_{s_n}]) = \sum_{i=1}^{n-1} \rm{AI}[X_{s_1},\compactldots,X_{s_i}] \wedge \rm{AI}[X_{s_{i+1}},\compactldots,X_{s_n}].\]
Secondly, $\sf{sder}_S^\vee$ is additively generated by symbols dual to the cyclic words in the symbols $\{X_s \mid s \in S\}$ or equivalently planar corollas with leaves labelled by these symbols. That is, a cyclic word $X_{s_0}\compactcdots X_{s_n}$ for $n \geq 0$ has a dual symbol $\rm{AC}[X_{s_0}: \compactcdots : X_{s_n}]$ and these are by definition cyclically symmetric: 
\[\rm{AC}[X_{s_0}: X_{s_1}:\compactcdots : X_{s_n}] = \rm{AC}[X_{s_1}:X_{s_2}: \compactcdots : X_{s_0}].\] 
The formula for the cobracket, which we will denote $\delta_\rm{cyc}$, is dual to the formula \eqref{eqn:kks-bracket} for the KKS bracket and is given by
\[\delta_\rm{cyc}(\rm{AC}[X_{s_0}: \compactcdots : X_{s_n}]) = \sum_{j=0}^n \sum_{i=1}^{n-1} \rm{AC}[X_{s_j}:\compactcdots : X_{s_{j+i}}] \wedge \rm{AC}[X_{s_j}:X_{s_{j+i+1}}:\compactcdots : X_{s_{j+n}}].\]

The surjective maps $\sf{as}_S^\vee \to \fr{lie}_S^\vee$ and $\sf{sder}_S^\vee \to \fr{sder}_S^\vee$ show that the target is additively generated by similar symbols $\rm{LI}[X_{s_1},\compactldots,X_{s_n}]$ and $\rm{LC}[X_{s_0}:\compactcdots: X_{s_n}]$ that are the images of the above symbols. These additionally satisfy the shuffle relations:
\begin{align*}\sum_{\sigma \in \rm{Sh}(n_1,n_2)} \rm{LI}[X_{s_{\sigma(1)}},\compactldots,X_{s_{\sigma(n_1+n_2)}}] &= 0 \quad \text{for $n=n_1+n_2$ with $n_1,n_2>0$,} \\
\sum_{\sigma \in \rm{Sh}(n_1,n_2)} \rm{LC}[X_{s_0}:X_{s_{\sigma(1)}}:\compactcdots:X_{s_{\sigma(n_1+n_2)}}] &= 0 \quad \text{for $n=n_1+n_2$ with $n_1,n_2>0$}.\end{align*}
Since the aforementioned surjective maps are maps of Lie coalgebras, the cobrackets are given by the same formulas; the additional information that we gain is that they are compatible with the shuffle relations.

\subsection{The Arnold algebra and Drinfeld--Kohno Lie coalgebra} Before defining our higher apartment classes we recall the \emph{Arnold algebra}, and its Koszul dual Lie coalgebra, the \emph{Drinfeld--Kohno Lie coalgebra}. We then relate the latter to special  derivations.

\subsubsection{Definitions of the Arnold algebra and Drinfeld--Kohno Lie coalgebra}

\begin{definition}For a finite set $S$ we define the \emph{Arnold algebra} to be the commutative algebra
\[\rm{Ar}_S \coloneq \frac{\free_\rm{Com}(\ds{Q}\{\omega_{ij} \mid i,j \in S,\,i \neq j\})}{(1){-}(2)} \in \Alg_\rm{Com}(\rm{GrMod}_\bb{Q}).\]
where the generators $\omega_{ij}$ have degree $1$ and the relations are given by
\begin{enumerate}[\indent (1)]
    \item $\omega_{ij} = \omega_{ji}$,
    \item $\omega_{ij}\omega_{jk}+\omega_{jk}\omega_{ki}+\omega_{ki}\omega_{ij} = 0$.
\end{enumerate}\end{definition}
For $S = \ul{n} \coloneq \{1,2,\ldots,n\}$, $\rm{Ar}_\ul{n}$ admits an action of the symmetric group $\fr{S}_n$. For subsets $S_1,S_2 \subseteq S$ there is a multiplication map $m \colon \rm{Ar}_{S_1} \otimes \rm{Ar}_{S_2} \to \rm{Ar}_S$. Define a graded vector space of (underived) indecomposables
\[\rm{Ar}^\indec_\ul{n} \coloneqq \rm{coker}\left(\bigoplus_{\substack{S = S_1 \sqcup S_2 \\ S_1,S_2 \neq \varnothing}} \rm{Ar}_{S_1} \otimes \rm{Ar}_{S_2} \overset{m}\lra \rm{Ar}_\ul{n}\right)\]
as the summand that cannot be written as product of elements from smaller subsets; this includes in particular all products of generators where not all $i \in \ul{n}$ appear as an index. It is concentrated in degree $n-1$ and in this degree the quotient $\rm{Ar}_\ul{n} \to \rm{Ar}^\indec_{\ul{n}}$ is an isomorphism. In fact, the degree $i$ part of $\rm{Ar}_\ul{n}$ can be expressed in terms of these as a sum over partitions of $\ul{n}$ into disjoint subsets $S_1,\ldots,S_r$ with $(|S_1|-1)+\cdots+(|S_r|-1)=i$ of terms $\rm{Ar}^\indec_{\ul{S}_1} \otimes \cdots \otimes \rm{Ar}^\indec_{\ul{S}_r}$ \cite[Section 6]{Cohen}.

It is well-known that the Arnold algebra is a Koszul algebra, as we explain now. For $S = \ul{n}$ the Arnold algebra in degree $r$ has a basis given by $\omega_{i_1j_1} \cdots \omega_{i_rj_r}$ where $i_1<j_1, \dots i_r<j_r$, and $j_1<j_2<\dots<j_r$ \cite[Corollary 3]{Arnold}. Secondly, with this order we get a PBW basis in the sense of \cite[Section 5.1]{Priddy}, cf.~\cite[Corollary 2.2]{Bezrukavnikov}, and deduce Koszulity as in \cite[Theorem 5.3]{Priddy}. Its Koszul Lie coalgebra is thus given by its quadratic dual Lie coalgebra; this will be degreewise finite-dimensional and may be more familiarly described through its linear dual Lie algebra:

\begin{definition}For the finite set $S$ we define the \emph{Drinfeld--Kohno Lie algebra} to be
\[\rm{DK}_S \coloneq \frac{\free_\rm{Lie}(\ds{Q}\{t_{ij} \mid i,j \in S,\,i \neq j\})}{(1){-}(3)}\]
where the generators $t_{ij}$ have degree $2$ and the relations are given by
\begin{enumerate}[\indent (1)]
    \item $t_{ij} = t_{ji}$,
    \item $[t_{ij},t_{kl}] = 0$ if $\{i,j\} \cap \{k,l\} = \varnothing$,
    \item $[t_{ij},t_{ik}+t_{jk}] = 0$.
\end{enumerate}\end{definition}

\begin{remark}We deviate slightly from the usual discussion: often the Arnold algebra is cohomologically graded with generators in degree 1 and thus the generators of the quadratic Drinfeld--Kohno Lie algebra are in degree 0. Here homological grading is more appropriate and the latter appear in degree 2 instead. This has no consequences apart from some degree shifts.\end{remark}

There are inclusions $\rm{DK}_{\ul{n}\setminus \{i\}} \to \rm{DK}_{\ul{n}}$ split by maps setting all generators with $i$ among their indices to zero. We define a graded vector space
\[\rm{tDK}_\ul{n} \coloneq \rm{coker}\left(\bigoplus_{i \in \ul{n}} \rm{DK}_{\ul{n}\setminus \{i\}} \overset{\inc}\lra \rm{DK}_\ul{n}\right) \overset{\cong}\longleftarrow \bigcap_{i \in \ul{n}} \ker\left(\rm{DK}_{\ul{n}} \to \rm{DK}_{\ul{n}\setminus \{i\}}\right)\]
with right isomorphism induced by the inclusion. It is the summand spanned by those bracketings of generators where each $i \in \ul{n}$ appears as an index. This is nonzero only in degree $\geq 2n-2$.

\begin{remark}Taking $S = \ul{n}$, assigning to $\omega_{ij}$ the form $\smash{\frac{1}{2\pi i}\frac{d(z_i-z_j)}{z_i-z_j}} \in \Omega^1_\rm{dR}(\rm{Conf}_n(\bb{C});\bb{C})$ induces an isomorphism \cite{Arnold}
\[\rm{Ar}_{\ul{n}} \otimes \bb{C} \overset{\cong}\lra H^*(\rm{Conf}_n(\bb{C});\bb{C}).\]
From this perspective, the degree completion of the Drinfeld--Kohno Lie algebra $\rm{DK}_\ul{n}$ agrees with the Malcev completion of the fundamental group of $\rm{Conf}_n(\bb{C})$. Alternatively, forgoing completion at the cost of degree shifts and sign representations, one can replace $\bb{C} = \bb{R}^2$ by $\bb{R}^d$ for $d>2$ and Malcev completion with the homotopy Lie algebra, cf.~\cite[Example 5.5]{BerglundKoszul}.
\end{remark}

\subsubsection{The Drinfeld--Kohno Lie algebra and special derivations}\label{sec:recol-dk-sder} To better understand the Drinfeld--Kohno Lie algebra $\rm{DK}_S$, one can think of it in terms of trees and special derivations. This is summarised by the following diagram of Lie algebras
\[\begin{tikzcd} & \rm{DK}_S \dar[hook] & \\[-5pt]
{|\fr{lie}_S|/{\textstyle \prod_{s \in S} |\fr{lie}_S|}} \rar{\cong} & \fr{tree}_S/{\textstyle \prod_{s \in S} \fr{tree}_s} \rar{\cong} & \fr{sder}_S,\end{tikzcd}\]
with definitions of the objects on the bottom row given in \cref{sec:lie-alg-trees-derivations}. Roughly, the vertical map has target a quotient of a Lie algebra of trees with external vertices labelled by $S$ and is determined as a map of Lie algebras by sending the generator $t_{ij}$ to an edge connecting the vertices $i$ and $j$. We will now provide details.

\smallskip

We first recall the description of the Drinfeld--Kohno Lie algebra in terms of graphs. For a finite set $S$, a \emph{graph with external vertices $S$} will be an unoriented finite graph $\Gamma$ with a linear order on its edges, a specified set $S$ of external vertices and possibly further internal vertices, satisfying
\begin{enumerate}
    \item it has no double edges,
    \item it has no simple loops,
    \item all internal vertices are at least trivalent,
    \item all internal vertices are connected by a path to an external vertex,
    \item the (open) graph obtained by deleting the external vertices is connected (``internally connected''),
\end{enumerate}
We let $\sf{CG}(S)$ \cite[Section 2, 3]{SeveraWillwacher} be the graded vector space spanned by graph with external vertices $S$, modulo the anti-symmetry relation $\Gamma^\sigma = (-1)^{\sigma} \Gamma$ where $\Gamma^\sigma$ is obtained from $\Gamma$ by permuting the edges by $\sigma$, and the grading of $\Gamma$ is $1-\# \text{edges}+2 \# \text{internal vertices}$. 

This admits the structure of an $L_\infty$-algebra whose differential (increasing degree) is given by splitting each of the vertices in all possible ways with new edge last in the order, and retaining only those terms that are internally connected. Its binary bracket $[\Gamma_1,\Gamma_2]$ is given by gluing $\Gamma_1$ and $\Gamma_2$ at $S$, concatenating the order of orders, and applying the differential; effectively, we sum over all $s \in S$ and all ways of gluing to a tripod both an edge in $\Gamma_1$ attached to $s$ and an edge in $\Gamma_2$ also attached to $s$:
\[\left[\begin{tikzpicture}[scale=.6,baseline]
   \foreach \i in {1,...,3} 
    {
    \node at (-360/3*\i:1.5) {$\i$};
    \node at (-360/3*\i:1) {$\bullet$};
    }
    \draw (-360/3*1:1) -- (-360/3*2:1);
\end{tikzpicture},
\begin{tikzpicture}[scale=.6,baseline]
   \foreach \i in {1,...,3} 
    {
    \node at (-360/3*\i:1.5) {$\i$};
    \node at (-360/3*\i:1) {$\bullet$};
    }
    \draw (-360/3*1:1) -- (-360/3*3:1);
\end{tikzpicture}\right] = \begin{tikzpicture}[scale=.6,baseline]
   \foreach \i in {1,...,3} 
    {
    \node at (-360/3*\i:1.5) {$\i$};
    \node at (-360/3*\i:1) {$\bullet$};
    \draw (-360/3*\i:1) -- (0,0);
    }
\end{tikzpicture}\]

The map sending $t_{ij}$ to the graph with a unique edge connecting the external vertices $i$ and $j$ induces an isomorphism \cite[Proposition 2]{SeveraWillwacher} (to deal with different conventions for encoding orientations, see \cite[Section 2.3.1]{ConantVogtmann})
\[\rm{DK}_S \overset{\cong}\lra H^0(\mathsf{CG}(S)).\]

We next define $\mathsf{CG}_\mathrm{tree}(S)$ as the quotient $\mathsf{CG}(S)/F^1 \mathsf{CG}(S)$ by the first step of the filtration by number of internal loops. The induced map of Lie algebras \cite[p.~185]{SeveraWillwacher}
\[H^0(\mathsf{CG}(S)) \lra \fr{tree}_S \cong H^0(\mathsf{CG}_\mathrm{tree}(S))\]
is injective. Here, as the notation suggests, the target can be identified as the Lie algebra of internally connected trivalent internal trees modulo the IHX relation. By construction, this map is uniquely determined by sending $t_{ij}$ to the tree with a unique edge connecting the external vertices $i$ and $j$, and bracket as described above.

\medskip

We use this to better understand the Drinfeld--Kohno Lie algebra. Recall from \cref{sec:lie-rep} that from the Lie operad we can extract $\bb{Z}[\fr{S}_n]$-modules $\rm{Lie}_n$ and $\rm{Lie}_n^\vee$, as well as a $\bb{Z}[\fr{S}_{n+1}]$-module $\rm{cycLie}_n$. They admit the following topological interpretation, whose statement uses the following objects: (i) the poset $\scr{P}(\ul{n})$ of proper nonempty partitions of the set $\ul{n}$ ordered by inclusion, which admits an $\fr{S}_n$-action by permuting the elements of $\ul{n}$ and whose geometric realisation is equivalent to a wedge of $(n-2)$-spheres, and (ii) the configuration space $\rm{Conf}_n(\bb{R}^k)$ of $n$ ordered points in $\bb{R}^k$ for $k \geq 2$.

\begin{lemma}\label{lem:reps-top-interpretation}\,
\begin{enumerate}[(i)]
\item \label{enum:reps-top-interpretation-i} $\rm{Lie}_n^\vee \otimes \bb{Z}_\rm{sign} \overset{\cong}\lra \widetilde{H}_{n-2}(\scr{P}(\ul{n});\bb{Z})$.
\item \label{enum:reps-top-interpretation-ii} $\rm{Lie}_n \otimes \bb{Z}_\rm{sign}^{\otimes (k-1)} \smash{\overset{\cong}\lra} \widetilde{H}_{(k-1)(n-1)}(\rm{Conf}_n(\bb{R}^k);\bb{Z})$.
\item \label{enum:reps-top-interpretation-iii} $\rm{Lie}_{n}^\vee \otimes \bb{Q}_\rm{sign} \smash{\overset{\cong}\lra} (\text{degree $n-1$ part of $\rm{Ar}^\indec_\ul{n}$})$.
\item \label{enum:reps-top-interpretation-iv} $\rm{cycLie}_n \otimes \ds{Q} \smash{\overset{\cong}\lra} (\text{degree $2n$ part of $\rm{tDK}_{[n]}$})$.
\end{enumerate}
\end{lemma}

\begin{proof}Part \eqref{enum:reps-top-interpretation-i} is \cite[Theorem 4.1]{Robinson}. Part \eqref{enum:reps-top-interpretation-ii} uses that $\rm{Conf}_n(\bb{R}^k)$ is equivalent to the space of $n$-ary operations in the $E_k$-operad. Using the identification of the homology of this operad in terms of a shifted Poisson operad \cite{Sinha}, the result then follows from \cite[Theorem 6.1]{Cohen}. Part \eqref{enum:reps-top-interpretation-iii} is obtained from this by taking $k=2$, linearly dualising, and recalling that the degree $n-1$ part of $\rm{Ar}_\ul{n}$ agrees with $\rm{Ar}^\indec_{\ul{n}}$. 

For part \eqref{enum:reps-top-interpretation-iv} one needs to trace through the isomorphisms of the previous subsection. We need to consider internally connected trees with $n$ internal edges, which must be spanning trees containing all $n+1$ external vertices; these in turn correspond to the summand of $\rm{cycLie}_n \otimes_{\fr{S}_{n+1}} (\ds{Q}\{X_0,\ldots,X_n\})^{\otimes (n+1)}$ where each $X_i$ appears exactly once, which is in turn isomorphic to $\rm{cycLie}_n \otimes \bb{Q}$.\end{proof}

\subsection{Higher apartment classes}\label{sec:higher-apts} We now give the precise definition of the higher apartment classes. We first recall a nonstandard construction of the apartment classes in the Steinberg module $\St(V)$. Fix an affine basis $\vec{v} = (v_0,\ldots,v_n)$ on $V$ (consequently $V$ is $n$-dimensional) and write $V_S \coloneq \rm{span}(v_i-v_j \mid i,j \in S)$ for a subset $S \subseteq [n] = \{0,1,\ldots,n\}$. Now consider the functor
\begin{align*}\Omega^{\vec{v}} \colon \rm{Sub}(V) &\lra \rm{GrMod}_\bb{Q} \\
W &\longmapsto \begin{cases}\bb{Q}\{\omega_{ij}\} & \text{if $W = V_{\{ij\}}$ with $i \neq j$,}\\
0 & \text{otherwise,}\end{cases}\end{align*}
where $\omega_{ij}$ has degree $1$ and by definition satisfies $\omega_{ij} = \omega_{ji}$. The free graded-commutative algebra $\free_\rm{Com}(\Omega^{\vec{v}})$ on this functor takes value $\bb{k}\{\omega_{ij}\omega_{jk},\omega_{jk}\omega_{ki},\omega_{ki}\omega_{ij}\}$ on $V_{\{ijk\}}$. Letting $R_{\Omega^{\vec{v}}}$ be the ideal generated by the elements $\omega_{ij}\omega_{jk}+\omega_{jk}\omega_{ki}+\omega_{ki}\omega_{ij}$, we then define
\[\rm{Ar}_{\vec{v}} \coloneq \frac{\free_\rm{Com}(\Omega^{\vec{v}})}{R_{\Omega^{\vec{v}}}} \in \rm{Alg}_\rm{Com}(\Fun(\rm{Sub}(V),\rm{GrMod}_\bb{Q})).\]
As the notation suggests, this is a lift of the Arnold algebra $\rm{Ar}_{[n]}$. Let us make this precise: as the unique functor $t \colon \rm{Sub}(V) \to \ast$ is lax symmetric promonoidal, it induces an oplax symmetric monoidal functor $t_! \colon \Fun(\rm{Sub}(V),\rm{GrMod}_\bb{Q}) \to \rm{GrMod}_\bb{Q}$ simply given by $F \mapsto \bigoplus_{W \subseteq V} F(W)$. Even though $t_!$ is \emph{not} strong symmetric monoidal, using that the $\omega_{ij}$ are in odd degree one may compute that sending generators to the image of generators under $t_!$ induces an isomorphism of graded vector spaces
\[\rm{Ar}_{[n]} \overset{\cong}\lra t_! \rm{Ar}_{\vec{v}}.\]
In fact, $\rm{Ar}_{\vec{v}}$ vanishes except when evaluated on $V_{S_1}+\cdots+V_{S_r}$ with $V_{S_i} \cap V_{S_j} = 0$ for a partition of $[n]$, in which case it is concentrated in degree $(|S_1|-1) + \cdots (|S_r|-1)$ and in that degree is isomorphic to $\rm{Ar}^\indec_{S_1} \otimes \cdots \otimes \rm{Ar}^\indec_{S_r}$.

There is a map of functors $\Omega_{\vec{v}} \to j_V^* \SSt$ given at $V_{\{ij\}}$ by sending $\omega_{ij}$ to the preferred generator of $\St(V_{\{ij\}}) \cong \Q$. The Bykovskii relations in $\St(V_{\{ijk\}})$ imply that the relations in $R_{\Omega^{\vec{v}}}$ go to zero, so there is an induced map of commutative algebras
\[\varpi_{\vec{v}} \colon \rm{Ar}_{\vec{v}} \lra j_V^* \SSt.\]
If we evaluate this on $V$ itself, we get a map $\rm{Ar}_{\vec{v}}(V) \to j_V^* \SSt(V)$ which in degree $n$, in terms of the isomorphism from \cref{lem:reps-top-interpretation} \eqref{enum:reps-top-interpretation-iii}, is given by
\begin{align*}(\text{degree $n$ part of $\rm{Ar}^\indec_{[n]}$}) \cong \rm{Lie}^\vee_{n+1} \otimes \ds{Q}_\rm{sign} &\lra \St(V) \\
\omega_{i_1j_1}  \cdots  \omega_{i_nj_n} &\longmapsto [V_{\{ i_1j_1\}}|\cdots |V_{\{i_nj_n\}}].\end{align*}
Thus we not only obtain the apartment classes, but also reveal some of their symmetries.

\begin{remark}A topological interpretation of this map, using \cref{lem:reps-top-interpretation} \eqref{enum:reps-top-interpretation-i}, is as follows. There is a map of posets $\scr{P}(\ul{n}) \to T(V)$ from the partition poset to the Tits building of $V$, sending a partition $\{S_1,\ldots,S_r\}$ of $[n]$ to the subspace $V_{S_1}+\cdots+V_{S_r} \subseteq V$. The induced map
\[\rm{Lie}^\vee_{n+1} \otimes \bb{Z}_\rm{sign} \cong \widetilde{H}_{n-2}(\scr{P}(\ul{n});\bb{Z}) \lra \widetilde{H}_{n-2}(T(V);\bb{Z}) \eqcolon \St(V)_{\bb{Z}}\]
agrees, after tensoring with $\ds{Q}$, with the one above.
\end{remark}

We next apply $\indec_{E_\infty^\rm{nu}}$, which can for example be computed by performing the Harrison complex construction in the category $\Fun(\rm{Sub}(V),\rm{Ch}_\bb{Q})$. The same proof as for the usual Arnold algebra---which constructs a Poincar\'e--Birkhoff--Witt basis and apply Priddy's criterion for Koszulity---yields the following Koszul duality result as it only involves symbolic manipulation of the generators and relations:

\begin{lemma}$\rm{Ar}_{\vec{v}}$ is Koszul.
\end{lemma}

The quadratic dual ${}^\vee \ol{\rm{DK}}_{\vec{v}}$ of $\rm{Ar}_{\vec{v}}$ is then a reduced variant of the predual of the Drinfeld--Kohno Lie algebra. Namely, its linear dual is the object of $\rm{Alg}_\rm{Lie}(\Fun(\rm{Sub}(V),\rm{GrMod}_\bb{Q}))$ given by
\[\ol{\rm{DK}}_{\vec{v}} \coloneq \frac{\free_{\rm{Lie}}(T^{\vec{v}})}{R_{T^{\vec{v}}}}\]
whose generators are given by
\begin{align*}T^{\vec{v}} \colon \rm{Sub}(V) &\lra \rm{GrMod}_\bb{Q} \\
W &\longmapsto \begin{cases}\Q\{t_{ij}\} & \text{if $W = V_{\{ij\}}$ with $i \neq j$,}\\
0 & \text{otherwise,}\end{cases}\end{align*}
where $t_{ij}$ has homological degree 2 and satisfies $t_{ij} = t_{ji}$, and $R_{T^{\vec{v}}}$ is the ideal generated by the elements $[t_{ij},t_{kl}]$ for $\{i,j\} \cap \{k,l\} = \varnothing$ as well as $[t_{ij},t_{ik}+t_{jk}]$. 

The oplax monoidality on $t_!$ induces an inclusion
\[t_! {}^\vee\overline{\rm{DK}}_{\vec{v}} \subseteq {}^\vee\rm{DK}_{[n]},\]
but this is generally \emph{not} an isomorphism. This is due to the tensor product on subspaces of $V$ being non-trivial only if they intersect in $0$. To describe its image observe that with notation as in \cref{sec:sder-dual}, there is a Lie sub-coalgebra
\[\ol{\fr{sder}}^\vee_{[n]} \subseteq \fr{sder}^\vee_{[n]}\]
where each symbol from $\{X_0,X_1,\ldots,X_n\}$ appears at most once. The linear dual of the inclusion $\rm{DK}_{[n]} \to \fr{sder}_{[n]}$ of \cref{sec:recol-dk-sder} yields a map of Lie coalgebras $\ol{\fr{sder}}^\vee_{[n]} \to {}^\vee \rm{DK}_{[n]}$ (identifying the target with its double dual) and this induces an isomorphism
\[t_! {}^\vee\ol{\rm{DK}}_{\vec{v}} \cong \rm{im}\big[\ol{\fr{sder}}^\vee_{[n]} \to {}^\vee \rm{DK}_{[n]}\big].\]

\medskip

Using that $\indec_\rm{Com}^\nil$ commutes with $j_V^*$, we thus obtain a map of Lie coalgebras
\[\tau_{\vec{v}} \colon {}^\vee \ol{\rm{DK}}_{\vec{v}} \lra j^*_V \SStL\]
which in degree $2n$, using \cref{lem:reps-top-interpretation} \eqref{enum:reps-top-interpretation-iv}, yields the \emph{higher apartment class} map in the following proposition. Here we use \cref{sec:sder-dual} to give generators for the domain.

\begin{proposition}\label{prop:apt-colie} There is a map
    \begin{align*}\rm{apt}_{\coLie} \colon \rm{cycLie}_n^\vee \otimes \ds{Q} &\lra \StL(V) \\
    \rm{LC}[X_0:X_1:\compactcdots:X_n] &\longmapsto \rm{C}[v_0:\compactcdots:v_n],\end{align*}
where the right side is as in \cref{def:steinberg-correlator}.
\end{proposition}

\begin{proof}It remains to justify the formula of $\rm{apt}_{\coLie}$. To do this, we use that $\tau_{\vec{v}}$ is a map of Lie coalgebras and both $n$-fold iterated cobrackets are injective by Koszul duality, we reduce to the case $n=1$. In this case it is true by construction.\end{proof}

\begin{remark}\label{rem:c-sder-cobracket}
Arguably \cref{prop:apt-colie} should have been used to \emph{define} $\rm{C}[v_0: \compactcdots : v_n]$. Let us explain how to deduce formula \eqref{eqn:stinfty-cobracket} if one were to define $\rm{C}[v_0:\compactcdots:v_n]$ as $\rm{apt}_{\coLie}(\rm{LC}[X_0:X_1:\compactcdots:X_n]) \in (j^*_V \StL)(V) \cong \StL(V)$. Recalling that $j_V^*$ is symmetric monoidal, it suffices to compute its cobracket in $j^*_V \StL$ and since $\tau_{\vec{v}}$ is a map of Lie coalgebras, it is the image under $\tau_{\vec{v}}$ of the cobracket in ${}^\vee \ol{\rm{DK}}_{\ul{v}}$. This is in turn determined by the cobracket of special derivations, using that $t_!$ sends no element to zero. The upshot is that the formula for the cobracket of $\rm{C}[v_0:\compactcdots:v_n]$ has the same formula as that for $\rm{LC}[X_0:X_1:\compactcdots:X_n]$.\end{remark}

\section{$E_\infty$-algebras of general linear groups and the definition of $\scr{G}(F)$} In this section we follow \cite{GKRW20} in defining a graded nonunital $E_\infty$-algebra $\BGLb(F)$ in terms of the general linear groups $\GL_n(F)$. We relate its indecomposables to Steinberg modules and define the Goncharov Lie coalgebra as in the introduction as
\[\scr{G}(F) \coloneq \bigoplus_{n \geq 1} \scr{G}_n(F) \qquad \text{with} \qquad \scr{G}_n(F) \coloneq H^{E_\infty}_{n,2n-1}(\BGLb(F)_\bb{Q}).\]
As \cref{fig:gltable1} indicates, these are those groups in each rank $n$ of lowest degree that can possibly be nonzero, except for $n=1$ (cf.~\cite[Figure 1]{GKRW20}). We explain these admit a Lie cobracket as well as an additional ``$\sigma$-component'' map, and we explain how to compute these in terms of infinite Steinberg modules. We follow \cref{conv:shorter-notation}.

\subsection{The $E_\infty$-algebras $\BGLb^+$ and $\BGLb$} \label{sec:bgl} We start by constructing the $E_\infty$-algebras $\BGLb^+$ and $\BGLb$, the former unital and the latter nonunital, as we will have use for both. We then use their $E_\infty$-homology to define the Goncharov Lie coalgebra.

Recall that Day convolution yields a symmetric monoidal structure on the category $\Fun(\Vect,\Spc)$ of functors from the symmetric monoidal groupoid of finite-dimensional vector spaces over $F$ with direct sum to the category of spaces with cartesian product, whose tensor product we denote by $\levi$. Similarly, there is a symmetric monoidal structure on $\Fun(\bb{N},\Spc)$ that we will also denote $\levi$. There is a symmetric monoidal functor $\dim \colon \Vect \to \bb{N}$ that assigns to a vector space its dimension, and it induces a symmetric monoidal functor $\dim_! \colon \Fun(\Vect,\Spc) \to \Fun(\bb{N},\Spc)$. This in turn induces functors on categories of algebras over operads in the domain and target. Let $\ul{\ast}$ denote the terminal object in $\Alg_{E^\rm{u}_\infty}(\Fun(\Vect,\Spc))$, and $\ul{\ast}_{>0} \in \Alg_{E^\rm{nu}_\infty}(\Fun(\Vect,\Spc))$ denote the nonunital $E_\infty$-algebra obtained by replacing the value on 0-dimensional vector spaces with $\varnothing$.

\begin{definition}We define
\begin{align*}&\BGLb^+ \coloneq \dim_!(\ul{\ast}) \in \Alg_{E^\rm{u}_\infty}(\Fun(\bb{N},\Spc)), \\
&\BGLb \coloneq \dim_!(\ul{\ast}_{>0}) \in \Alg_{E^\rm{nu}_\infty}(\Fun(\bb{N},\Spc)).\end{align*}
\end{definition}

\begin{example}The underlying objects have values at nonnegative numbers $n \in \bb{N}$ given by
\begin{align*}\BGLb^+(n) &\coloneqq n^* \fgt_{E_\infty^\rm{u}}\, \BGLb^+ \simeq \BGL_n \quad \text{for $n \geq 0$, } \\
\BGLb(n) &\coloneqq n^* \fgt_{E_\infty^\rm{nu}}\, \BGLb \simeq \begin{cases} \BGL_n & \text{for $n \geq 1$}, \\
\varnothing & \text{for $n=0$}.\end{cases}\end{align*} 
Under these equivalences, the multiplication  $\BGL_n \times \BGL_m \to \BGL_{n+m}$ of the $E_\infty$-algebra structures is induced by block sum of matrices.\end{example} 

Post-composing with the map induced by the symmetric monoidal functor $C_*(-;\ds{Q}) \colon \Spc \to \DQ$ we obtain
\[\BGLb^+_\ds{Q} \in \Alg_{E_\infty^\rm{u}}(\Fun(\bb{N},\DQ)) \quad \text{and} \quad \BGLb_\ds{Q} \in \Alg_{E_\infty^\rm{nu}}(\Fun(\bb{N},\DQ)).\]        
These are related as follows: the unique map $\ul{\ast}_{>0} \to \ul{\ast}$ in $\Alg_{E_\infty^\rm{nu}}(\Fun(\Vect,\Spc))$ induces a map in $\Alg_{E_\infty^\rm{nu}}(\Fun(\Vect,\DQ))$
\[\ul{\ds{Q}}_{>0} \coloneq C_*(\ul{\ast}_{>0};\ds{Q}) \lra  C_*(\ul{\ast};\ds{Q}) \eqcolon \ul{\ds{Q}},\]
whose target is reduced, so admits a canonical augmentation, and this map exhibits the domain as the augmentation ideal of the target. Applying $\dim_!$ we obtain from this a map $\BGLb_{\ds{Q}} \to \BGLb^+_{\ds{Q}}$ and this exhibits the domain as the augmentation ideal of the target (see \cref{sec:augmentation-ideals}).

\subsection{$E_\infty$-homology and the Goncharov Lie coalgebra} \label{sec:ek-homology}
We can consider $\BGLb_\bb{Q}$ as a nonunital $E_k$-algebra by restriction along the map of operads $E^\rm{nu}_k \to E^\rm{nu}_\infty$ and define \emph{$E_k$-homology} in terms of its cotangent complex from \cref{sec:operads-algebras}, (see also \cite[Section 10.1.6]{GKRW18}, called $E_k$-indecomposables there)
\[H^{E_k}_{n,d}(\BGLb_\bb{Q}) \coloneq H_d(\cot_{E_k^\rm{nu}}(\BGLb_\bb{Q})(n)).\]
The following describes these groups in terms of Steinberg modules---implicitly with rational coefficients---and their variants \cite[Section 6]{GKRW20}. 

\begin{theorem}[Galatius-Kupers--Randal-Williams] \label{thm:ekhomology-steinberg} We have isomorphisms
    \[H^{E_k}_{*,*}(\BGLb_\bb{Q}) \cong \begin{cases} H_{*,*}(\GL;\Sigma^{-1} \ISSt) & \text{if $k=1$,} \\
    H_{*,*}(\GL;\Sigma^{-2} \ISStH) & \text{if $k=2$,} \\
    H_{*,*}(\GL;\Sigma^{-2} \SStL) & \text{if $k=\infty$.}\end{cases}\]
\end{theorem} 

\begin{example}\label{exam:ekhomology-steinberg-gradings} Recall that our grading convention places $\St(V)$ in degree $\dim(V)$, and $\StH(V)$ and $\StL(V)$ in degree $2\dim(V)$. Thus we have \cite[Theorems 6.2, 6.5]{GKRW20}
\begin{align*}H_{n,d}(\GL;\Sigma^{-1} \ISSt) &\cong H_{d-n+1}(\GL_n;\St_n) \\
H_{n,d}(\GL;\Sigma^{-2} \ISStH) &\cong H_{d-2n+2}(\GL_n;\StH_n) \\
H_{n,d}(\GL;\Sigma^{-2} \SStL) &\cong H_{d-2n+2}(\GL_n;\StL_n).\end{align*} \end{example}

Koszul duality in the guise of \cref{lem:bialgebra-structure} endows these objects with additional structure, each in the category $\Fun(\bb{N},\rm{GrVect}_\bb{Q})$ with symmetric monoidal structure given by Day convolution with respect to addition on $\bb{N}$ and the graded tensor product (with Koszul sign rule) on $\rm{GrVect}_\bb{Q}$:

\begin{theorem}\label{thm:gl-st-structure} We have that
    \begin{align*} &H_{*,*}(\GL;\SSt) \text{ is a commutative bialgebra,} \\
    &H_{*,*}(\GL;\SStH) \text{ is a commutative cocommutative bialgebra,} \\
    & H_{*,*}(\GL;\Sigma^{-1} \SStL) \text{ is a Lie coalgebra}.\end{align*}
\end{theorem}

\begin{remark}This may be a bit surprising, as the coproduct on $\SStH$ is \emph{not} cocommutative. However, the Nesterenko--Suslin property as in \cref{sec:nesterenko-suslin} gives an isomorphism between $H_{*,*}(\GL;\SStH)$ and homology of general linear groups with a ``split'' variant of $\SStH$, which does have cocommutative coproduct.
\end{remark}

The last statement of \cref{thm:gl-st-structure} concretely says $H_{*,*}(\GL;\SStL)$ is a shifted Lie coalgebra: it has a cobracket of degree 1, i.e.~with components
\[\delta_{k,n-k} \colon H_d(\GL_n;\StL_n) \lra \bigoplus_{d'+d''=d+1} H_{d'}(\GL_k; \StL_k) \otimes H_{d''}(\GL_{n-k}; \StL_{n-k}).\]

We now recall some vanishing results for the $E_\infty$-homology groups of $\BGLb_\bb{Q}$. We start with the observation that since $\SStL$ are concentrated in bidegrees of the form $(n,2n)$, we must have that $\smash{H^{E_\infty}_{n,d}(\BGLb_\bb{Q})} = 0$ for $d<2n-2$. The following was established in \cite[Theorem B]{GKRW20} and reproved in \cite[Corollary 38]{CharltonRadchenkoRudenko}: by computing coinvariants we see that $H_0(\GL_n;\St^\infty(F^n)) = 0$ for $n \geq 2$ and $H_0(\GL_1;\St^\infty(F^1)) \cong \bb{Q}$ and obtain one additional degree of vanishing.

\begin{theorem}\label{thm:bgl-critical-line-vanishing} We have
\[H^{E_\infty}_{n,2n-2}(\BGLb_\bb{Q}) \cong \begin{cases} 0 & \text{for $n \geq 2$,} \\
\bb{Q} & \text{for $n=1$.}\end{cases}\]
\end{theorem}

Using \cref{thm:ekhomology-steinberg} we can rephrase these observations as a vanishing result for the Lie coalgebra $H_{*,*}(\GL;\Sigma^{-1}\SStL)$ (see \cref{fig:gltable1}). 

For $n=1$ we have $\GL_1 \cong F^\times$ and $\StL \cong \bb{Q}$, which should be thought of as being in degree $2$, so we get an isomorphism $\smash{H_{1,d}(\GL;\Sigma^{-1}\SStL) \cong \Lambda^{d-1} F^\times_\bb{Q}}$ \cite[Section 9.1.2]{GKRW20}.
We will occasionally write $H_{1,1}(\GL;\Sigma^{-1} \SStL) \cong \bb{Q}\{\sigma\}$. This may conflict with our notation for permutations and symmetry isomorphisms, but we believe there is no risk of confusion.

For $n \geq 2$ we can use the vanishing result \cref{thm:bgl-critical-line-vanishing} to get 
\[H_{n,d}(\GL;\Sigma^{-1} \SStL) \cong \begin{cases} \scr{G}_n & \text{if $d=2n$,} \\
0 & \text{if $d<2n$,}\end{cases} \qquad \text{with }
\scr{G}_n \coloneq H^{E_\infty}_{n,2n-1}(\BGLb_\bb{Q}) \cong H_1(\GL_n;\StL_n).\]
It follows that the cobracket on $H_{*,*}(\GL;\Sigma^{-1} \SStL)$ induces on $\scr{G} = \bigoplus_{n=1} \scr{G}_n$ not one but two structures:
\begin{enumerate}[(i)]
\item a \emph{cobracket} $\delta \colon \scr{G} \to \Lambda^2\scr{G}$,
\item a \emph{$\sigma$-component} $\delta_\sigma \colon \scr{G}_n\to H_2(\GL_{n-1};\StL_{n-1})$.
\end{enumerate}

\subsection{Vanishing results for $E_1$-homology} We will now explain how to bootstrap this to a vanishing result for $E_1$-homology, improving rationally on \cite[Theorem 10.2]{GKRW20} and \cite[Theorem A, B]{MPPII}, and independently observed by Randal-Williams and Galatius:

\begin{theorem}\label{thm:steinberg-homology-improved-vanishing} $H_{n,d}(\GL;\SSt) = 0$ when $d < 2n-2$.
\end{theorem}

\begin{example}Making grading conventions concrete, this says that $H_d(\GL_n;\St_n) = 0$ for $d<n-2$.\end{example} 

\cref{thm:steinberg-homology-improved-vanishing} will be a consequence of a more general result, phrased in terms of the augmented $E_\infty^\rm{u}$-algebra $\Bar(\BGLb^+_\bb{Q})$, as the proof of \cref{thm:ekhomology-steinberg} provides an isomorphism
\begin{equation}\label{eqn:bar-bgl-st} H_{*,*}(\Bar(\bf{BGL}_\bb{Q}^+)) \overset{\cong}\lra H_{*,*}(\GL;\SSt).\end{equation}
It will be useful to consider the more general case of an augmented $E_\infty^\rm{u}$-algebra $\bf{R}^+$ in $\Fun(\bb{N},\DQ)$ and the associated sequence of augmented $E_\infty^\rm{u}$-algebras $\Bar^k(\bf{R}^+)$.

In the following, $S^*$ denotes the free commutative algebra, with Koszul sign with respect to the second ``homological'' grading but no Koszul sign with respect to the first ``rank'' grading, and $c_0$ is the abstract connectivity given by $c_0(0) = 0$ and $c_0(n) = -\infty$ for $n \geq 1$.

\begin{proposition}\label{prop:homology-of-bar-rat} For an augmented $E_\infty^\rm{u}$-algebra $\bf{R}^+$ in $\Fun(\bb{N},\DQ)$ with $c_0$-connected augmentation ideal and $k \geq 1$, there is an isomorphism
\[S^*(\Sigma^k H_{*,*}^{E_\infty}(\bf{R})) \overset{\cong}\lra H_{*,*}(\Bar^k(\bf{R}^+)).\]
\end{proposition}

\begin{proof}We first observe that $H_{*,*}(\Bar^k(\bf{R}^+))$ is free as a graded-commutative algebra. As it is a $c_0$-connected commutative graded bialgebra by \cref{lem:bialgebra-structure}, this follows from \cite[Theorem 3.8.3, Remark 3.8.2]{Cartier} (observing that $c_0$-connected is sufficient in the bigraded setting). 

We next claim that if $\bf{S}^+$ is an augmented $E_\infty^\rm{u}$-algebra in $\Fun(\bb{N},\DQ)$ so that $H_{*,*}(\bf{S}^+)$ is a free graded-commutative algebra, then we have an isomorphism 
\[S^*(H^{E_\infty}_{*,*}(\bf{S})) \overset{\cong}\lra H_{*,*}(\bf{S}^+).\]
To see this, pick lifts of generators $\{x_i\}_{i \in I}$ of the graded-commutative algebra $H_{*,*}(\bf{S}^+)$ to cycles and use these to construct a map
\[\free_{E^\rm{u}_\infty}(\bb{Q}\{x_i\}_{i \in I}) \lra \bf{S}^+\]
of $E^\rm{u}_\infty$-algebras. This is an equivalence, because by construction it induces an equivalence on homology, using that we are working with rational coefficients. Taking augmentation ideals and considering the induced map on $E_\infty$-homology
\[\bb{Q}\{x_i\}_{i \in I} \cong H_{*,*}^{E_\infty}(\free_{E^\rm{nu}_\infty}(\bb{Q}\{x_i\}_{i \in I}))  \lra H_{*,*}^{E_\infty}(\bf{S})\]
must then also be an isomorphism, implying the claim.

We finally observe that if $\bf{S}^+$ is an augmented $E_\infty^\rm{u}$-algebra in $\Fun(\bb{N},\DQ)$ with $c_0$-connected augmentation ideal, then the $E_\infty$-homology of $\Bar^k(\bf{S}^+)$ is isomorphic to $\Sigma^k H^{E_\infty}_{*,*}(\bf{S})$ by the description of the $E_\infty$-homology in terms of iterated bar constructions \cite[Section 13.7]{GKRW18}.
\end{proof}

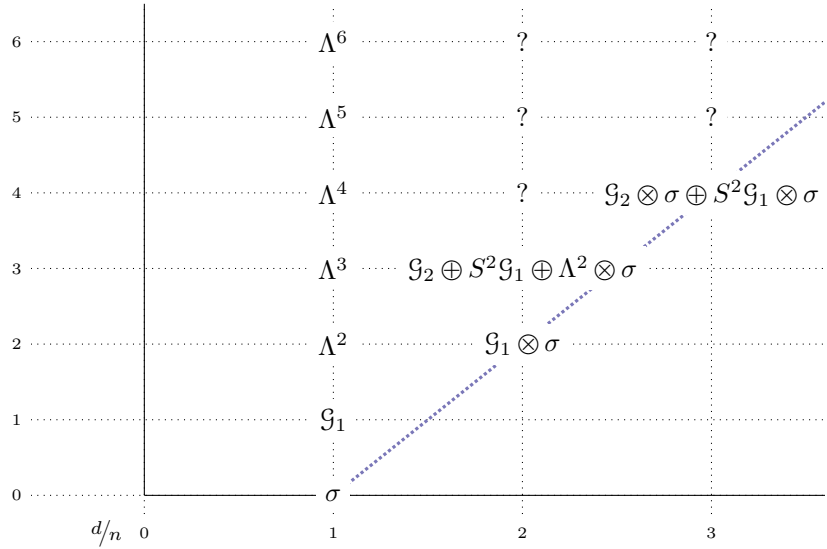
\begin{figure}[ht]
	\begin{tikzpicture}
	\begin{scope}
	\clip (-2,-1) rectangle ({2.5*4-1},6.5);
	\draw (0,0)--(10.5,0);
	\draw (0,0) -- (0,6.5);
	\foreach \s in {0,...,6}
	{
		\draw [dotted] (-1.5,\s)--(10.5,\s);
		\node [fill=white] at (-1.5,\s) [left] {\tiny $\s$};
	}
	\foreach \s in {0,...,4}
	{
		\draw [dotted] ({2.5*\s},-0.5)--({2.5*\s},6.5);
		\node [fill=white] at ({2.5*\s},-.5) {\tiny $\s$};
	}
	\draw [very thick,Periwinkle,densely dotted] (2.5,0) -- (10,6);

	\node [fill=white] at (2.5,0) {$\sigma$};
	\node [fill=white] at (2.5,1) {$\scr{G}_1$};
	\node [fill=white] at (2.5,2) {$\Lambda^2$};
	\node [fill=white] at (2.5,3) {$\Lambda^3$};
	\node [fill=white] at (2.5,4) {$\Lambda^4$};
	\node [fill=white] at (2.5,5) {$\Lambda^5$};
	\node [fill=white] at (2.5,6) {$\Lambda^6$};

	\node [fill=white] at (5,2) {$\scr{G}_1 \otimes \sigma$}; 
	\node [fill=white] at (5,3) {$\scr{G}_2 \oplus S^2 \scr{G}_1 \oplus \Lambda^2 \otimes \sigma$};
	\node [fill=white] at (5,4) {?};
	\node [fill=white] at (5,5) {?};
	\node [fill=white] at (5,6) {?};

	\node [fill=white] at (7.5,4) {$\scr{G}_2 \otimes \sigma \oplus S^2 \scr{G}_1 \otimes \sigma$};
	\node [fill=white] at (7.5,5) {?};
	\node [fill=white] at (7.5,6) {?};
	
	\node at (-.5,-.5) {$\nicefrac{d}{n}$};
	\end{scope}
	\end{tikzpicture}
	\caption{The $E_1$-homology of $\BGLb(F)_\ds{Q}$ is, up to a shift by $-1$, the free graded-commutative algebra on the $E_\infty$-homology, and thus vanishes for $d < 2n-2$ as indicated by the dashed line. We use the abbreviations $\Lambda^k \coloneqq \Lambda^k F^\times_\bb{Q}$ and $\sigma \coloneqq \bb{Q}\{\sigma\}$.}
	\label{fig:gltable1e1}
\end{figure}

\begin{proof}[Proof of \cref{thm:steinberg-homology-improved-vanishing}] Returning to the $E_\infty^\rm{u}$-algebra $\Bar(\BGLb^+_\bb{Q})$, \cref{prop:homology-of-bar-rat} yields an isomorphism
\[S^*(\Sigma H^{E_\infty}_{*,*}(\BGLb_\bb{Q})) \overset{\cong}\lra H_{*,*}(\Bar(\BGLb^+_\bb{Q}))\]
and combining this with \cref{thm:bgl-critical-line-vanishing} and the isomorphism \eqref{eqn:bar-bgl-st} we get that $H_{*,*}(\GL;\SSt)$ is isomorphic to a free graded-commutative algebra with a single generator $\sigma$ in bidegree $(1,1)$ and all remaining generators in bidegrees $(n,d)$ with $d \geq 2n$. Thus it is concentrated in bidegrees $(n,d)$ with $d \geq 2n-1$, and using \cref{thm:ekhomology-steinberg} involves a degree shift by $-1$, $H^{E_1}_{n,d}(\BGLb_\bb{Q})$ is concentrated in degrees $(n,d)$ with $d \geq 2n-2$. (see \cref{fig:gltable1e1}).
\end{proof}

\subsection{Recollections from proof of \cref{thm:ekhomology-steinberg}} We will need some ingredients of the proof of \cref{thm:ekhomology-steinberg} from \cite{GKRW20}. We recall those here.

\subsubsection{Steinberg modules and buildings}\label{sec:buildings} We first recall the connection between Steinberg modules and buildings and justify the Koszulity hypotheses \eqref{hyp:ast-koszul} and \eqref{hyp:st-koszul} for fields. 

For each $k \geq 1$, \cite[Section 5.1]{GKRW20} defines a \emph{$k$-fold building} $D^k(V)$. This is a pointed space with $\GL(V)$-action, arising as the geometric realisation of a $k$-fold pointed simplicial set of ``lattices'' of flags in $V$. We can assemble these to $D^k \in \Fun(\Vect,\Spc_*)$, which satisfy the following property, implicit in \cite[Lemma 1.3]{MPW23}:

\begin{proposition}\label{prop:bar-is-building} There are equivalences in $\Fun(\Vect,\Spc_*)$
\[\Bar^k_\levi(\Bar_{\para}(\ul{\ast}_+)) \simeq D^{k+1} \qquad \text{for all $k \geq 0$}.\]
\end{proposition}

\begin{corollary}\label{cor:koszul-st} Hypotheses \eqref{hyp:ast-koszul} and \eqref{hyp:st-koszul} hold.
\end{corollary}

\begin{proof}By applying $\widetilde{C}_*(-;\bb{Q})$, we reduce hypotheses \eqref{hyp:ast-koszul} and \eqref{hyp:st-koszul} to verifying connectivity results for $D^1$ and $D^2$. Firstly, $D^1(V)$ is isomorphic to a two-fold simplicial suspension of the Tits building $T(V)$ \cite[Lemma 6.1]{GKRW20}, the nerve of the poset of proper nonempty subspaces of $V$ ordered by inclusion. The Solomon--Tits theorem says that $T(V)$ is equivalent to a wedge of $(\dim(V)-2)$-spheres, see \cite[Theorem 2.2]{GKRW20} and the references there. Hence $D^1(V)$ is equivalent to a wedge of $\dim(V)$-spheres and we get
\[\SSt = H_*(\Bar_{\para}(\ul{\bb{Q}})(V)) \cong \widetilde{H}_*(D^1(V);\bb{Q}) = 0 \text{ unless $* = \dim(V)$.}\]
Secondly, there are $\GL(V)$-equivariant pointed \emph{sum maps}
    \[D^k(V) \lra \underbrace{D^1(V) \wedge \cdots \wedge D^1(V)}_k\]
which for $k=2$ yields an equivalence \cite[Proposition 6.3]{GKRW20} \begin{equation}\label{eqn:d2-is-d1-wedge-d1} D^2(V) \xrightarrow{\simeq} D^1(V) \wedge D^1(V).\end{equation} As the latter is equivalent to a wedge of $2\dim(V)$-spheres, we get
\[\SStH = H_*(\Bar_\levi(\SSt)(V)) \cong \widetilde{H}_*(D^2(V);\bb{Q}) = 0 \text{ unless $* = 2\dim(V)$.}\qedhere\]
\end{proof}

\begin{remark}
It is also possible to give an interpretation of the infinite Steinberg modules in terms of buildings. The $k$-fold buildings come with $\GL(V)$-equivariant pointed \emph{suspension maps} \cite[Section 5.1]{GKRW20}
    \[S^1 \wedge D^k \lra D^{k+1}\]
using which we can define a \emph{stable building} $D^\infty \coloneq \{D^k\}_{k \geq 1} \in \Fun(\Vect,\Sp)$ and we have
\[\SStL = \colim_{k \to \infty} H_{*-k+1}(\Bar^k_{\levi}(\Bar_{\para}(\ul{\ast}_+));\ds{Q}) \cong H_*(D^\infty;\ds{Q}).\]
By \cite[Proposition 5.3]{MPW23} the spectrum $D^\infty(V)$ is in fact equivalent to the suspension spectrum of Rognes' common basis complex \cite[Definition 14.5]{Rognes} and this is in turn equivalent to Br\"uck--Piterman--Welker's partial decomposition poset \cite[Section 1]{BPW} by Corollary 1.1 loc.cit..
\end{remark}

\begin{proof}[Proof of \cref{prop:bar-is-building}] The proof requires a generalisation of the $k$-fold buildings due to Miller--Patzt--Wilson \cite[Definition 3.5]{MPW23}: for each $k,\ell \geq 0$ they define for a finite-dimensional vector space $V$ a pointed space $D^{k,\ell}(V)$ with $\GL(V)$-action, arising as the thick geometric realisation of a $(k+\ell)$-fold pointed semisimplicial set with $k$ directions of flags and $\ell$ directions of splittings, all satisfying the common basis property. For $k=0=\ell$ this is simply the constant functor $\ul{\ast}_+ \in \Fun(\Vect,\Spc_*)$.

The latter is an $E_1^\rm{u}$-algebra with respect to $\para$, and admits a ``canonical'' augmentation by the map that is the identity in dimension $0$ and the map to the basepoint otherwise, and we can form $\Bar_{\para}(\ul{\ast}_+)$ with respect to this augmentation. We claim that there is an equivalence in $\Fun(\Vect,\Spc_*)$
\[\Bar_{\para}(\ul{\ast}_+) \simeq D^1.\]
Using the standard simplicial model for the bar construction, we see that $\rm{Bar}_{\para}(\ul{\ast}_+)$ is given by the geometric realisation of
\[[p] \longmapsto 1 \para (\ul{\ast}_+)^{\para p} \para 1\]
where $1$ is the monoidal unit of $\Fun(\Vect,\Spc_*)$. The inner face maps are induced by the multiplication of $\ul{\ast}$ and outer face maps are induced by the canonical augmentation. Using the formula for Day convolution and the description of the profunctor $\obackslash_k$, this is given
\[(1 \para (\ul{\ast}_+)^{\para p} \para 1)(V) \simeq \frac{\left\{\parbox[c]{6cm}{\centering flags $0 = V_0 \subseteq V_1 \subseteq \cdots \subseteq V_{p+1} = V$}\right\}}{\left\{\parbox[c]{6cm}{\centering flags $0 = V_0 \subseteq V_1 \subseteq \cdots \subseteq V_{p+1} = V$ with $V_1/V_0 \neq 0$ or $V_{p+1}/V_p \neq 0$}\right\}}\]
where now all face maps forget a term in a flag. Comparing to \cite[Definition 5.4]{GKRW20} we recognise this as $D^1(V)$, and this equivalence is evidently $\GL$-equivariant.

There is a remaining $E_\infty$-algebra structure on the bar construction and its product \emph{is} simplicial for the standard simplicial model for the bar construction we used above, given by levelwise direct sum. That is, it corresponds to the product given by the case $k=1$ of the pointed $\rm{GL}(V)\times \rm{GL}(V')$-equivariant \emph{sum maps} \cite[Section 6.3]{GKRW20}
\[D^1(V) \wedge D^1(V') \lra D^1(V \oplus V').\]
Arguing as above using the tensor product $\levi$ instead of $\para$, one proves 
\[\Bar_\levi(\Bar_{\para}(\ul{\ast}_+)) \simeq \Bar_\levi(D^1) \simeq D^{1,1}.\]
There are pointed $\GL(V)$-equivariant forgetful maps $D^{1,1}(V) \lra D^{2,0}(V)$ turning a direct sum decomposition direction into a flag one and the corresponding map
\[D^{1,1} \lra D^{2,0}\]
is an equivalence by \cite[Theorem 3.15]{MPW23}. Iterating this argument, one proves that
\[\Bar^k_\levi(\Bar_{\para}(\ul{\ast}_+)) \simeq \Bar_\levi(D^{k,0}) \simeq D^{k,1} \simeq D^{k+1,0}.\qedhere\]
\end{proof}

\begin{remark} Bar--cobar duality also yields an $E_1$-coalgebra structure on $D^1$, which is more difficult to determine than the $E_1$-algebra structure as it does \emph{not} arise from simplicial maps; see \cref{rem:coprod-ad-hoc}.
\end{remark}

\subsubsection{The Nesterenko--Suslin property} \label{sec:nesterenko-suslin} We next recall a result comparing the homology of general linear groups preserving a splitting with that of groups preserving a flag, going back to Nesterenko and Suslin \cite[Section 1]{NesterenkoSuslin}.

For $V = U_1 \oplus U_2$, let $\rm{P}_{U_1}(V) \subseteq \GL(V)$ denote the parabolic subgroup preserving the flag $U_1 \subseteq V$. The restriction map $\rm{P}_{U_1}(V) \to \GL(U_1) \times \GL(U_2)$ admits a section and we have \cite[Theorem 1.1]{NesterenkoSuslin} \cite[Definition 5.12]{GKRW20}:

\begin{lemma}[Nesterenko--Suslin] \label{lem:nesterenko-suslin} For all nonzero $U_1,U_2$ both of the maps 
\[\GL(U_1) \times \GL(U_2) \to G_{U_1,U_2} \to \GL(U_1) \times \GL(U_2)\] induce an isomorphism on $H_*(-;\ds{Q})$.
\end{lemma} 

This allows us to identify the tensor products $\levi$ and $\para$ after applying $\dim_!$. Recall that the projective general linear groups $\PGL(V)$ are defined as the quotient of $\GL(V)$ by the subgroup $F^\times$ of linear automorphisms given by scaling, so in particular we can consider a $\ds{Q}[\PGL(V)]$-module as a $\ds{Q}[\GL(V)]$-module by restriction along the quotient homomorphism $\GL(V) \to \PGL(V)$. The image of $\rm{P}_{U_1}(V)$ in the projective group is denoted $\rm{P}(\rm{P}_{U_1}(V)) \subseteq \PGL(V)$. The following is deduced from \cref{lem:nesterenko-suslin} using Hochschild--Serre spectral sequences (cf.~the proof of \cite[Theorem 1.1]{NesterenkoSuslin}:

\begin{lemma}\label{lem:para-levi-comparison-pgl} If $M_1$ is a $\ds{Q}[\PGL(U_1)]$-module and $M_2$ is a $\ds{Q}[\PGL(U_2)]$-module, the following maps are isomorphisms
\begin{align*}H_*(\GL(U_1) \times \GL(U_2);M_1 \otimes M_2) &\overset{\cong}\lra H_*(\rm{P}_{U_1}(V);M_1 \otimes M_2) \\
H_*(\rm{P}(\GL(U_1) \times \GL(U_2));M_1 \otimes M_2) &\overset{\cong}\lra H_*(\rm{P}(\rm{P}_{U_1}(V));M_1 \otimes M_2).\end{align*}
\end{lemma}

Combining this with the formulas from \cref{exam:levi-para-explicit-formulas} and Shapiro's lemma, we deduce:

\begin{lemma}\label{lem:dim-levi-para} If $M,N \in \Fun(\Vect,\DQ)$ are objectwise pulled back from the projective general linear groups, then the following induces an isomorphism on homology
\[\dim_!(M \levi N) \lra \dim_!(M \para N).\]
\end{lemma}

\subsubsection{$E_k$-homology and split buildings} In \cref{sec:buildings} we saw that Steinberg modules and their variants arise from the homology of buildings. These are obtained from collections of flags, and in this section we will consider versions obtained from splittings, i.e.~direct sum decompositions.

The \emph{split buildings} $\widetilde{D}^k(V) \in \Fun(B\GL(V);\Spc_*)$ were defined in \cite[Section 5.2]{GKRW20} as the thick geometric realisation of a $k$-fold pointed semisimplicial set $\tilde{D}^k_{\bullet,\ldots,\bullet}(V)$ of $k$-dimensional grids of direct sum decompositions of $V$; in the proof of \cref{prop:bar-is-building} these were denoted $D^{0,k}(V)$. We can assemble these to $\widetilde{D}^k \in \Fun(\Vect,\Spc_*)$ and this object arises in a similar manner as the nonsplit buildings \cite[Theorem 5.20]{GKRW20}:

\begin{lemma}\label{lem:bar-is-split-building} There are equivalences in $\Fun(\Vect,\Spc_*)$
    \[\Bar_\levi^k(\ul{\ast}_+) \simeq \widetilde{D}^k \qquad \text{for all $k \geq 1$.}\]
\end{lemma}

In particular, it has a preferred lift to an object of $\coAlg_{E_k^\rm{u}}(\Alg_{E_\infty^\rm{u}}(\Fun(\Vect,\Spc_*)))$, where both the coalgebra and algebra structures are with respect to $\levi$.  Sending direct sum decompositions to flags induces \emph{forgetful maps} \cite[(5.6)]{GKRW20} in $\Fun(\Vect,\Spc_*)$
\[\widetilde{D}^k \lra D^k.\]
These are not equivalences, but do induce isomorphisms on homology after applying $\dim_!$ using \cref{lem:nesterenko-suslin}, see \cite[Theorem 5.18]{GKRW20}. We now use this to recall why $E_k$-homology is computed by homology with coefficients in Steinberg modules:

\begin{proof}[Proof of \cref{thm:ekhomology-steinberg}] We will establish isomorphisms
\[H_{*,*}(\Sigma^k \indec_{E_k}(\BGLb_\ds{Q})^+) \overset{\cong}\lra H_{*,*}(\Bar^k(\BGLb^+_\ds{Q})) \overset{\cong}{\lra} \widetilde{H}_{*,*}(\dim_! D^{k+1};\ds{Q})\]
Unwinding the definition of the right side gives the result for $k=1,2$, and for $k=\infty$ one uses iterated bar spectral sequences to identify it with the quotient of the case $k=2$ by the decomposables, as in the proof of \cite[Corollary 6.12]{GKRW20}. The left map is an instance of \cref{thm:indec-is-bar}, so it remains to establish that the right map is. This follows by considering
\[\Bar^k(\dim_!(\ul{\ast}_+)) \simeq \dim_!(\Bar^k_\levi(\ul{\ast}_+)) \simeq \dim_! \widetilde{D}^k \lra \dim_! D^k\]
and applying $\widetilde{C}_*(-;\ds{Q})$. Here the left-most map is an equivalence as $\dim_!$ is a symmetric monoidal left adjoint, the middle map is an equivalence by \cref{lem:bar-is-split-building}, and the right map becomes an isomorphism as a consequence of \cref{lem:nesterenko-suslin}.\end{proof}

We need to understand how the previous proof interacts with coalgebraic structures. 

\begin{lemma}\label{lem:gl-st-coproduct} The isomorphism 
\[H_{*,*}(\Bar(\BGLb^+_\bb{Q})) \cong H_{*,*}(\GL;\SSt)\]    
is one of coalgebras, where the coproduct on the left is induced by the bar construction and the one on the right is induced by that on $\SSt$.
\end{lemma}

\begin{proof}Consider the commutative diagram of promonoidal categories
\[\begin{tikzcd} (\Vect,\oplus) \arrow{rr}{\id} \arrow{rd}[swap]{\dim} & & (\Vect,\obackslash) \arrow{ld}{\dim} \\[-5pt]
& (\bb{N},+) & \end{tikzcd}\]
where the diagonal maps are monoidal and the horizontal map is lax monoidal. Applying the Day convolution construction, we obtain a commutative diagram of monoidal categories
\[\begin{tikzcd} (\Fun(\Vect,\Spc_*),\levi) \arrow{rr}{\id_! = \id} \arrow{rd}[swap]{\dim_!} & & (\Fun(\Vect,\Spc_*),\obackslash) \arrow{ld}{\dim_!} \\[-5pt]
& (\Fun(\bb{N},\Spc_*),\levi) & \end{tikzcd}\]
where the left diagonal map is monoidal and the other maps are oplax monoidal. It thus suffices to show that there is a map
\[\Bar_{\levi}(\ul{\ast}_+) \lra \Bar_{\para}(\ul{\ast}_+).\]
of $E^u_1$-coalgebras and identify it on underlying objects with the map $\widetilde{D}^1 \to D^1$.

Given an adjunction $L \dashv R$ where $R \colon \scr{D} \to \scr{C}$ is lax monoidal, $L \colon \scr{C} \to \scr{D}$ is oplax monoidal using the mate correspondence, and $\scr{C}$ and $\scr{D}$ are sufficiently nice, \cref{sec:bar-cobar-naturality-general} provides a natural transformation of functors $L^\coAlg \Bar_\scr{D} R^\Alg \Rightarrow \Bar_\scr{C}$. On underlying objects, for an augmented algebra $\bf{A}$ this is the map induced by the maps $L(R(\bf{A})^{\otimes p}) \to \bf{A}^{\otimes p}$ of $p$-simplices given by the oplax monoidality of $L$ and the counit. When we apply this result to $L = \id_! = R$, we get the desired map: levelwise it is simply the oplax monoidality of $\id_!$ and recalling equivalences between bar constructions and building from \cref{prop:bar-is-building} and \cref{lem:bar-is-split-building}, we recognise this as the desired map.\end{proof}

\begin{remark}\label{rem:coprod-ad-hoc} Let us outline a more direct argument for \cref{lem:gl-st-coproduct} avoiding \cite{BlansBlomKupers}. \cref{exam:e1-coproduct} describes the homotopy class of the coproduct on the bar construction in terms of relative tensor products as
\[1 \otimes_\bf{A} 1 \overset{\simeq}\longleftarrow 1 \otimes_\bf{A} \bf{A} \otimes_\bf{A} 1 \lra 1 \otimes_\bf{A} 1 \otimes_\bf{A} 1.\]
In functor categories with values in pointed spaces the left equivalence is modelled by the pointed homeomorphism $|\rm{Bar}_\bullet(1,\bf{A},1)| \cong |\rm{esd}\,\rm{Bar}_\bullet(1,\bf{A},1)|$ (here geometric realisations are taken in pointed spaces) using B\"okstedt--Hsiang--Madsen's variant \cite[Section 1]{BHM} of the edgewise subdivision construction of Segal \cite[Appendix 1]{SegalConfiguration} and the right map by the map on $p$-simplices
\[1 \otimes \bf{A}^{\otimes p} \otimes \bf{A} \otimes \bf{A}^{\otimes p} \otimes 1 \lra 1 \otimes \bf{A}^{\otimes p} \otimes 1 \otimes \bf{A}^{\otimes p} \otimes 1\]
given by the augmentation on the middle term and identity elsewhere. Passing to the case of interest, we need to prove that the following square of pointed sets commutes for each $p$
\[\begin{tikzcd}\frac{\left\{\parbox[c]{5.5cm}{\centering splittings $U_0 \oplus \cdots \oplus U_{2p+1} = V$}\right\}}{\left\{\parbox[c]{5.5cm}{\centering splittings $U_0 \oplus \cdots \oplus U_{2p+1} = V$ with $U_0 \neq 0$ or $U_{2p+1} \neq 0$}\right\}} \rar \dar & 
\frac{\left\{\parbox[c]{5.5cm}{\centering splittings $U_0 \oplus \cdots \oplus U_{2p+1} = V$}\right\}}{\left\{\parbox[c]{5.5cm}{\centering splittings $U_0 \oplus \cdots \oplus U_{2p+1} = V$ with $U_0 \neq 0$, $U_p \neq 0$, or $U_{2p+1} \neq 0$}\right\}} \dar \\
\frac{\left\{\parbox[c]{5.8cm}{\centering flags $0 = V_{-1} \subseteq V_0 \subseteq \cdots \subseteq V_{2p+1} = V$}\right\}}{\left\{\parbox[c]{5.8cm}{\centering flags $0 = V_{-1} \subseteq V_1 \subseteq \cdots \subseteq V_{2p+1} = V$ with $V_0/V_{-1} \neq 0$ or $V_{2p+1}/V_{2p} \neq 0$}\right\}} \rar & 
\frac{\left\{\parbox[c]{5.8cm}{\centering flags $0 = V_{-1} \subseteq V_1 \subseteq \cdots \subseteq V_{2p+1} = V$}\right\}}{\left\{\parbox[c]{5.8cm}{\centering flags $0 = V_{-1} \subseteq V_1 \subseteq \cdots \subseteq V_{2p+1} = V$ with $V_0/V_{-1} \neq 0$, $V_p/V_{p-1} \neq 0$, or $V_{2p+1}/V_{2p} \neq 0$}\right\}}\end{tikzcd}\]
where the vertical maps sum up splittings to flags and the horizontal maps are projections; it visibly does. This should agree up to a suspension with the coproduct of Campbell--Zakharevich \cite[Section 2.2]{CampbellZakharevich} (see also \cite[Section 4.2]{KKMMW}) and, after passing to homology, that of Brown--Chan--Galatius--Payne \cite[Section 3.1]{BCGP}.
\end{remark}

\subsection{The cobracket on the $E_\infty$-homology of $\BGLb_\ds{Q}$}
As we have discussed, by Koszul duality, the indecomposables $\Sigma \,\cot_{E^\rm{nu}_\infty}(\BGLb_\ds{Q})$ admit the structure of a Lie coalgebra in $\Fun(\bb{N},\DQ)$.  It is a crucial point that this structure is \emph{not} induced by the corresponding structure on $\SStL$ under the isomorphism of \cref{thm:ekhomology-steinberg}. This can be seen by a simple comparison of degrees: the Lie coalgebra structure on $\SStL$ translates to a Lie coalgebra structure whose \emph{wrong} cobracket on $H_{*,*}(\GL;\SStL)$ has degree $0$, but the \emph{correct} cobracket has degree $1$. 

The explanation is that $\SStL$ is not only a Lie coalgebra with respect to $\levi$, but also has a remaining coproduct with respect to $\para$ which is compatible with the cobracket. The two tensor products $\levi$ and $\para$ become identified upon applying $\dim_!$ by \cref{lem:dim-levi-para} and via the same phenomenon as the additivity theorem or the Eckmann--Hilton argument, this leads to a vanishing result for the induced cobracket and a secondary cobracket of one degree higher. Our goal in this section is to give a strategy for the computation of this correct secondary cobracket using the lift of the original coproduct to the cobar complexes in \cref{prop:coproduct-on-cobar}, which will be implemented in \cref{sec:cobracket}. 

\medskip

Our starting point is to recall that since
\[\rm{H}^+ \coloneq H_{*,*}(\GL;\SSt)\]
arises as the homology of the bar construction $\Bar(\BGLb_\bb{Q}^+)$, it comes with the structure of an augmented commutative bialgebra; the coproduct here arises from that on the bar construction but agrees with that induced by the coproduct on Steinberg modules by \cref{lem:gl-st-coproduct}. The bialgebra $\rm{H}^+$ is concentrated in bidegrees $(n,d)$ for $d \geq 2n-2$ by \cref{thm:steinberg-homology-improved-vanishing} and it is related to the Lie coalgebra $H_{*,*}(\GL;\Sigma^{-1} \SStL)$ as follows:

\begin{proposition}\label{prop:cobracket-via-coproduct} There is an isomorphism
\[H_{*,*}(\GL;\Sigma^{-1} \SStL) \cong \rm{H}/\rm{H}^2\]
where $\rm{H} \coloneq I(\rm{H}^+) \subset \rm{H}^+$ is the augmentation ideal. Under this isomorphism, the cobracket is induced by the antisymmetrisation $\ol{\Delta}-\sigma \circ \ol{\Delta}$ of the reduced coproduct. \end{proposition}

\begin{proof}The statement about the isomorphism follows from \cref{prop:homology-of-bar-rat}. For the statement about the cobracket, we first observe that the inclusion $\bb{R} \to \bb{R}^\infty$ induces a map $i \colon E_1^\rm{nu} \to E_\infty^\rm{nu}$ of operads, with Koszul dual map $Bi \colon BE_1^\rm{nu} \to BE_\infty^\rm{nu}$ of cooperads. There is then a Beck--Chevalley transformation
\[(Bi)_* \indec^\nil_{E_1^\rm{nu}} i^* \lra \indec^\nil_{E_\infty^\rm{nu}}.\]
Working rationally, the map $i$ is equivalent to the map $\rm{As}^\rm{nu} \to \rm{Com}^\rm{nu}$ forgetting the product is commutative, and its Koszul dual is equivalent to the operadic suspension of the map $\rm{coAs}^\rm{nu} \to \coLie$ which considers a coassociative coalgebra with (necessarily reduced) coproduct $\ol{\Delta}$ as a Lie coalgebra with cobracket $\delta = \ol{\Delta}-\sigma \circ \ol{\Delta}$ given by the antisymmetrisation of the coproduct. Passing to homology, we see that the (shifted) cobracket on $E_\infty$-homology is induced by antisymmetrisation of the (shifted) coproduct on $E_1$-homology.\end{proof}

\begin{remark}Since we have $\ol{\Delta}-\sigma \circ \ol{\Delta} = \Delta-\sigma \circ \Delta$ on $\rm{H}^+$, one can also use the antisymmetrisation of the coproduct itself.\end{remark}

More precisely, \cref{prop:cobracket-via-coproduct} says that the cobracket is given by the dashed map in
\[\begin{tikzcd} \rm{H} \dar[swap]{\overline{\Delta}-\sigma \circ \overline{\Delta}} \rar[two heads] & \rm{H}/\rm{H}^2 \dar[dashed]{\delta} \\[-5pt]
\rm{H} \otimes \rm{H} \rar[two heads] & \rm{H}/\rm{H}^2 \otimes \rm{H}/\rm{H}^2 & \end{tikzcd}\]
obtained by lifting along the top surjective map and applying the left-bottom composition; the result is independent of the choice of lift using
\[\ol{\Delta}(xy) \equiv x \otimes y+(-1)^{|x||y|} y \otimes x \pmod{ \rm{H} \otimes \rm{H}^2+\rm{H}^2 \otimes \rm{H}}.\]
Thus our first task is to understand how to compute this coproduct. By \cref{lem:gl-st-coproduct} it is induced by the coproduct $\Delta$ on $\SSt$ with respect to the tensor product $\para$. \cref{prop:coproduct-on-cobar} provides a lift of the Koszul duality equivalence
\[\Sigma^{-1}\Omega^\coLie(\SStL) \overset{\simeq}\lra \SSt_{>0}\]
to one of counital coassociative coalgebras with respect to the tensor product $\para$. We will explain how to obtain the quotients $\rm{H}/\rm{H}^2$ from this perspective.

Our starting point is a concrete instance of the canonical multiplicative filtration from \cite[Section 5.4]{GKRW18}, obtained as the functor $(-1)^\alg_!$ on algebras obtained from $(-1)_! \colon \scr{C} \to \Fun(\bb{Z}_\leq,\scr{C})$ by the construction of \cref{sec:alg-nat-left-adjoint}. Rectifying using \cref{sec:rect-dg} we work in the category $\Fun(\rm{A},\rm{Ch}_\bb{Q})$ for a symmetric monoidal $1$-category $\rm{C}$. 

Recall from \cref{def:cobar-colie} that the cobar construction of a Lie dg-coalgebra $\bf{L}$ in $\Fun(\rm{A},\rm{Ch}_\bb{Q})$ has underlying object given by $\Sigma(\rm{Com}^\rm{nu} \circ \Sigma^{-1} \bf{L})$, where $\circ$ denotes the composition product of symmetric sequences. This can be lifted to a filtered object
\[\rm{fil}_\rm{can} \Omega^\coLie(\ul{L}) \coloneq \Sigma(\rm{Com}^\rm{nu} \circ \Sigma^{-1} (-1)_! \bf{L})\]
by putting $\bf{L}$ in filtration degree $-1$; when $\bf{L}= B^\rm{Com}(\bf{R})$ for a nonunital dg-commutative algebra $\bf{R}$ this recovers the canonical multiplicative filtration on $\bf{R}$. Moreover, in the pointed setting the functor $a_!$ for $a \in \bb{Z}$ admits a further left adjoint $a^\dagger$ \cite[Section 5.2.2]{GKRW18} (called $a^!$ there), given by sending a filtered object $X$ to $\colim(X)/X(a-1)$. The natural transformation $(n+1)_! \to n_!$ induces natural transformations $n^\dagger \to (n+1)^\dagger$ and this allows to extract from the canonical filtration a second filtered object
\[\rm{fil}_\rm{pow} \Omega^\coLie(\bf{L}) \coloneq \big(a \mapsto (-a)^\dagger(\rm{fil}_\rm{can} \Omega^\coLie(\bf))\big).\]
We call this the \emph{power filtration}, justified by the following result in the special case $\rm{C} = \bb{N}$:

\begin{lemma}\label{lem:fil-pow-free} Suppose that $\bf{R}^+$ is an $E^\rm{u}_\infty$-algebra in $\Fun(\bb{N},\DQ)$ so that $\rm{H}^+ \coloneq H_{*,*}(\bf{R}^+)$ is free as a graded-commutative algebra and has a $c_0$-connected augmentation ideal $H$. Then there is an isomorphism
\[H_*(\rm{fil}_\rm{pow}(\Omega^{\coLie}(B^\rm{Com} \bf{R}))) \cong \rm{H}/\rm{H}^{\bullet},\]
where the right side is the filtered object obtained from $\bf{H}$ by taking quotients by powers.
\end{lemma}

\begin{proof}By the proof of \cref{prop:homology-of-bar-rat} there is an equivalence 
\[S^*(H^{E_\infty}_{*,*}(\bf{R})) \overset{\simeq}\lra \bf{R}^+,\]
so without loss of generality we may replace $\bf{R}^+$ with the left side. Since this is free,  $B^\rm{Com}(S^*(H^{E_\infty}_{*,*}(\bf{R})))$ is equivalent to the trivial Lie coalgebra $\rm{cotriv}_\rm{coLie}(H^{E_\infty}_{*,*}(\bf{R}))$, and the canonical filtration on its cobar construction is split, isomorphic to the free nonunital graded-commutative algebra on $H^{E_\infty}_{*,*}(\bf{R})$ with filtration degree corresponding to the number of products.
\end{proof}

Returning to the case at hand, we take $L = \SStL$ in $\Fun(\Vect,\rm{GrMod}_\bb{Q}) \subset \Fun(\Vect,\rm{Ch}_\bb{Q})$. Then the coLie cobar construction $\Omega^\coLie(\SStL)$ takes the form
\[\SStL \to \Lambda^2 \SStL \to \Lambda^3 \SStL \to \cdots,\]
where on $F^n \in \Vect$ the left-most term is given by $\StL_n$ in degree $2n$. The power filtration on this is given in filtration degree $-r$ for $r \geq 0$ by having only the first $r$ terms starting from the left. Applying the construction of \cref{prop:coproduct-on-cobar} in filtered chain complexes lifts it to a filtered map
\[\fil_\rm{can} \overline{\Delta} \colon \fil_\rm{can}(\Sigma^{-1}\Omega^\coLie(\SStL)) \lra \fil_\rm{can}(\Sigma^{-1}\Omega^\coLie(\SStL)) \para \fil_\rm{can}(\Sigma^{-1}\Omega^\coLie(\SStL)).\]
Plotting the power filtration horizontally, we get maps of chain complexes
\[\begin{tikzcd}\big[\Sigma^{-1}\Omega^\coLie(\SStL)\big] \dar \rar &[-10pt] \cdots \rar &[-10pt] \big[\SStL \to \Lambda^2 \SStL\big] \rar\dar &[-10pt]\SStL \dar \\[-5pt]
\big[\Sigma^{-1}\Omega^\coLie(\SStL) \para \Sigma^{-1}\Omega^\coLie(\SStL)\big] \rar & \cdots \rar & \big[0 \to \SStL \para \SStL] \rar & 0 \end{tikzcd}\]
where the left-most vertical map is equivalent to the reduced coproduct $\SSt_{>0} \to \SSt_{>0} \para \SSt_{>0}$. Applying $\dim_!$, which as symmetric monoidal left adjoint commutes with the constructions, and taking homology, we obtain from this by \cref{lem:fil-pow-free} a commutative diagram
\[\begin{tikzcd} \rm{H} \rar[two heads] \dar{\ol{\Delta}} & \cdots \rar[two heads] \dar & \rm{H}/\rm{H}^3 \rar[two heads] \dar & \rm{H}/\rm{H}^2 \dar \\[-5pt]
\rm{H} \otimes \rm{H} \rar & \cdots \rar & \rm{H}/\rm{H}^2 \otimes \rm{H}/\rm{H}^2 \rar & 0.\end{tikzcd}\]
This conclusion is the following procedure, the computational implementation of which appears in \cref{sec:cobracket}.

\begin{proposition}\label{prop:construction-of-cobracket} The cobracket is obtained by 
\begin{enumerate}[(1)]
\item applying $\dim_!$ to the zigzag of chain complexes
\[\SStL \longleftarrow [\SStL \to \Lambda^2 \SStL] \longrightarrow [0 \to \SStL \para \SStL]\]
where the left map is projection to $\SStL$ and the right map is projection to $\Lambda^2 \SStL$ followed by an instance of $\zeta^\alt$ from \eqref{eqn:zeta-alt}, and
\item antisymmetrising.
\end{enumerate}
\end{proposition}

\section{The cobracket of $\scr{G}(F)$} \label{sec:cobracket} In this section we describe how to obtain a presentation of $\scr{G}_n(F)$ close to that in \cref{thm:polyl-presentation-additive} (it has the same generators but the relations are phrased differently) and compute the cobracket
\[\delta \colon \scr{G}(F) \lra \Lambda^2 \scr{G}(F) \qquad \text{where $\scr{G}_n(F) = H_1(\GL_n(F);\StL_n)$}\]
to obtain \cref{thm:polyl-presentation-cobracket}. We follow \cref{conv:shorter-notation}.

\subsection{Resolving infinite Steinberg modules by formal correlators} We will use formal analogues of the generators and relations in \cref{sec:pres-form-stl} and decomposition operators \eqref{eqn:decomposition-operator} to construct resolutions of the infinite Steinberg modules $\StL(V)$ as in \cref{sec:intro-presentation-resolution}.

\subsubsection{Formal correlators} Our starting point is a Lie coalgebra $\FC \in \Fun(\Vect,\rm{GrMod}) \subset \Fun(\Vect,\DQ)$ of ``formal'' analogues of the Steinberg correlators. Its value on a vector space $V$ of positive dimension $n$ will be generated by symbols $\FC[u_0:\compactcdots:u_n]$ that we call \emph{formal correlators} for affine bases $u_0,\ldots,u_n$, satisfying the analogues of \eqref{enum:stl-relations-i}--\eqref{enum:stl-relations-iii}:
\begin{enumerate}[(1)]
\item \label{enum:fc-relations-i} They are homogeneous: $\rm{FC}[u_0:\compactcdots:u_n] = \rm{FC}[u_0-u:\compactcdots:u_n-u]$ for any $u \in V$.
\item \label{enum:fc-relations-ii} They are cyclically symmetric: $\rm{FC}[u_0:u_1:\compactcdots:u_n] = \rm{FC}[u_1:u_2:\compactcdots:u_0]$.
\item \label{enum:fc-relations-iii} They satisfy the shuffle relations: 
\[\sum_{\sigma \in \rm{Sh}(n_1,n_2)} \rm{FC}[u_0:u_{\sigma(1)}:\compactcdots:u_{\sigma(n_1+n_2)}] =0 \quad \text{for $n=n_1+n_2$ with $n_1,n_2>0$.}\]
\end{enumerate}

\begin{definition}If $V$ is of dimension $n$, then we define
\[\FC(V) \coloneq \begin{cases} \displaystyle \frac{\bb{Q}\{\FC[u_0:\compactcdots:u_n] \text{ for affine bases $u_0,\compactldots,u_n$}\}}{\text{\eqref{enum:fc-relations-i}--\eqref{enum:fc-relations-iii}}} & \text{if $n>0$}, \\
0 & \text{if $n=0$}.\end{cases}\]
\end{definition}

The action of $\GL(V)$ on affine bases induces a $\GL(V)$-action on $\FC(V)$, and we can assemble all to an object $\FC \in \Fun(\Vect,\rm{GrMod}_\bb{Q}) \subset \Fun(\Vect,\DQ)$. We will endow this with a Lie coalgebra structure by repeating the formula \eqref{eqn:stinfty-cobracket} of \cref{prop:stl-explicit-cobracket}, cf.~\cite[p.~436]{Goncharov01} and \cite[Lemma 2, Proposition 15]{CMRR24}.

\begin{definition}The \emph{cyclic cobracket} $\delta_\rm{cyc} \colon \rm{FC} \to \Lambda^2 \rm{FC}$ is given by the formula
\[\delta_\rm{cyc}(\FC[u_0:\compactcdots:u_n]) = \sum_{j=0}^n \sum_{i=1}^{n-1} \FC[u_j:u_{j+1}:\compactcdots:u_{j+i}] \wedge \FC[u_j:u_{j+i+1}:\compactcdots:u_{j+n}]\]
where the indices are to be interpreted cyclically.
\end{definition}

\begin{remark}\label{rem:delta-cyc-from-sder} One can use \cref{sec:lie-alg-trees-derivations} to justify that $\delta_\rm{cyc}$ is a well-defined Lie cobracket (compare to \cref{rem:c-sder-cobracket}): that is, $\delta_\rm{cyc}$ is compatible with the relations \eqref{enum:fc-relations-i}--\eqref{enum:fc-relations-iii}, and satisfies the anti-symmetry and co-Jacobi relations. It suffices to verify these properties on a generator $\FC[u_0 : \compactcdots : u_n]$. On such a generator it agrees---after applying the functor $t_!$ induced by $t \colon \Vect \to \ast$, which sums up all values on all vector spaces and hence sends no elements to zero---with the same-named Lie cobracket from \cref{sec:sder-dual} on the generator $\rm{LC}[X_0: \compactcdots: X_n] \in \fr{sder}_S^\vee$ where $S = [n] = \{0,1,\ldots,n\}$ and this has the aforementioned properties.\end{remark}

The vanishing result of \cref{lem:stl-cobracket-vanishing} holds also for this cobracket. Note that it implies the corresponding result in $\SStL$, since the map $\FC \to \SStL$ is a morphism of Lie coalgebras, where $\zeta^\alt$ is $x \wedge y \mapsto \frac{1}{2}(x \para \ol{y}-y \para \ol{x})$ as in \eqref{eqn:zeta-alt}:

\begin{lemma}\label{lem:fc-cobracket-symmetry} The cobracket on $\FC$ has the following vanishing property:
\[\zeta^\alt \circ \delta_\rm{cyc} = 0.\]
\end{lemma}

\begin{proof}The terms cancel pairwise: given a marked cut $i,(j+i,j+i+1)$ where the marked polygon is counterclockwise from the cut vertex, we can construct another marked cut $j+i,(i-1,i)$ where the marked polygon is clockwise from the cut vertex, giving a pairing between the two sets. For example, the following marked cuts are paired where the first is $j=2$, $i=2$ (markings are indicated by shading):
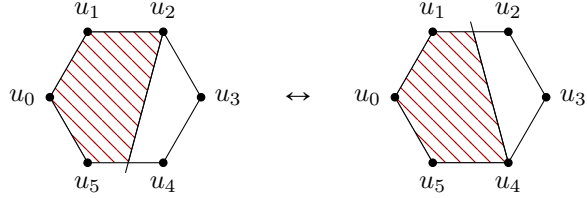
\[\begin{tikzpicture}[baseline={([yshift=-.5ex]current bounding box.center)}]
   \begin{scope}
       \clip (60:1)--(120:1)--(180:1)--(240:1)--(300:1)--(360:1)--cycle;
       \draw[pattern={Lines[angle=-45,distance=4pt]},pattern color=Mahogany] (60:1)--(270:1)--(240:1) -- (180:1) -- (120:1) -- cycle;
   \end{scope}
   \draw (0:1) \foreach \x in {60,120,...,360} {  -- (\x:1) };
   \foreach \x/\l/\p in
     { 60/{$u_2$}/above,
      120/{$u_1$}/above,
      180/{$u_0$}/left,
      240/{$u_5$}/below,
      300/{$u_4$}/below,
      360/{$u_3$}/right
     }
     \node[inner sep=1pt,circle,draw,fill,label={\p:\l}] at (\x:1) {};
     \draw (60:1) -- (270:1);
\end{tikzpicture}
\quad \leftrightarrow \quad
\begin{tikzpicture}[baseline={([yshift=-.5ex]current bounding box.center)}]
    \begin{scope}
       \clip (60:1)--(120:1)--(180:1)--(240:1)--(300:1)--(360:1)--cycle;
       \draw[pattern={Lines[angle=-45,distance=4pt]},pattern color=Mahogany] (90:1)--(300:1) -- (240:1) -- (180:1) -- (120:1) -- cycle;
   \end{scope}
   \draw (0:1) \foreach \x in {60,120,...,360} {  -- (\x:1) };
   \foreach \x/\l/\p in
     { 60/{$u_2$}/above,
      120/{$u_1$}/above,
      180/{$u_0$}/left,
      240/{$u_5$}/below,
      300/{$u_4$}/below,
      360/{$u_3$}/right
     }
     \node[inner sep=1pt,circle,draw,fill,label={\p:\l}] at (\x:1) {};
     \draw (300:1) -- (90:1);
\end{tikzpicture}\]
By construction, these contribute (up to a factor of $\frac{1}{2}$)
\begin{align*}&\rm{FC}[u_i:u_{i+1}:\cdots:u_{i+j}] \para \overline{\rm{FC}[u_j:u_{j+i+1}:\cdots:u_{j+n}]}\\
&\qquad -\rm{FC}[u_{i+j}:u_i:\cdots:u_{i+j-1}] \para \overline{\rm{FC}[u_{i+j}:u_{i+j+1}:\cdots:u_{j+n}]}.\end{align*}
We now observe that using the cyclic symmetry of formal correlators the first parts of the two terms agree, and that upon applying naturality in taking the quotient by $\rm{span}(u_i-u_{i+j},\ldots,u_{i+j-1}-u_{i+j})$ in particular $u_i$ and $u_{i+j}$ get identified and also the second parts of the two agree; thus they cancel.
\end{proof}

\begin{example}\label{exam:zeta-alt-d-cyc} For $V = F^2$ with affine basis $e_1+e_2,e_1,e_2$, the composition $\zeta^\alt \circ \delta_\rm{cyc}$ has six terms (up to a factor of $\tfrac{1}{2})$:
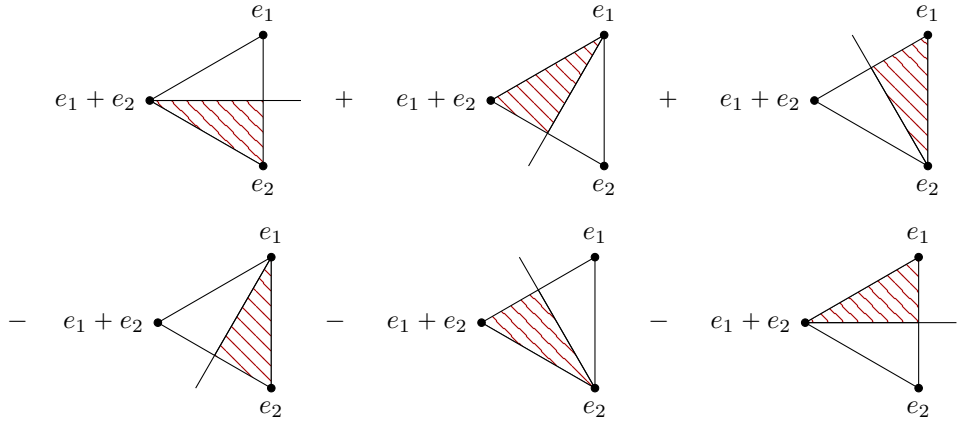
\[\, \quad \begin{tikzpicture}[baseline={([yshift=-.5ex]current bounding box.center)}]
   \begin{scope}
       \clip (180:1) -- (300:1) -- (60:1) -- cycle;
       \draw[pattern={Lines[angle=-45,distance=4pt]},pattern color=Mahogany] (180:1) -- (0:1) -- (-60:1) -- cycle;
   \end{scope}
   \draw (180:1) -- (300:1) -- (60:1) -- cycle;
   \node[inner sep=1pt,circle,draw,fill,label={left:$e_1+e_2$}] at (180:1) {};
   \node[inner sep=1pt,circle,draw,fill,label={above:$e_1$}] at (60:1) {};
   \node[inner sep=1pt,circle,draw,fill,label={below:$e_2$}] at (300:1) {};
   \draw (180:1) -- (0:1);
\end{tikzpicture}
\quad + \quad
\begin{tikzpicture}[baseline={([yshift=-.5ex]current bounding box.center)}]
   \begin{scope}
       \clip (180:1) -- (300:1) -- (60:1) -- cycle;
       \draw[pattern={Lines[angle=-45,distance=4pt]},pattern color=Mahogany] (60:1) -- (180:1) -- (-120:1) -- cycle;
   \end{scope}
   \draw (180:1) -- (300:1) -- (60:1) -- cycle;
   \node[inner sep=1pt,circle,draw,fill,label={left:$e_1+e_2$}] at (180:1) {};
   \node[inner sep=1pt,circle,draw,fill,label={above:$e_1$}] at (60:1) {};
   \node[inner sep=1pt,circle,draw,fill,label={below:$e_2$}] at (300:1) {};
   \draw (60:1) -- (-120:1);
\end{tikzpicture}
\quad + \quad
\begin{tikzpicture}[baseline={([yshift=-.5ex]current bounding box.center)}]
   \begin{scope}
       \clip (180:1) -- (300:1) -- (60:1) -- cycle;
       \draw[pattern={Lines[angle=-45,distance=4pt]},pattern color=Mahogany] (120:1) -- (60:1) -- (-60:1) -- cycle;
   \end{scope}
   \draw (180:1) -- (300:1) -- (60:1) -- cycle;
   \node[inner sep=1pt,circle,draw,fill,label={left:$e_1+e_2$}] at (180:1) {};
   \node[inner sep=1pt,circle,draw,fill,label={above:$e_1$}] at (60:1) {};
   \node[inner sep=1pt,circle,draw,fill,label={below:$e_2$}] at (300:1) {};
   \draw (120:1) -- (-60:1);
\end{tikzpicture}
\]
\[-\quad \begin{tikzpicture}[baseline={([yshift=-.5ex]current bounding box.center)}]
   \begin{scope}
       \clip (180:1) -- (300:1) -- (60:1) -- cycle;
       \draw[pattern={Lines[angle=-45,distance=4pt]},pattern color=Mahogany] (60:1) -- (-60:1) -- (-120:1) -- cycle;
   \end{scope}
   \draw (180:1) -- (300:1) -- (60:1) -- cycle;
   \node[inner sep=1pt,circle,draw,fill,label={left:$e_1+e_2$}] at (180:1) {};
   \node[inner sep=1pt,circle,draw,fill,label={above:$e_1$}] at (60:1) {};
   \node[inner sep=1pt,circle,draw,fill,label={below:$e_2$}] at (300:1) {};
   \draw (60:1) -- (-120:1);
\end{tikzpicture}
\quad - \quad
\begin{tikzpicture}[baseline={([yshift=-.5ex]current bounding box.center)}]
   \begin{scope}
       \clip (180:1) -- (300:1) -- (60:1) -- cycle;
       \draw[pattern={Lines[angle=-45,distance=4pt]},pattern color=Mahogany] (120:1) -- (180:1) -- (-60:1) -- cycle;
   \end{scope}
   \draw (-180:1) -- (300:1) -- (60:1) -- cycle;
   \node[inner sep=1pt,circle,draw,fill,label={left:$e_1+e_2$}] at (180:1) {};
   \node[inner sep=1pt,circle,draw,fill,label={above:$e_1$}] at (60:1) {};
   \node[inner sep=1pt,circle,draw,fill,label={below:$e_2$}] at (300:1) {};
   \draw (120:1) -- (-60:1);
\end{tikzpicture}
\quad - \quad
\begin{tikzpicture}[baseline={([yshift=-.5ex]current bounding box.center)}]
   \begin{scope}
       \clip (180:1) -- (300:1) -- (60:1) -- cycle;
       \draw[pattern={Lines[angle=-45,distance=4pt]},pattern color=Mahogany] (180:1) -- (0:1) -- (60:1) -- cycle;
   \end{scope}
   \draw (180:1) -- (300:1) -- (60:1) -- cycle;
   \node[inner sep=1pt,circle,draw,fill,label={left:$e_1+e_2$}] at (180:1) {};
   \node[inner sep=1pt,circle,draw,fill,label={above:$e_1$}] at (60:1) {};
   \node[inner sep=1pt,circle,draw,fill,label={below:$e_2$}] at (300:1) {};
   \draw (180:1) -- (0:1);
\end{tikzpicture}
\]
Each term on top cancels against the one below it, as the reader should verify from
\begin{align*}&\rm{FC}[e_1{+}e_2{:}e_1] \wedge \ol{\rm{FC}[e_1{+}e_2{:}e_2]}+\rm{FC}[e_1{:}e_2] \wedge \ol{\rm{FC}[e_1{:}e_1{+}e_2]}+\rm{FC}[e_2{:}e_1{+}e_2] \wedge \ol{\rm{FC}[e_2{:}e_1]} \\
&-\rm{FC}[e_1{:}e_1{+}e_2] \wedge \ol{\rm{FC}[e_1{:}e_2]}-\rm{FC}[e_2{:}e_1] \wedge \ol{\rm{FC}[e_2{:}e_1{+}e_2]}-\rm{FC}[e_1{+}e_2{:}e_2] \wedge \ol{\rm{FC}[e_1{+}e_2{:}e_1]}\end{align*}
using cyclic symmetry and identifying vectors in the quotient.
\end{example}

Since the relations \eqref{enum:fc-relations-i}--\eqref{enum:fc-relations-iii} hold for the Steinberg correlators in the infinite Steinberg module and its cobracket is given by the analogous formula, there is a map
\begin{align*}\pr^\St \colon \FC &\lra \SStL \\
\FC[u_0 : \compactcdots : u_n] &\longmapsto \rm{C}[u_0: \compactcdots : u_n]\end{align*}
of Lie coalgebras. It is clear from its construction that the map $\rm{C}^\St_h$ admits a lift to a map $\rm{C}^\FC_h \colon (\rm{B}^\rm{Com} \SSt)^H_n(V) \to \FC(V)$, allowing us to define a decomposition operator
\begin{equation}\label{eqn:dec-fc} D^\FC_h \coloneq (\rm{C}^\FC_h \circ s_H \circ \pr^\St) \colon \FC(V) \lra \FC(V).\end{equation}
Unlike for $\StL$, this is generally \emph{not} the identity. However, it does satisfy 
\[D^\FC_{h_1} D^\FC_{h_2} = D^\FC_{h_1}\]
since $\pr^\St$ annihilates the image of $\id - D^\FC_{h_2}$ as $\pr^\St D^\FC_{h_2} = D^\St_{h_2}$ and this is the identity by the discussion following \eqref{eqn:decomposition-operator}.

\subsubsection{Formal iterated integrals} It is convenient to next introduce a formal analogue of the Steinberg iterated integrals \cite[Definition 28]{CMRR}, or rather their images under the projection $\SStH \to \SStL$. Its value on a vector space $V$ of positive dimension $n$ will be generated by symbols $\FI[v_1,\compactldots,v_n]$ that we call \emph{formal iterated integrals} for bases $v_1,\ldots,v_n$, satisfying only the analogue of \eqref{enum:stl-relations-iii}:
\begin{enumerate}[(1)] \setcounter{enumi}{2}
\item \label{enum:fi-relations-iii} They satisfy the shuffle relations: 
\[\sum_{\sigma \in \rm{Sh}(n_1,n_2)} \rm{FI}[v_{\sigma(1)},\compactldots,v_{\sigma(n_1+n_2)}] =0 \quad \text{for $n=n_1+n_2$ with $n_1,n_2>0$.}\]
\end{enumerate}

\begin{definition}If $V$ is of dimension $n$, then we define
\[\FI(V) \coloneq \begin{cases} \displaystyle \frac{\bb{Q}\{\FI[v_1,\compactldots,v_n] \text{ for bases $v_1,\ldots,v_n$}\}}{\text{\eqref{enum:fi-relations-iii}}} & \text{if $n>0$,} \\
0 & \text{if $n=0$.}\end{cases}\]
\end{definition}

The action of $\GL(V)$ on bases induces a $\GL(V)$-action on $\FI(V)$, and we can assemble all to an object $\FI \in \Fun(\Vect,\rm{GrMod}_\bb{Q}) \subset \Fun(\Vect,\DQ)$. We endow this with a Lie coalgebra structure using the deconcatenation cobracket:

\begin{definition}The \emph{deconcatenation cobracket} $\delta_\rm{dec} \colon \rm{FI} \to \Lambda^2 \rm{FI}$ is given by the formula
\[\delta_\rm{dec}(\rm{FI}[v_1,\compactldots,v_n]) = \sum_{i=1}^{n-1} \rm{FI}[v_1,\compactldots,v_i] \wedge \rm{FI}[v_{i+1},\compactldots,v_n].\]
\end{definition}

There is a more economical construction of this Lie coalgebra: it is isomorphic to the cofree Lie coalgebra on the functor $X \colon \Vect \to \rm{GrMod}_\bb{Q} \subset \DQ$ given by
\[V \longmapsto \begin{cases} \bb{Q}\{\text{nonzero $v \in V$}\} & \text{if $\dim(V)=1$} \\
0 & \text{otherwise.}\end{cases}\]
Indeed, there is an obvious map from underlying object of $\FI$ to $X$ and the adjoint map $\FI \to \cofree_{\coLie}(X)$ is an isomorphism, once we recall from \cref{sec:lie-rep} that the $\coLie$-representations are obtained by taking the quotient of the permutation representations by shuffles. This justifies that the deconcatenation cobracket is well-defined.

\begin{remark}\label{rem:delta-dec-from-sder} One can also use \cref{sec:lie-alg-trees-derivations} to justify that $\delta_\rm{dec}$ is a well-defined Lie cobracket. Indeed, by the analogous procedure to \cref{rem:delta-cyc-from-sder}, it corresponds to the Lie cobracket on $\fr{lie}_{\ul{n}}^\vee$.
\end{remark}

There is a map in $\Fun(\Vect,\DQ)$
\begin{equation}\label{eqn:qfi-to-fc} \begin{aligned} \pr^\FC \colon \rm{FI} &\lra \rm{FC} \\
\FI[v_1,\compactldots,v_n] &\longmapsto (-1)^n\FC[0:v_1:\compactcdots:v_n]\end{aligned} \end{equation}
which is \emph{not} a map of Lie coalgebras, as the deconcatenation cobracket on the cofree Lie coalgebra only involves those terms in the cyclic cobracket that ``cut through the zeroth vertex''. We will return to the relationship between these in \cref{sec:cyclic-coaction}.

It is clear from the construction that the map $\rm{C}^\FC_h$ admits a further lift to a map $\rm{C}^\FI_h \colon (\rm{B}^\rm{Com} \SSt)^H_n(V) \to \FI(V)$, with associated decomposition operator
\[D^\FI_h \coloneq (\rm{C}^\FI_h \circ s_H \circ \pr^\St \circ \pr^\FC) \colon \FI(V) \lra \FI(V),\]
which satisfies $\pr^\FC \circ D^\FI_h = D^\FC_h$ and $D^\FI_{h_2} D^\FI_{h_1} = D^\FI_{h_2}$.

\subsubsection{A resolution of infinite Steinberg modules}\label{sec: reso;ution inf st} We now construct a resolution of $\StL(V)$ by projective  $\bb{Q}[\GL(V)]$-modules. The starting point is the bar construction for a semigroup (i.e.~nonunital monoid) $S$ with a left action on a module $M$ in a category left tensored over $\DQ$; this is simply the bar construction of $M$ as a left $\bb{Q}[S]$-module and explicitly given by
\[B_k(S,M) \coloneq \bb{Q}[S]^{\otimes k} \otimes M\]
with differential 
\begin{align*}d([s_1|\compactcdots|s_k] \otimes m) &= [s_2|\compactcdots|s_k] \otimes m + \sum_{i=1}^{k-1} (-1)^i [s_1|\compactcdots|s_is_{i+1}|\compactcdots|s_k] \otimes m \\
&\qquad + (-1)^k [s_1|\compactcdots|s_{k-1}] \otimes s_km.\end{align*}
We will apply this in the category of $\bb{Q}[\GL(V)]$-modules and to the following semigroup:

\begin{definition}$\Dec_V$ is the semigroup whose elements are nonzero linear functionals $h \in V^\vee$ with multiplication given by $h_1h_2 = h_1$.\end{definition}

It has been constructed so that letting $h$ act on $\FC(V)$ or $\FI(V)$ by decomposition operators gives an action of $\Dec_V$ on these $\bb{Q}[\GL(V)]$-modules. There is also an action of $\Dec_V$ on $\bb{Q}$ by \emph{zero maps}, and we denote the associated $\bb{Q}[\Dec_V]$-module by $\bb{Q}_0$:

\begin{lemma}\label{lem:decv-q0-acyclic} $H_*(\Dec_V;\bb{Q}_0)=0$ for all $* \in \bb{Z}$.\end{lemma}

\begin{proof}This is the chain complex given by
\[B_k(\Dec_V;\bb{Q}_0) = \bb{Q}\{[h_1|\compactldots|h_k] \mid h_1,\ldots,h_k \in V^\vee \setminus \{0\}\}\]
with differential given by $d([h_1|\compactldots|h_k]) = \sum_{i=1}^{k} (-1)^{i+1} [h_1|\compactldots|\widehat{h}_i|\compactldots|h_k]$ and fixing some $h \in V^\vee \setminus \{0\}$ we can give a chain null-homotopy of the identity by $H([h_1|\compactldots|h_k]) = [h|h_1|\compactldots|h_k]$, verified by the computation
\begin{align*}(dH+Hd)([h_1|\compactldots|h_k]) &= \sum_{i=0}^{k} (-1)^i [h|h_1|\compactldots|\widehat{h}_i|\compactldots|h_k] -\sum_{i=1}^k (-1)^i[h|h_1|\compactldots|\widehat{h}_i|\compactldots|h_k] \\
&= [h_1|\compactldots|h_k]\end{align*}
where in the first sum in the term $i=0$ we delete $h$.
\end{proof}

\begin{lemma}\label{lem:stl-fi-resolution} We have
\[H_k(\Dec_V;\FI(V)) \cong \begin{cases} \StL(V) & \text{if $k=0$,} \\
0 & \text{otherwise,}\end{cases}\]
with isomorphism induced by the map $\pr^\St \circ \pr^\FC \colon \FI(V) \to \StL(V)$.\end{lemma}

\begin{proof}We claim that the evidently well-defined map
\begin{align*}\alpha \colon \bb{Q}[\Dec_V] \otimes \StL(V) &\lra \FI(V) \\
[h] \otimes \rm{C}[0:v_1:\compactcdots:v_n] &\longmapsto (\rm{C}^\FI_h \circ s_H)(\rm{C}[0:v_1:\compactcdots:v_n])\end{align*}
is an isomorphism. Its inverse is given by taking $\FI[v_1,\compactldots,v_n]$ to $(-1)^n [h] \otimes \rm{C}[0:v_1:\compactcdots:v_n]$ where $h$ is the unique linear functional so that $h(v_1) = \cdots = h(v_n)=1$, which is well-defined because $h$ does not depend on the order of the $v_i$ and the Steinberg correlators satisfy the shuffle relations. Then analogous to the proof of \cref{lem:ch-inverse-of-symbol}, when computing $\alpha([h] \otimes \rm{C}[0:v_1:\compactcdots:v_n])$ taking the projection of the symbol onto the summands where none of the lines is contained in $H = \ker(h)$ will annihilate all terms in the symbol of $\rm{C}[0:v_1:\compactcdots:v_n]$ except $[P_1|\compactcdots|P_n]$ with $P_i = \rm{span}(v_i)$, which $C^\FI_h$ maps to the formal iterated integral $(-1)^n \FI[v_1,\compactldots,v_n]$; since formal iterated integrals satisfy the shuffle relations this is well-defined. More generally, we have an isomorphism
\begin{align*} \alpha_k \colon B_{k+1}(\Dec_V;\StL(V)) &\overset{\cong}\lra B_k(\Dec_V;\FI(V)) \\
[h_1|\compactldots|h_{k+1}] \otimes \rm{C}[0:v_1:\compactcdots:v_n] &\longmapsto [h_1|\compactldots|h_k] \otimes \alpha([h_{k+1}] \otimes \rm{C}[0:v_1:\compactcdots:v_n]).\end{align*}
Augmenting $B_*(\Dec_V;\FI(V))$ by $B_{-1}(\Dec_V;\FI(V)) \coloneq \StL(V)$ using the map $\pr^\St \circ \pr^\FC$, and letting $\Dec_V$ act on $\StL(V)$ by zero maps, the $\alpha_k$ assemble to an isomorphism of chain complexes as we have
\[\alpha([h_k] \otimes \rm{C}[0:v_1:\compactcdots:v_n]) = D^\FI_{h_k}(\alpha([h_{k+1}] \otimes \rm{C}[0:v_1:\compactcdots:v_n]))\]
as a consequence of the equation $\rm{C}^\FI_{h_k} = D^\FI_{h_k} \rm{C}^{\FI}_{h_{k+1}}$ of maps $\StL(V) \to \FI(V)$. We conclude that $B_*(\Dec_V;\FI(V))$ augmented by $\StL(V)$ is acyclic, and the result follows.
\end{proof}

The following result in particular implies \cref{prop:stl-explicit-pres}:

\begin{proposition}\label{prop:stl-resolution} We have
\[H_p(\Dec_V;\FC(V)) \cong \begin{cases} \StL(V) & \text{if $p =0$,} \\
0 & \text{otherwise,}\end{cases}\]
with isomorphism induced by the map $\pr^\St \colon \FC(V) \to \StL(V)$.
\end{proposition}

\begin{proof}As $\pr^\FC \colon \FI(V) \to \FC(V)$ is compatible with the $\Dec_V$-action, we may define $K(V)$ as $\bb{Q}[\GL(V)]$-module with $\Dec_V$-action through the short exact sequence
\[0 \lra K(V) \lra \FI(V) \overset{\pr^\FC}\lra \FC(V) \lra 0.\]
Note that each $D^\FI_h$ annihilates $K(V)$, and thus there is an isomorphism of chain complexes $B_*(\Dec_V;K(V)) \cong B_*(\Dec_V;\bb{Q}_0) \otimes K(V)$ and the latter is acyclic by \cref{lem:decv-q0-acyclic}. Using the long exact sequence we conclude that $\pr^\FC$ induces an isomorphism
\[H_*(\Dec_V;\FI(V)) \overset{\cong}\lra H_*(\Dec_V;\FC(V))\]
which is compatible with the map to $\StL(V)$.
\end{proof}

It remains to establish the quality of the terms in our resolutions of $\StL(V)$. We will repeatedly make use of the following fact: if $R$ is a $\bb{Q}$-algebra, $M$ is a projective $R$-module with $\fr{S}_k$-action, and $S$ is a finite-dimensional $\bb{Q}[\fr{S}_k]$-module, then $M \otimes_{\bb{Q}[\fr{S}_k]} S$ is also a projective $R$-module. To see this, use that a summand of a projective module is projective and every $\bb{Q}[\fr{S}_k]$-module $S$ is a summand of a free $\bb{Q}[\fr{S}_k]$-module so $M \otimes_{\bb{Q}[\fr{S}_k]} S$ is also a summand of a free $R$-module.

\begin{lemma} \label{lemma: resolution is projective} Each term in $B_*(\Dec_V;\FI(V))$ and $B_*(\Dec_V;\FC(V))$ is a projective $\bb{Q}[\GL(V)]$-module.\end{lemma}

\begin{proof} If $V$ is of dimension $n$, $\bb{Q}[\Dec_V]^{\otimes k} \otimes \FI(V)$ is isomorphic to $F \otimes_{\bb{Q}[\fr{S}_n]} \rm{coLie}_n$ for the free module with basis as a $\bb{Q}[\GL(V)]$-module given by $[h_1|\compactldots|h_k] \otimes [e_1,\compactldots,e_n]$ for nonzero functionals $h_i$ on $V$ and a fixed basis $e_1,\ldots,e_n$ of $V$, and hence is projective. Similarly, $\bb{Q}[\Dec_V]^{\otimes k} \otimes \FC(V)$ is isomorphic to $F' \otimes_{\bb{Q}[\fr{S}_{n+1}]} \rm{cycLie}_n^{\vee}$ for the free module $F'$ with basis as a $\bb{Q}[\GL(V)]$-module given by $[h_1|\compactldots|h_p] \otimes [0:e_1:\compactldots:e_n]$ for nonzero functionals $h_i$ on $V$ and a fixed affine basis $0,e_1,\ldots,e_n$ of $V$.
\end{proof}

This tells us that the higher homology groups of $\GL(V)$ acting on the terms vanish. We will also need the following computations of coinvariants:

\begin{lemma}\label{lem:fi-fc-coinvariants} We have that
\[\FI(V)_{\GL(V)} \cong \begin{cases} \bb{Q} & \text{if $\dim(V) = 1$,} \\
0 & \text{else,}\end{cases} \quad \text{and} \quad \FC(V)_{\GL(V)} \cong \begin{cases} \bb{Q} & \text{if $\dim(V) = 1$,} \\
0 & \text{else,}\end{cases}\]
and the map $\FI(V)_{\GL(V)} \to \FC(V)_{\GL(V)}$ induced by $\pr^\FC$ is an isomorphism.
\end{lemma}

\begin{proof}In the case $\dim(V)=1$, it follows from the definitions that the surjective map $\FI(V) \to \FC(V)$ is an isomorphism, as upon picking an identification $V \cong F$, $\FI(V)$ is isomorphic to $\bb{Q}[F^\times]$ with standard $\bb{Q}[F^\times]$-action while $\FC(V)$ is isomorphic to its quotient $\bb{Q}[F^\times/\{\pm 1\}]$. In the case $\dim(V)>1$, it suffices to prove that $\FI(V)_{\GL(V)} = 0$. This follows because $\GL(V)$ acts transitively on bases, implying the coinvariants are generated by the class of $\FI[e_1,\compactldots,e_n]$ for a fixed basis $e_1,\compactldots,e_n$, and then the shuffle relations imply that 
\[{n_1+n_2 \choose n_1} \FI[e_1,\compactldots,e_n] = 0 \qquad \text{for $n=n_1+n_2$ with $n_1,n_2>0$}\]
in the coinvariants, so that $\FI[e_1,\compactldots,e_n] = 0$ since we are working over the rationals.
\end{proof}

\subsection{A presentation for $\scr{G}$}\label{sec:generators of G}  We use the results of the previous subsection to give a presentation for $\scr{G}$ and then make concrete how to obtain elements in it.

\subsubsection{The presentation} Our starting point is the short exact sequence
\[0 \lra \FCR(V) \lra \FC(V) \overset{\pr^\St}\lra \StL(V) \lra 0\]
defining the left term as the kernel of the map $\pr^\St$: it is simply given by elements of $\FC(V)$ that represent \emph{relations} in $\StL(V)$. From \cref{lemma: resolution is projective} we know $\FC(V)$ is projective, and from \cref{thm:bgl-critical-line-vanishing} and \cref{lem:fi-fc-coinvariants} that the coinvariants of $\StL(V)$ and $\FC(V)$ vanish if $\dim(V) \neq 1$ and are isomorphic to $\bb{Q}$ otherwise. Thus the connecting homomorphism yields an isomorphism
\[\scr{G}(V) \coloneq H_1(\GL(V);\StL(V)) \overset{\cong}\lra H_0(\GL(V);\FCR(V)).\]
To understand the right side, we use that \cref{prop:stl-resolution} provides a resolution by projective $\GL(V)$-modules by \cref{lemma: resolution is projective}
\[B_2(\Dec_V;\FC(V)) \lra B_1(\Dec_V;\FC(V)) \lra \FCR(V),\]
inducing an exact sequence
\[H_0(\GL(V);B_2(\Dec_V;\FC(V))) \lra H_0(\GL(V);B_1(\Dec_V;\FC(V))) \lra \scr{G}(V) \lra 0.\]
Spelling out the definitions, this gives a presentation for $\scr{G}(V)$: if $V$ is of dimension $n$ then it has generators $[h] \otimes \FC[u_0:\compactcdots:u_n]$ for affine bases $u_0,\ldots,u_n$ and nonzero linear functionals $h$. These satisfy in addition to the relations of formal correlators \eqref{enum:fc-relations-i}--\eqref{enum:fc-relations-iii} the following two relations:
\begin{enumerate}[\noindent (1)] \setcounter{enumi}{3}
    \item \label{enum:gd-relation-iv} They satisfy the coinvariant relations
    \[[h] \otimes \FC[u_0:\compactcdots:u_n] = [g^*h] \otimes \FC[gu_0:\compactcdots:gu_n] \qquad \text{for $g \in \GL(V)$}.\]
    \item \label{enum:gd-relation-v} They satisfy the decomposition relations 
    \[\qquad [h_2] \otimes \FC[u_0:\compactcdots:u_n]-[h_1] \otimes \FC[u_0:\compactcdots:u_n]+[h_1] \otimes D^\FC_{h_2}(\FC[u_0:\compactcdots:u_n])=0.\]
\end{enumerate}
Let us record this:

\begin{proposition}\label{prop:presentation for G(F)} Suppose that $V$ is of dimension $n$, then there is an isomorphism
\[\frac{\bb{Q}\{[h] \otimes \FC[u_0:\compactcdots:u_n] \text{ for nonzero functionals $h$ and affine bases $u_0,\ldots,u_n$}\}}{\text{\eqref{enum:fc-relations-i}--\eqref{enum:gd-relation-v}}} \overset{\cong}\lra \scr{G}(V).\] 
\end{proposition}

If $h(u_i) = x_i \in F$, we denote the image of $[h] \otimes \FC[u_0:u_1:\compactcdots : u_n]$ by
\[\Cor^\scr{G}(x_0,x_1,\compactldots,x_n) \in \scr{G}(V),\] 
justified using \eqref{enum:fc-relations-i} and \eqref{enum:gd-relation-iv}, which imply that the image only depends on the elements $x_i$.

\begin{remark}Using \cref{lem:stl-fi-resolution} instead of \cref{prop:stl-resolution} we could have equally well given a presentation using $[h] \otimes \FI[v_1,\compactldots,v_n]$ instead.
\end{remark}

\subsubsection{Projection from $\FCR(V)$ to $\PolyL$}\label{section: projection from FC to G} We explain an explicit formula for the projection $\FCR(V)\to H_0(\GL(V);\FCR(V))$. Given $V$ of dimension $n$ and nonzero functional $h\in V^{\vee}$, consider the map
\begin{align*} E_h\colon \FC(V) &\lra \PolyL(V) \\
\FC[u_0:u_1:\compactcdots:u_n] &\longmapsto \CorG(h(u_0),h(u_1),\compactldots,h(u_n)).
\end{align*}

\begin{proposition}\label{prop: formula for projection to coinvariants}
The projection $\FCR(V)\lra H_0(\GL(V);\FCR(V))\cong \PolyL(V)$ coincides with the composition
\[
\FCR(V) \overset{\inc}\lra \FC(V) \overset{E_h}{\lra} \PolyL(V).
\]
\end{proposition} 
\begin{proof}
Consider an element $
x= \sum_i a_i \FC[u_0^{(i)}:\compactcdots:u_n^{(i)}]\in  \FCR(V)$. Since the decomposition operator $D_h^{\FC}\colon \FC(V)\to \FC(V)$  factors through $\StL(V)$, it has to vanish on $\FCR(V)$, so we have
\[
 \sum a_i D_h^{\FC}\left( \FC[u_0^{(i)}:\compactcdots:u_n^{(i)}]\right)=0.
\]
Thus we can write
\[
x= \sum a_i \left(\FC[u_0^{(i)}:\compactcdots:u_n^{(i)}]-D_h^{\FC}\FC[u_0^{(i)}:\compactcdots:u_n^{(i)}]\right).
\]
As the element $\FC[u_0^{(i)}:\compactcdots:u_n^{(i)}]-D_h^{\FC}(\FC[u_0^{(i)}:\compactcdots:u_n^{(i)}])\in \FCR(V)$ projects to $\CorG(h(u_0^{(i)}),\compactcdots, h(u_n^{(i)}))$, the projection of $x$ to $\PolyL(V)$ equals $E_h(x)$. 
\end{proof}

\cref{prop: formula for projection to coinvariants} implies that the value of the map $E_h$ on an element in $\FCR(V)$ is independent of a nonzero linear functional $h$. This statement can be interpreted as a functional equation in $\PolyL(V)$, yielding \cref{thm:polyl-relations-from-stl-relations}:

\begin{customthm}{F}
Suppose that an identity 
\[
\sum a_i \rm{C}\bigl[u_0^{(i)}:\compactcdots:u_n^{(i)}\bigr]=0\]
for $a_i\in \Q$ and affine bases $u_0^{(i)},\dots,u_n^{(i)}$ of $V$ holds in $\StL(V)$. Then the element 
\[
\sum a_i \CorG\Bigl(h\bigl(u_0^{(i)}\bigr),\compactcdots, h\bigl(u_n^{(i)}\bigr)\Bigr)\in \PolyL_n(V)
\]
is independent of a nonzero linear functional $h$ on $V$.
\end{customthm}

\subsection{Computing the cobracket on $\scr{G}$}

\begin{notation}For the sake of readability, we will write the direct sum $\oplus$ as $+$ in this subsection.
\end{notation}

We will now explicitly compute the cobracket on $\scr{G}$ in terms of the presentation from \cref{prop:presentation for G(F)}.

\subsubsection{Computing the cobracket on $\scr{G}$: outline} \label{sec:cobracket-outline} As we have explained in \cref{sec:ek-homology}, the shifted Lie coalgebra structure on $E_\infty$-indecomposables yields in particular a Lie coalgebra structure on $\scr{G}$ with cobracket $\delta \colon \scr{G} = H_1(\GL;\SStL) \to \Lambda^2 \scr{G}= \Lambda^2 H_1(\GL;\SStL)$. Recall the following is \cref{thm:polyl-presentation-cobracket}:

\begin{customthm}{C.b} The cobracket $\delta \colon \scr{G} \to \Lambda^2 \scr{G}$ is given by the formula
\[\delta(\CorG(x_0,x_1,\compactldots,x_n)) = \sum_{j=0}^n \sum_{i=1}^{n-1} \CorG(x_j,x_{j+1},\compactldots,x_{j+i}) \wedge \CorG(x_j,x_{j+i+1},\compactldots,x_{j+n})\]
as long as $x_0,\ldots,x_n$ are all distinct.
\end{customthm}

\begin{remark}This is \emph{not} the entirety of the structure induced on $\scr{G}$, as the nonvanishing of the coinvariants $H_0(\GL_1,\SStL(F^1)) \cong \bb{Q}$ also yields a ``$\sigma$-component''
\[\delta_\sigma \colon \scr{G}_n =  H_1(\GL_n;\StL_n) \lra  H_2(\GL_{n-1},\StL_{n-1}).\]
This will be discussed in more detail in \cref{sec:sigma-component}.\end{remark}

Before outlining the proof, we recall that \cref{proposition: generic correlators} justifies that $\scr{G}_n(F)$ is generated by generic correlators, i.e.~those with $x_0,\ldots,x_n$ all distinct, so the above formula determines the cobracket completely. We expect the same formula holds in general, though, with the caveat that any term of the form $\CorG(x,\compactldots,x)$ must be set to zero; we intend to return to this in future work. We will explicitly state below when we use the genericity property.

To prove \cref{thm:polyl-presentation-cobracket}, we recall from \cref{prop:cobracket-via-coproduct} the recipe for computing the cobracket on $\scr{G}$. We start applying $\dim_!$ to the zigzag in $\Fun(\Vect,\DQ)$
\[[\St^\infty \to 0] \longleftarrow [\St^\infty \to \Lambda^2 \St^\infty] \lra [0 \to \St^\infty \para \St^\infty].\]
Evaluating on the vector space $F^n$, we take homology in degree $2n$: on the left this yields $H_1(\GL_n;\StL_n)$ and on the right---after applying the K\"unneth theorem and discarding terms corresponding to $\sigma$-component $\delta_\sigma$---this yields terms $\bigoplus_{d=0}^n H_1(\GL_d;\StL_d) \otimes H_1(\GL_{n-d},\StL_{n-d})$. The left map is surjective, so we pick a lift and apply the right map; after antisymmetrising we get the cobracket, which is independent of the choice of lift.

We will perform this computation in the 1-category $\rm{Ch}_\bb{Q}$ of rational chain complexes, by resolving the zigzag by maps of chain complexes of $\bb{Q}[\GL(V)]$-modules which have sufficient vanishing properties to reduce a computation on $H_1$ to a computation on $H_0$. More precisely, we will construct chain complexes $A_*$ \eqref{eqn:res-cobracket-i}, $B_*$ \eqref{eqn:res-cobracket-ii}, $C_*$ \eqref{eqn:res-cobracket-iii}, and maps of chain complexes fitting in a commutative diagram
\[\begin{tikzcd} A_* \rar{\simeq} & {[\StL(V) \to 0]} \\[-5pt]
B_* \rar{\simeq} \uar \dar & {[\StL(V) \to (\Lambda^2 \StL)(V)]} \uar \dar \\[-5pt]
C_* \rar{\simeq} \dar & {[0 \to (\StL \para \StL)(V)]} \\[-5pt]
(\FCR \para \FCR)(V), \end{tikzcd}\]
whose horizontal maps are quasi-isomorphisms so induce an isomorphism on $H_*(\GL(V);-)$ for all $* \geq 0$. We also understand the effect of the left-bottom zig-zag on $H_*(\GL(V);-)$.

The main work will lie in the following: 
\begin{enumerate}[(a)]
    \item \label{enum:technical-cobracket-i} The construction of $B_*$ will amount to prescribing relations that give a reason that $B_1(\Dec_V;\FI(V))$ maps to zero in $(\Lambda^2 \StL)(V)$.
    \item \label{enum:technical-cobracket-ii} The explicit evaluation of the map $B_* \to C_*$ will amount to a vanishing result along the lines of \cref{lem:fc-cobracket-symmetry}.
\end{enumerate}

\subsubsection{The cyclic coaction}\label{sec:cyclic-coaction} We start by writing the cyclic cobracket as a sum of the deconcatenation cobracket and a correction term that captures the remaining terms; the \emph{cyclic coaction}. The latter is given by letting the first index start at $j=1$ rather than $j=0$ in \cref{prop:stl-explicit-cobracket}, and interpreting the unique term containing the entry $u_0 = 0$ as a formal iterated integral, moving it to the front with an appropriate sign if necessary:

\begin{definition}The \emph{cyclic coaction} $\delta_\rm{coact} \colon \rm{FI} \to \rm{FI} \otimes \rm{FC}$ is given by the formula
\begin{align*}
\delta_{\mathrm{coact}}&(\FI[v_1,\ldots,v_n])
\\
=&
\sum_{j=1}^{n}\ \ \: \sum_{i=1}^{n-j}
(-1)^{i+1}
\FI[v_1,\ldots,v_j,v_{j+i+1},\ldots,v_n]
\otimes
\FC[v_j:\cdots:v_{j+i}]\\
+&
\sum_{j=1}^{n}\sum_{i=n+1-j}^{n-1}
(-1)^{n-i}
\FI[v_1,\ldots,v_{j+i-n-1},v_j,\ldots,v_n]
\otimes
\FC[v_j:v_{j+i+1}:\cdots:v_{j+n}].
\end{align*}
where $v_0$ is put to be $0$ and the indices are to be interpreted cyclically.
\end{definition}

We need to justify the well-definedness of $\delta_\rm{coact}$, and though we will not need it, will also make precise some of its properties; this can be done by direct computation but we opt to take an approach along the lines of \cref{rem:delta-cyc-from-sder,rem:delta-dec-from-sder}. In the following, we consider exterior powers as subsets of tensor powers via $x \wedge y \mapsto \frac{1}{2}(x \otimes y - y \otimes x)$.

\begin{lemma}The map $\delta_\rm{coact}$ is well-defined and has the property that the following diagrams commute
\[\begin{tikzcd} \FI \rar{\delta_\rm{dec}} \dar{\delta_\rm{coact}} &[20pt] \FI \otimes \FI \dar{(\id \otimes \sigma) \circ (\delta_\rm{coact} \otimes \id) + \id \otimes \delta_\rm{coact}}\\[-5pt]
\FI \otimes \FC \rar{\delta_\rm{dec} \otimes \id} & \FI \otimes \FI \otimes \FC, \\[-5pt]
\FI \rar{\delta_\rm{coact}} \dar{\delta_\rm{coact}} &[20pt] \FI \otimes \FC \dar{(\id-\sigma \otimes \id) \circ (\delta_\rm{coact} \otimes \id)} \\[-5pt]
\FI \otimes \FC \rar{\id \otimes \delta_\rm{cyc}} & \FI \otimes \FC \otimes \FC.
\end{tikzcd}\]
\end{lemma}

\begin{proof}It suffices to verify this on a generator $\FI[v_1,\compactldots,v_n]$. As we need to prove certain equations hold, it suffices to do this after applying the functor that sums up the values at all vector spaces. We now refer to \cref{sec:lie-alg-trees-derivations} for the Lie coalgebras $\fr{sder}_S^\vee$ and $\fr{lie}_S^\vee$ for $S = [n]$. Upon identifying $\FI[v_{i_1},\compactldots,v_{i_j}]$ for $\{i_1,\ldots,i_j\} \subseteq \ul{n}$ with the element $\rm{LI}[X_{i_1},\compactldots,X_{i_j}]$ of $\fr{lie}_S^\vee$ and $\FC[0:v_{j_1}:\compactcdots:v_{j_k}]$ with the element $\rm{LC}[X_0:X_{j_1}:\compactcdots:X_{j_k}]$, we see that $\delta_\rm{cyc}$ and $\delta_\rm{dec}$ are given by the same formulas. Moreover, in this setting $\delta_\rm{coact}$ admits an interpretation as the coaction of $\fr{sder}_S^\vee$ on $\fr{lie}_S^\vee$ that is dual to the action of the special derivations $\fr{sder}_S$ on the free Lie algebra $\fr{lie}_S$. This implies the statement, as this action is well-defined and the diagrams in the statement of the lemma are dual to those saying it is an action. 
\end{proof}

The following is a direct consequence of the definitions, where we recall that $\smash{\pr^\FC} \colon \FI \to \FC$ was defined in \eqref{eqn:qfi-to-fc} and $\alt \colon \FC \otimes \FC \to \Lambda^2 \FC$ is given by $x \otimes y \mapsto x \wedge y$.

\begin{lemma}\label{lem:coact-correction} There is an equation of maps $\FI \to \Lambda^2 \FC$
\[\delta_\rm{cyc} \circ \pr^\FC - \Lambda^2 \pr^\FC \circ \delta_\rm{dec} = \alt \circ (\pr^\FC \otimes \id) \circ \delta_\rm{coact}.\]
\end{lemma}

Recall that the map $\rm{C}^\FC_h$ admits a further lift to a map $\rm{C}^\FI_h \colon (\rm{B}^\rm{Com} \SSt)^H_n(V) \to \FI(V)$, with associated decomposition operator $D^\FI_h = (\rm{C}^\FI_h \circ s_H \circ \pr^\St \circ \pr^\FC) \colon \FI(V) \to \FI(V)$, which satisfies
\[D^\FI_{h_1} \circ D^\FI_{h_2} = D^\FI_{h_1} \quad \text{and} \quad \pr^\FC \circ D^\FI_h = D^\FC_h \circ \pr^\FC.\] 
There is a more subtle compatibility between the decomposition operators and cobrackets:

\begin{lemma}\label{lem:cyc-dec-h} Given nonzero linear functional $h \in V^\vee$, there is an equation of maps $\FI(V) \to \Lambda^2 \FC(V)$
\[\Lambda^2 \pr^\FC \circ \delta_\rm{dec} \circ D_h^{\FI} = \Lambda^2 D_h^\FC \circ \delta_\rm{cyc} \circ \pr^\FC\]
where $D_h^{\FC} \colon \FC(U) \to \FC(U)$ for $U \subseteq V$ is defined to be $0$ if $U \subseteq \ker(h)$ and $D_{h|_U}^{\FC}$ otherwise.
\end{lemma}

\begin{proof}Writing $H = \ker(h)$, consider the following diagram
\[\hspace{-.15cm}\begin{tikzcd} \FC(V) \rar{\pr^\St} \dar{\delta_\rm{cyc}} &[-7pt] \StL(V) \rar{s} \dar{\delta} &[-6pt] (\rm{B}^{\rm{Com}} \SSt)_n(V) \rar{\pi_H} &[-5pt] (\rm{B}^{\rm{Com}} \SSt)^H_n(V) \rar{C_h^\rm{FI}} &[-5pt] \FI(V) \dar{\delta_\rm{dec}} \\
(\Lambda^2 \FC)(V) \rar{\Lambda^2 \pr^\St} & (\Lambda^2 \StL)(V) \rar{\Lambda^2 s} & (\Lambda^2 \rm{B}^{\rm{Com}} \SSt)_n(V) \rar{\Lambda^2 \pi_H} & (\Lambda^2 \rm{B}^{\rm{Com}} \SSt)^H_n(V) \rar{\Lambda^2 \rm{C}^\FI_h} & (\Lambda^2 \FI)(V),\end{tikzcd}\]
where the map $\Lambda^2 \pi_H$ on the bottom-right needs some further explanation as $H$ may not be a subspace of a summand of $V$: its domain $(\Lambda^2 \rm{B}^{\rm{Com}} \SSt)_n(V)$ is a direct sum indexed by an (unordered) pair of sets $\{P_1,\ldots,P_i\}$ and $\{P_{i+1},\ldots,P_n\}$ of lines $P_1,\ldots,P_n \subseteq V$ and we define $\Lambda^2 \pi_H$ as projection onto the summand $(\Lambda^2 \rm{B}^{\rm{Com}} \SSt)^H_n(V)$ of those terms where none of these lines are contained in $H$. 

Given this definition, the maps in the statement of this lemma are obtained by taking the top-right and bottom-left composites in the diagram, precomposing with $\pr^\FC$, and postcomposing with $\Lambda^2 \pr^\FC$. It hence suffices to prove that this diagram commutes. 

\medskip

The left square commutes because $\pr^\St \colon \FC \to \SStL$ is a map of Lie coalgebras. For the right square, we use that $\St^\infty(V)$ is generated by $\rm{C}[0:v_1:\compactcdots:v_n]$ where $h(v_1) = \cdots = h(v_n) = 1$. We first compute the value of top-right composition on such a $\rm{C}[0:v_1:\compactcdots:v_n]$. As in the proof of \cref{lem:ch-inverse-of-symbol} the latter condition implies that $\smash{\rm{C}^\FI_h} \circ \pi_H \circ s$ sends this to $(-1)^n\FI[v_1,\compactldots,v_n]$ and in turn $\delta_\rm{dec}$ sends this to $(-1)^n\smash{\sum_{i=1}^{n-1}}\FI[v_1,\compactldots,v_i] \wedge \FI[v_{i+1},\compactldots,v_n]$. 

To see the right square commutes, we next compute the value of the left-bottom composition on $\rm{C}[0:v_1:\compactcdots:v_n]$ and start with the terms in the formula for $\delta$ where $j > 0$. All of these terms are concentrated on pairs of subspaces $\{V_I,V_J\}$ where at least one $V_J$ is contained in $H$ and hence are annihilated by $\Lambda^2 \pi_H$. It remains to consider the terms in the formula $\delta$ where $j = 0$, given by $\sum_{i=1}^{n-1} \rm{C}[0:v_1:\compactcdots:v_i] \wedge \rm{C}[0:v_{i+1}:\compactcdots:v_n]$ concentrated at pairs of subspaces $\{V_I,V_J\}$ both not contained in $H$ and arguing as above $\Lambda^2 \rm{C}^\FI_h \circ \Lambda^2 \pi_H \circ \Lambda^2 s$ sends this to $\smash{\sum_{i=1}^{n-1}} (-1)^i \FI[v_1,\compactldots,v_i] \wedge (-1)^{n-i} \FI[v_{i+1},\compactldots,v_n]$. We conclude the square commutes.
\end{proof}

\subsubsection{Choosing relations} \label{sec:cobracket-choosing-relations} We now make precise the first of the technical inputs, \eqref{enum:technical-cobracket-i}, by constructing a dashed map
\begin{equation}\label{eqn:choosing-relations-diag}\begin{tikzcd} B_1(\Dec_V;\rm{FI}(V)) \rar \dar[dashed]{\rho} & \FI(V)  \dar{\delta_\rm{cyc} \circ \pr^\FC} \\
(\FCR \otimes \FC+\FC \otimes \FCR)(V) \rar & (\Lambda^2\FC)(V)  \end{tikzcd}\end{equation}
making the diagram commute, where the top map is part of the resolution in \cref{lem:stl-fi-resolution} and the bottom map is induced by the inclusion $\FCR(V) \to \rm{FC}(V)$ followed by anti-symmetrising. When restricted to generic pairs $[h] \otimes \FI[v_1,\compactldots,v_n]$, i.e.~those where $h(v_1),\ldots,h(v_n)$ are all distinct and nonzero, it will have the following two properties:
\begin{enumerate}[(i)]
    \item $\rho$ is a difference of two terms $\rho^\rm{cyc}$ and $\rho^\rm{coact}$,
    \item one of the components of $\rho^\rm{coact}$ under the map
    \[\zeta^\alt \circ \alt \colon (\FCR \otimes \FC+\FC \otimes \FCR)(V) \lra (\FCR \para \FC+\FC \para \FCR)(V)\]
    defined in \cref{sec:map-b-to-c}, vanishes.
\end{enumerate}
We will not attempt to specify the value of $\rho$ on nongeneric pairs $[h] \otimes \FI[v_1,\compactldots,v_n]$: there is a direct sum decomposition
\[ B_1(\Dec_V;\rm{FI}(V)) \cong  B_1(\Dec_V;\rm{FI}(V))^\rm{gen}+B_1(\Dec_V;\rm{FI}(V))^\rm{nongen}\]
into the subspace spanned by generic pairs and the subspace spanned by nongeneric pairs. We will define $\rho$ explicitly on the first term. On the second term we merely prove $\rho$ exists making \eqref{eqn:choosing-relations-diag} commute, using that (a) by the proof of \cref{lemma: resolution is projective}, $B_1(\Dec_V;\rm{FI}(V))^\rm{nongen}$ spanned by nongeneric pairs is projective, and (b) we know that the composition $B_1(\Dec_V;\rm{FI}(V)) \to \FI(V) \to (\Lambda^2 \FC)(V) \to (\Lambda^2 \StL)(V)$ vanishes.

The top-left corner of \eqref{eqn:choosing-relations-diag} is spanned by pairs $[h] \otimes \FI[v_1,\compactldots,v_n]$, which get sent by the top horizontal map to $(\id-D_h^{\FI})(\FI[v_1,\compactldots,v_n]) \in \FI(V)$. When we apply the right vertical map $\delta_\rm{cyc} \circ \pr^\FC$ to this, by \cref{lem:coact-correction,lem:cyc-dec-h} we get a sum of two terms
\[\delta_\rm{cyc} \circ \pr^\FC \circ (\id-D_h^{\FI}) = \delta^\rm{cyc}_h-\delta^\rm{coact}_h, \qquad \text{with} \qquad
\begin{aligned}&\delta^\rm{cyc}_h \coloneq (\id-\Lambda^2 D_h) \circ \delta_\rm{cyc} \circ \pr^\FC \\
&\delta^\rm{coact}_h \coloneq \overline{\delta}_\rm{coact} \circ D^{\FI}_h\end{aligned}\]
where we abbreviate 
\[\ol{\delta}_\rm{coact} \coloneq \alt \circ (\pr^\FC \otimes \id) \circ \delta_\rm{coact} \qquad \text{and} \qquad D_h \coloneqq D^\FC_h.\]
We will find an explicit lift of $\delta^\rm{cyc}_h(\FI[v_1,\compactldots,v_n])$ to the bottom-left term, and then use this to prove the existence of a nonexplicit lift of $\delta^\rm{coact}_h(\FI[v_1,\compactldots,v_n])$ with certain properties.

\subsubsection{Choosing relations for the ``cyclic cobracket'' term} Assume that $[h] \otimes \FI[v_1,\compactldots,v_n]$ is generic. We start with the construction of the explicit lift of $\delta^\rm{cyc}_h(\FI[v_1,\compactldots,v_n])$, viewing $\Lambda^2 \FC$ as a summand of $\FC \otimes \FC$ so that $\id - \Lambda^2 D_h$ is the restriction of $\id \otimes \id-D_h \otimes D_h$. Here the meaning of $D_h$ is as in \cref{lem:cyc-dec-h}: $D_h \colon \FC(U) \to \FC(U)$ for $U \subseteq V$ is defined to be $0$ if $U \subseteq \ker(h)$ and $D_{h|_U}^{\FC}$ otherwise; in the generic case one only encounters the latter. 

We can write this as
\[\id \otimes \id - D_h \otimes D_h =  \id \otimes (\id - D_h) + (\id - D_h) \otimes \id-(D_h - \id) \otimes (D_h -\id).\]
When restricting to generic pairs, there are then lifts of the right terms to $\FC \otimes \FCR$, $\FCR \otimes \FC$, and $\FCR \otimes \FC$ respectively (note that we have made a choice in the last term):

\begin{definition}\label{def:rho-cyc} The element $\rho^\rm{cyc}([h] \otimes \FI[v_1,\compactldots,v_n]) \in (\FCR \otimes \FC+\FC \otimes \FCR)(V)$ is defined as follows. Write $x = (\delta_\rm{cyc} \circ \pr^\FC)(\FI[v_1,\compactldots,v_n])$, and take the sum of terms
\[(\id \otimes [\id - D_h])(x)+([\id - D_h] \otimes \id)(x)-([D_h - \id] \otimes (D_h -\id))(x),\]
where the square brackets indicate an element is to be considered as a relation.
\end{definition}

Note we crucially use that $[h] \otimes \FI[v_1,\compactldots,v_n]$ is generic: if not, on some terms $\id-D_h$ (or $D_h - \id$) could equal $\id$ (or $-\id$) and would not give a relation.

\subsubsection{Choosing relations for the ``coaction'' term} Assume still that $[h] \otimes \FI[v_1,\compactldots,v_n]$ is generic. For the nonexplicit lift of $\delta^\rm{coact}_h(\FI[v_1,\compactldots,v_n])$, we start by recalling that the map $\pr^\St \colon \FC \to \SStL$ is one of Lie coalgebras, and hence $\delta_\rm{cyc} \circ \pr^\FC \circ (\id - D_h^{\FI})$ maps $\FI[v_1,\compactldots,v_n]$ to zero in $(\Lambda^2 \SStL)(V)$. The element $\delta^\rm{cyc}_h(\FI[v_1,\compactldots,v_n])$ also maps to zero in $(\Lambda^2 \SStL)(V)$ since it admits a lift to $(\FCR \otimes \FC+\FC \otimes \FCR)(V)$, namely as in \cref{def:rho-cyc}. 

Thus we must have that the same is true for the remaining terms:
\begin{equation}\label{eqn:coact-lift-1} \text{$\delta^\rm{coact}_h(\FI[v_1,\compactldots,v_n]) \in (\Lambda^2 \FC)(V)$ maps to zero in $(\Lambda^2 \SStL)(V)$.}\end{equation}
We next separate the different components of $\delta^\rm{coact}_h(\FI[v_1,\compactldots,v_n])$. To do so, recall a nonzero linear functional $h$ is fixed and write $v_0 \coloneq 0$, $H = \ker(h) \subset V$.

\begin{lemma}The element
\[(\overline{\delta}_\rm{coact} \circ D_h^{\FI})(\FI[v_1,\compactldots,v_n])\]
is supported in the summands of $(\Lambda^2 \FC)(V)$ indexed by $\{V_1,V_2\}$ where $V_1 \not \subseteq H$ and $V_2 \subseteq H$.
\end{lemma}

\begin{proof}By the iterative formula \eqref{eqn:symbol-iterated-integral} for the symbol, $D^{\FI}_h(\FI[v_1,\compactldots,v_n])$ is a sum of terms $\FI[w_1,\compactldots,w_n]$ for bases $w_1,\ldots,w_n$ where each $w_i$ is of the form $\frac{v_i-v_j}{h_i-h_j}$ for $v_i-v_j \not \in \ker(h)$. The coaction is given by taking only certain terms of $\delta_\rm{cyc}(\FC[0:w_1:\compactldots:w_n])$, e.g.~in
\[\begin{tikzpicture}
   \draw (0:1) \foreach \x in {60,120,...,360} {  -- (\x:1) };
   \foreach \x/\l/\p in
     { 60/{$w_2$}/above,
      120/{$w_1$}/above,
      180/{$0$}/left,
      240/{$w_5$}/below,
      300/{$w_4$}/below,
      360/{$w_3$}/right
     }
     \node[inner sep=1pt,circle,draw,fill,label={\p:\l}] at (\x:1) {};
\end{tikzpicture}\]
we cut along lines that do \emph{not} pass through zero. The component containing the $0$th vertex will lie in a term $V_1$ spanned by $w_i$ and since each of these is not contained in $H$, we have $V_1 \not \subseteq H$. The component not containing the $0$th vertex lies in a term $V_2$ spanned by $w_i-w_j$ and since each of these lies in $H$ as $h(w_i) = 1 = h(w_j)$, we have $V_2 \subseteq H$.
\end{proof}

We now recall that if a subspace $V_1 \not \subseteq H$ of $V$ is of dimension $d \leq n$ then the map $\StL(V_1) \to (B^\rm{Com} \SSt)_d^H(V_1)$ (the target is as in the proof of \cref{lem:cyc-dec-h}) is an isomorphism with inverse $\rm{C}^\St_h$ given by sending $[P_1|\compactldots|P_d]$ to the Steinberg correlator $\rm{C}[0:v_1:\compactcdots:v_d]$ where $v_i \in P_i$ is the unique vector so that $h(v_i) = 1$. In particular, such $\rm{C}[0:v_1:\compactcdots:v_d]$ almost form a basis of $\StL(V_1)$; the only relations between them are the shuffle relations. We will refer to these as well as the formal correlators $\rm{FC}[0:v_1:\compactcdots:v_d]$ that map to them, as \emph{basic}. The crucial observation---a direct consequence of this discussion---is
\begin{equation}\label{eqn:coact-lift-2}\text{independent basic formal correlators in $\FC(V_1)$ remain independent in $\StL(V_1)$.}\end{equation}

In the previous lemma, the terms supported at $\{V_1,V_2\}$ where $V_1 \not \subseteq H$ and $V_2 \subseteq H$ are linear combinations of elements of the form $\rm{FC}[0:v_1,\ldots,v_d] \otimes \rm{FC}[v'_1:\ldots:v'_{n-d}]$ where the first term is basic. Using facts \eqref{eqn:coact-lift-1} and \eqref{eqn:coact-lift-2}, we can thus express the component of $\delta^\rm{coact}_h(\FI[v_1,\ldots,v_d])$ in the term indexed by $\{U,W\}$, where we must have $U \not \subseteq H$ and $W \subseteq H$ by the above lemma, as $\sum_\alpha r_\alpha^{U,W} \otimes x_\alpha^U$ where the $r_\alpha^{U,W}$ are in $\FCR(W)$ and the $x_\alpha^U \in \FC(U)$ are basic. Doing this for all components we obtain the desired choice of lift:

\begin{definition}\label{def:rho-coact}The element $\rho^\rm{coact}([h] \otimes \FI[v_1,\ldots,v_n]) \in (\FCR \otimes \FC)(V) \subset (\FCR \otimes \FC+\FC \otimes \FCR)(V)$ is defined as follows. Write $x = \delta^\rm{coact}_h(\FI[v_1,\compactldots,v_n])$ and express it as
\[\sum_{W,U} \sum_\alpha [r_\alpha^{U,W}] \otimes x_\alpha^U,\]
where the square brackets indicate an element is to be considered as a relation.  
\end{definition}

Finally, we define $\rho$ as
\begin{align*}\rho \colon B_1(\Dec_V;\StL(V)) &\lra (\FCR \otimes \FC+\FC \otimes \FCR)(V)\\
[h] \otimes \FI[v_1,\compactldots,v_n] &\longmapsto \begin{cases} (\rho^\rm{cyc} - \rho^\rm{coact})([h] \otimes \FI[v_1,\compactldots,v_n]) & \text{if the pair is generic,} \\
\text{some unspecified choice} & \text{else,}\end{cases}\end{align*}
where the second case is justified in \cref{sec:cobracket-choosing-relations}.

\begin{example}\label{exam:rho-coact} Let us compute $\rho^\rm{coact}$ on $\FI[v_1,v_2]$ with $v_1,v_2 \notin H = \ker(h)$. By the proof of \cref{prop: G_2=B_2} we have that $D^\FI_h(\FI[v_1,v_2])$ is given by 
\[\FI[\tfrac{v_1}{h_1},\tfrac{v_2}{h_2}]-\FI[\tfrac{v_1}{h_1},\tfrac{v_{2}-v_1}{h_2-h_1}]+\FI[\tfrac{v_2}{h_2},\tfrac{v_{2}-v_{1}}{h_2-h_1}].\]
The map $\ol{\delta}_\rm{coact}$ sends this to
\begin{align*}&-\FC[0:\tfrac{v_1}{h_1}]\wedge \FC[\tfrac{v_1}{h_1}:\tfrac{v_2}{h_2}]+\FC[0:\tfrac{v_2}{h_2}]\wedge \FC[\tfrac{v_1}{h_1}:\tfrac{v_2}{h_2}]\\
&+\FC[0:\tfrac{v_1}{h_1}]\wedge \FC[\tfrac{v_1}{h_1}:\tfrac{v_{2}-v_1}{h_2-h_1}]-\FC[0:\tfrac{v_{2}-v_1}{h_2-h_1}]\wedge \FC[\tfrac{v_1}{h_1}:\tfrac{v_{2}-v_1}{h_2-h_1}]\\
&-\FC[0:\tfrac{v_2}{h_2}] \wedge \FC[\tfrac{v_2}{h_2}:\tfrac{v_{2}-v_{1}}{h_2-h_1}]+\FC[0:\tfrac{v_{2}-v_{1}}{h_2-h_1}] \wedge \FC[\tfrac{v_2}{h_2}:\tfrac{v_{2}-v_{1}}{h_2-h_1}]\end{align*}
where all of the first terms in the wedge products are supported at a $U \not \subseteq H$ and the second terms at $H$. We can collect the second terms as relations, in this case rather easily because $\StL(V) \cong \bb{Q}$ if $V$ is 1-dimensional. Denoting these first in square brackets and then tensoring them with basic elements second, we get that $\rho^\rm{coact}(\FI[v_1,v_2]) \in (\FCR \otimes \FC)(V)$ is given by
\begin{align*}&\big[\FC[\tfrac{v_1}{h_1},\tfrac{v_2}{h_2}]-\FC[\tfrac{v_1}{h_1}:\tfrac{v_{2}-v_1}{h_2-h_1}]\big] \otimes \FC[0:\tfrac{v_1}{h_1}] \\
&+\big[\FC[\tfrac{v_2}{h_2}:\tfrac{v_{2}-v_{1}}{h_2-h_1}]-\FC[\tfrac{v_1}{h_1}:\tfrac{v_2}{h_2}]\big] \otimes \FC[0:\tfrac{v_2}{h_2}] \\
&+\big[\FC[\tfrac{v_1}{h_1}:\tfrac{v_{2}-v_1}{h_2-h_1}]-\FC[\tfrac{v_2}{h_2}:\tfrac{v_{2}-v_{1}}{h_2-h_1}]\big] \otimes \FC[0:\tfrac{v_{2}-v_1}{h_2-h_1}].
\end{align*}\end{example}

\subsubsection{The map $B_* \to A_*$} In this subsection we define the chain complexes that are partial resolutions of the first two chain complexes that appear in the computation of the cobracket
\[A_* \lra [\StL(V) \to 0] \qquad \text{and} \qquad B_* \lra [\StL(V) \to (\Lambda^2 \SStL)(V)]\]
and construct a map from the latter to the former.

\medskip

We first construct $A_*$. Recall we have the exact sequence
\[\cdots \lra B_2(\Dec_V;\FI(V)) \lra B_1(\Dec_V;\FI(V)) \lra \FI(V) \lra \StL(V) \lra 0,\]
and think of this as a resolution of $\StL(V)$ by a chain complex 
\begin{equation}\label{eqn:res-cobracket-i} A_* \coloneq [\cdots \to B_2(\Dec_V;\FI(V)) \to B_1(\Dec_V;\FI(V)) \to \FI(V)],\end{equation}
where we place the entry $\FI(V)$ in degree $0$. As a consequence the map
\[A_* \overset{\simeq}\lra [\StL(V) \to 0]\]
whose target is thought of as a chain complex with entry $\StL(V)$ placed in degree $0$, is a quasi-isomorphism.

\medskip

We next construct $B_*$. This will arise from a double complex obtained from a first row $A_*$ by adding a second row given by the second exterior power of the resolution $0 \to \FCR(V) \to \FC(V) \to \StL(V) \to 0$. For this we recall that more generally if 
\[0 \lra A \overset{f}\lra B \overset{g}\lra C \lra 0\]
is exact then so is
\[0 \lra \rm{Sym}^2 A \lra \tfrac{A \otimes B+B \otimes A}{\{a \otimes b+b \otimes a\}} \lra \Lambda^2 B \lra \Lambda^2 C \lra 0\]
with maps given respectively by $aa' \mapsto a \otimes f(a')-f(a) \otimes a'$, by $a \otimes b+b'\otimes a' \mapsto f(a) \wedge b+b' \wedge f(a')$, and by $b \wedge b' \mapsto g(b) \wedge g(b')$. Thus both rows in the following double complex are exact:
\[\begin{tikzcd}[column sep=.5cm] \cdots \rar & B_2(\Dec_V;\FI(V)) \rar \dar{\rho'} & B_1(\Dec_V;\FI(V)) \rar \dar{\alt \circ \rho} & \FI(V) \rar \dar{\delta_\rm{cyc} \circ \pr^\FC} & \StL(V) \dar{\delta}\\[-5pt]
0 \rar & (\rm{Sym}^2 \FCR)(V) \rar & \frac{(\FCR \otimes \FC+\FC \otimes \FCR)}{\{a \otimes b + b \otimes a\}}(V) \rar & (\Lambda^2 \FC)(V) \rar & (\Lambda^2 \StL)(V)
\end{tikzcd}\]
where the quotient on the bottom identifies $a \otimes b \in \FCR(U) \otimes \FC(W)$ for $U \oplus W = V$ with $b \otimes a \in \FC(W) \otimes \FCR(U)$ for a direct sum decomposition $U \oplus W = V$. The right square commutes because $\pr^\St$ is a map of Lie coalgebras, that the middle square commutes follows from \eqref{eqn:choosing-relations-diag} as its bottom map factors over the antisymmetrisation of the domain, and there exists a map $\rho'$ making the right square commute since its domain is part of a projective resolution. Taking total complexes we define
\begin{equation}\label{eqn:res-cobracket-ii} B_* \coloneq \left[\compactcdots \to \parbox{2.8cm}{\centering$(\rm{Sym}^2 \FCR)(V)$ \\[-2pt]
$+B_1(\Dec_V;\FI(V))$} \to \tfrac{(\FCR \otimes \FC{+}\FC \otimes \FCR)}{\{a \otimes b {+} b \otimes a\}}(V){+}\FI(V) \to (\Lambda^2 \FC)(V)\right]
\end{equation} 
with a quasi-isomorphism $B_* \overset{\simeq}\lra [\StL(V) \to (\Lambda^2\StL)(V)]$.

\medskip

We finally construct the map $B_* \to A_*$ of chain complexes of $\bb{Q}[\GL(V)]$-modules as that induced by projection onto the first row of the double complex.

\subsubsection{The map $B_* \to C_*$} \label{sec:map-b-to-c} We now define $C_*$. By taking $\para$-tensor products of two copies of the resolution $0 \to \FCR \to \rm{FC} \to \SStL \to 0$ and evaluating at $V$, we get a chain complex
\begin{equation}\label{eqn:res-cobracket-iii}C_* = \big[(\FCR \para \FCR)(V) \to (\FCR \para \FC+\FC \para \FCR)(V) \to (\FC \para \FC)(V)\big]\end{equation}
with a quasi-isomorphism $C_* \overset{\simeq}\lra [0 \to (\StL \para \StL)(V)]$. For this we observe that if 
\[0 \lra A \overset{f}\lra B \overset{g}\lra C \lra 0\]
is exact then so is
\[0 \lra A \para A \lra A \para B + B \para A \lra B \para B \lra C \para C \lra 0\]
with maps given respectively by $a \para a' \mapsto a \otimes f(a')-f(a) \para a'$, by $a \para b+b'\para a' \mapsto f(a) \para b+b' \para f(a')$, and by $b \para b' \mapsto g(b) \para g(b')$. 

\medskip

To construct the map $B_* \to C_*$, we recall that there are maps to the parabolic tensor product:
\[\FCR \otimes \FC+\FC \otimes \FCR \overset{\alt}\lra \frac{\FCR \otimes \FC+\FC \otimes \FCR}{\{a \otimes b + b\otimes a\}} \xrightarrow{\zeta^\alt} \FCR \para \FC+\FC \para \FCR,\]
where $\alt$ is the projection and we have 
\begin{align*}\zeta^\alt \coloneq \tfrac{1}{2}(\ol{\zeta}-\ul{\zeta}) \colon \frac{X \otimes Y+Y \otimes X}{\{x \otimes y+y \otimes x\}}&\lra X \para Y + Y \para X \\
x \otimes y + y' \otimes x' &\longmapsto \tfrac{1}{2}(x \para y-x' \para y' +y'\para x'-y \para x).\end{align*}

Observe that the inclusion $\FCR \to \rm{FC}$ induces a canonical map $\FCR \para \FC+\FC \para \FCR \to \FC \para \FC$, which fits into a commutative diagram
\[\begin{tikzcd} \FCR \otimes \FC+\FC \otimes \FCR \dar[swap]{\zeta^\alt \circ \alt} \rar & \Lambda^2 \FC \dar{\zeta^\alt} \\[-5pt]
\FCR \para \FC+\FC \para \FCR \rar & \FC \para \FC.\end{tikzcd}\]

\begin{lemma}The image of $(\zeta^\alt \circ \alt \circ \rho)(h \otimes \FI[v_1,\ldots,v_n]) \in (\FCR \para \FC+\FC \para \FCR)(V)$ in $\FC \para \FC$ vanishes.
\end{lemma}

\begin{proof}This follows from the commutative diagram
\[\begin{tikzcd} B_1(\Dec_V,\rm{FI}(V)) \rar \dar[swap]{\rho} & \FI(V)  \dar{\delta_\rm{cyc} \circ \pr^\FC} \\[-5pt]
(\FCR \otimes \FC+\FC \otimes \FCR)(V) \rar \dar[swap]{\zeta^\alt \circ \alt} & (\Lambda^2\FC)(V)  \dar{\zeta^\alt} \\[-5pt]
(\FCR \para \FC+\FC \para \FCR)(V) \rar & (\FC \para \FC)(V)\end{tikzcd}\]
obtained from \eqref{eqn:choosing-relations-diag} and the commutative square preceding the statement, and the fact that $\zeta^\alt \circ \delta_\rm{cyc} = 0$ by \cref{lem:fc-cobracket-symmetry}.
\end{proof}

We now define the map $B_* \to C_*$ on terms coming from each row of the double complex defining $B_*$ separately. For the bottom row it is 
\[\begin{tikzcd}0 \rar & (\rm{Sym}^2 \FCR)(V) \rar \dar{\zeta^\sym} & \frac{(\FCR \otimes \FC+\FC \otimes \FCR)}{\{a \otimes b + b \otimes a\}}(V) \rar \dar{\zeta^\alt} & (\Lambda^2 \FC)(V) \dar{\zeta^\alt} \\[-5pt]
0 \rar & (\FCR \para \FCR)(V) \rar  & (\FCR \para \FC+\FC \para \FCR)(V) \rar & (\FC \para \FC)(V) 
\end{tikzcd}\]
where we have
\begin{align*} \zeta^\sym \coloneq \tfrac{1}{2}(\ol{\zeta}+\ul{\zeta}) \colon \rm{Sym}^2 X &\lra X \para X \\[-5pt]
xy &\longmapsto \tfrac{1}{2}(x \para y + y \para x).\end{align*}
For the top row it is 
\[\begin{tikzcd}\cdots \rar & B_1(\Dec_V;\FI(V)) \rar \dar{\zeta^\alt \circ \alt \circ \rho} & \FI(V)  \dar{0} \rar & 0 \dar{0} \\[-5pt]
0 \rar & (\FCR \para \FCR)(V) \rar  & (\FCR \para \FC+\FC \para \FCR)(V) \rar & (\FC \para \FC)(V) 
\end{tikzcd}\]
using that $\zeta^\alt \circ \alt \circ \rho$ has a unique lift into $(\FCR \para \FCR)(V)$ using the previous lemma. We need to check this induces a map of chain complexes: compatibility with differentials starting at $B_1(\Dec_V;\FI(V))$ is evident, and compatibility with differentials starting at $B_2(\Dec_V;\FI(V))$ is a consequence of the definition of $\rho'$.

\subsubsection{Projecting away the $\sigma$-component} \label{sec:map-c-to-fcr} By construction of $C_*$ as in \eqref{eqn:res-cobracket-iii}, there is a map of chain complexes projecting onto the term $(\FCR \para \FCR)(V)$. The short exact sequence $0 \to \FCR(V) \to \FC(V) \to \StL(V) \to 0$ induces a connecting homomorphism
\[H_*(\GL(V);\StL(V)) \overset{\partial}\lra H_{*-1}(\GL(V);\FCR(V)),\]
which for $\dim(V) \neq 0$ is an isomorphism and for $\dim(V) = 1$ exactly projects to zero the term $\bb{Q}\{\sigma\}$ in degree $*=0$. The equivalence $[0 \to \StL(V)] \simeq [\FCR(V) \to \FC(V)]$ induces an isomorphism on homology and, tracing through the definitions, following this by projection onto $\FCR(V)$ exactly implements the construction of the connecting homomorphism. Thus, using the K\"unneth theorem and Shapiro's lemma, the maps
\[\bigoplus_{d'+d''=2} H_{d'}(\GL_k;\StL_k) \otimes H_{d''}(\GL_{n-k};\StL_{n-k}) \lra \scr{G}_k \otimes \scr{G}_{n-k}.\]
induced by $C_* \to (\FCR \para \FCR)(V)$, are also given by projecting terms $\bb{Q}\{\sigma\}$ to zero.

\subsubsection{A vanishing result and computation} Assume that $[h] \otimes \FI[v_1,\ldots,v_n]$ is generic. We start with making precise the second technical input \eqref{enum:technical-cobracket-ii}, by understanding the unique lift of the element $(\zeta^\alt \circ \alt \circ \rho)([h] \otimes \FI[v_1,\ldots,v_n])$ to $(\FCR \para \FCR)(V)$. This is done by combining the injectivity of the map $\FCR \para \FCR \to \FCR \para \FC$ with the following vanishing result:

\begin{lemma}\label{lem:xi-rhoact-cancellation} The image of $(\zeta^\alt \circ \alt \circ \rho^\rm{coact})([h] \otimes \FI[v_1,\ldots,v_n])$ in the component $(\FCR \para \FC)(V)$ vanishes.
\end{lemma}

\begin{proof}It suffices to verify this after composition with the injection $(\FCR \para \FC)(V) \to (\FC \para \FC)(V)$. Looking at \cref{def:rho-coact} and the definition of $\zeta^\alt \circ \alt$, the image will be supported at flags $0 \subseteq W \subseteq V$ with $W \subseteq H$ so there exists $U \not\subseteq H$ with $U \smash{\overset{\cong}\lra} V/W$. Moreover, up to a factor of $\tfrac{1}{2}$, there it is given by image in $\FC(W) \otimes \FC(V/W)$ under $\ol{\zeta}$ of the contributions to $\delta^\rm{coact}_h(\FI[v_1,\compactldots,v_n])$ at $W$ and such $U$. In particular, it suffices to prove that contributions to $\delta^\rm{coact}_h(\FI[v_1,\compactldots,v_n])$ at a given $\{U,W\}$ cancel pairwise when mapped to $\FC(W) \otimes \FC(V/W)$. 

To see this, observe the contributions are of two types: ``vertex first'' ones go through a vertex $v$ ($1 \leq v \leq n$) and an edge $(e,e+1)$ ($v+1 \leq e \leq n)$ and ``vertex last'' ones go through a vertex $v'$ ($1 \leq v' \leq n$) and edge $(e',e'+1)$ ($0 \leq e' \leq v'-2$). We pair the ``vertex first'' contribution $\{v,(e,e+1)\}$ to the ``vertex last'' contribution $\{v',(e',e'+1)\} = \{e,(v-1,v)\}$; indeed, $v+1 \leq e$ in the former is equivalent to $e' = v-1 \leq e-2 = v'-2$ in the latter; see the following examples of pairs:
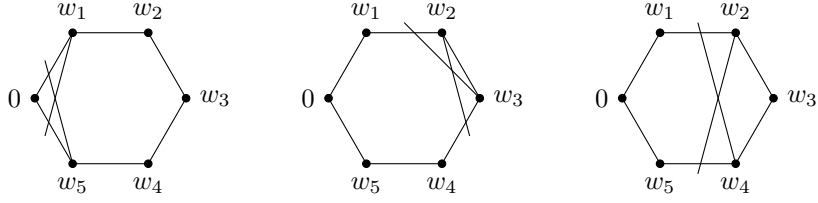
\[\begin{tikzpicture}
   \draw (0:1) \foreach \x in {60,120,...,360} {  -- (\x:1) };
   \foreach \x/\l/\p in
     { 60/{$w_2$}/above,
      120/{$w_1$}/above,
      180/{$0$}/left,
      240/{$w_5$}/below,
      300/{$w_4$}/below,
      360/{$w_3$}/right
     }
     \node[inner sep=1pt,circle,draw,fill,label={\p:\l}] at (\x:1) {};
    \draw (120:1) -- (210:1);
    \draw (150:1) -- (240:1);
\end{tikzpicture} \qquad 
\begin{tikzpicture}
   \draw (0:1) \foreach \x in {60,120,...,360} {  -- (\x:1) };
   \foreach \x/\l/\p in
     { 60/{$w_2$}/above,
      120/{$w_1$}/above,
      180/{$0$}/left,
      240/{$w_5$}/below,
      300/{$w_4$}/below,
      360/{$w_3$}/right
     }
     \node[inner sep=1pt,circle,draw,fill,label={\p:\l}] at (\x:1) {};
    \draw (60:1) -- (330:1);
    \draw (360:1) -- (90:1);
\end{tikzpicture}
\qquad 
\begin{tikzpicture}
   \draw (0:1) \foreach \x in {60,120,...,360} {  -- (\x:1) };
   \foreach \x/\l/\p in
     { 60/{$w_2$}/above,
      120/{$w_1$}/above,
      180/{$0$}/left,
      240/{$w_5$}/below,
      300/{$w_4$}/below,
      360/{$w_3$}/right
     }
     \node[inner sep=1pt,circle,draw,fill,label={\p:\l}] at (\x:1) {};
    \draw (60:1) -- (270:1);
    \draw (300:1) -- (90:1);
\end{tikzpicture}\]
So let $\{v,(e,e+1)\}$ and $\{e,(v-1,v)\}$ be a paired vertex-first and vertex-last contribution. Both have the same associated $W \subseteq H$ and the associated $U,U' \not \subseteq H$ have the same image in $V/W$, as $w_v \equiv w_e$ there. Moreover, the formal correlators from the term not containing $0$ are literally the same and those from the term containing $0$ become the same after projection to $V/W$ up to a single sign (this uses that putting first the element clockwise from the vertex has a positive sign, and putting it second has a negative sign).
\end{proof}

\begin{example}We work out the right-most example: we have $W = \rm{span}(w_2-w_3,w_3-w_4)$, $U = \rm{span}(w_1,w_2,w_5)$, $U' = \rm{span}(w_1,w_4,w_5)$, and the contributions to the cobracket at $\{W,U\}$ and $\{W,U'\}$ are respectively
\[\FC[w_2:w_3:w_4] \wedge \FC[w_2:w_5:0:w_1] \qquad \FC[w_4:w_5:0:w_1] \wedge \FC[w_4:w_2:w_3].\]
As $w_2$ and $w_4$ map to the same vector in $U/W$ and formal correlators are cyclically invariant, these cancel.
\end{example}

It will not in general be the case that the image of $(\zeta^\alt \circ \alt \circ \rho^\rm{coact})([h] \otimes \FI[v_1,\ldots,v_n])$ in the component $(\FC \para \FCR)(V)$ vanishes:

\begin{example}We continue \cref{exam:rho-coact}. Applying $\zeta^\alt \circ \alt$ to $\rho^\rm{coact}(\FI[v_1,v_2])$, we obtain (up to sign and a factor of $\tfrac{1}{2}$) terms in both $(\FCR \para \FC)(V)$ and $(\FC \para \FCR)(V)$. The former are given by
\begin{align*}&\big[\FC[\tfrac{v_1}{h_1}:\tfrac{v_2}{h_2}]-\FC[\tfrac{v_1}{h_1}:\tfrac{v_{2}-v_1}{h_2-h_1}]\big] \otimes \FC[0:\ol{\tfrac{v_1}{h_1}}] \\
&+\big[\FC[\tfrac{v_2}{h_2}:\tfrac{v_{2}-v_{1}}{h_2-h_1}]-\FC[\tfrac{v_1}{h_1}:\tfrac{v_2}{h_2}]\big] \otimes \FC[0:\ol{\tfrac{v_2}{h_2}}] \\
&+\big[\FC[\tfrac{v_1}{h_1}:\tfrac{v_{2}-v_1}{h_2-h_1}]-\FC[\tfrac{v_2}{h_2}:\tfrac{v_{2}-v_{1}}{h_2-h_1}]\big] \otimes \FC[0:\ol{\tfrac{v_{2}-v_1}{h_2-h_1}}]
\end{align*}
where the overlines mean we apply the map to the quotient; as $\ol{\tfrac{v_1}{h_1}} = \ol{\tfrac{v_2}{h_2}} = \ol{\tfrac{v_{2}-v_1}{h_2-h_1}}$ we can collect the terms of the left and cancel them pairwise. This is the cancellation used in \cref{lem:xi-rhoact-cancellation}. 

The latter are given (up to sign and a factor of $\tfrac{1}{2}$) by 
\begin{align*}& \FC[0:\tfrac{v_1}{h_1}] \otimes \big[\FC[\ol{\tfrac{v_1}{h_1}}:\ol{\tfrac{v_2}{h_2}}]-\FC[\ol{\tfrac{v_1}{h_1}}:\ol{\tfrac{v_{2}-v_1}{h_2-h_1}}]\big]\\
&+\FC[0:\tfrac{v_2}{h_2}] \otimes \big[\FC[\ol{\tfrac{v_2}{h_2}}:\ol{\tfrac{v_{2}-v_{1}}{h_2-h_1}}]-\FC[\ol{\tfrac{v_1}{h_1}}:\ol{\tfrac{v_2}{h_2}}]\big] \\
&+\FC[0:\tfrac{v_{2}-v_1}{h_2-h_1}] \otimes \big[\FC[\ol{\tfrac{v_1}{h_1}}:\ol{\tfrac{v_{2}-v_1}{h_2-h_1}}]-\FC[\ol{\tfrac{v_2}{h_2}}:\ol{\tfrac{v_{2}-v_{1}}{h_2-h_1}}]\big]
\end{align*}
and these do \emph{not} cancel, e.g.~the top-right term is equal to $\FC[0:\ol{\tfrac{v_2}{h_2}}-\ol{\tfrac{v_1}{h_1}}]-\FC[0:\ol{\tfrac{v_{2}}{h_2-h_1}}-\ol{\tfrac{v_1}{h_1}}]$.
\end{example}

To understand $(\zeta^\alt \circ \alt \circ \rho)([h] \otimes \FI[v_1,\ldots,v_n])$ we may look at its image in the $(\FCR \para \FC)(V)$-component. By \cref{lem:xi-rhoact-cancellation} this agrees with the image of $(\zeta^\alt \circ \alt \circ \rho^\rm{cyc})([h] \otimes \FI[v_1,\ldots,v_n])$ in that component. Looking at \cref{def:rho-cyc} and the definition of $\zeta^\alt \circ \alt$, it is given by a sum of four terms (writing $v_0 \coloneq 0$):
\begin{equation}\label{eqn:xi-rho-cyc} \begin{aligned}&\tfrac{1}{2}\sum_{j=0}^n \sum_{i=1}^{n-1} \big[(\id-D_h)(\FC[v_j:\compactcdots:v_{j+i}])\big] \para \FC[\ol{v}_j:\ol{v}_{j+i+1}:\compactcdots:\ol{v}_{j+n}] \\
&\quad -\tfrac{1}{2}\sum_{j=0}^n \sum_{i=1}^{n-1} \big[(\id-D_h)(\FC[v_j:v_{j+i+1}:\compactcdots:v_{j+n}])\big] \para \FC[\ol{v}_j:\compactcdots:\ol{v}_{j+i}] \\
&\quad +\tfrac{1}{2}\sum_{j=0}^n \sum_{i=1}^{n-1} \big[(\id-D_h)(\FC[v_j\compactcdots:v_{j+i}])\big] \para \ol{(\id-D_h)(\FC[v_j:v_{j+i+1}:\compactcdots:v_{j+n}])}\\
&\quad -\tfrac{1}{2}\sum_{j=0}^n \sum_{i=1}^{n-1} \big[(\id-D_h)(\FC[v_j:v_{j+i+1}:\compactcdots:v_{j+n}])\big] \para \ol{(\id-D_h)(\FC[v_j:\compactcdots:v_{j+i}])}\end{aligned}\end{equation}
where the terms in square brackets are interpret as relations, and the overlines indicate we apply naturality in the isomorphism $W \to V \to V/U$ for a splitting $U \oplus W \overset{\cong}\lra V$. Including \eqref{eqn:xi-rho-cyc} further in $(\FC \para \FC)(V)$, we note that first two terms in \eqref{eqn:xi-rho-cyc} are $((\id-D_h) \otimes \id) \circ \overline{\zeta} \circ \delta_\rm{cyc}$ and we already know $\overline{\zeta} \circ \delta_\rm{cyc} = 0$ by \cref{lem:fc-cobracket-symmetry}. The result is:

\begin{proposition}\label{prop:para-after-rho-computation} The image of $[h] \otimes \FI[v_1:\compactldots:v_n]$ in $(\FCR \para \FCR)(V)$ is given by
\begin{align*} &\tfrac{1}{2} \sum_{j=0}^n \sum_{i=1}^{n-1} \big[(\id-D_h)(\FC[v_j:\compactcdots:v_{j+i}])\big] \otimes \big[\ol{(\id-D_h)(\FC[v_j:v_{j+i+1}:\compactcdots:v_{j+n}])}\big] \\
&-\tfrac{1}{2} \sum_{j=0}^n \sum_{i=1}^{n-1} \big[(\id-D_h)(\FC[v_j:v_{j+i+1}:\compactcdots:v_{j+n}])\big] \para \big[\ol{(\id-D_h)(\FC[v_j:\compactcdots:v_{j+i}])}\big].\end{align*}
\end{proposition}

\subsubsection{Finishing the computation}
This allows us to finish the proof of \cref{thm:polyl-presentation-cobracket}. At this point we know that the cobracket is induced by the zigzag
\[H_1(\GL(V);A_*) \longleftarrow H_1(\GL(V);B_*) \lra H_1(\GL(V);C_*),\]
which we know is isomorphic to the restriction of the zigzag
\[\rm{H}/\rm{H}^2 \longleftarrow \rm{H}/\rm{H}^3 \lra \rm{H}/\rm{H}^2 \otimes \rm{H}/\rm{H}^2\]
to rank $n = \dim(V)$ and degree $2n+1$, followed by discarding terms corresponding to $\sigma$. In particular, the left map is surjective and upon passing to the quotient $\Lambda^2(\rm{H}/\rm{H}^2)$ of the right term by antisymmetrising, the result is independent of a choice of lift. It remains to explicitly perform the following four steps in explicit models for group homology: 
\begin{enumerate}[(A)]
    \item \label{enum:cobracket-a} find representatives in the left term, 
    \item \label{enum:cobracket-b} lift them to the middle term, 
    \item \label{enum:cobracket-c} apply the map to the right term,
    \item \label{enum:cobracket-d} discard the components corresponding to $\sigma$.
\end{enumerate}

\begin{proof}[Proof of \cref{thm:polyl-presentation-cobracket}] Denote the canonical bar resolution $\bb{Q}[\GL(V)]^{\otimes \bullet}$ of $\bb{Q}$ by $P_*$, then we choose our explicit model to be total chain complexes of the zigzag of double complexes
\[A_* \otimes_{\bb{Q}[\GL(V)]} P_* \longleftarrow B_* \otimes_{\bb{Q}[\GL(V)]} P_* \lra C_* \otimes_{\bb{Q}[\GL(V)]} P_*.\]
We will denote the differential in the first term by $d$ and that in the second term by $\partial$.

\smallskip

For step \eqref{enum:cobracket-a}, recall that our goal is to compute the cobracket of a correlator $\CorG(x_0,\compactldots,x_n) \in \scr{G}(V)$ with  all arguments $x_0,\ldots,x_n$ all distinct, and by homogeneity we may assume $x_0 = 0$. To do so, we must represent it by an 1-cycle in $A_* \otimes_{\bb{Q}[\GL(V)]} P_*$, the relevant part of which is
\[\begin{tikzcd}[column sep=.2cm, row sep=.25cm] & \FI(V) \otimes \bb{Q}[\GL(V)] \arrow{rd} & \\
B_1(\Dec_V;\FI(V)) \arrow{rr} & & \FI(V). \end{tikzcd}\]
We choose to the pair of the necessarily generic element $x=[h] \otimes \FI[v_1,\compactldots,v_n] \in B_1(\Dec_V;\FI(V))$ satisfying $h(v_i) = x_i$ and an $z \in \FI(V) \otimes \bb{Q}[\GL(V)]$ satisfying $d(x) + \partial(z) = 0 \in \FI(V)$, which must exist because the coinvariants of $\FI(V)$ vanish.

\smallskip

For step \eqref{enum:cobracket-b}, we must extend this to a 1-cycle in $B_* \otimes_{\bb{Q}[\GL(V)]} P_*$, the relevant part of which is
\[\begin{tikzcd}[column sep=-1cm, row sep=.25cm] & & & \FI(V) \otimes \bb{Q}[\GL(V)] \arrow{rd} & \\
& & B_1(\Dec_V;\FI(V)) \arrow{rr} \arrow{dd} & & \FI(V) \\ 
& \frac{(\FCR \otimes \FC+\FC \otimes \FCR)}{\{a \otimes b + b \otimes a\}}(V) \otimes \bb{Q}[\GL(V)] \arrow{rd} & & & \\
(\rm{Sym}^2 \FCR)(V) \arrow{rr} & & \frac{(\FCR \otimes \FC+\FC \otimes \FCR)}{\{a \otimes b + b \otimes a\}}(V). & & \end{tikzcd}\]
This means finding elements $y \in (\rm{Sym}^2 \FCR)(V)$ and $w \in (\tfrac{\FCR \otimes \FC + \FC \otimes \FCR}{a \otimes b + b \otimes a})(V) \otimes \bb{Q}[\GL(V)]$ so that $d(x,y)-\partial(w) = 0$, which is always possible by surjectivity of the map $\rm{H}/\rm{H}^3 \to \rm{H}/\rm{H}^2$.

\smallskip

For step \eqref{enum:cobracket-c}, we apply the map $B_* \to C_*$ to get a 1-cycle in $C_* \otimes_{\bb{Q}[\GL(V)]} P_*$, the relevant part of which is
\[\begin{tikzcd}[column sep=-.5cm, row sep=.25cm] & (\FCR \para \FC+\FC\para \FCR)(V) \otimes \bb{Q}[\GL(V)] \arrow{rd} & \\
(\FCR \para \FCR)(V) \arrow{rr} & & (\FCR \para \FC+\FC\para \FCR)(V). \end{tikzcd}\]
The result is an element $(\zeta^\alt \circ \alt \circ \rho)(x)+\zeta^\sym(y) \in (\FCR\para \FCR)(V)$ where this first term was computed in \cref{prop:para-after-rho-computation}, as well as the image $\zeta^\alt(w) \in (\FCR \para \FC + \FC \para \FCR)(V) \otimes \bb{Q}[\GL(V)]$. 

\smallskip

For step \eqref{enum:cobracket-d}, to discard the components corresponding to $\sigma$ we simply project to the term $(\FCR\para \FCR)(V)$, according to the discussion in \cref{sec:map-c-to-fcr}. (The element $\zeta^\alt(w)$ in fact yields the $\sigma$-component $\delta_\sigma$). 

\smallskip

Passing to homology, we can make an identification 
\[H_0(\GL(V);(\FCR \para \FCR)(V)) \cong (\scr{G} \otimes \scr{G})(V) \cong (\rm{Sym}^2 \scr{G})(V)+(\Lambda^2 \scr{G})(V)\]
by splitting the coefficients $(\FCR \para \FCR)(V)$ into a symmetric and antisymmetric part. The element $\zeta^\sym(y)$ contributes only to the former, while $(\zeta^\alt \circ \alt \circ \rho)(x)$ contributes only to the latter. Upon antisymmetrising only the latter will remain and is exactly given by the element in the formula of \cref{thm:polyl-presentation-cobracket}. This completes the computation.
\end{proof}

\section{The $\sigma$-component of the cobracket} \label{sec:sigma-component} In this section we compute the $\sigma$-component
\[\delta_\sigma \colon \scr{G}(F) \lra H_{2}(\GL_{n-1}(F);\StL_{n-1}(F))\]
in terms of the cobracket, or at least a part of it. More precisely, we will prove using a scaling action that the target splits into two summands and describe the projection of the $\sigma$-component to the easier summand. We follow \cref{conv:shorter-notation}.

\subsection{The scaling action} \label{sec:action-by-scaling} If we let $\rm{CMon}(\rm{Cat}_1)$ denote the $(2,1)$-category whose objects are symmetric monoidal 1-categories, whose morphisms are symmetric monoidal functors, and whose (invertible) 2-morphisms are symmetric monoidal natural isomorphisms, then the nerve construction identifies this as a full subcategory of the category of symmetric monoidal categories \cite[Section 8]{GepnerHaugsengNikolaus}. The restriction to the symmetric monoidal 1-groupoids lands in $\rm{CMon}(\Spc) \subset \rm{CMon}(\rm{Cat})$ and if we can identify the latter with $\Alg_{E_\infty^\rm{u}}(\Spc)$, the nerve is given by classifying space construction
\[B \colon \rm{CMon}(\rm{Gpd}_1) \overset{\subset}\lra \Alg_{E_\infty^\rm{u}}(\Spc).\]
To add gradings, one uses that the commutative monoid $\bb{N}$ of natural numbers under addition can be thought of as a symmetric monoidal 1-groupoid where all morphisms are the identity with associated $E_\infty^\rm{u}$-algebra in spaces denoted $\bf{N}$, and use that
\[\rm{CMon}(\rm{Gpd}_1)_{/\bb{N}} \overset{B}{\lra} \Alg_{E_\infty^\rm{u}}(\Spc)_{/\bf{N}} \simeq \Alg_{E_\infty^\rm{u}}(\Fun(\bb{N},\Spc)).\]

\medskip

From this perspective, $\BGLb(F)^+$ arises as the classifying space of the symmetric monoidal groupoid $\rm{Vect}$ of finite-dimensional vector spaces over $F$ under direct sum, with the symmetric monoidal functor to $\bb{N}$ given by $\dim$. Every $\lambda \in F^\times$ gives rise to a symmetric monoidal natural isomorphism $\phi_\lambda \colon \id_\rm{Vect} \to \id_\rm{Vect}$ with components at an object $V \in \Vect$ given by $\lambda \cdot \id_V \colon V \to V$. These satisfy $\phi_{\lambda} \phi_\mu = \phi_{\lambda \mu}$, and hence lift $\Vect$ to a functor
\[\rm{Vect} \colon \rm{B}^2 F^\times \lra \rm{CMon}(\rm{Cat}_1)_{/\bb{N}}.\]
whose domain has a unique object $\ast$ and 1-morphism $\id_\ast$, with 2-endomorphisms of $\id_\ast$ given by $F^\times$ with composition given by multiplication. Taking classifying spaces we thus lift $\BGLb(F)^+$ to a functor
\[\BGLb(F)^+ \colon \rm{B}^2 F^\times \lra \Alg_{E_\infty^\rm{u}}(\Fun(\bb{N},\Spc)).\]
Unwinding the definitions, in rank $n$ this is the following map: for an $n$-dimensional vector space $V$, the inclusion of the orbit groupoid $*{\sslash}\GL(V) \to \Vect$ induces an equivalence $\BGL(V) \to \BGLb(F)^+(n)$. On mapping spaces, the above functor then induces a map $\rm{B}F^\times \to \Map_{\Spc}(\BGL(V),\BGL(V))$ with adjoint map 
\[\rm{B}F^\times \times \BGL(V) \to \BGL(V)\]
that is induced by the group homomorphism
\begin{align*}\alpha \colon F^\times \times \GL(V) &\lra \GL(V) \\
(\lambda,A) &\longmapsto \lambda \, A.\end{align*}

Given the construction of this action on $\BGLb^+$, it is inherited by any object obtained naturally from it. In particular, postcomposing with the functor induced on $E_\infty^\rm{u}$-algebras in $\Fun(\bb{N},\Spc)$ by rationalisation $\Spc \to \DQ$, passing to the augmentation ideal of the canonical augmentation, and taking $E_\infty^\rm{nu}$-indecomposables, we obtain a functor
\[\cot_{E_\infty^\rm{nu}}(\BGLb) \colon \rm{B}^2 F^\times \lra \Fun(\bb{N},\DQ)\]
that induces a $\Lambda^* F^\times$-action on the $E_\infty$-homology.

To understand this more explicitly, we take a slightly different perspective. The symmetric monoidal $\rm{B}F^\times$-action on $\Vect$ induces one on $\Fun(\Vect,\scr{C})$ for any presentable symmetric monoidal category $\scr{C}$, and hence induces an action on $\Alg_{E_\infty^\rm{u}}(\Fun(\Vect,\scr{C}))$. Taking $\scr{C} = \Spc$, the terminal $E_\infty^\rm{u}$-algebra $\ul{\ast}$ is necessarily a fixed point of this action, and since the action is by symmetric monoidal equivalences so are the iterated bar construction of $\ul{\ast}_+$. By naturality of the identification 
\[H_{n,d}^{E_\infty}(\BGLb(F)) \cong H_{n,d-2n+2}(\GL_n;\StL_n),\]
we see that the action on the right side is given by that induced by $\alpha$, using the canonical identification $\alpha^* \StL_n \cong \pr_2^* \StL_n$. 

\subsection{Splittings} A useful feature of Steinberg modules and their variants is that the action of general linear groups on them factors over the projective general linear groups. Though it is easy to see it is a consequence of the fact that $\ul{*}$ and its iterated bar constructions are fixed points for the $BF^\times$-action. We take advantage of this to produce splittings, compatible with the scaling action, that in particular hold with trivial coefficients and coefficients in infinite Steinberg modules.

\medskip

Recall the projective general linear groups over a field $F$ are defined as follows: there is an injective homomorphism as diagonal matrices
\begin{align*} F^\times &\lra \GL_n \\
\lambda & \longmapsto \lambda\,\id\end{align*}
taking values in the centre, and we identify $F^\times$ with its image. Then we set
\[\PGL_n \coloneq \GL_n/F^\times.\]
Thus there is a canonical central extension
\[1 \lra F^\times \lra \GL_n \overset{\pr}\lra \PGL_n \lra 1.\]
Though the determinant $\det \colon \GL_n \to F^\times$ does \emph{not} split this, it does include a splitting on rational homology:

\begin{lemma}The homomorphism $(\pr,\det) \colon \GL_n \to \PGL_n(F) \times F^\times$ induces for all $\bb{Q}[\PGL_n]$-modules $M$ an isomorphism.
\begin{equation}\label{eqn:splitting} \varpi_n \colon H_*(\GL_n;\pr^* M) \overset{\cong}\lra  H_*(\PGL_n;M) \otimes \Lambda^* F^\times.\end{equation}
\end{lemma}

\begin{proof}There is a map of short exact sequences of groups
\[\begin{tikzcd} 1 \rar & F^\times \dar{(-)^n} \rar{\inc} & \GL_n \dar{(\pr,\det)} \rar{\pr} & \PGL_n \dar{\id} \rar & 1 \\
1 \rar & F^\times \rar{i_2} & \PGL_n \times F^\times \rar{\pi_1} & \PGL_n \rar & 1 \end{tikzcd}\]
and the induced map on Hochschild--Lyndon--Serre spectral sequences for homology with coefficients in $M$, is an isomorphism on the $E^2$-page, giving the result as the second spectral sequence collapses. To see this, we observe (a) that the top sequence is central, so the K\"unneth theorem provides an isomorphism $H_*(F^\times,\rm{incl}^* \rm{pr}^* M) \cong H_*(F^\times;\bb{Q}) \otimes M$ as $\Q[\PGL_n]$-modules, where the first factor has trivial action, and (b) that the maps $H_*(F^\times;\bb{Q}) \to H_*(F^\times;\bb{Q})$ induced by $\lambda \mapsto \lambda^n$, are isomorphisms as they act by multiplication with $n^k$ in degree $k$.

To see the latter, recall that if $A$ is any abelian group then $H_1(A;\bb{Q}) \cong A \otimes \Q$ and $H_*(A;\bb{Q})$ is a graded-commutative algebra with product induced by the $E_\infty$-structure on $BA$ (induced by the multiplication on $A$ viewed as a group homomorphism $A \times A \to A$). Thus, there is a canonical algebra map $\Lambda^*_\bb{Q}(A \otimes \Q) \to H_*(A;\bb{Q})$ as the domain is free. This is an isomorphism for any finitely generated abelian group $A$ by direct computation, and as both functors $\Lambda_\bb{Q}^*(- \otimes \Q)$ and $H_*(-;\Q)$ commute with filtered colimits of abelian groups, and any abelian group is the filtered colimit of its finitely generated subgroups, it is an isomorphism for all abelian groups. As $F$ is commutative this discussion applies to $A=F^\times$, giving $H_*(F^\times;\bb{Q}) \cong \Lambda^* F^\times_\bb{Q}$, so the map induced by $\lambda \mapsto \lambda^n$ acts by multiplication by $n^k$ in degree $k$.\end{proof}

The following describes the scaling action by $\Lambda^* F^\times_\bb{Q}$ under this splitting; this in particular applies to $M = \StL_n$.

\begin{lemma}Under the isomorphism $\varpi_n \colon H_*(\GL_n(F);\pr^* M) \xrightarrow{\cong} H_*(\PGL_n(F);M) \otimes \Lambda^* F^\times$ for a $\bb{Q}[\PGL_n]$-module $M$, the action of $\Lambda^* F^\times$ on $H_*(\GL_n(F);\pr^* M)$ is given by
\[x \cdot (y \otimes z) = (-1)^{|x||y|} n^{|x|} y \otimes (x \cdot z)\]
where $\cdot$ denotes the multiplication on $\Lambda^* F^\times$.\end{lemma}

\begin{proof}The diagram 
\[\begin{tikzcd} F^\times \times \GL_n \rar{\alpha} \dar[swap]{\id \times (\pr,\det)} & \GL_n \dar{(\pr,\det)} \\[-5pt]
F^\times \times \PGL_n \times F^\times \rar & \PGL_n \times F^\times \end{tikzcd}\]
commutes if the bottom group homomorphism is given by $(\lambda,[A],\mu) \longmapsto ([A],\lambda^n \mu)$. The result follows upon passing to homology.
\end{proof} 

\subsection{Components of the cobracket} These splittings allow us to decompose the operations---products, coproducts, cobrackets---on the homology groups in \cref{thm:ekhomology-steinberg} into different components. In particular, recall from \cref{sec:ek-homology} that on $\scr{G}$ there are a cobracket and a $\sigma$-component
\[\delta \colon \scr{G} \lra \Lambda^2 \scr{G}, \quad  \text{and} \quad \delta_\sigma \colon \scr{G} \lra H_2(\GL_{n-1};\StL_{n-1}).\] 
The target of the $\sigma$-component can be described more explicitly:
\begin{enumerate}[(1)]
    \item For $n=2$, we know the target is given by $\Lambda^2 F^\times$.
    \item For $n \geq 3$, \cref{thm:bgl-critical-line-vanishing} implies there is an isomorphism
    \[\quad H_2(\GL_n;\StL_n) \underset{\varpi_n}{\overset{\cong}\lra}H_1(\PGL_n;\StL_n) \otimes F^\times \oplus H_2(\PGL_n;\StL_n)\]
which we may further simplify using $H_1(\PGL_n;\StL_n) \cong H_1(\GL_n;\StL_n) = \scr{G}_n$.
\end{enumerate}

Our goal for the remainder of this section is to explain how to compute $\delta_\sigma$ in case (1) and the first component $\rm{pr}_1 \delta_\sigma$ in case (2).

\begin{proposition}\label{prop:cobr-vs-sigma-component} \,
\begin{enumerate}[(i)]
    \item For $n=2$, the cobracket and $\sigma$-component (both maps $\scr{G}_2 \to \Lambda^2 F^\times$) agree up to a sign: 
    \[\delta = -\delta_\sigma.\]
    \item For $n \geq 3$, the $(n-1,1)$-component of cobracket and the first term of the $\sigma$-component (both maps $\scr{G}_n \lra \scr{G}_{n-1}\otimes F^\times$) agree up to a sign:
\[\delta_{1,n-1} = -\pr_1 \delta_\sigma.\]
\end{enumerate}
\end{proposition}

\subsubsection{A projective cobracket} Up to an application of the K\"unneth isomorphism, the cobracket on $H_{*,*}(\GL;\StL)$ has components given by the dashed maps
\[\begin{tikzcd} H_*(\GL_n;\StL_n) \rar[dashed]{\Delta_{k,n-k}-\sigma \circ \Delta_{n-k,k}} &[5pt] H_{*+1}(\GL_k \times \GL_{n-k};\StL_k \otimes \StL_{n-k}) \arrow{dd}{\cong}\\[-5pt] 
H_*(\GL_n;[\StL_n \to (\SStL \para \SStL)(F^n)]) \dar \uar[two heads]& \\[-5pt] 
H_{*+1}(\GL_n;(\SStL \para \SStL)(F^n))  \rar & H_{*+1}(\rm{P}_{n,n-k};\StL_k \otimes \StL_{n-k})\end{tikzcd}\]
obtained using zigzag from \cref{prop:construction-of-cobracket}, Shapiro's lemma followed by projection onto the $(k,n-k)$-term, the Nesterenko--Suslin property in the guise of \cref{lem:para-levi-comparison-pgl}, and antisymmetrising. As the action of $\GL_n$ on $\StL_n$ factors over $\PGL_n$, we can construct a projective version, though this is \emph{not} a cobracket given the form of its domain and target:
\[\begin{tikzcd} H_*(\PGL_n;\StL_n) \rar[dashed]{P\Delta_{k,n-k}-\sigma \circ P\Delta_{n-k,k}} &[5pt] H_{*+1}(\rm{P}(\GL_k \times \GL_{n-k});\StL_k \otimes \StL_{n-k}) \arrow{dd}{\cong} \\[-5pt]
H_{*}(\PGL_n;[\StL_n \to (\SStL \para \SStL)(F^n)]) \uar[two heads] \dar& \\[-5pt]
H_{*+1}(\PGL_n;(\SStL \para \SStL)(R^n))  \rar & H_{*+1}(\rm{P}(\rm{P}_{n,n-k});\StL_k \otimes \StL_{n-k}).\end{tikzcd}\] 
By construction of the splitting the following diagram commutes
\begin{equation}\label{eqn:splitting-compatibility-1} \begin{tikzcd} H_*(\GL_n;\StL_n) \rar{\Delta_{k,n-k}-\sigma \circ \Delta_{n-k,k}} \dar{\cong}[swap]{\varpi_n} & H_{*+1}(\GL_k \times \GL_{n-k};\StL_k \otimes \StL_{n-k}) \dar{(\pr,\det)_*} \\[-5pt]
H_*(\PGL_n;\StL_n) \otimes \Lambda^* F^\times \rar & H_{*+1}(\rm{P}(\GL_k \times \GL_{n-k});\StL_k \otimes \StL_{n-k}) \otimes \Lambda^* F^\times \end{tikzcd}\end{equation}
where the bottom map is $(P\Delta_{k,n-k}-\sigma \circ P\Delta_{n-k,k}) \otimes \id$. We may apply the K\"unneth theorem in the top-right term of this square and apply the splittings to each term, resulting in the top map of the zigzag
\begin{equation}\label{eqn:splitting-compatibility-2}\begin{tikzcd} H_*(\GL_k \times \GL_{n-k};\StL_k \otimes \StL_{n-k}) \dar{(\pr,\det)_*} \rar{\varpi_{k} \otimes \varpi_{n-k}} &[-10pt] \parbox{4.8cm} {\centering $H_*(\PGL_k;\StL_k) \otimes \Lambda^* F^\times$ \\ $\otimes H_*(\PGL_{n-k};\StL_{n-k}) \otimes \Lambda^* F^ \times$} \\[-5pt]
H_*(\rm{P}(\GL_k \times \GL_{n-k}),\StL_k \otimes \StL_{n-k}) \otimes \Lambda^* F^\times.\end{tikzcd}\end{equation}
In the case $k=n-1$, there is an isomorphism from the bottom-left to top-right term of \eqref{eqn:splitting-compatibility-2} making the diagram obtained pasting \eqref{eqn:splitting-compatibility-1} and \eqref{eqn:splitting-compatibility-2} side-by-side commute:

\begin{lemma}\label{lem:delta-vs-pdelta} For $n \geq 2$ the following diagram commutes
\[\begin{tikzcd} H_{*}(\GL_n;\StL_n) \dar{\cong}[swap]{\varpi_n} \rar{\Delta_{n-1,1}-\sigma \circ \Delta_{1,n-1}} &[40pt] H_{*+1}(\GL_{n-1};\StL_{n-1}) \otimes \Lambda^* F^\times \dar{\cong}[swap]{\varpi_{n-1} \otimes \varpi_1} \\[-5pt]
H_*(\PGL_n;\StL_n) \otimes \Lambda^* F^\times \rar & H_{*+1}(\PGL_{n-1},\StL_{n-1}) \otimes \Lambda^* F^\times \otimes \Lambda^* F^\times\end{tikzcd}\]
where the bottom map is $(\id \otimes \rm{sh}_*)^{-1} \circ ((\rm{P}\Delta_{n-1,1}-\sigma \circ P\Delta_{1,n-1}) \otimes \id)$
and $\rm{sh} \colon F^\times \times F^\times \to F^\times \times F^\times$ is the ``shear'' homomorphism given by $(\mu,\lambda) \mapsto (\mu\lambda^{-(n-1)},\mu \lambda)$.
\end{lemma}

\begin{proof}As explained above, it suffices to construct an isomorphism in the zigzag \eqref{eqn:splitting-compatibility-2}, and to do so we consider the zigzag of groups
\[\begin{tikzcd} & \arrow{ld}[swap]{(\pr,\det)} \GL_k \times \GL_{n-k} \arrow{rd}{((\pr,\det),(\pr,\det))} & \\[-7pt]
\rm{P}(\GL_k \times \GL_{n-k})\times F^\times  & &  \PGL_k \times F^\times \times \PGL_{n-k} \times F^\times\end{tikzcd}\]
in the case $k=n-1$. The composition $\pr \circ \inc \colon \GL_{n-1} \to \GL_{n-1} \times \GL_1 \to \rm{P}(\GL_{n-1} \times \GL_1)$ is an isomorphism, so there is a map $\phi \colon \rm{P}(\GL_{n-1} \times \GL_1) \cong \GL_{n-1} \to \PGL_{n-1} \times F^\times$ where the second map is $(\pr,\det)$. We use this to produce a diagram
\[\begin{tikzcd} & \arrow{ld}[swap]{(\pr,\det)} \GL_{n-1} \times \GL_1 \arrow{rd}{((\pr,\det),(\pr,\det))} & \\[-7pt]
\rm{P}(\GL_{n-1} \times \GL_1) \times F^\times \dar{\phi \times \id} & &  \PGL_{n-1} \times F^\times \times \PGL_1 \times F^\times \dar{\cong} \\[-2pt]
\PGL_{n-1} \times F^\times \times F^\times & & \PGL_{n-1} \times F^\times \times F^\times \arrow[dashed]{ll}[swap]{\id \times \rm{sh}}
\end{tikzcd}\]
where the left composition is given by $(A,\lambda) \mapsto ([A],\lambda^{-(n-1)}\det(A),\lambda\det(A))$ and the right composition is given by $(A,\lambda) \mapsto ([A],\det(A),\lambda)$. Thus we can make it commute by setting the dashed map to be $([A],\mu,\lambda) \mapsto ([A],\lambda^{-(n-1)}\mu,\lambda\mu)$.\end{proof}

We will now use this to establish \cref{prop:cobr-vs-sigma-component} comparing the cobracket to the $\sigma$-component.

\begin{proof}[Proof of \cref{prop:cobr-vs-sigma-component}] We start by recalling from \cref{prop:cobracket-via-coproduct} that the cobracket is induced by the antisymmetrisation of the reduced coproduct as $\delta = \ol{\Delta}-\sigma \circ \ol{\Delta}$; here $\delta_{n-1,1} = \Delta_{n-1,1} - \sigma \circ \Delta_{1,n-1}$. We specialise the commutative square of \cref{lem:delta-vs-pdelta} to $*=1$:
\[\begin{tikzcd} H_1(\GL_n;\StL_n) \dar{\cong}[swap]{\varpi_n} \rar{\Delta_{n-1,1} - \sigma \circ \Delta_{1,n-1}} &[20pt] \big(H_{*+1}(\GL_{n-1};\StL_{n-1}) \otimes \Lambda^* F^\times\big)_1 \dar{\cong}[swap]{\varpi_{n-1} \otimes \id} \\[-5pt]
H_1(\PGL_n;\StL_n) \otimes \Lambda^0 F^\times \rar{} & \big(H_{*+1}(\PGL_{n-1};\StL_{n-1}) \otimes \Lambda^* F^\times \otimes \Lambda^* F^\times\big)_1\end{tikzcd}\]
where the bottom horizontal map is
\[(\id \otimes \rm{sh}_*)^{-1} \circ ((\rm{P}\Delta_{n-1,1}-\sigma \circ \rm{P}\Delta_{1,n-1}) \otimes \id).\]
Here we have used the assumption that $n\geq 2$ and hence the coinvariants of $\StL_n$ vanish, to see that the left vertical map reduces to an isomorphism $\varpi_n \colon H_1(\GL_n;\StL_n) \smash{\overset{\cong}\lra} H_1(\PGL_n;\StL_n)$.

\medskip

For $n=2$, $\PGL_1$ is trivial and applying $\rm{P}\Delta_{1,1}-\sigma \circ \rm{P}\Delta_{1,1}$ we land in the summand $\bb{Q} \otimes \Lambda^2 F^\times \otimes \Lambda^0 F^\times$ of $\bb{Q} \otimes (\Lambda^* F \otimes \Lambda^* F)_2$. The map $(\id \otimes \rm{sh}_*)^{-1}$ sends this into the subspace spanned by elements of the form $1 \otimes (x \wedge y) \otimes 1 + 1 \otimes 1 \otimes (x \wedge y)-1 \otimes x \wedge y$, where we consider $\Lambda^2 F^\times$ as summand of $F^\times \otimes F^\times$. Here we have used that if we write $F^\times$ additively, we can identify the homomorphisms $A$ and $A^{-1}$ with the matrices
\[A = \begin{bmatrix} 1 & -(n-1) \\ 1 & 1 \end{bmatrix} \qquad \text{and} \qquad A^{-1} = \frac{1}{n} \begin{bmatrix} 1 & n-1 \\ -1 & 1 \end{bmatrix}.
\] 
Moreover, $\GL_1 = F^\times$, $\StL_1 = \bb{Q}$, and $\varpi_1$ is the usual identification. By definition, the projection of $\Delta_{1,1}-\sigma \circ \Delta_{1,1}$ onto $\bb{Q} \otimes \Lambda^2 F^\times \cong \Lambda^2 F^\times \otimes \bb{Q}$ yields the $\sigma$-component $\delta_\sigma$ and projection to $\bb{Q} \otimes \Lambda^2 F^\times$ yields the cobracket $\delta_{1,1}$. The restrictions we found above on the image thus imply $\delta_\sigma = -\delta_{1,1}$.

\medskip

For $n \geq 3$ we argue similarly, but now applying $\rm{P}\Delta_{n-1,1}-\sigma \circ \rm{P}\Delta_{1,n-1}$ we land in the summands $H_1(\PGL_{n-1};\St_{n-1}) \otimes \Lambda^1 F^\times \otimes \Lambda^0 F^\times$ and $H_2(\PGL_{n-1};\St_{n-1}) \otimes \Lambda^0 F^\times \otimes \Lambda^0 F^\times$. The map $(\id \otimes \rm{sh}_*)^{-1}$ fixes the latter pointwise, but sends the former into the subspace spanned by $x \otimes y \otimes 1-x \otimes 1 \otimes y$. By definition the projection of $(\varpi_{n-1} \otimes \varpi_1) \circ (\Delta_{n-1,1}-\sigma \circ \Delta_{n-1,1})$ to $H_1(\PGL_{n-1};\St_{n-1}) \otimes \Lambda^1 F^\times \otimes \Lambda^0 F^\times$ is the first part of the $\sigma$-component $\rm{pr}_1 \delta_\sigma$ and its projection to $H_1(\PGL_{n-1};\St_{n-1}) \otimes \Lambda^0 F^\times \otimes \Lambda^1 F^\times$ is the cobracket $\delta_{n-1,1}$. The restriction found above on the image thus implies $\pr_1 \delta_\sigma = -\delta_{n-1,1}$.
\end{proof}

\subsection{The rank $2$ case.} In this subsection we describe a more elementary approach to think of the cobracket in rank $2$, 
\[\delta \colon H_{2,d}^{E_\infty}(\BGL(F)_\Q) \lra (\Lambda^2 H_{1,*}^{E_\infty}(\BGL(F)_\Q))_{d+1}\]
and connect it to previous computations in the literature. This is not used in the remainder of this paper.

\medskip

Recall that the counit $\fgt_{E^\rm{nu}_\infty} \rm{triv}_{E^\rm{nu}_\infty} \simeq \rm{id}$ is an equivalence of endofunctors of $\Fun(\N,\DQ)$. Precomposing it with $\rm{cot}_{E^\rm{nu}_\infty}$ and using the unit of the adjunction $\cot_{E^\rm{nu}_\infty} \dashv \triv_{E^\rm{nu}_\infty}$ yields a natural transformation of functors $\Alg_{E_\infty^\rm{nu}}(\Fun(\N,\DQ)) \to \Fun(\N,\DQ)$
\[\rm{pr}_{E^\rm{nu}_\infty} \colon \fgt_{E^\rm{nu}_\infty} \lra \cot_{E^\rm{nu}_\infty}.\]
We use this to define the \emph{decomposables} functor 
\begin{align*}\rm{dec}_{E^\rm{nu}_\infty} \colon  \Alg_{E_\infty^\rm{nu}}(\Fun(\N,\DQ)) &\lra \Fun(\N,\DQ) \\
\bf{R} &\longmapsto \rm{fib}\big[\rm{pr}_{E^\rm{nu}_\infty} \colon \fgt_{E_\infty^\rm{nu}}(\bf{R}) \to \cot_{E_\infty^\rm{nu}}(\bf{R})\big].\end{align*}

To compute with this, we use the rectification results of \cref{sec:rect-dg} to assume that $\bf{R} \in \Alg_{E^\rm{nu}_\infty}(\Fun(\N,\DQ))_{\geq 1}$ is modelled by a dg-commutative algebra with additional rank grading. In this case, the decomposables admit a description in terms of the commutative bar construction of \cref{def:bar-comm}. 

\begin{lemma}
If $\bf{R} \in \Alg_{\rm{Com}^\rm{nu}}(\Fun(\N,\rm{Ch}_\bb{Q}))_{\geq 1}$ then there is a natural equivalence 
   \[\rm{dec}_{E^\rm{nu}_\infty}(\bf{R}) \simeq  \Sigma^{-1} B^\rm{Com}(\bf{R})_{\ge 2} = \Sigma^{-2} (\rm{coLie}_{\ge 2} \circ \Sigma \bf{R}, d_\bf{R}+d_B)\]
so that the canonical map $\rm{dec}_{E^\rm{nu}_\infty}(\bf{R}) \to \fgt_{E^\rm{nu}_\infty}\bf{R}$ is induced by the multiplication on $R$. In particular, there is a natural map
   \[\inc \colon S^2(\bf{R}) \simeq  \Sigma^{-2} \rm{coLie}(2) \otimes_{\fr{S}_2} (\Sigma \bf{R})^{\otimes 2} \overset{\subset}\lra \Sigma^{-1} B^\rm{Com}(\bf{R})_{\ge 2} \overset{\simeq}\lra \rm{dec}_{E^\rm{nu}_\infty}(\bf{R}).\]
\end{lemma}

\begin{proof}
    The equivalence $B^\rm{Com}(\bf{R}) \simeq \cot_{E^\rm{nu}_\infty}(\bf{R})$ is such that the canonical map $\rm{pr}_{E^\rm{nu}_\infty}$ can be identified with the inclusion $\bf{R} \cong B^\rm{Com}(\bf{R})_{\le 1} \hookrightarrow B^\rm{Com}(\bf{R})$. The result then follows by taking the mapping cone. 
\end{proof}

The next lemma will be the key to computing the cobracket in rank $2$ and requires us to introduce one further construction. For $\bf{R} \in \Alg_{\rm{Com}^\rm{nu}}(\Fun(\N,\rm{Ch}_\bb{Q}))_{\geq 1}$ we introduce the following cofibre in $\Fun(\N,\rm{Ch}_\Q)$
\[Q(\bf{R}) \coloneq \rm{cofib}\big[S^2(\bf{R}) \xrightarrow{m}  \fgt_{E^\rm{nu}_\infty}\bf{R}\big]\]
where $m$ denotes the multiplication. There is a map of cofibre sequences in $\Fun(\N,\rm{Ch}_\bb{Q})$
\[\begin{tikzcd} S^2(\bf{R}) \rar{m} \dar{\inc} & \fgt_{E^\rm{nu}_\infty} \bf{R} \dar[equal] \rar & Q(\bf{R}) \dar[dashed] \\[-5pt]
\rm{dec}_{E_\infty^\rm{nu}}(\bf{R}) \rar & \fgt_{E_\infty^\rm{nu}}(\bf{R}) \rar{\pr_{E_\infty^\rm{nu}}} & \cot_{E_\infty^\rm{nu}}(\bf{R})\end{tikzcd}\]
where the dashed map is induced.

\begin{lemma}
   Let $\bf{R} \in \Alg_{\rm{Com}^\rm{nu}}(\Fun(\N,\rm{Ch}_\bb{Q}))$, then the following diagram commutes after taking homology
  \[
\begin{tikzcd}
Q(\bf{R})
  \rar
  \dar{\partial}
&[20pt]
\cot_{E_\infty}(\bf{R})
  \dar{\Sigma^{-1}\delta}
\\[-5pt]
\Sigma S^2 \bf{R}
  \rar[swap]{\Sigma S^2 \rm{pr}_{E_\infty}}  &
  \Sigma S^2 \cot_{E_\infty}(\bf{R})
\end{tikzcd}
\]
   where $\delta$ is obtained from the Lie cobracket
   \[\Sigma \cot_{E_\infty}(\bf{R}) \overset{\delta}\lra \rm{coLie}(2) \otimes_{S_2}(\Sigma \cot_{E_\infty}(\bf{R}))^{\otimes 2} \simeq \Sigma^2 S^2 \cot_{E_\infty}(\bf{R})\]
   and $\partial$ the connecting homomorphism of the cofibre sequence $S^2 \bf{R} \to \fgt_{E_\infty} \bf{R} \to Q(\bf{R})$. 
\end{lemma}

\begin{proof}
First we prove this for $\bf{R}=\rm{triv}_{\rm{E_\infty}}(A)$ for some $A \in \Fun(\N,\DQ)$. In this case we have (1) $m \colon S^2\bf{R} \to \fgt_{E_\infty}\bf{R}$ is the zero map $S^2 A \to A$ and hence $Q(\bf{R}) \simeq A \oplus \Sigma S^2 A$, (2) $\cot_{E_\infty}(\bf{R}) \simeq B^\rm{Com}(\bf{R}) \simeq \Sigma^{-1}\coLie(\Sigma A)$ and the map $Q(\bf{R}) \to \cot_{E_\infty}(\bf{R})$ is just the inclusion of the first two summands, (3) $\rm{pr}_{E_\infty} \colon\fgt_{E_\infty}(\bf{R}) \simeq A \to \cot_{E_\infty}(\bf{R}) \simeq  \Sigma^{-1}\coLie(\Sigma A)$ is the inclusion of the first term, and (4) the cobracket is that of the cofree Lie coalgebra under the equivalence $\cot_{E_\infty}(\bf{R}) \simeq  \Sigma^{-1}\coLie(\Sigma A)$. Then the result follows since the following diagram commutes
 \[
\begin{tikzcd}
A \oplus \Sigma S^2 A
  \rar{\inc}
  \dar{\pr_2}
&[20pt] 
\Sigma^{-1}\coLie(\Sigma A)
  \dar{\Sigma^{-1}\delta_{\mathrm{cofree}}}
\\[-5pt]
\Sigma S^2 A
  \rar{\Sigma S^2\mathrm{incl}} & \Sigma S^2(\Sigma^{-1}\coLie(\Sigma A)).
\end{tikzcd}
\]

In the general case, for $\bf{R} \in \Alg_{\rm{Com}^\rm{nu}}(\Fun(\N,\DQ))$ there is a canonical algebra map 
\[\nu \colon \bf{R} \lra \rm{triv}_{\rm{E_\infty}} \cot_{E_\infty} \bf{R}.\]
Let us denote $A \coloneq \cot_{E_\infty}(\bf{R})$ and observe $\cot_{E_\infty} \nu$ is the inclusion $A \hookrightarrow \Sigma^{-1} \coLie(\Sigma A)$. Thus the map 
\[\Sigma S^2 \cot_{E_\infty} \bf{R} \xrightarrow{ \Sigma S^2 \cot_{E_\infty} \nu} \Sigma S^2 (\Sigma^{-1} \coLie(\Sigma A))\] is injective in homology. To prove the result it suffices to prove it after post-composing with the map $\Sigma S^2 \cot_{E_\infty} \nu$ and taking homology, reducing it to the previous case by naturality. 
\end{proof}

Recall that if $\bf{R}(0) \simeq 0$ we say it is \emph{reduced}. In this case the canonical maps $S^2 (\bf{R}(1)) \to (S^2 \bf{R})(2) \to \rm{dec}_{E_\infty}(\bf{R})(2)$ are both equivalences. Thus, both of the following maps are also equivalences
\[\frac{\bf{R}(2)}{S^2 (\bf{R}(1))} \xrightarrow{\simeq }Q(\bf{R})(2) \xrightarrow{\simeq} \cot_{E_\infty}(\bf{R})(2).\]
Moreover, it also follows from the formula for $B^\rm{Com}$ that $ \fgt_{E_\infty}\bf{R}(1) \xrightarrow{\rm{pr}_{E_\infty}} \cot_{E_\infty}(\bf{R})(1)$ is an equivalence. Thus, we have a commutative diagram
  \[
\begin{tikzcd}
H_d(\bf{R}(2),S^2(\bf{R}(1)))
  \rar{\simeq}
  \dar{\partial}
&
H_{2,d}^{E_\infty}(\bf{R})
  \dar{\delta}
\\[-5pt]
(S^2 H_{1,*}(\bf{R}))_{d-1}
  \rar{\simeq} & (S^2 H_{1,*}^{E_\infty}(\bf{R}))_{d-1}
\end{tikzcd}
\]

\begin{corollary}
 If $\bf{R} \in \Alg_{\rm{Com}^\rm{nu}}(\Fun(\N,\rm{Ch}_\bb{Q}))_{\geq 1}$ is reduced then the cobracket on the $E_\infty$-homology groups in rank $n=2$ is computed by 
 \[\partial \colon H_d(\bf{R}(2),S^2(\bf{R}(1))) \lra H_{2,d-1}(S^2 \bf{R}) \simeq (S^2 H_{1,*}(\bf{R}))_{d-1}\]
 under the canonical equivalences $H_d(\bf{R}(2),S^2(\bf{R}(1))) \simeq H_{2,d}^{E_\infty}(\bf{R})$ and $(S^2 H_{1,*}(\bf{R}))_{d-1} \simeq (S^2 H_{1,*}^{E_\infty}(\bf{R}))_{d-1}$
\end{corollary}

Let us apply this to $\bf{R}=\BGLb_\Q$, implicitly fixing a field $F$, and $d=3$. First, one can use the proof of \cite[Theorem 9.1]{GKRW20} to build an explicit inverse map 
\[B_2(F) \cong H_{2,3}^{E_\infty}(\BGLb_\Q) \overset{\cong}\lra H_3(\BGLb_\Q(2), S^2 \BGLb_\Q(1)) \cong  H_3(\GL_2,\rm{GM}_2;\Q)\]
and then \cite[(9.3)]{GKRW20} identifies $\partial$ as follows: 
as $H_{1,*}(\BGL(F)_\Q) \cong \Lambda^* F^\times_\Q$ canonically, we have a preferred isomorphism 
\begin{align*}(S^2 H_{1,*}(\BGL(F)_\Q)_2 &\cong \Lambda^2 H_{1,1}(\BGL(F)_\Q) \oplus H_{1,0}(\BGL(F)_\Q \otimes H_{1,2}(\BGL(F)_\Q) \\
&\cong \Lambda^2 F^\times_\Q \oplus \Q\{\sigma\} \otimes (\Lambda^2 F^\times_\Q)\end{align*}
and then the map $\partial$ is given by (using work of Suslin \cite{Sus90} and Mirzaii \cite{Mirzaii,MirzaiiErratum})
\begin{align*}\partial \colon B_2(F) &\lra \Lambda^2 F^\times_\Q \oplus \Q\{\sigma\} \otimes (\Lambda^2 F^\times_\Q) \\
{x}_2 &\longmapsto ((x) \wedge (1-x), \sigma \otimes (- (x) \wedge (1-x))),\end{align*}
which both recovers our computation of the cobracket $B_2(F) \cong \mathcal{G}_2(F) \to \Lambda^2 \mathcal{G}_1(F)) \cong \Lambda^2 F^\times_\Q$ and shows in this case explicitly the relationship $\delta_\sigma=-\delta$ between the $\sigma$-component and the cobracket.

\section{The relationship of $\PolyL(F)$ to polylogarithms}
In this section we discuss the Lie coalgebra $\PolyL(F)$ in more details. We start by spelling out how the results of \cref{sec:cobracket} yield \cref{thm:polyl-presentation-additive}. We next identify $\PolyL_n(F)$ for $n \leq 3$ in terms of more classical definitions proving \cref{thm:polyl-identification}. This requires the construction of several families of relations in $\PolyL_n(F)$ that may be of independent interest. We finally explain how \cref{thm:polyl-presentation-additive} implies the existence of a map to the Lie coalgebra of formal polylogarithms of Charlton, Matveiakin, Radchenko, and Rudenko \cite{CMRR24}, yielding Hodge and motivic realisations and thus proving Theorems \ref{thm:motivic-realisation} and \ref{thm:hodge-realisation}.

\begin{convention}\label{conv:less-short} We fix a field $F$ but do not suppress it from the notation. We work rationally and suppress this from the notation unless there is a risk of confusion, e.g.~write $F^\times$ for $F^\times_\bb{Q}$.
\end{convention}

\subsection{Generators and relations}
We now prove \cref{thm:polyl-presentation-additive}: 

\begin{customthm}{C.a} 
The Lie coalgebra $\PolyL(F)$ is generated by correlators 
    \[\Cor^\scr{G}(x_0,x_1,\compactldots,x_n) \in \scr{G}_n(F) \qquad \text{for $x_0,\ldots,x_n \in F$ not all equal}\]  
    subject to the following relations:
    \begin{enumerate}[\noindent (1)]
\item  Homogeneity: $\Cor^\scr{G}(x_0,x_1,\compactldots,x_n)=\Cor^\scr{G}(x_0+b,x_1+b,\compactldots,x_n+b)$ for $b\in F$.
\item  Cyclic symmetry: $\Cor^\scr{G}(x_0,x_1,\compactldots,x_n)=\Cor^\scr{G}(x_1,x_2,\compactldots,x_0)$.
\item Shuffle relations: 
\[\sum_{\sigma \in \rm{Sh}(n_1,n_2)} \Cor^\scr{G}(x_0,x_{\sigma(1)},\cdots,x_{\sigma(n_1+n_2)}) =0 \quad 
\text{for $n=n_1+n_2$, $n_1,n_2>0$.}\]
\item Decomposition relations:
\begin{align*}
    &\CorG(x_0,\compactldots,x_n)-\CorG(y_0,\compactldots,y_n)\\
     &=\sum_{\iota=((i_1,j_1),\dots,(i_{n},j_n))\in T(n)}\rm{sign}(\iota)\,
     \CorG\left(0,\frac{x_{i_1}-x_{j_1}}{y_{i_1}-y_{j_1}},\compactldots,\frac{x_{i_n}-x_{j_n}}{y_{i_n}-y_{j_n}}\right),
\end{align*}
where we omit terms with $y_{i_k}=y_{j_k}$ for some k, and the set $T(n)$ as well as the sign $\rm{sign}(\iota)$ are as in \cref{prop:universal-symbol-combinatorics}.
\end{enumerate}
\end{customthm}

\begin{remark}These relations are redundant, e.g.~homogeneity follows from the decomposition relations. We believe cyclic symmetry also follows from the decomposition relations, but this is left an exercise to the interested reader.
\end{remark}

\begin{proof}[Proof of \cref{thm:polyl-presentation-additive}]
Recall the presentation of the Lie coalgebra $\PolyL(F)$ from \cref{prop:presentation for G(F)}: for a vector space $V$  of dimension $n \geq 1$ over $F$ we have
\[
\scr{G}(V)\overset{\cong}{\lra} \frac{\bb{Q}\{[h] \otimes \FC[u_0:\compactcdots:u_n] \text{ for nonzero functionals $h$ and affine bases $u_0,\ldots,u_n$}\}}{\text{\eqref{enum:fc-relations-i}--\eqref{enum:gd-relation-v}}},
\]
with the following relations:
\begin{enumerate}[\noindent (1)] 
\item  Homogeneity: $\rm{FC}[u_0:\cdots:u_n] = \rm{FC}[u_0-u:\cdots:u_n-u]$ for any $u \in V$.
\item  Cyclic symmetry: $\rm{FC}[u_0:u_1:\cdots:u_n] = \rm{FC}[u_1:u_2:\cdots:u_0]$.
\item  Shuffle relations: 
\[\sum_{\sigma \in \rm{Sh}(n_1,n_2)} \rm{FC}[u_0:u_{\sigma(1)}:\cdots:u_{\sigma(n_1+n_2)}] =0 \quad \text{for $n=n_1+n_2$ with $n_1,n_2>0$.}\]
\item  Coinvariant relations:
\[[h] \otimes \FC[u_0:\compactcdots:u_n] = [g^*h] \otimes \FC[gu_0:\compactcdots:gv_n] \qquad \text{for $g \in \GL(V)$}.\]
\item Decomposition relations:
\[\qquad [h_2] \otimes \FC[u_0:\compactcdots:u_n]-[h_1] \otimes \FC[u_0:\compactcdots:u_n]+[h_1] \otimes D^\FC_{h_2}(\FC[u_0:\compactcdots:u_n])=0.\]
\end{enumerate}

To get the generating set of correlators, we define
\[\Cor^\scr{G}(x_0,x_1,\compactldots,x_n) \in \scr{G}(V)\] 
as the image of an element $[h] \otimes \FC[u_0:\compactcdots : u_n]$ for $h(u_i) = x_i \in F$. 

To get the relations, we reformulate \cref{prop:presentation for G(F)} more explicitly. To do so, we first make the decomposition operator more explicit. Recall from \eqref{eqn:dec-fc} the definition of the decomposition operator as $D^\FC_h = \rm{C}^\FC_h \circ s_H \circ \pr^\St$ in terms of four steps: it is given by 
\begin{enumerate}[(i)]
\item \label{eqn:decomposition-steps-i} we map a formal correlator to a Steinberg correlator, 
\item \label{eqn:decomposition-steps-ii} take its symbol,
\item \label{eqn:decomposition-steps-iii} project onto those terms where no line in contained in $H = \ker(h)$, and 
\item \label{eqn:decomposition-steps-iv} take formal correlators with entries the unique vectors in the lines on which the functional $h$ takes the value $1$.
\end{enumerate} This yields that the decomposition operator has the following explicit form: using the universal formula for the symbol from \cref{prop:universal-symbol-combinatorics} to implement step \cref{eqn:decomposition-steps-ii}, there exists unique set $T(n)$ of pairs of indices and a unique sign function sending $\iota \in T(n)$ to $\rm{sign}(\iota)\in\{\pm 1\}$ such that
\begin{equation}\label{eqn:dec-explicit}\begin{aligned}
&D^\FC_h(\FC[u_0:\compactcdots:u_n])\\
&=\sum_{\iota=((i_1,j_1),\dots,(i_{n},j_n))\in T(n)}\rm{sign}(\iota)
\FC\left[0:\frac{u_{i_1}-u_{j_1}}{h(u_{i_1})-h(u_{j_1})}:\compactcdots:\frac{u_{i_n}-u_{j_n}}{h(u_{i_n})-h(u_{j_n})}\right],
\end{aligned}\end{equation}
where we omit the terms $\iota$ with $h(u_{i_k})=h(u_{j_k})$ for some $k$. Here the omission of certain term occurs in the step \eqref{eqn:decomposition-steps-iii} when projecting.

Then, relations (1), (2) and (3) above give the respective homogeneity \eqref{enum:rel-goncharov-1}, cyclic symmetry \eqref{enum:rel-goncharov-2}, and shuffle relations \eqref{enum:rel-goncharov-3} for correlators in \cref{thm:polyl-presentation-additive}. The coinvariant relations (4) above show that the elements $\CorG(x_0,\compactldots,x_n) \in \PolyL(F)$ are well-defined. Finally, the decomposition relations (5) above give the decomposition relations for correlators \eqref{enum:rel-goncharov-4} by evaluating \eqref{eqn:dec-explicit} for $h_1, h_2$ such that $h_1(u_i)=x_i$ and $h_2(u_i)=y_i$.\end{proof}

We next note consequences of these relations. The shuffle relation \eqref{enum:rel-goncharov-3} implies that for any presentation $n=n_1+n_2$ we have
\[
{\textstyle{n_1+n_2 \choose n_1}}\,\CorG(x_0,\underbrace{x_1,\compactldots,x_1}_n)=0 \quad \text{for } n\geq 2,
\]
and since $\PolyL(F)$ is a $\Q$-vector space, it follows that
\begin{equation}\label{eqn:cor vanishing}
\CorG(x_0,\underbrace{x_1,\compactldots,x_1}_n)=0 \quad \text{for } n\geq 2.
\end{equation}

\begin{lemma} \label{lemma: multiplicative homogenity relation for correlators}  The following identity holds for $n \geq 2$:
    \[
    \CorG(ax_0,\compactldots,ax_n)=\CorG(x_0,\compactldots, x_n) \qquad \text{ for $a \in F^{\times}$}.
    \]
\end{lemma}

\begin{proof}
We apply the decomposition relation \eqref{enum:rel-goncharov-4} with $y_i=a x_i$. For any $\iota \in T(n)$ we have
\[
 \CorG\left(0,\frac{x_{i_1}-x_{j_1}}{y_{i_1}-y_{j_1}},\compactldots,\frac{x_{i_n}-x_{j_n}}{y_{i_n}-y_{j_n}}\right)=\CorG\Bigl(0,\underbrace{\frac{1}{a},\compactldots,\frac{1}{a}}_n\Bigr)=0,
\]
with right equality by \eqref{eqn:cor vanishing}, which yields the result.
\end{proof}

\begin{lemma}\label{lemma: reversal relation for correlators}The following identity holds for $n\geq 1$: 
    \[
    \CorG(x_0,x_1,\compactldots,x_n)=(-1)^{n-1}\CorG(x_0,x_n,\compactldots,x_1).
    \]
\end{lemma}

\begin{proof} This is a consequence of the shuffle relations \eqref{enum:rel-goncharov-3}. To see that, recall that the cofree Lie coalgebra is the Lie coalgebra of indecomposable elements of the Hopf algebra $\bigoplus_{n\geq 0} V^{\otimes n}$ with shuffle product and deconcatenation coproduct. The antipode is given by the formula \cite[Example 1.6.3]{GrinbergReiner}
\[
S(v_1\otimes \dots \otimes v_n)=(-1)^n v_n\otimes \dots\otimes v_1
\]
and thus acts by $(-1)$ on the Lie coalgebra of indecomposable elements. It follows that in the cofree Lie coalgebra the projections of elements $v_1\otimes \dots \otimes v_n$ and $(-1)^{n-1} v_n\otimes \dots\otimes v_1$ coincide and from here the statement follows.
\end{proof}

\begin{proposition}\label{proposition: generic correlators}
    The vector space $\PolyL_n(F)$ is spanned by elements $\CorG(x_0,\compactldots,x_n)$ such that  $x_i\neq x_j$ for $i\neq j$.
\end{proposition}

\begin{proof} If the field $F$ is finite, $\PolyL_n(F)=0$ and there is nothing to prove, so we may assume that $F$ is infinite. Consider a correlator $\CorG(y_0,\compactldots,y_n)$ where some of the arguments may coincide. Since $F$ is infinite, we can choose a tuple of distinct elements $x_0,\dots,x_n \in F$ such that 
\[
\frac{x_i-x_j}{y_i-y_j} \neq \frac{x_k-x_l}{y_k-y_l}
\]
for all $i,j,k,l$ with $y_i\neq y_j$ and $y_k\neq y_l$. The decomposition relation \eqref{enum:rel-goncharov-4} gives a presentation of the correlator $\CorG(y_0,\compactldots,y_n)$ as a linear combination of correlators with distinct arguments.
\end{proof}

\begin{proposition} \label{prop: injectivity Goncharov complex transcendental extensions}
    The maps induced by the inclusion $F \to F(t)$
    \[\PolyL_n(F) \lra \PolyL_n(F(t))\]
    are injective for any field $F$ and $n \ge 1$. 
\end{proposition}

\begin{proof}
The statement is trivial for finite $F$ since $\PolyL(F)=0$, so let us assume that $F$ is infinite. It suffices to show that if element $\alpha \in \PolyL_n(F)$ vanishes $\PolyL(F(t))$ then $\alpha$ must be zero. First, write $\alpha$ as a rational linear combination of correlators with entries in $F$. Define $\mathcal{R}_n(F(t)) \subseteq \Q \{ \CorG(x_0,\compactldots,x_n) \mid x_0,\dots,x_n \in F(t)\}$ to be the subspace spanned by relations \eqref{enum:rel-goncharov-1}--\eqref{enum:rel-goncharov-4} in \cref{thm:polyl-presentation-additive} for the field $F(t)$. That $\alpha$ vanishes in $\PolyL_n(F(t))$ means that $\alpha \in \mathcal{R}(F(t))$ so we can find some $N \in \N$, rational functions $f_i^j \in F(t)$ for $0 \le i \le n$ and $1 \le j \le N$, and rational numbers $\lambda_j$ for $1 \le j \le N$ such that 
    \[\alpha = \sum_{j=1}^N  \lambda_j \CorG(f_0^j,\compactldots,f_n^j) \in \mathcal{R}_n(F(t)).\]
Now use that $F$ is infinite to pick $t_0 \in F$ such that all the $f_i^j(t_0)$ are distinct and nonzero. By construction of the relations, evaluating at $t_0$ sends relations in $\mathcal{R}_n(F(t))$ to relations in $\mathcal{R}(F)$, so we get
    \[\alpha = \sum_{j=1}^N  \lambda_j \CorG(f_0^j(t_0),\compactldots,f_n^j(t_0)) \in \mathcal{R}_n(F),\]
showing that the original class vanishes. 
\end{proof}

\begin{remark}
    Similar arguments as in \cref{prop: injectivity Goncharov complex transcendental extensions} can be used to show that the maps
    \[H_{n,d}^{E_\infty}(\BGLb(F)_\Q) \lra H_{n,d}^{E_\infty}(\BGLb(F(t))_\Q) \quad \text{and} \quad H^*(\PolyL(F)) \lra H^*(\PolyL(F(t))) \]
on $E_\infty$-homology, resp.~Chevalley--Eilenberg homology, are always injective.
\end{remark}

\subsection{Explicit description of  $\PolyL_1$ and  $\PolyL_2$}\label{sec: G1 and G_2} We now identify $\PolyL_1(F)$ and $\PolyL_2(F)$ in classical terms, proving two-thirds of \cref{thm:polyl-identification}.

\subsubsection{Explicit description of  $\PolyL_1$}
First we show that the weight one component of the Lie coalgebra $\PolyL(F)$ is isomorphic to the rationalisation of the group of units $F^\times$. From the preferred isomorphism $\StL_1(F) \cong \bb{Q}$, we obtain a preferred isomorphism
\begin{equation}\label{eqn:isomorphism in weight 1}
    \PolyL_1(F) = H_1(\GL_1(F);\StL_1(F))\cong H_1(\GL_1(F);\bb{Q}) \underset{\cong}{\overset{\rm{ab}}\lra} F^{\times}_{\Q}.
\end{equation}

\begin{lemma} \label{LemmaCorrelatorWeightOne} Under the isomorphism  \eqref{eqn:isomorphism in weight 1} an element  $(x_1-x_0)\in F^{\times}_{\Q}$ corresponds to $\CorG(x_0, x_1)\in \PolyL_1(F)$. 
\end{lemma}
\begin{proof}
Using homogeneity \eqref{enum:rel-goncharov-1}, it is sufficient to prove the statement for $x_0=0$ and $x_1=a\in F^{\times}$. We analyse the exact sequence 
\[0 \lra \FCR(V) \lra \FC(V) \overset{p}\lra \StL(V) \lra 0\]
for a one-dimensional vector space $V$. Choose a nonzero vector $e\in V$ and a functional $h\in V^{\vee}$ such that $h(e)=1$.  We have an isomorphism $\FC(V)\smash{\stackrel{\cong}{\lra}} \Q[F^{\times}/\{\pm 1\}]$ induced by the map sending $\FC[x_1e,x_2e]$ to $[x_1-x_2]\in \Q[F^{\times}/\{\pm 1\}]$. 

The tautological exact sequence $0\to I\to \Q[F^{\times}]\to \Q \to 0$ defining the augmentation ideal induces an isomorphism $H_1(\GL_1(F);\bb{Q}) \cong I_{\GL_1(F)}$, so an element in $H_1(\GL_1(F);\bb{Q})$ corresponding to $a\in F^{\times}$ is represented by $([a]-[1])\in I$.
The natural chain map 
\[\begin{tikzcd} 0 \rar & I \rar \dar{}& \Q[F^{\times}] \rar\dar{} & \Q \rar\dar{\text{id}} & 0 \\[-5pt]
0 \rar & \FCR(V)  \rar & \Q[F^{\times}/\{\pm 1\}] \rar &  \Q \rar & 0 \end{tikzcd}\]
induces an isomorphism $I_{\GL_1(F)} \cong \FCR(V)_{\GL_1(F)}$. The map $I \to \FCR(V)$ sends $[a]-[1]\in I$ to $\FC[0,ae]-\FC[0,e]$. The image of $\FC[0,ae]-\FC[0,e]$ in coinvariants coincides with $[h] \otimes \FC[0,ae]$ and thus equals $\CorG(0,a)$.  
\end{proof}

In other words, we have that there is an isomorphism
\begin{align*} F^\times_\bb{Q} &\overset{\cong}\lra \scr{G}_1(F) \\
a &\longmapsto \CorG(0,a).\end{align*}

\begin{remark}Alternatively, one may set $\{a\}_1 \coloneqq (1-a)^{-1}$ and then this isomorphism takes the form $\{a\}_1 \mapsto -\CorG(1,a)$ (the elements $(1-a)^{-1}$ and $(a-1)^{-1}$ agree in $F^\times_\bb{Q}$), which is more in line with the other formulas in \cref{thm:polyl-identification} and is supposed to call to mind that $\rm{Li}_1(z) = -\log(1-z)$.\end{remark}

\subsubsection{Explicit description of  $\PolyL_2$} Next we show that the weight two component of the Lie coalgebra $\PolyL(F)$ is isomorphic to the Bloch group. 

\begin{definition}The \emph{Bloch group} $B_2(F)$ is defined as a quotient of the group $\Z[F^{\times}\setminus \{1\}]$ by the subgroup spanned by elements 
\begin{equation}\label{Equation5term}
R_2(a,b)=[a]-[b]+\left [ \frac{b}{a}\right ] - \left [ \frac{1-a^{-1}}{1-b^{-1}}\right ] + \left [ \frac{1-a}{1-b}\right ]
\end{equation}
for $a\neq b \in F^{\times}\setminus \{1\}$. We denote the projection of $[a]$ to $B_2(F)$ by $\{a\}_2$.\end{definition}

\begin{remark}Suslin called this the \emph{pre-Bloch group} in  \cite{Sus90}, and denoted it $\fr{p}(F)$. He also showed for any $a\in F^{\times}\setminus \{1\}$ elements $\{a\}_2+\{a^{-1}\}_2$ and $\{a\}_2+\{1-a\}_2$ are $6$-torsion.\end{remark}

\begin{proposition} \label{prop: G_2=B_2} There is a well-defined isomorphism
\begin{align*}
B_2(F)_{\Q} &\stackrel{\cong}{\lra} \PolyL_2(F) \\
\{a\}_2 &\longmapsto -\CorG(1,0,a) = \CorG(0,1,a).
\end{align*}
\end{proposition} 
\begin{proof}The decomposition operator is given by
\begin{align*}
D^\FC_h(\FC[u_0:u_1:u_2])
=&\,\FC\left[0:\frac{u_{1}-u_{0}}{h(u_{1})-h(u_0)}:\frac{u_{2}-u_{0}}{h(u_{2})-h(u_0)}\right]\\
&-\FC\left[0:\frac{u_{1}-u_{0}}{h(u_{1})-h(u_0)}:\frac{u_{2}-u_{1}}{h(u_{2})-h(u_1)}\right]\\
&+\FC\left[0:\frac{u_{2}-u_{0}}{h(u_{2})-h(u_0)}:\frac{u_{2}-u_{1}}{h(u_{2})-h(u_1)}\right],
\end{align*}
leading to the following identity for correlators:
\begin{align*}
    &\CorG(x_0,x_1,x_2)-\CorG(y_0, y_1,y_2)\\
     &=\CorG\Bigl(0,\frac{x_1-x_0}{y_1-y_0},\frac{x_2-x_0}{y_2-y_0}\Bigr)-\CorG\bigl(0,\frac{x_1-x_0}{y_1-y_0},\frac{x_2-x_1}{y_2-y_1}\Bigr)+\CorG\Bigl(0,\frac{x_2-x_0}{y_2-y_0},\frac{x_2-x_1}{y_2-y_1}\Bigr).
\end{align*}
Specializing to $(x_0,x_1,x_2)=(0,1,a)$ and $(y_0,y_1,y_2)=(0,1,b)$ we obtain an identity
\begin{align*}
    &\CorG(0,1,a)-\CorG(0, 1,b)=\CorG\Bigl(0,1,\frac{a}{b}\Bigr)-\CorG\bigl(0,1,\frac{a-1}{b-1}\Bigr)+\CorG\Bigl(0,\frac{a}{b},\frac{a-1}{b-1}\Bigr).
\end{align*}
We have 
\[\CorG\Bigl(0,1,\frac{a}{b}\Bigr)=\CorG\Bigl(0,b,a\Bigr)=-\CorG\Bigl(0,1,\frac{b}{a}\Bigr) \qquad \text{and}\]
\[\CorG\Bigl(0,\frac{a}{b},\frac{a-1}{b-1}\Bigr)=\CorG\Bigl(0,1,\frac{1-a^{-1}}{1-b^{-1}}\Bigr).
\]
Thus we obtain that
\begin{align*}
    &\CorG(0,1,a)-\CorG(0, 1,b)+\CorG\Bigl(0,1,\frac{b}{a}\Bigr)+\CorG\bigl(0,1,\frac{a-1}{b-1}\Bigr)-\CorG\Bigl(0,1,\frac{1-a^{-1}}{1-b^{-1}}\Bigr)=0.
\end{align*}
This shows that the map from $B_2(F)$ to $\PolyL_2(F)$ is well-defined. To show it is an isomorphism, we use that the formula 
\[
\CorG(x_0, x_1,x_2) \longmapsto
\begin{cases}
\left\{\dfrac{x_2-x_0}{x_1-x_0}\right\}_2 & \text{ if } x_0,x_1,x_2 \text{ are distinct},\\
0& \text{ otherwise,}
 \end{cases}
\]
defines an inverse. To see it is well-defined, the only nontrivial observation needed is one made above: the decomposition relation corresponds to the 5-term relation.
\end{proof}

\begin{example}As mentioned in the introduction, let us spell out that under this isomorphism the weight 2 component of the Chevalley--Eilenberg complex computing the homology of the Lie coalgebra $\PolyL(F)$ is the rationalisation of the two-term Bloch complex \cite{Sus90}
\[
B_2(F)\lra \Lambda^2 F^{\times}
\]
with the differential given by $\{a\}_2$ to $a\wedge (1-a)$. This was stated in the introduction and follows from the computation of the cobracket $\delta(\CorG(x_0, x_1,x_2))$ as
\[\CorG(x_0, x_1)\wedge \CorG(x_0, x_2)+\CorG(x_1, x_2)\wedge \CorG(x_1, x_0)+\CorG(x_2, x_0)\wedge \CorG(x_2, x_1),\]
which leads to $\delta(\CorG(0,1,a))=a\wedge (1-a)$.
\end{example}

\subsubsection{Iterated integrals and multiple polylogarithms}\label{sec:iterated-integrals-and-multiple-polylogarithms}
We start by giving an analogue of the definition \cite[Definition 16]{CMRR24} in the setting of formal polylogarithms.

\begin{definition}For $x_0,\dots,x_{n+1} \in F$ and $n\geq 1$ we define the \emph{iterated integral} by the formula:
\[
			\ItG(x_0; x_1,\compactldots,x_n; x_{n+1}) \coloneq  \CorG(x_1,x_2,\compactldots,x_{n+1}) - \CorG(x_0, x_1,\compactldots, x_n) \in \PolyL_n(F).
\]
By convention, $\ItG(x_0; x_1) = 0$.\end{definition}

\cref{LemmaCorrelatorWeightOne} implies that under the isomorphism \eqref{eqn:isomorphism in weight 1} the iterated integral $\ItG(x_0; x_1; x_{2})$ corresponds to $\tfrac{x_2-x_1}{x_0-x_1}\in F^{\times}$. The properties of correlators immediately imply the shuffle relations for iterated integrals: 
\[\sum_{\sigma \in \rm{Sh}(n_1,n_2)} \ItG(x_0;x_{\sigma(1)},\compactldots,x_{\sigma(n_1+n_2)},x_{n+1}) =0.\]
Furthermore, we have the following formula for the cobracket of iterated integrals:
\begin{align*}
\delta & \ItG(x_0; x_1,\compactldots,x_n; x_{n+1}) \\
&=\sum_{0 \leq i < j \leq n+1} \!\!\! \ItG(x_0; x_1,\compactldots,x_i, x_j,\ldots, x_n; x_{n+1}) \wedge \ItG(x_i; x_{i+1}, \compactldots, x_{j-1}; x_j)  \,.
\end{align*}
The proof is identical to that of \cite[Proposition 19]{CMRR24}.

\medskip

We next give an analogue of the definition \cite[(23)]{CMRR24}.

\begin{definition}For an integer $n_0\geq 0$, $k\ge1$, positive integers $n_1,\dots,n_k$, and elements $a_1,\dots,a_k\in F^\times$ we define the \emph{multiple polylogarithm}
\begin{align*}
&\LiG_{n_0; n_1,\compactldots,n_k}(a_1,a_2,\compactldots,a_k)\\
& \coloneq (-1)^{k}\ItG(0;\underbrace{0,\compactldots,0,1}_{n_0+1},\underbrace{0,\compactldots,0,a_1}_{n_1},\compactldots,\underbrace{0,\compactldots,0,a_1a_2\compactldots a_{k-1}}_{n_{k-1}},\underbrace{0,\compactldots,0;a_1a_2\compactldots a_{k}}_{n_k})\,.
\end{align*}
This is an element of $\PolyL_n(F)$ for $n=n_0+n_1+n_2+\dots+n_k$.\end{definition}

If $n_0=0$, we omit it from the notation:
\[
\LiG_{n_1,\compactldots,n_k}(a_1,a_2,\compactldots,a_k) \coloneq \LiG_{0;n_1,\compactldots,n_k}(a_1,a_2,\compactldots,a_k)\,.
\]
We call elements $\LiG_{n}(a)$ \emph{classical polylogarithms}, and more specifically, call $\LiG_2(a)$ the \emph{dilogarithm} and $\LiG_3(a)$ the \emph{trilogarithm}. We will need that they satisfy the following property:

\begin{lemma}[Inversion relation for classical polylogarithms] \label{lemma: inversion for classical polylogs} For $n\geq 2$ we have
\[
\LiG_n\left(\frac{1}{a}\right)=(-1)^{n-1}\LiG_n(a) \quad \text{for $a\in F^{\times}$.}
\]
\end{lemma}

\begin{proof}
By definition, we have
\[
\LiG_n(a)=-\ItG(0;1,\underbrace{0,\compactldots,0}_{n-1};a)=-\CorG(1,\underbrace{0,\compactldots,0}_{n-1},a)+\CorG(0,1,\underbrace{0,\compactldots,0}_{n-1}).
\]
By \eqref{eqn:cor vanishing} and the cyclic symmetry \eqref{enum:rel-goncharov-2} we have 
\[
\LiG_n(a)=-\CorG(\underbrace{0,\compactldots,0}_{n-1},a,1).
\]
By \cref{lemma: multiplicative homogenity relation for correlators}, \cref{lemma: reversal relation for correlators}, and the cyclic symmetry \eqref{enum:rel-goncharov-2} we have 
\begin{align*}
\LiG_n\left(\frac{1}{a}\right)&=-\CorG(\underbrace{0,\compactldots,0}_{n-1},\frac{1}{a},1)=-\CorG(\underbrace{0,\compactldots,0}_{n-1},1,a)=-\CorG(a,\underbrace{0,\compactldots,0}_{n-1},1)\\
&=(-1)^{n}\CorG(a,1,\underbrace{0,\compactldots,0}_{n-1})=(-1)^{n}\CorG(\underbrace{0,\compactldots,0}_{n-1},a,1)=(-1)^{n-1}\LiG_n(a).
\end{align*}
\end{proof}

\subsubsection{The depth filtration}\label{sec:depth-filtration} Like multiple polylogarithms, the Goncharov Lie coalgebra admits an increasing \emph{depth filtration}. We start by observing $B_1(\Dec_V;\FC(V))$ admits an increasing 
\[\mathcal{D}_\bullet B_1(\Dec_V;\FC(V))\] 
defined as follows: the subspace $\mathcal{D}_k B_1(\Dec_V;\FC(V))$ is spanned by those elements $[h] \otimes \FC[u_0:\compactcdots:u_n]$ with at most $k+1$ indices $i$ satisfying $h(u_i)\neq h(u_{0})$. Note that we have $\mathcal{D}_{n-1} B_1(\Dec_V;\FC(V))=B_1(\Dec_V;\FC(V))$.  

\begin{definition}The \emph{depth filtration} on $\PolyL_n(F)$ is defined as
\[\mathcal{D}_\bullet\PolyL_n(F) \coloneq \rm{im}\big[\mathcal{D}_\bullet B_1(\Dec_V;\FC(V)) \to \PolyL_n(F)\big].\]
\end{definition}

By definition, the space $\mathcal{D}_k \PolyL_n(F)$ is spanned by those elements $\CorG(x_0,\compactldots,x_n)$ with at most $k+1$ indices $i$ satisfying $x_{i}\neq x_0$. The multiple polylogarithm  
$\smash{\LiG_{n_0; n_1,\compactldots,n_k}(a_1,a_2,\compactldots,a_k)}$ lies in $\mathcal{D}_k\PolyL_n(F)$. The following statement can be deduced from shuffle relations similarly to \cite[Corollary 30]{CMRR24}.

\begin{lemma}\label{lemma: polylogs span} The space $\mathcal{D}_k\PolyL_n(F)$ is spanned by multiple polylogarithms
\[\LiG_{n_1,\compactldots,n_d}(a_1,a_2,\compactldots,a_d)\] 
for $d\leq k$ and $a_1,\dots,a_d\in F^{\times}$.
\end{lemma}

Inspired by \cite{Gon95}, we introduce the following definition:

\begin{definition}The \emph{Bloch group $B_n^{\PolyL}(F)$ of weight $n$} is the depth one subspace $\mathcal{D}_1\PolyL_n(F)$.\end{definition} 

\cref{lemma: polylogs span} implies that $B_n^{\PolyL}(F)$ is spanned by classical polylogarithms. Note that $B_n^\PolyL(F) = \PolyL_n(F)$ for $n=1,2,3$, but in general this is no longer the case for $n \geq 4$.

\subsection{Explicit description of $\PolyL_3$} \label{sec: bloch group weight 3}
In this section we show $\PolyL_3(F)$ is spanned by trilogarithms and that the classical trilogarithmic identities hold---the $3$-term relation, Kummer’s relation, and Goncharov’s $22$-term relation. We prove that all relations in $\PolyL_3(F)$ follow from Goncharov’s $22$-term relation and use this to construct an isomorphism between $\PolyL_3(F)$ and Goncharov's $B_3(F)$. Together with the results of \cref{sec: G1 and G_2}, this finishes the identification of $\PolyL_n(F)$ for $n \leq 3$ in more classical terms.

\subsubsection{The Bloch group $B_3$ and equations for trilogarithm} 
Goncharov defined the Bloch group $B_3(F)$ as a conjectural symbolic description of the weight $3$ component of the motivic Lie coalgebra, whose existence is itself conjectural. In this section we will prove it agrees with $\scr{G}_3(F)$, proving the remaining one-third of \cref{thm:polyl-identification}. There are slight variations in the definition of $B_3(F)$ in the literature and to avoid potential subtleties related to torsion, we work with the rationalized Bloch group $B_3(F)_{\Q}$ defined as follows:

\begin{definition}The group $B_3(F)_{\Q}$  is a quotient of $\Q[F^{\times}]$ by the subspace spanned by elements 
\begin{equation} \label{eqn: 2 term in B3}
[a]-\left[a^{-1}\right] \text{ for }a\in F^{\times},
\end{equation}
the \emph{$3$-term relation} 
\begin{equation}\label{eqn: 3 term in B3}
T_3(a)= [a]+[1-a]+\left[1 - a^{-1}\right]-[1]  \text{ for }a\in F^{\times}\setminus \{1\}
\end{equation}
and the \emph{Goncharov  22-term relation} 
\begin{equation}\label{eqn: 22 term in B3}
\begin{aligned}
R_3(a,b,c) = {} & [ca-a+1]+[ab-b+1]+[bc-c+1]\\
&+ \left[\frac{ca-a+1}{ca}\right]+\left[\frac{ab-b+1}{ab}\right]+\left[\frac{bc-c+1}{bc}\right]\\
&+ \left[\frac{bc-c+1}{(ca-a+1)b}\right]+\left[\frac{ca-a+1}{(ab-b+1)c}\right]+\left[\frac{ab-b+1}{(bc-c+1)a}\right]\\
&- \left[\frac{ca-a+1}{c}\right]-\left[\frac{ab-b+1}{a}\right]-\left[\frac{bc-c+1}{b}\right]\\
&+ \left[-\frac{(bc-c+1)a}{(ca-a+1)}\right]+\left[-\frac{(ca-a+1)b}{(ab-b+1)}\right]+\left[-\frac{(ab-b+1)c}{(bc-c+1)}\right]\\
&-\left[\frac{(bc-c+1)}{(ca-a+1)bc}\right]-\left[\frac{(ca-a+1)}{(ab-b+1)ca}\right]-\left[\frac{(ab-b+1)}{(bc-c+1)ab}\right]\\
& + [a]+[b]+[c]+[-abc]-3 \, [1]
\end{aligned}
\end{equation}
for $a, b, c \in F^{\times}$ such that $ca-a+1,\: ab-b+1,\: bc-c+1\neq 0$. We denote the image of an element $[a]$ in $B_3(F)_{\Q}$ by $\{a\}_3$.\end{definition}

\begin{remark}
The above definition is a slight variation of the definition \cite[p.30]{Gon95b}, in which only elements $[a]$ with $a\neq 1$ are allowed. After rationalising, this definition yields the same group, as is shown in \cite[\S 5.3]{Gon95b}.
\end{remark}

\begin{remark}
The restriction of Goncharov's 22-term relation to $a=1$ gives  \emph{Kummer's 9-term relation}:
\begin{equation}\label{eqn: 9 term in B3}
\begin{aligned}
 0&=-\left\{\frac{bc-c+1}{b}\right\}_3-\left\{\frac{(bc-c+1)}{bc^2}\right\}_3-\left\{(bc-c+1)b\right\}_3+ 2\left\{-\frac{(bc-c+1)}{c}\right\}_3 \\ 
 &+2\{bc-c+1\}_3+2\left\{\frac{bc-c+1}{bc}\right\}_3
+2\{b\}_3+2\{c\}_3+2\{-bc\}_3-2\{1\}_3.
\end{aligned}
\end{equation}
\end{remark}

\subsubsection{Explicit form of the decomposition relations in weight $3$}\label{sec:decomp-rels-wt-3}
We begin by writing the equations defining $\PolyL_3(F)$ explicitly. To this end, recall that we have the following formula for the symbol of a Steinberg iterated integral:
\begin{align*}
s(\rm{I}[e_1,e_2,e_3]) = {} & 
 -s(\rm{I}[e_1,e_2])\otimes[e_3]
+\big(s(\rm{I}[e_1,e_3])
-s(\rm{I}[e_2,e_3])\big)\otimes[e_2{-}e_1] \\
& 
{} +\big(s(\rm{I}[e_1,e_2])-s(\rm{I}[e_1,e_3])\big)\otimes[e_3{-}e_2]\\[1ex]
= {} & - [e_1 \smid e_2 \smid e_3]+[e_1 \smid e_2{-}e_1 \smid e_3]-[e_2 \smid e_2{-}e_1 \smid e_3]
\\
 &+[e_1 \smid e_3 \smid e_2{-}e_1]-[e_1 \smid e_3{-}e_1 \smid e_2{-}e_1]+[e_3 \smid e_3{-}e_1 \smid e_2{-}e_1]\\
 &-[e_2 \smid e_3 \smid e_2{-}e_1]+[e_2 \smid e_3{-}e_2 \smid e_2{-}e_1]-[e_3 \smid e_3{-}e_2 \smid e_2{-}e_1]\\
 &+[e_1 \smid e_2 \smid e_3{-}e_2]-[e_1 \smid e_2{-}e_1 \smid e_3{-}e_2]+[e_2 \smid e_2{-}e_1 \smid e_3{-}e_2]\\
 &-[e_1 \smid e_3 \smid e_3{-}e_2]+[e_1 \smid e_3{-}e_1 \smid e_3{-}e_2]-[e_3 \smid e_3{-}e_1 \smid e_3{-}e_2]\,.
\end{align*}
The following shuffle relations hold in $\BSt_3(V)\otimes_{\SGroup_3}\rm{coLie}(3)$:
\begin{align*}
&[e_1 \smid e_2{-}e_1 \smid e_3]+[e_1 \smid e_3 \smid e_2{-}e_1]=-[e_3\smid e_1\smid e_2-e_1],\\
&[e_2 \smid e_2{-}e_1 \smid e_3]+[e_2 \smid e_3 \smid e_2{-}e_1]=-[e_3\smid e_2 \smid e_2-e_1],\\
&[e_2 \smid e_2{-}e_1 \smid e_3{-}e_2]+[e_2 \smid e_3{-}e_2 \smid e_2{-}e_1]=-
[e_3{-}e_2 \smid e_2 \smid e_2{-}e_1].
\end{align*}
Using them, we get
\begin{align*}
s(\rm{I}[e_1,e_2,e_3]) = 
 {} & - [e_1 \smid e_2 \smid e_3]-[e_3\smid e_1\smid e_2-e_1]+[e_3\smid e_2 \smid e_2-e_1]
\\
 &-[e_1 \smid e_3{-}e_1 \smid e_2{-}e_1]+[e_3 \smid e_3{-}e_1 \smid e_2{-}e_1]\\
 &-[e_3{-}e_2 \smid e_2 \smid e_2{-}e_1]-[e_3 \smid e_3{-}e_2 \smid e_2{-}e_1]\\
 &+[e_1 \smid e_2 \smid e_3{-}e_2]-[e_1 \smid e_2{-}e_1 \smid e_3{-}e_2]\\
 &-[e_1 \smid e_3 \smid e_3{-}e_2]+[e_1 \smid e_3{-}e_1 \smid e_3{-}e_2]-[e_3 \smid e_3{-}e_1 \smid e_3{-}e_2]\,.
\end{align*}

Next, consider an affine basis $u_0, u_1, u_2, u_3$ and a nonzero functional $h\in V^{\times}$ with  $h_i \coloneq h(u_i)$. If $h_0,\ h_1,\ h_2,\ h_3$ are pairwise distinct, we deduce that
\begin{align*}
D^\FC_h  &(\FC[u_0:u_1:u_2:u_3]) \\  
&=\FC\Bigl[0:\frac{u_1-u_0}{h_1-h_0} :\frac{u_2-u_0}{h_2-h_0} : \frac{u_3-u_0}{h_3-h_0}\Bigr]
+\FC\Bigl[0:\frac{u_3-u_0}{h_3-h_0}: \frac{u_1-u_0}{h_1-h_0}:\frac{u_2-u_1}{h_2-h_1}\Bigr]\\
&-\FC\Bigl[0:\frac{u_3-u_0}{h_3-h_0}: \frac{u_2-u_0}{h_2-h_0} : \frac{u_2-u_1}{h_2-h_1}\Bigr]
+\FC\Bigl[0:\frac{u_1-u_0}{h_1-h_0} : \frac{u_3-u_1}{h_3-h_1} : \frac{u_2-u_1}{h_2-h_1}\Bigr]\\
 &-\FC\Bigl[0:\frac{u_3-u_0}{h_3-h_0} : \frac{u_3-u_1}{h_3-h_1} : \frac{u_2-u_1}{h_2-h_1}\Bigr]
 +\FC\Bigl[0:\frac{u_3-u_2}{h_3-h_2} : \frac{u_2-u_0}{h_2-h_0} : \frac{u_2-u_1}{h_2-h_1}\Bigl]\\
& +\FC\Bigl[0:\frac{u_3-u_0}{h_3-h_0} : \frac{u_3-u_2}{h_3-h_2} : \frac{u_2-u_1}{h_2-h_1}\Bigr]
 -\FC\Bigl[0:\frac{u_1-u_0}{h_1-h_0} :\frac{u_2-u_0}{h_2-h_0} : \frac{u_3-u_2}{h_3-h_2}\Bigr]\\
 &+\FC\Bigl[0:\frac{u_1-u_0}{h_1-h_0} : \frac{u_2-u_1}{h_2-h_1} : \frac{u_3-u_2}{h_3-h_2}\Bigr]
 +\FC\Bigl[0:\frac{u_1-u_0}{h_1-h_0} : \frac{u_3-u_0}{h_3-h_0} : \frac{u_3-u_2}{h_3-h_2}\Bigr]\\
 &-\FC\Bigl[0:\frac{u_1-u_0}{h_1-h_0} : \frac{u_3-u_1}{h_3-h_1} : \frac{u_3-u_2}{h_3-h_2}\Bigr]
+\FC\Bigl[0:\frac{u_3-u_0}{h_3-h_0} : \frac{u_3-u_1}{h_3-h_1} : \frac{u_3-u_2}{h_3-h_2}\Bigr]\,.
\end{align*}
If some of the elements $h_0, h_1, h_2, h_3$ coincide, the above expression simplifies; all terms in which one of the denominators appearing in the entries vanishes are omitted.

Next, we consider a pair of nonzero functionals $h_1, h_2$ and denote $h_1(e_i):=x_i, h_2(e_i)=y_i$. If $y_0, y_1, y_2, y_3$ are distinct, we get the decomposition relation
\begin{align*}
\CorG&\bigl(x_0,x_1 ,x_2 , x_3\bigr)-\CorG\bigl(y_0, y_1, y_2, y_3\bigr)\\
&=\CorG\Bigl(0,\frac{x_1-x_0}{y_1-y_0} ,\frac{x_2-x_0}{y_2-y_0} , \frac{x_3-x_0}{y_3-y_0}\Bigr)
+\CorG\Bigl(0,\frac{x_3-x_0}{y_3-y_0}, \frac{x_1-x_0}{y_1-y_0},\frac{x_2-x_1}{y_2-y_1}\Bigr)\\
&-\CorG\Bigl(0,\frac{x_3-x_0}{y_3-y_0}, \frac{x_2-x_0}{y_2-y_0}, \frac{x_2-x_1}{y_2-y_1}\Bigr)
+\CorG\Bigl(0,\frac{x_1-x_0}{y_1-y_0} , \frac{x_3-x_1}{y_3-y_1} , \frac{x_2-x_1}{y_2-y_1}\Bigr)\\
&-\CorG\Bigl(0,\frac{x_3-x_0}{y_3-y_0} , \frac{x_3-x_1}{y_3-y_1} , \frac{x_2-x_1}{y_2-y_1}\Bigr)
+\CorG\Bigl(0,\frac{x_3-x_2}{y_3-y_2},\frac{x_2-x_0}{y_2-y_0}, \frac{x_2-x_1}{y_2-y_1}\Bigr)\\
&+\CorG\Bigl(0,\frac{x_3-x_0}{y_3-y_0} , \frac{x_3-x_2}{y_3-y_2} , \frac{x_2-x_1}{y_2-y_1}\Bigr)
-\CorG\Bigl(0,\frac{x_1-x_0}{y_1-y_0} ,\frac{x_2-x_0}{y_2-y_0} , \frac{x_3-x_2}{y_3-y_2}\Bigr)\\
&+\CorG\Bigl(0,\frac{x_1-x_0}{y_1-y_0} , \frac{x_2-x_1}{y_2-y_1} , \frac{x_3-x_2}{y_3-y_2}\Bigr)
 +\CorG\Bigl(0,\frac{x_1-x_0}{y_1-y_0} , \frac{x_3-x_0}{y_3-y_0} , \frac{x_3-x_2}{y_3-y_2}\Bigr)\\
&-\CorG\Bigl(0,\frac{x_1-x_0}{y_1-y_0} , \frac{x_3-x_1}{y_3-y_1} , \frac{x_3-x_2}{y_3-y_2}\Bigr)
 +\CorG\Bigl(0,\frac{x_3-x_0}{y_3-y_0} , \frac{x_3-x_1}{y_3-y_1} , \frac{x_3-x_2}{y_3-y_2}\Bigr)\,.
\end{align*}
If some of the elements $y_0, y_1, y_2, y_3$ coincide, we omit all correlators containing a fraction with a vanishing denominator.  

\subsubsection{Classical trilogarithms} \label{section: trilog} Recall that the trilogarithm is defined by the formula 
\[\LiG_3(a)=-\CorG(0,0,a,1) = -\CorG(1,0,0,a)\] 
for $a\in F$. By \cref{lemma: inversion for classical polylogs}, $\LiG_3(\frac{1}{a})=\LiG_3(a)$.  We start with the following statement:
\begin{lemma}[3-term relation]\label{lem 3-term relation}
    For any $a\in F^{\times}$ we have
    \begin{equation}
    \LiG_3\left(a\right)+\LiG_3(1-a)+\LiG_3\left(1-\frac{1}{a}\right)=\LiG_3(1)\label{EqTrilog2}.
    \end{equation}
\end{lemma}
\begin{proof} We consider a degenerate version of the decomposition relation with $y_0=y_3$:
\begin{align*}
\CorG&\bigl(x_0,x_1 ,x_2 , x_3\bigr)-\CorG\bigl(y_0, y_1, y_2, y_3\bigr)\\
&=\CorG\Bigl(0,\frac{x_1-x_0}{y_1-y_0} , \frac{x_3-x_1}{y_3-y_1} , \frac{x_2-x_1}{y_2-y_1}\Bigr)\\
&+\CorG\Bigl(0,\frac{x_3-x_2}{y_3-y_2},\frac{x_2-x_0}{y_2-y_0}, \frac{x_2-x_1}{y_2-y_1}\Bigr)\\
&-\CorG\Bigl(0,\frac{x_1-x_0}{y_1-y_0} ,\frac{x_2-x_0}{y_2-y_0} , \frac{x_3-x_2}{y_3-y_2}\Bigr)\\
&+\CorG\Bigl(0,\frac{x_1-x_0}{y_1-y_0} , \frac{x_2-x_1}{y_2-y_1} , \frac{x_3-x_2}{y_3-y_2}\Bigr)\\
&-\CorG\Bigl(0,\frac{x_1-x_0}{y_1-y_0} , \frac{x_3-x_1}{y_3-y_1} , \frac{x_3-x_2}{y_3-y_2}\Bigr).
\end{align*}
Evaluating it at $x_0=x_2=x_3=y_0=y_3=0,\: x_1=y_1=1,\: y_2=a$ we get
\begin{align*}
\CorG\bigl(0,1 ,0 ,0 \bigr)&-\CorG\bigl(0, 1,a, 0\bigr)=\CorG\Bigl(0,1, 1, \frac{-1}{a-1}\Bigr)+\CorG\Bigl(0,0,0, \frac{-1}{a-1}\Bigr)\\
&-\CorG\Bigl(0,1 ,0 , 0\Bigr)+\CorG\Bigl(0,1 , \frac{-1}{a-1} ,0\Bigr)-\CorG\Bigl(0,1 , 1 , 0\Bigr).
\end{align*}
Using \eqref{eqn:cor vanishing}, we can simplify this to 
\begin{equation}\label{eqn:3-term relatio proof}
-\CorG\bigl(0, 1,a, 0\bigr)=\CorG\Bigl(0,1, 1, \frac{-1}{a-1}\Bigr)+\CorG\Bigl(0,1 , \frac{-1}{a-1} ,0\Bigr)-\CorG\Bigl(0,1 , 1 , 0\Bigr).
\end{equation}
Notice that by homogeneity, cyclic symmetry, \cref{lemma: multiplicative homogenity relation for correlators}, and \cref{lemma: inversion for classical polylogs} we have
\begin{align*}
&\CorG\Bigl(0,1, 1, \frac{-1}{a-1}\Bigr)=\CorG\Bigl(-1,0,0, \frac{a}{1-a}\Bigr)=\CorG\Bigl(0,0, \frac{a}{1-a},-1\Bigr)\\
&=
\CorG\Bigl(0,0, \frac{a}{a-1},1\Bigr)=-\LiG\left(\frac{a}{a-1}\right)=-\LiG\left(1-\frac{1}{a}\right).
\end{align*}
Similarly,
\begin{align*}
    &\CorG\bigl(0, 1,a, 0\bigr)=-\LiG_3(a), \\
    &\CorG\Bigl(0,1 , \frac{-1}{a-1} ,0\Bigr)=-\LiG_3(1-a)\\
    &\CorG\Bigl(0,1 , 1 , 0\Bigr)=-\LiG_3(1).
\end{align*}
Thus \eqref{eqn:3-term relatio proof} is equivalent to the $3$-term relation.
\end{proof}

\begin{proposition}\label{prop: G3 spanned by trilogs} The space $\PolyL_3(F)$ is spanned by classical trilogarithms: we have $\PolyL_3(F)=B_3^{\PolyL}(F)$.
\end{proposition}
\begin{proof}
We consider the degenerate version of the decomposition relation with $y_0=y_3$ and specialize to $x_0=y_0=y_3=0,\: x_1=y_1=1 \: x_2=x_3=a, \: y_2=b$:
\begin{align*}
\CorG&\bigl(0,1 ,a ,a\bigr)-\CorG\bigl(0, 1, b, 0\bigr)\\
&=\CorG\Bigl(0,1 , \frac{a-1}{-1} , \frac{a-1}{b-1}\Bigr)+\CorG\Bigl(0,0,\frac{a}{b}, \frac{a-1}{b-1}\Bigr)\\
&-\CorG\Bigl(0,1 ,\frac{a}{b} , 0\Bigr)+\CorG\Bigl(0,1 , \frac{a-1}{b-1} , 0\Bigr)-\CorG\Bigl(0,1, \frac{a-1}{-1} , 0\Bigr).
\end{align*}
This identity can be rewritten as follows:
\begin{align*}
\CorG&\Bigl(0,1, 1-a , \frac{a-1}{b-1}\Bigr)\\
&=-\LiG_3\left(1-\frac{1}{a}\right)
+\LiG_3(b)
+\LiG_3\left(\frac{b(a-1)}{a(b-1)}\right)
-\LiG_3\left(\frac{a}{b}\right)
+\LiG_3\left(\frac{a-1}{b-1}\right)
-\LiG_3\Bigl(1-a \Bigr).
\end{align*}
Using the properties of the correlator and the $3$-term relation, we obtain the following identity:
\begin{equation}\label{eqn: correlator via trilogs}
\begin{aligned}
&\CorG(x_0,x_1, x_2, x_3) =+\LiG_3\left(\frac{x_3 - x_2}{x_3 - x_0}\right)+\LiG_3\left(\frac{x_3 - x_2}{x_1 - x_2}\right) \\
&-\LiG_3\left(\frac{(x_1 - x_2)(x_3 - x_0)}{(x_1 - x_0)(x_3 - x_2)}\right) + 
   \LiG_3\left(\frac{x_3 - x_0}{x_1 - x_0}\right)+\LiG_3\left(\frac{x_2 - x_1}{x_0 - x_1}\right)- \LiG_3(1).
\end{aligned}
\end{equation}
\end{proof}

\subsubsection{The isomorphism $\PolyL_3\cong B_3$}
Consider the map 
\begin{equation} \label{eqn: map from B_3 to G_3}
\begin{aligned}
  \tilde{L} \colon\Q[F^{\times}] &\lra \PolyL_3(F) \\
  [a] &\longmapsto \LiG_3(a).\end{aligned}
\end{equation}

\begin{lemma}The map   $\tilde{L}$ vanishes on relations \eqref{eqn: 2 term in B3}, \eqref{eqn: 3 term in B3}, \eqref{eqn: 22 term in B3} and thus defines a map $L$ from $B_3(F)$ to $\PolyL_3(F)$. In particular, the Goncharov $22$-term relation holds in $\PolyL_3(F)$.
\end{lemma}
\begin{proof}
By \cref{lemma: inversion for classical polylogs} the map $\tilde{L}$ vanishes on relation \eqref{eqn: 2 term in B3}. By \cref{lem 3-term relation} the map $\tilde{L}$ vanishes on relation \eqref{eqn: 3 term in B3}. It remains to show that $\tilde{L}$ vanishes on the 22-term relation. 

We start with a degenerate version of the decomposition relation with $y_0=y_3$. We put $x_0=y_0=y_3=0, \: x_1=y_1=1, \: x_2=a, x_3=b, y_2=c$:
\begin{align*}
&\CorG(0,1,a,b) - \CorG(0,1,c,0) = 
\CorG\left(0,\,1,\,1{-}b,\,\frac{a{-}1}{c{-}1}\right)  + \CorG\left(0,\,\frac{a{-}b}{c},\,\frac{a}{c},\,\frac{a{-}1}{c{-}1}\right) \\
&- \CorG\left(0,\,1,\,\frac{a}{c},\,\frac{a{-}b}{c}\right)  + \CorG\left(0,\,1,\,\frac{a{-}1}{c{-}1},\,\frac{a{-}b}{c}\right) - \CorG\left(0,\,1,\,1{-}b,\,\frac{a{-}b}{c}\right).
\end{align*}
We rewrite each term in terms of trilogarithms, and, after applying the $3$-term relation several times, obtain the following identity:
\begin{align*}
0=&\LiG_3\!\left(\frac{b-a}{b}\right)
+ \LiG_3\!\left(\frac{b-a}{1-a}\right)
- \LiG_3\!\left(\frac{(1-a)b}{b-a}\right)
+ \LiG_3(1-a)
+ \LiG_3(c)
+ \LiG_3(b)\\
&\quad
+ \LiG_3\!\left(\frac{a-b}{c}\right)
+ \LiG_3\!\left(\frac{c-a}{c}\right)
+ \LiG_3\!\left(\frac{c-a}{1-a}\right)
- \LiG_3\!\left(\frac{a-b-c+bc}{-1+a}\right) \\
&\quad
+ \LiG_3\!\left(\frac{(-1+a)b}{a-b-c+bc}\right)
- \LiG_3\!\left(\frac{a-c}{(-1+a)c}\right)
+ \LiG_3\!\left(-\frac{(-1+a) b c}{(a-b)(a-c)}\right)
+ \LiG_3\!\left(-\frac{b}{c-a}\right) \\
&\quad
- \LiG_3\!\left(\frac{(a-b)(a-c)}{b c}\right)
- \LiG_3\!\left(\frac{a-b-c+bc}{(a-c)c}\right)
+ \LiG_3\!\left(\frac{(a-b)(a-c)}{a-b-c+bc}\right) \\
&\quad
+ \LiG_3\!\left(-\frac{a-b-c+bc}{-a+b}\right)
+ \LiG_3\!\left(\frac{a-b-c+bc}{b c}\right)
- \LiG_3\!\left(\frac{(a-b)b}{a-b-c+bc}\right) \\
&\quad
+ \LiG_3\!\left(\frac{(-1+a)c}{a-b-c+bc}\right)
+ \LiG_3\!\left(-\frac{-a+c}{a-b-c+bc}\right)- 3\,\LiG_3(1) .
\end{align*}
This identity can be rewritten in the form
\[
\tilde{L}\left(R_3\Bigl(\frac{1}{1-a},-\frac{1-a}{a-b},c\Bigr)\right)=0.
\]
Thus the map $\tilde{L}$ vanishes on the $22$-term relation.
\end{proof}

We obtain a map 
\[
L\colon B_3(F)_\bb{Q}\lra \PolyL_3(F),
\]
which is surjective by \cref{prop: G3 spanned by trilogs}. To show injectivity, we construct a map $M$ in the opposite direction. We define it on generators $\CorG(x_0,x_1, x_2, x_3)$ as follows. If  the arguments $x_i$ are pairwise distinct, we put
\begin{align*}
&M(\CorG(x_0,x_1, x_2, x_3)\} =\left\{\frac{x_3 - x_2}{x_3 - x_0}\right\}_3+\left\{\frac{x_3 - x_2}{x_1 - x_2}\right\}_3 \\
&-\left\{\frac{(x_1 - x_2)(x_3 - x_0)}{(x_1 - x_0)(x_3 - x_2)}\right\}_3 + 
   \left\{\frac{x_3 - x_0}{x_1 - x_0}\right\}_3+\left\{\frac{x_2 - x_1}{x_0 - x_1}\right\}_3- \{1\}_3.
\end{align*}
If a certain pair of arguments coincide, we use the same formula omitting the terms which are not in $F^{\times}$. If there exist three arguments equal to each other, the corresponding element of $\PolyL_3(F)$ vanishes. Finally, in the case when there are two pairs of equal arguments, we put 
\[
    M(\CorG(x_0,x_1, x_2, x_3))=
    \begin{cases}
 -\{1\}_3 & \text{ if } x_0=x_1,\: x_2=x_3,\\
  2\{1\}_3 & \text{ if } x_0=x_2,\: x_1=x_3,\\
   -\{1\}_3 & \text{ if } x_0=x_3,\: x_1=x_2.
    \end{cases}
\]
The check that $M$ is well defined is rather tedious and we omit it. The key observation is that the decomposition relation follows from the $22$-term relation. The corresponding substitution was found by Steven Charlton via a computer-assisted search.
It is easy to see that $L$ and $M$ are mutually inverse. For instance, for the distinct arguments $x_i$, we have
\[
 L(M(\CorG(x_0,x_1, x_2, x_3)))=\CorG(x_0,x_1, x_2, x_3)
\]
by \eqref{eqn: correlator via trilogs}. This finishes the proof of the following result of \cref{thm:polyl-identification}:

\begin{theorem}There is a well-defined isomorphism
\begin{align*} B_3(F)_\bb{Q} &\overset{\cong}\lra  \PolyL_3(F)\\
\{a\}_3 &\longmapsto \LiG_3(a) = -\CorG(1,0,0,a).\end{align*}
\end{theorem}

\subsection{Realisations}
The Lie coalgebra of formal multiple polylogarithms $\PolyLF(F)$ was constructed in \cite{CMRR24}. The goal of this section is to construct a map of Lie coalgebras from $\PolyL(F)$ to $\scr{L}^\rm{f}(F)$. Using the results of \cite[\S 5]{CMRR24}, we then construct Hodge and motivic realisations of $\PolyL(F)$. 

\subsubsection{Lie coalgebra of formal polylogarithms}\label{sec:formal-polylogarithms}
We recall the main properties of the Lie coalgebra $\PolyLF(F)$ of formal polylogarithms from \cite[\S2]{CMRR24}. For any field $F$, there exists a graded Lie coalgebra $\PolyLF(F)$, which is generated, as a $\Q$-vector space, by \emph{formal correlators}:
\[
\CorF(x_0,\compactldots,x_n)\in \PolyLF_n(F)\quad  \text{for}\quad x_0,\dots,x_n\in F.
\]
The cobracket of formal correlators is given by the following formula:
\begin{equation} \label{eqn:form-cor-cobracket}
\delta\CorF (x_0,\compactldots, x_n )=\sum_{j=0}^n\sum_{i=1}^{n-1}\CorF( x_{j}, x_{j+1}, \compactldots, x_{j+i}) \wedge \CorF( x_{j},  x_{j+i+1}, \compactldots, x_{j+n})\,.
\end{equation}
For $n=1$, we have an isomorphism $u\colon \PolyLF_1(F)\cong F^{\times}_{\Q}$ given by sending $\CorF(x_0,x_1)$ to $(x_1-x_0)\in F^{\times}_{\Q}$. The formal correlators satisfy the following identities, amongst others:
\begin{align}
&\CorF(x_0+b,\compactldots,x_n+b)=\CorF(x_0,\compactldots,x_n),\label{eqn:form-cor-add}\\
&\CorF(x_0,\compactldots,x_{n-1},x_n)=\CorF(x_1,\compactldots,x_n,x_0), \label{eqn:form-cor-cyc-sym}\\
&\CorF(m x_0,\compactldots,m x_n)=\CorF(x_0,\compactldots,x_n) \text{ for } m\in F^{\times} \text{ and } n\geq 2, \label{eqn:form-cor-homog}\\
& \CorF(x_0,\underbrace{x_1, \compactldots, x_1}_{n-1})=0 \text{ for } n\geq 2\label{eqn:form-cor-norm}\\
&\sum_{\sigma \in \rm{Sh}(n_1,n_2)} \CorF(x_0,x_{\sigma(1)},\compactldots,x_{\sigma(n_1+n_2)})=0\,,
\end{align}

The inclusion $F\hookrightarrow F(t)$ induces an injective map $i\colon \PolyLF(F)\lra \PolyLF(F(t))$. For any $t_0 \in F$, there is a specialisation homomorphism
\[
\rm{Sp}_{t\to t_0} \colon \PolyLF(F(t))\lra  \PolyLF(F) 
\]
defined as follows. For a tuple of functions $f_1(t),\dots,f_n(t)\in F(t)$ consider the smallest integer $d$ such that functions $u_i(t)=(t-t_0)^d f_i(t)$ are regular at $t=t_0$. We put 
\[
\rm{Sp}_{t\to t_0}\CorF(0,f_1(t),\compactldots,f_n(t))=\CorF(0, u_1(t_0),\compactldots,u_n(t_0)).
\]
For any point $t_0\in F$ we have $\rm{Sp}_{t\to t_0}\circ i=\id$. The following lemma allows one to prove identities in $\PolyLF(F)$ by induction. 

\begin{lemma}\label{lem:specialisation in formal polylogs} Consider $R\in \PolyLF_n(F(t))$ such that $\delta(R)=0$ and $n\geq 2$. Then for any $t_0, t_1\in F$ we have
\[
\Sp_{t\to t_0}R=\Sp_{t\to t_1}R.
\] 
\end{lemma}
\begin{proof}
By \cite[Corollary 11]{CMRR24}, the inclusion  $i\colon \PolyLF(F)\to \PolyLF(F(t))$ induces an isomorphism $H^1(\PolyLF(F))\cong H^1(\PolyLF(F(t)))$ on Chevalley--Eilenberg homology. Since $H^1(\PolyLF(F(t)))$ can be identified with the kernel of the cobracket, the element $R$ lies in the image of $i$. Moreover, 
\[
\rm{Sp}_{t\to t_0}\circ i=\rm{Sp}_{t\to t_1}\circ i=\id,
\]
so $\Sp_{t\to t_0}R=\Sp_{t\to t_1}R.$
\end{proof}

\subsubsection{Formal, motivic, and Hodge realisations}
The goal of this section is to construct a morphism from $\PolyL(F)$ to $\PolyLF(F)$ sending correlators to formal correlators. The key step is to prove that decomposition relations \eqref{enum:rel-goncharov-4} hold in $\PolyLF(F)$. We prove this by induction on weight, using the fact that the cobracket on $\PolyL(F)$ is well-defined, as we believe this statement is difficult to verify by a direct computation.

\begin{theorem} \label{prop: formal realization} There exists a unique morphism of graded Lie coalgebras
\begin{align*}
r^{\rm{f}}\colon \PolyL(F) &\lra \PolyLF(F) \\ \CorG(x_0,\compactldots,x_n) &\longmapsto \CorF(x_0,\compactldots,x_n).\end{align*}
\end{theorem}
\begin{proof}
We will show that the map $r^{\rm{f}}$ is well-defined by induction on the  weight. We start with the base case $n=1$. By \cref{LemmaCorrelatorWeightOne}, there is an isomorphism
$\PolyL_1(F)\cong F^{\times}$ sending $\smash{\CorG(x_0,x_1)}$ to $(x_0-x_1)\in  F^{\times}$. Moreover, by \cite[Lemma 1]{CMRR24}, we have an  isomorphism $\PolyLF_1(F)\cong F^{\times}$ sending $\smash{\CorF(x_0,x_1)}$ to  $(x_0-x_1)\in  F^{\times}$. This proves the base case.

Assume that the map $r^{\mathrm{f}}$ is well defined for any field in weights less than $n$. The properties \eqref{eqn:form-cor-add}, \eqref{eqn:form-cor-cyc-sym}, and \eqref{eqn:form-cor-norm} imply that the map $r^{\mathrm{f}}_n$ vanishes on the homogeneity relation \eqref{enum:rel-goncharov-1}, the cyclic symmetry relation \eqref{enum:rel-goncharov-2}, and the shuffle relations \eqref{enum:rel-goncharov-3}. It remains to show that $r^{\mathrm{f}}$ kills the decomposition relation \eqref{enum:rel-goncharov-4} in weight $n$. 

First consider a field $L=F(X_0,\dots,X_n,Y_0,\dots,Y_n)$.
The decomposition relation for the tuples $tX_0+(1-t)Y_0,\dots, tX_n+(1-t)Y_n$ and  $Y_0,\dots, Y_n$ yields the following identity in $\PolyL_n(L(t))$:
\begin{align*}
  &\CorG(tX_0+(1-t)Y_0,\compactldots,tX_n+(1-t)Y_n)-\CorG(Y_0,\compactldots,Y_n)\\
     &=\sum_{\iota=((i_1,j_1),\compactldots,(i_{n},j_n))\in T(n)}\rm{sign}(\iota)
     \CorG\left(0,(1-t)+t\frac{X_{i_1}-X_{j_1}}{Y_{i_1}-Y_{j_1}},\compactldots,(1-t)+t\frac{X_{i_n}-X_{j_n}}{Y_{i_n}-Y_{j_n}}\right).
\end{align*}
Applying the cobracket, we obtain the identity 
\begin{align*}
  &\delta\CorG(tX_0+(1-t)Y_0,\compactldots,tX_n+(1-t)Y_n)-\delta\CorG(Y_0,\compactldots,Y_n)\\
     &=\sum_{\iota=((i_1,j_1),\compactldots,(i_{n},j_n))\in T(n)}\rm{sign}(\iota)
     \delta\CorG\left(0,(1-t)+t\frac{X_{i_1}-X_{j_1}}{Y_{i_1}-Y_{j_1}},\compactldots,(1-t)+t\frac{X_{i_n}-X_{j_n}}{Y_{i_n}-Y_{j_n}}\right)
\end{align*}
which holds in the weight $n$ component of $\Lambda^2 \PolyL(L(t))$. By the induction hypothesis, the map
$\Lambda^2r^{\mathrm{f}}\colon \Lambda^2 \PolyL(L(t)) \lra \Lambda^2 \PolyLF(L(t))$ is well-defined on the weight $n$ component. All terms of the relation above have distinct arguments, so the cobracket of each term can be computed using \cref{thm:polyl-presentation-cobracket}. Since the cobracket in $\PolyLF(L(t))$ is given by the analogous formula \eqref{eqn:form-cor-cobracket},  the element 
\begin{align*}
R=&\CorF(tX_0+(1-t)Y_0,\compactldots,tX_n+(1-t)Y_n)-\CorF(Y_0,\compactldots,Y_n)\\
     &-\sum_{\iota=((i_1,j_1),\compactldots,(i_{n},j_n))\in T(n)}\rm{sign}(\iota)
   \CorF\left(0,(1-t)+t\frac{X_{i_1}-X_{j_1}}{Y_{i_1}-Y_{j_1}},\compactldots,(1-t)+t\frac{X_{i_n}-X_{j_n}}{Y_{i_n}-Y_{j_n}}\right)
\end{align*}
satisfies the identity $\delta(R)=0$. Recall the specialisation homomorphisms from \cref{sec:formal-polylogarithms}. Using \eqref{eqn:form-cor-norm}, we obtain on the one hand an equation
\begin{align*}
\Sp_{t\to 0}(R)= &\CorF(Y_0,\compactldots,Y_n)-\CorF(Y_0,\compactldots,Y_n)-\sum_{\iota\in T(n)}\rm{sign}(\iota)
   \CorF\left(0,1,\compactldots,1\right)=0.
\end{align*}
On the other hand, \cref{lem:specialisation in formal polylogs} yields that 
\begin{align*}
 0=\Sp_{t\to 1}(R)= &\CorF(X_0,\compactldots,X_n)-\CorF(Y_0,\compactldots,Y_n)\\
     &-\sum_{\iota=((i_1,j_1),\dots,(i_{n},j_n))\in T(n)}\rm{sign}(\iota)
     \CorF\left(0,\frac{X_{i_1}-X_{j_1}}{Y_{i_1}-Y_{j_1}},\compactldots,\frac{X_{i_n}-X_{j_n}}{Y_{i_n}-Y_{j_n}}\right).
\end{align*}
This establishes the decomposition relation for tuples $(x_0,\dots,x_n)$ and $(y_0,\dots,y_n)$ such that points $y_0,\dots,y_n$ are distinct after specialisation of $X_i$ to $x_i$ and $Y_j$ to $y_j$. We call this relation the \emph{generic decomposition relation}. 

It remains to prove the decomposition relation for an arbitrary tuple $y_0,\dots,y_n$. Points $y_1,\dots,y_n$ are not all the equal, and we may assume $y_0\neq y_n$, by cyclic symmetry. Consider variables $C_0,C_1,\ldots,C_n$ and $t$. The generic decomposition relation implies that 
\begin{align*}
 0= &\CorF(x_0,\compactldots,x_n)-\CorF(y_0+C_0 T,\compactldots,y_n+C_n T)\\
     &-\sum_{\iota\in T(n)}\rm{sign}(\iota)
     \CorF\left(0,\frac{X_{i_1}-X_{j_1}}{(y_{i_1}-y_{j_1})+(C_{i_1}-C_{j_1})t},\compactldots,\frac{X_{i_n}-X_{j_n}}{(y_{i_n}-y_{j_n})+(C_{i_n}-C_{j_n})t}\right).
\end{align*}

Consider an equivalence relation $\sim$ on the set $\{0,\dots,n\}$ defined by $i\sim j$ when $y_i=y_j$; note that $0 \not \sim n$. We call a term 
\[
\CorF\left(0,\frac{X_{i_1}-X_{j_1}}{(y_{i_1}-y_{j_1})+(C_{i_1}-C_{j_1})t},\compactldots,\frac{X_{i_n}-X_{j_n}}{(y_{i_n}-y_{j_n})+(C_{i_n}-C_{j_n})t}\right)
\]
\emph{regular} if $i_1 \not \sim j_1$, \dots,  $i_n \not \sim j_n$; the other terms are called \emph{singular}. We claim that the sum of singular terms vanishes under the specialisation map
\[
\Sp_{t\to 0}\colon  \PolyLF(F(t,X_1,\dots,X_n,C_0,C_1,\dots,C_n))\lra  \PolyLF(F(X_1,\dots,X_n,C_0,C_1,\dots,C_n)).
\]
This would imply the statement of this proposition, after applying further specialisations sending $C_i$ to $0$ and $X_i$ to $x_i$. To prove the claim note that specialisation of a singular term can be described by making following substitutions:
\[\frac{X_{i}-X_{j}}{(y_{i}-y_{j})+(C_{i}-C_{j})t} \longmapsto \begin{cases} 0 & \text{if $i \not \sim j$,} \\
\frac{X_{i}-X_{j}}{C_{i}-C_{j}} & \text{if $i \sim j$.}\end{cases}\]
For example, specialisation of the term
\[
\CorF\left(0,\frac{X_{0}-X_{1}}{y_{0}-y_{1}+(C_{0}-C_{1})t},\frac{X_{1}-X_{2}}{(y_{1}-y_{2})+(C_{1}-C_{2})t},\frac{X_{2}-X_{3}}{(y_{2}-y_{3})+(C_{2}-C_{3})t}\right)
\]
for $y_0\neq y_1=y_2=y_3$ equals 
\[
\CorF\left(0,0,\frac{X_{1}-X_{2}}{(C_{1}-C_{2})},\frac{X_{2}-X_{3}}{(C_{2}-C_{3})}\right).
\]
\cref{lem:symbol and equivalence relation} applied to the set $S=\{0, \dots,n\}$ and equivalence relation $\sim$ implies that the specialisations of singular terms cancel out.
\end{proof}

The first part of the following result is \cref{thm:motivic-realisation} and the second part is \cref{thm:hodge-realisation}.

\begin{theorem}\,
\begin{enumerate}
    \item If $F$ is a number field there is a unique \emph{motivic realisation} map of Lie coalgebras
    \begin{align*} r^{\rm{MTM}}\colon \PolyL(F) &\lra \MotCoLie(F) \\
    \CorG(x_0,\compactldots,x_n) &\longmapsto \Cor^\rm{MTM}(x_0,\compactldots,x_n). \end{align*}
    \item For an embedding $\sigma \colon F \to \bb{C}$ there is a unique \emph{Hodge realisation} map of Lie coalgebras
    \begin{align*} r^\rm{Hod}_\sigma \colon \PolyL(F) &\lra \cal{L}^\rm{Hod}(F) \\
    \CorG(x_0,\compactldots,x_n) &\longmapsto \Cor^\rm{Hod}(\sigma(x_0),\compactldots,\sigma(x_n)). \end{align*}
\end{enumerate}
\end{theorem}

\begin{proof}The first part is a consequence of \cref{prop: formal realization} and \cite[Proposition 50]{CMRR24}, and the second is a consequence of \cref{prop: formal realization} and \cite[Proposition 48]{CMRR24}. The uniqueness follows since the correlators generate $\PolyL_n(F)$.\end{proof}

\begin{remark}The formal realisation $r^\rm{f} \colon \scr{G}_n(F) \to \scr{L}^\rm{f}_n(F)$ is an isomorphism for $n \leq 2$ by \cite[\S 4.4]{CMRR24}. One might wonder whether it is an isomorphism for $n \geq 3$ as well.
\end{remark}

As an application, we use formal realisation to give a noncomputational proof of the cobracket of classical trilogarithm, though a computational proof is also possible albeit lengthy:

\begin{lemma}The cobracket of the classical trilogarithm is given by the  formula 
\begin{equation}\label{eqn:cobracket of classical polylog}
    \delta \bigl( \LiG_3(a)\bigr)=\LiG_2(a)\otimes a.
\end{equation}
\end{lemma}

\begin{proof} By  \cite[Proposition 45]{CMRR24} the map $\PolyL_2(F)\to \PolyLF_2(F)$ is an isomorphism, so the cobracket of $\LiG_3(a)$ agrees with the cobracket of $\mathrm{Li}^{\mathrm{f}}_3(a)$ which is given by \eqref{eqn:cobracket of classical polylog}. 
\end{proof}

\section{The duality involution on $\scr{G}(F)$}\label{section duality involution} In this section we investigate the duality involution on the Goncharov Lie coalgebra. More precisely, taking duals is an involution of the category $\Vect$ of finite-dimensional vector spaces over $F$. This induces an involution on $\scr{G}(F)$ and we determine in this section that it acts by $(-1)^n$ in rank $n$. We follow \cref{conv:shorter-notation}.

\subsection{The duality involution} 
The symmetric monoidal groupoid $\Vect$ of finite-dimensional vector spaces over $F$ admits a symmetric monoidal automorphism
\[\vee \colon \Vect \lra \Vect\]
that takes an object $V$ to its linear dual $V^\vee$ and a linear isomorphism $A \colon V \to W$ to the inverse of its linear dual $(A^\vee)^{-1} \colon V^\vee \to W^\vee$, in terms of matrices given by taking the transpose inverse. This is an involution up to the natural equivalence $\rm{ev} \colon V \to (V^\vee)^\vee$ given by sending $v$ to the functional that evaluates linear functionals at $v$. This fits in a commutative diagram with $\dim \colon \Vect \to \bb{N}$ and the identity automorphism of $\bb{N}$. Taking classifying spaces as in \cref{sec:action-by-scaling}, this lifts $\BGLb^+$ to a functor
\[\BGLb^+ \colon \rm{B}C_2 \lra \Alg_{E_\infty^\rm{u}}(\Fun(\N,\Spc)).\] 
Unwinding the definitions, in rank $n$ this is the following map. For each $n$-dimensional vector space $V$, the inclusion the orbit groupoid $*{\sslash} \GL(V) \to \Vect$ induces an equivalence $\BGL(V) \to \BGLb^+(n)$ and the duality involution restricts to the map
\[\rm{B}{\vee} \colon \BGL(V) \lra \BGL(V^\vee)\]
induced by the isomorphism of groups given by $A \mapsto (A^\vee)^{-1}$. Upon picking an isomorphism $V^\vee \cong V$, we can identify the domain and target; any two isomorphisms induce homotopic identifications on classifying spaces. In the special case that $V = F^n$, the canonical basis allows us to identify $\GL(V) = \GL_n = \GL(V^\vee)$ and $\rm{B}{\vee}$ is induced by the automorphism of the group $\GL_n$ given by $A \mapsto (A^t)^{-1}$.

Any object obtained naturally from $\BGLb^+$ as a unital $E_\infty$-algebra in spaces inherits an involution. In particular, postcomposing with the rationalisation $\Spc \to \DQ$, passing to the augmentation ideal of the canonical augmentation, and taking $E_\infty^\rm{nu}$-indecomposables, we get a functor
\[\cot_{E_\infty^\rm{nu}}(\BGLb_\bb{Q}) \colon \rm{B}C_2 \lra \Fun(\bb{N},\DQ)\]
that induces an involution on the $E_\infty$-homology.

\subsubsection{Description in terms of buildings} Our next goal is describe this involution more concretely in terms of the homotopy orbits of buildings. To do so, we take a slightly different perspective. First, note that $\vee$ induces a symmetric monoidal automorphism $\vee_!$ of $\Fun(\Vect,\scr{C})$ for any presentable symmetric monoidal category $\scr{C}$, and by naturality of the Day convolution an automorphism of $\Alg_{E_\infty^\rm{u}}(\Fun(\Vect,\scr{C}))$. Taking $\scr{C} = \Spc$, for the terminal $E_\infty^\rm{u}$-algebra $\ul{\ast}$ there is an essentially unique equivalence $\tau \colon  \ul{\ast} \simeq \vee_!(\ul{\ast})$. Using the natural equivalence $\psi \colon \dim_! \vee_! \simeq \dim_!$, we obtain the above automorphism as the composition
\[\BGLb^+ = \dim_!(\ul{\ast}) \underset{\simeq}{\xrightarrow{\dim_! \tau_!}} \dim_! \vee_!(\ul{\ast}) \underset{\simeq}{\xrightarrow{\psi}} \dim_!(\ul{\ast}) = \BGLb^+.\]

Since $\vee_!$ is a symmetric monoidal equivalence, on iterated bar constructions we obtain a natural equivalence $\Bar^k_\levi \vee_! \simeq \vee_! \Bar^k_\levi$ of functors $\Alg^\aug_{E_k^\rm{u}}(\Fun(\Vect,\Spc_*)) \to \Fun(\Vect,\Spc_*)$ and from this we obtain induced equivalences
\[\tau^k_! \colon \Bar^k_\levi(\ul{\ast}_+) \simeq \Bar^k_\levi(\vee_! \ul{\ast}_+) \simeq \vee_! \Bar^k_\levi(\ul{\ast}_+)\]
for $k \geq 1$. This can be understood more concretely using the identification \cite[Section 5.4]{GKRW20} explained in \cref{sec:buildings}
\[\Bar^k_\levi(\ul{\ast}_+)(V) \simeq \widetilde{D}^k(V)\]
of the bar construction with the split buildings of \cite[Definition 5.9]{GKRW20}. Recall the split building is a $k$-fold pointed simplicial set of direct sum decompositions of $V$ encoded by $k$-dimensional grids of summands (see loc.cit.~for details) and under the aforementioned identification, the involution is given by the map $\smash{\widetilde{D}^k(V)} \to \smash{\widetilde{D}^k(V^\vee)}$ taking a direct sum decomposition of $V$ to the dual one of $V^\vee$.

Applying $\dim_!$ we get an equivalence
\[\dim_! \Bar^k_\levi(\ul{\ast}_+) \underset{\simeq}{\xrightarrow{\dim_! \tau^k_!}} \dim_! \vee_! \Bar^k_\levi(\ul{\ast}_+) \underset{\simeq}{\xrightarrow{\psi_!}}  \dim_! \Bar^k_\levi(\ul{\ast}_+)\]
that in rank $n$ is as follows. For each $n$-dimensional vector space $V$, the inclusion induces an equivalence $\smash{\widetilde{D}^k(V)}_{\GL(V)} \to (\dim_! \Bar^k_\levi(\ul{\ast}_+))(n)$ on orbits and the dualisation involution restricts along these to the map
\[\smash{\widetilde{D}^k(V)}_{\GL(V)} \lra \smash{\widetilde{D}^k(V^\vee)}_{\GL(V^\vee)}\]
induced on orbits by the dualisation map.

To relate this to the nonsplit buildings that appear in \cref{sec:buildings}, recall the map $\widetilde{D}^k(V) \to D^k(V)$ forgetting from splitting to flags \cite[(5.6)]{GKRW20} induces a homology isomorphism \cite[Theorem 5.18]{GKRW20}
\[\dim_! \widetilde{D}^k(V) \lra \dim_! D^k(V)\]
and there is an evidently commutative square
\[\begin{tikzcd} \widetilde{D}^k(V) \rar{\vee} \dar & \widetilde{D}^k(V^\vee) \dar \\[-5pt]
D^k(V) \rar{\vee} & D^k(V^\vee) \end{tikzcd}\]
where the bottom map taking flags of $V$ to the dual flags of $V^\vee$. The nonsplit buildings are related to $E_k$-homology (and hence $E_\infty$-homology, by taking infinite bar spectra) by passing to rational chains and applying the natural equivalence $\Bar^k \simeq (\Sigma^k \indec_{E_k^\rm{nu}})^+$ of \cref{thm:indec-is-bar}, so we can understand the duality involution on $E_\infty$-homology by computing the isomorphisms $\StL(V) \to \StL(V^\vee)$ via the action on buildings and taking the induced maps on homology
\[H_*(\GL(V);\StL(V)) \lra H_*(\GL(V^\vee);\StL(V^\vee)).\]
The same holds for $E_1$- and $E_2$-homology, in terms of the Steinberg modules $\St(-)$ and double Steinberg modules $\StH(-)$. Upon picking an isomorphism $V \cong V^\vee$, we can identify the domain and target; this identification is independent of our choice of isomorphism, since the map induced by any other choice differs by an inner automorphism and acts as the identity on group homology \cite[III.8.1]{Brown}.

\subsubsection{The duality involution on Steinberg modules} We now implement this. We read off from the induced action on buildings that the involution on $\Vect$ induces an isomorphism of Steinberg modules
\begin{align*} \vee_* \colon \St(V) &\lra \St(V^\vee) \\
[v_1,\compactldots,v_n] &\longmapsto [v_1^\vee,\compactldots,v_n^\vee],\end{align*}
sending the apartment corresponding to a basis to the apartment corresponding to a dual basis. Similarly, we read off from the induced action on double buildings that the involution on $\Vect$ induces an isomorphism of double Steinberg modules
\begin{equation}\label{eqn:sth-duality-vee}\begin{aligned}\vee_* \colon \StH(V) &\lra \StH(V^\vee) \\
[v_1,\compactldots,v_n] \otimes [w_1,\compactldots,w_n] &\longmapsto [v_1^\vee,\compactldots,v_n^\vee] \otimes [w_1^\vee,\compactldots,w_n^\vee],\end{aligned}\end{equation}
sending the apartment corresponding to a basis to the apartment corresponding to a dual basis. This is \emph{not} the duality isomorphism $D$ considered in \cite[(30)]{CharltonRadchenkoRudenko}, which is rather given by
\begin{equation}\label{eqn:sth-duality-d} \begin{aligned} D \colon \StH(V) &\lra \StH(V^\vee) \\
[v_1,\compactldots,v_n] \otimes [w_1,\compactldots,w_n] &\longmapsto [w_1^\vee,\compactldots,w_n^\vee] \otimes [v_1^\vee,\compactldots,v_n^\vee].\end{aligned}\end{equation} 
Finally, the map $\vee_* \colon \StL(V) \to \StL(V^\vee)$ is induced by map on $\SStH$. The duality isomorphism $D$ of \eqref{eqn:sth-duality-d} also is compatible with the algebra structure and also induces a map $D \colon \StL(V) \to \StL(V^\vee)$. They are related as follows:

\begin{lemma}\label{lem:duality-on-stl-sign} As maps $\StL(V) \to \StL(V^\vee)$ we have
\[\vee_* = (-1)^{\dim(V)-1} D.\]
In particular, on $H_*(\GL_n;\StL_n)$ we have $\vee_* = (-1)^{n-1} D$.
\end{lemma}

\begin{proof}Comparing the formulas \eqref{eqn:sth-duality-d} and \eqref{eqn:sth-duality-vee}, we see that they differ by composition with an instance of the swap map
\begin{align*}\rm{rev} \colon \StH(W) & \lra \StH(W) \\
[v_1,\compactldots,v_n] \otimes [w_1,\compactldots,w_n] &\longmapsto [w_1,\compactldots,w_n] \otimes [v_1,\compactldots,v_n].\end{align*}
(One may be tempted to call it the Poincar\'e or Verdier duality involution, following \cite[2.1.6]{BGSV90}.) Recall that the symbol maps fit into a commutative diagram
\[\begin{tikzcd} \StH(V) \dar[two heads] \rar[hook]{s^\rm{As}} &[20pt]  (\rm{B}^\rm{As}\SSt)_n(V) \dar[two heads]\\[-5pt]
\StL(V)  \rar[hook]{s^\rm{Com}}  & (\rm{B}^\rm{Com}\SSt)_n(V)\end{tikzcd}\]
with horizontal maps injective, surjective left vertical map the canonical projection, and surjective right vertical map the quotient by shuffles (see \cref{sec:symbol-maps}). We now use the formula of \cite[Lemma 19]{CharltonRadchenkoRudenko}:
\[s([v_1,\compactldots,v_n] \otimes [w_1,\compactldots,w_n]) = \sum_{\sigma,\tau \in \fr{S}_n} (-1)^\sigma (-1)^\tau [F^\sigma_1 \cap G^\tau_n|F^\sigma_2 \cap G^\tau_{n-1}|\compactldots|F^\sigma_n \cap G^\tau_1]\]
where $F^\sigma_i = \rm{span}(v_{\sigma(1)},\ldots,v_{\sigma(i)})$ and $G^\tau_j = \rm{span}(w_{\tau(n-j+1)},\ldots,w_{\tau(n)})$, and a term in the sum is zero if the two flags are not in general position. Letting $\rho \in \fr{S}_n$ denote the involution with $\rho(i) = n-i+1$, we see that $s([w_1,\ldots,w_n] \otimes [v_1,\ldots,v_n])$ is given instead by 
\[\sum_{\sigma,\tau \in \fr{S}_n} (-1)^{\rho \sigma \rho^{-1}} (-1)^{\rho \tau \rho^{-1}} [G_1^{\rho \tau \rho^{-1}} \cap F_n^{\rho \sigma \rho^{-1}}|G_2^{\rho \tau \rho^{-1}} \cap F_{n-1}^{\rho \sigma \rho^{-1}}|\compactldots|G_n^{\rho \tau \rho^{-1}} \cap F_1^{\rho \sigma \rho^{-1}}].\]
Up to shuffles we have that $[v_1|\compactldots|v_n] = (-1)^{n-1} [v_n|\compactldots|v_1]$ so this is equal to 
\[(-1)^{n-1} \sum_{\sigma,\tau \in \fr{S}_n} (-1)^{\rho \sigma \rho^{-1}} (-1)^{\rho \tau \rho^{-1}} [G_n^{\rho \tau \rho^{-1}} \cap F_1^{\rho \sigma \rho^{-1}}|G_{n-1}^{\rho \tau \rho^{-1}} \cap F_{2}^{\rho \sigma \rho^{-1}}|\compactldots|G_1^{\rho \tau \rho^{-1}} \cap F_n^{\rho \sigma \rho^{-1}}]\]
and reindexing the sum, we get the desired equality.
\end{proof}

\subsubsection{Compatibility between the duality involution and the scaling action} Finally, we comment on the compatibility between the duality involution and the scaling action. This can be seen by combining both the scaling action and duality involution into a functor
\[\BGLb^+ \colon B(C_2 \ltimes BF^\times) \lra \Alg_{E_\infty^\rm{u}}(\Fun(\bb{N},\Spc))\]
where the domain is obtained from the $2$-category with unique object $\ast$, $1$-morphisms given by $\id$ and $\vee$ with $\vee^2 = \id$, and $2$-morphisms given by $F^\times$, where the ``whiskering'' of $\vee$ with $\lambda$ is given by $\lambda^{-1}$.

\begin{lemma}
Under the splittings of \eqref{eqn:splitting}
\[H_{*}(\GL_n,\StL_n) \cong H_*(\PGL_n,\StL_n) \otimes \Lambda^*F^\times\]
the $C_2$-action splits and acts by $(-1)^*$ on the second factor.
\end{lemma}

\subsection{Duality involution} \label{section: duality computation}
In this section, we compute the action of the duality automorphism $\vee_*$ on $\PolyL_n$ using our presentation:

\begin{theorem}\label{theorem: duality} The map $\vee_*$ acts by $(-1)^n$ on $\PolyL_n$.
\end{theorem}

By \cref{lem:duality-on-stl-sign}, it suffices to show that the map $D_*$ acts by $(-1)$ on $H_1(\GL_n;\StL_n)$.

\subsubsection{An explicit description of the map $D_*$}
Recall from \cref{def:steinberg-iterated-integral} that for an ordered basis $v_1,\dots,v_n$ the corresponding Steinberg iterated integral is given by
\[
\rm{I}[v_1,\compactldots,v_n]=(-1)^n [v_n,v_{n-1},\compactldots,v_1]\otimes [v_n,v_{n-1}-v_n,\compactldots,v_1-v_2]\in \StH(V).
\]
Denote by $v^1,\dots,v^n$ the dual basis, so that \eqref{eqn:sth-duality-d} implies that we have 
\begin{align*}
    D(\rm{I}[v_1,\compactldots,v_n])&=(-1)^n [v^1+\compactldots+v^n,v^1+\compactldots+v^{n-1},\compactldots,v^1]\otimes [v^n,v^{n-1},\compactldots,v^1]\\
    &=(-1)^n [v^1,\compactldots,v^1+\compactldots+v^{n-1}, v^1+\compactldots+v^n]\otimes [v^1,\compactldots,v^{n-1},v^n]\\
    &=\rm{I}[v^1+\compactldots+v^n,v^1+\compactldots+v^{n-1},\compactldots,v^1].
\end{align*}
The Steinberg correlator $\rm{C}[0:v_1:\compactcdots:v_n]$ is the projection of the element $(-1)^n\rm{I}[v_1,\compactldots,v_n]$ from $\StH(V)$ to $\StL(V)$, so
\begin{equation}\label{eqn: duality for correlators}
    D(\rm{C}[0:v_1:\compactcdots:v_n])= \rm{C}[0:v^1+\compactldots+v^n:v^1+\compactldots+v^{n-1}:\compactcdots:v_1]
\end{equation}
as in \cite[Proposition 32]{CMRR24}. Recall from \cref{sec:generators of G} that we have an exact sequence 
\[
0  \lra \FCR(V) \lra  \FC(V) \lra  \StL(V) \lra 0 
\]
and $\PolyL_n(F)$ is isomorphic to $H_0(\GL(V);\FCR(V))$. The map $D$ does not admit a natural extension to the projective module $\FC(V)$, so we will use a different resolution instead. 

Let $\rm{Bas}_V$ be the set of bases of $V$ and define
\begin{align*}\rm{Bas} \colon \Vect &\lra \rm{GrMod}_{\bb{Q}} \\
V &\longmapsto \Q[\rm{Bas}_V]\end{align*}
with the free $\bb{Q}$-vector space on ordered bases of $V$, a free $\GL(V)$-module concentrated in degree $0$. Concatenation of bases lifts $\rm{Bas} \in \Fun(\Vect,\rm{GrMod}_{\bb{Q}})$ to a $E_1$-algebra with respect to the tensor product $\levi$. 

There is a surjective morphism 
\begin{align*}
\rm{Bas}(V) &\lra \StL(V) \\
[(v_1,\dots,v_n)] &\longmapsto \rm{C}[0:v_1:\compactcdots:v_n],\end{align*}
and we denote its kernel by $\widetilde{\FCR}(V)$. Since $\rm{Bas}(V)$ is a free $\GL(V)$-module with coinvariants $\bb{Q}$, the exact sequence $0 \to \widetilde{\FCR}(V) \to \rm{Bas}(V) \to \StL(V) \to 0$ induces an exact sequence of group homology:
\begin{equation} \label{eqn: exact sequence for tilde FR}
    0 \lra \PolyL_n(F)\lra H_0(\GL(V);\widetilde{\FCR}(V)) \lra \Q \lra 0.
\end{equation}

Consider now the surjective map
\begin{align*} \rm{pr}\colon \rm{Bas}(V) &\lra \FC(V) \\
[(v_1,\compactldots,v_n)] &\longmapsto \FC[0:v_1:\compactcdots:v_n],\end{align*} 
fitting in a morphism of exact sequences
\[\begin{tikzcd}0\rar & \widetilde{\FCR}(V)\rar \dar & \rm{Bas}(V) \rar \dar{\rm{pr}} & \StL(V) \dar{\cong} \rar & 0\\[-5pt]
0 \rar & \FCR(V) \rar &  \FC(V) \rar & \StL(V) \rar & 0.
\end{tikzcd}\]
The induced map $H_0(\GL(V);\widetilde{\FCR}(V))\to H_0(\GL(V);\FCR(V))\cong \PolyL(V)$ on coinvariants gives a splitting of the exact sequence \eqref{eqn: exact sequence for tilde FR}.

The formula \eqref{eqn: duality for correlators} for the action of $D$ on Steinberg correlators implies that the map $D\colon \StL(V)\to  \StL(V)$ extends to a map of short exact sequences
\[\begin{tikzcd}0\rar &  \widetilde{\FCR}(V)\rar \dar{D} & \rm{Bas}(V)\rar \dar{D} & \StL(V) \dar{D} \rar & 0\\[-5pt]
0 \rar & \widetilde{\FCR}(V^{\vee}) \rar &  \rm{Bas}(V^{\vee})\rar & \StL(V) \rar & 0
\end{tikzcd}\]
where $D$ acts on $\rm{Bas}(V)$ by the formula
\[
D([(v_1,\compactldots,v_n)])=[(v^1+\compactldots+v^n,v^1+\compactldots+v^{n-1},\compactldots,v^1)].
\]
The map of short exact sequences induces the map 
\[\begin{tikzcd}0\rar &  \PolyL_n \rar \dar{D_*} & H_0(\GL(V);\rm{Bas}(V))\rar \dar{D_*} & \Q \dar{D_*} \rar & 0\\[-5pt]
0 \rar & \PolyL_n \rar &  H_0(\GL(V);\rm{Bas}(V^{\vee}))\rar & \Q \rar & 0
\end{tikzcd}\]
which we will use to compute the action of $D$ on $\PolyL_n$. Here we have used that since $\scr{G}_n$ is defined as the coinvariants $\FCR(V)_{\GL(V)}$, the isomorphism
\[\FCR(V)_{\GL(V)} \overset{\cong}\lra \FCR(F^n)_{\GL_n} = \scr{G}_n\]
induced by a linear isomorphism $V \cong F^n$ is in fact independent of the choice of this linear isomorphism.

\begin{example}To illustrate the discussion above, we take $V = F^2$ and compute the action of $D_*$ on $\CorG(0,x_1,x_2)$. We can assume that $x_1,x_2,x_2-x_1\neq 0$ by \cref{proposition: generic correlators}. The element $\CorG(0,x_1,x_2)$ equals the projection of the element
\[
x=\FC[0:e_1:e_2]-D^\FC_h(\FC[0:e_1:e_2])\in \FCR(V)
\]
with $h(e_1)=x_1$ and $h(e_2)=x_2$, where   
\begin{align*}
&D^\FC_h(\FC[0:e_1:e_2])=\FC\left[0:\frac{e_{1}}{x_1}:\frac{e_{2}}{x_2}\right]-\FC\left[0:\frac{e_{1}}{x_1}:\frac{e_{2}-e_{1}}{x_2-x_1}\right]+\FC\left[0:\frac{e_{2}}{x_2}:\frac{e_{2}-e_{1}}{x_2-x_1}\right].
\end{align*}
Consider an element
\[
\tilde{x}=[(e_1,e_2)]-
\left[\Bigl ( \frac{e_{1}}{x_1},\frac{e_{2}}{x_2}\Bigr)\right]+\left[\Bigl ( \frac{e_{1}}{x_1},\frac{e_{2}-e_{1}}{x_2-x_1}\Bigr)\right]
-\left[\Bigl ( \frac{e_{2}}{x_2},\frac{e_{2}-e_{1}}{x_2-x_1}\Bigr)\right]\in \rm{Bas}(V)
\]
which maps to $x$ by the projection $\pr \colon \rm{Bas}(V)\to \FC(V)$. We have
\begin{align*}
D(\tilde{x})&=[(e^1+e^2,e^2)]-
\left[(x_1e^{1}+x_2 e^{2},x_1 e^{1})\right]\\
&+\left[(x_1e^{1}+x_2 e^{2},x_1e^{1}+x_1e^2)\right]
-\left[(x_1e^{1}+x_2 e^{2}, x_2 e^{1}+x_2 e^{2})\right].
\end{align*}
The projection $\FC(V^{\vee})$ to $H_0(\GL(V^{\vee}),\FC(V^{\vee}))$ can be computed using \cref{prop: formula for projection to coinvariants}. Any vector $v\in V$ can be viewed as a functional on $V^{\vee}$; the composition 
\[
\rm{Bas}(V^\vee)\stackrel{\rm{pr}}{\lra} \FC(V^{\vee})\stackrel{E_{v}}\lra \PolyL_2
\]
sends $[(f_1,f_2)]$ to $\CorG(0, f_1(v),f_2(v))$. Using $v=e_1$, we obtain
\begin{align*}
D_*(\CorG(0,x_1,x_2))&=\CorG(0,1,0)-\CorG(0,x_1,x_1)+\CorG(0,x_1,x_1)-\CorG(0,x_1,x_2)\\ 
&=-\CorG(0,x_1,x_2).
\end{align*}
This shows that $D_*$ acts by $-1$ on $\PolyL_2$.\end{example}

\subsubsection{Proof of \cref{theorem: duality}} 
By \cref{proposition: generic correlators}, it is sufficient to prove that 
\[
D_*(\CorG(0,x_1,\compactldots,x_n))=-\CorG(0,x_1,\compactldots,x_n)
\]
for distinct arguments $x_i\in F^{\times}$. The element $\CorG(0,x_1,\compactldots,x_n)\in \PolyL_n$  is represented by the element
\[
x=\FC[0:e_1:\compactcdots:e_n]-D^\FC_h(\FC[0:e_1:\compactcdots:e_n])\in \FCR_n(F)
\]
for a basis $e_1,\dots,e_n$ of $V$ and a functional $h\in V^*$ such that $h(v_i)=x_i$. 

\smallskip

Consider a basis $v_1,\dots,v_n$ of $V$ and a functional $h$ such that $0, h(v_1),\dots,h(v_n)\in F$ are distinct and write $h_i \coloneq h(v_i)$. We then define inductively define elements
\[
f_h(v_1,\compactldots,v_n)\in \rm{Bas}(V)
\]
as follows: For $n=1$, we put $f_h(v_1)=[(\tfrac{v_1}{h_1})]$ and for $n\geq 2$, we define
\begin{equation*}\begin{aligned}
f_h(v_1,\compactldots,v_n) &= f_h(v_1,\compactldots,v_{n-1})\otimes \left[\Bigl(\frac{v_n}{h_n}\Bigr)\right] \\
&\quad +\sum_{i=1}^{n-1} \big(f_h(v_1,\compactldots,\widehat{v}_i,\compactldots,v_n)-f_h(v_1,\compactldots,\widehat{v}_{i+1},\compactldots,v_n)\big)\levi \left[\Bigl(\frac{v_{i+1}-v_i}{h_{i+1}-h_i}\Bigr)\right].
\end{aligned}\end{equation*} 

\begin{lemma}\label{lemma: lift of x} The element $f_h(v_1,\compactldots,v_n)\in \rm{Bas}(V)$ maps to $D^\FC_h[0:v_1:\compactcdots:v_n]$ under the projection $\rm{Bas}(V) \to \FC(V)$.
\end{lemma}

\begin{proof}
Comparing the formula for $f_h$ with the inductive definition \eqref{eqn:symbol-correlator} for the symbol of Steinberg correlator used in the formula for the decomposition operator.
\end{proof}

By \cref{lemma: lift of x}, the element 
\[
\tilde{x}=[(e_1,\compactldots,e_n)]-f_h(e_1,\compactldots,e_n)\in \rm{Bas}(V)
\]
gives a lift of $x\in \FC(V)$.  Our next goal is to compute the image of $\tilde{x}$ under the composition
\[
\rm{Bas}(V) \stackrel{D}{\lra} \rm{Bas}(V^{\vee})\stackrel{\rm{pr}}{\lra} \FC(V^{\vee}) \stackrel{E_u}{\lra} \PolyL_n(F),
\]
where $u$ is a nonzero vector in $V$ and $E_u$ is the map defined in \cref{section: projection from FC to G}. This composition can be computed in a different way: if we define maps
\begin{align*}\Psi_u\colon \rm{Bas}(V) &\lra \Q[F]^{\otimes n} \\
\Psi_u([(v_1,\compactldots,v_n)])&\longmapsto [(v^1+\compactldots+ v^n)(u)]\otimes [(v^1+\compactldots+ v^{n-1})(u)]\otimes \compactcdots \otimes [v^1(u)] \\
\CorG \colon  \Q[F]^{\otimes n} &\lra \PolyL_n \\
[x_1]\otimes  \compactcdots \otimes [x_n] &\longmapsto \CorG(0,x_1,\compactldots,x_n),\end{align*}
then we have an equation $E_u\circ \rm{pr} \circ D = \CorG \circ \Psi_u$, and thus
\[
D_*(\CorG(0,x_1,\compactldots,x_n))=E_u(\rm{pr}(D(\tilde{x})))= \CorG(\Psi_u(\tilde{x})).
\]
To evaluate the right side of this equation, we will use the following properties of $\Psi_u$:

\begin{lemma}\label{lemma: decomposition} Suppose we have a direct sum decomposition $V=V_1\oplus V_2$ and bases $v_1,\dots, v_{n_1}$  of $V_1$ and $v_{n_1+1},\dots, v_{n_1+n_2}$ of $V_2$.
For a vector $u\in V_1$ we have
\[
\Psi_u([(v_1,\dots,v_n)])=[(v^1+\dots+v^n)(u)]^{\otimes n_2} \otimes \Psi_u([(v_{1},\dots,v_{n_1})]).
\]
\end{lemma}

\begin{proof} 
Since $u\in V_1$, we have $v^{n_1+1}(u)=\dots=v^{n_1+n_2}(u)=0$ and so for any $m\geq n_1$ we have $(v^1+\dots+v^m)(u)=(v^1+\dots+v^n)(u)$, from which the statement follows.
\end{proof}

\begin{lemma} Suppose we have a basis $v_1,\dots,v_n$ of $V$ and a nonzero functional $h$ such that $0,h(v_1),\dots,h(v_n)$ are distinct. Then we have 
\[
\Psi_{v_1}(f_h(v_1,\compactldots,v_n))=[h(v_1)]\otimes \dots \otimes \bigl[h(v_n)\bigr].
\]
\end{lemma}
\begin{proof}
We prove the statement by induction on $n$, using the abbreviation $h_i \coloneq h(v_i)$. For $n=1$, we have 
$f_h(v_1)=\big[(\tfrac{v_1}{h_1})\big]$ and so 
\[
\Psi_{v_1}(f_h(v_1))=[h_1v^1(v_1)]=[h_1].
\]
Assume now that the statement holds in dimensions less than $n$. Using the inductive definition of $f$ one may verify that every term $\pm[w_1,\dots,w_n]$ in $f_h(v_1,\compactldots,v_n)$ we have $w^1+\dots+w^n=h$, so in particular we have $(w^1+\dots+w^n)(v_1)=h_1$. By \cref{lemma: decomposition} applied to $V_1=\rm{span}(v_1,\dots, v_{n-1})$ and  $V_2=\rm{span}(\tfrac{v_n}{h_n})$ and the induction hypothesis, we then see that
\begin{align*}
\Psi_{v_1}\Bigl(f_h(v_1,\compactldots,v_{n-1})\otimes \Bigl[\frac{v_n}{h_n}\Bigr] \Bigr)&=[h_1]\otimes \Psi_{v_1}(f_h(v_1,\compactldots,v_{n-1})) \\
&=[h_1]\otimes ([h_1] \otimes \dots \otimes [h_{n-1}]).\end{align*}
By the same  argument, for $2\leq i \leq n-1$ we have 
\[
\Psi_{v_1}\Bigl(f_h(v_1,\compactldots,\widehat{v}_{i},\compactldots,v_n)\otimes \Bigl[\Bigl(\frac{v_{i+1}-v_i}{h_{i+1}-h_i}\Bigr)\Bigr]\Bigr)=[h_1] \otimes \bigl([h_1]\otimes \dots \otimes [\widehat{h_{i}}] \otimes \dots \otimes [h_n]\bigr).
\]
Similarly, for $1\leq i \leq n-1$ we have
\[\Psi_{v_1}\Bigl(f_h(v_1,\compactldots,\widehat{v}_{i+1},\compactldots,v_n)\otimes \Bigl[\Bigl(\frac{v_{i+1}-v_i}{h_{i+1}-h_i}\Bigr)\Bigr]\Bigr)=[h_1] \otimes \bigl([h_1]\otimes \dots \otimes \widehat{[h_{i+1}]} \otimes \dots \otimes [h_n]\bigr).
\]

The remaining term $\Psi_{v_1}\big(f_h(v_2,\compactldots,v_n)\otimes \left[(\tfrac{v_{2}-v_1}{h_{2}-h_1})\right]\big)$ has a different structure, as $v_1$ does not lie in the span of the vectors $v_2,\dots,v_n$. For each term $[(w_1,\dots,w_{n-1})]$ of $f_h(v_2,\compactldots,v_n)$ consider the dual basis $u^1,\dots,u^{n}$ to the basis 
\[
w_1,\dots,w_{n-1}, \frac{v_2-v_1}{h_{2}-h_1}
\]
of $V$ (this is a basis because $w_1,\ldots,w_{n-1}$ is a basis of $\rm{span}(v_2,\ldots,v_n)$). Then $u^k(v_1)=u^k(v_2)$ for $k \leq n-1$ by choice $\tfrac{v_2-v_1}{h_{2}-h_1}$ of the last basis vector, and we see 
\[
(u^1+\dots+u^k)(v_1)=(u^1+\dots+u^k)(v_2) \qquad \text{for $k \leq n-1$.}
\]
Since $(u^1+\dots+u^n)(v_1)=h_1$ and functionals $u^1,\dots, u^{n-1}$ restrict to the dual basis to $w_1,\dots,w_{n-1}$, we obtain 
\[
\Psi_{v_1}\Bigl(f_h(v_2,\compactldots,v_n)\otimes \Bigl[\Bigl(\frac{v_{2}-v_1}{h_{2}-h_1}\Bigr)\Bigr]\Bigr)=[h_1] \otimes \Psi_{v_2}(f_h(v_2,\compactldots,v_n))=[h_1] \otimes \dots \otimes [h_n].
\]
To prove the result we now perform the following computation
\begin{equation*}\begin{aligned}
\Psi_{v_1}&\Bigl(f_h(v_1,\compactldots,v_n)\Bigr)\\
&=\Psi_{v_1}\Bigl ( f_h(v_1,\compactldots,v_{n-1})\otimes \Bigl[\Bigl(\frac{v_n}{h_n}\Bigr)\Bigr]\Bigr) +\Psi_{v_1}\Bigl (f_h(v_2,\compactldots,v_n)\otimes \Bigl[\Bigl(\frac{v_{2}-v_1}{h_{2}-h_1}\Bigr)\Bigr]\Bigr)\\
&\quad +\sum_{i=2}^{n-1} \Psi_{v_1}\Bigl(f_h(v_1,\compactldots,\widehat{v}_i,\compactldots,v_n)\otimes \Bigl[\Bigl(\frac{v_{i+1}-v_i}{h_{i+1}-h_i}\Bigr)\Bigr]\Bigr)\\
&\quad - \sum_{i=1}^{n-1}\Psi_{v_1}\Bigl(f_h(v_1,\compactldots,\widehat{v}_{i+1},\compactldots,v_n)\otimes \left[\Bigl(\frac{v_{i+1}-v_i}{h_{i+1}-h_i}\Bigr)\right]\Bigr)\\
&= [h_1]\otimes ( [h_1] \otimes \dots \otimes [h_{n-1}])+[h_1] \otimes [h_2]\otimes \dots \otimes [h_n]\\
&\quad +\sum_{i=2}^{n-1} [h_1] \otimes \bigl([h_1]\otimes \dots \otimes [\widehat{h_{i}}] \otimes \dots \otimes [h_n]\bigr)\\
&\quad - \sum_{i=1}^{n-1}[h_1] \otimes \bigl([h_1]\otimes \dots \otimes \widehat{[h_{i+1}]} \otimes \dots \otimes [h_n]\bigr)\\
&=[h_1] \otimes [h_2]\otimes \dots \otimes [h_n].
\end{aligned}\end{equation*} 
\end{proof}

We are ready to finish the proof of \cref{theorem: duality}. At this point we have shown that 
\[D_*(\CorG(0,x_1,\compactldots,x_n) = \CorG(\Psi_u([(e_1,\compactldots,e_n)-f_h(e_1,\compactldots,e_n)]))\] and the key observation is that the choice $u=e_1$ makes the computation particularly simple: we have 
\[
\CorG \Bigl(\Psi_{e_1}([(e_1,\compactldots,e_n)])\Bigr)=\CorG([(e^1+\dots+e^n)(e_1)]\otimes \dots \otimes [e^1(e_1)])=\CorG(0,1,\compactldots,1)=0
\]
and thus
\[
D_*(\CorG(0,x_1,\compactldots,x_n))=-\CorG \Bigl( \Psi_{e_1}(f_h(e_1,\compactldots,e_n)\Bigr)= -\CorG(0,x_1,\compactldots,x_n).
\]

\section{The Rognes rank spectral sequence relating $K(F)$ and $\scr{G}(F)$} In this section we describe a rank spectral sequence, which is inspired by the one obtained from Rognes' spectrum-level rank filtration \cite{Rognes} and is an instance of the group completion spectral sequence of Galatius, Kupers, and Randal-Williams \cite{GKRW18}. It plays a major role in the next section and in particular yields an edge homomorphism mapping the rationalised algebraic $K$-theory groups to the Goncharov Lie coalgebra. Here we develop only what is needed for those applications, postponing a more extensive discussion to \cite{KRS2}.

\subsection{The rank spectral sequence} The following is the main result of this section. We suggest a reader interested in applications skip its proof and move to the next section. 

\begin{theorem}\label{thm:rank-ss-omnibus} There exists a strongly convergent \emph{rank spectral sequence}
\[E^1_{n,d} \cong H^{E_\infty}_{n,d}(\BGLb(F)_\ds{Q}) \Longrightarrow \pi_d(K(F))_\ds{Q}\]
with $d^r$-differentials have bidegree $(-r,-1)$ and the following properties:
\begin{enumerate}
    \item The $d^1$-differential $d^1 \colon E^1_{n,2n-1} \to E^1_{n-1,2n-2}$ agrees with the $\sigma$-component $\delta_\sigma$.    \item It is a spectral sequence of $\Lambda^* F^\times_\bb{Q}$-modules, compatible with the actions on the $E^1$-page on abutment (see \cref{sec:rank-ss-action} for details).
    \item It is a spectral sequence with involution, compatible with the duality involutions on the $E^1$-page and abutment (see \cref{sec:rank-ss-involution} for details).
\end{enumerate}
\end{theorem}

Given this theorem, we may define:

\begin{definition}The maps 
\[\rm{edge}_n \colon K_{2n-1}(F)_\ds{Q} \lra \scr{G}_n(F)\]
are  defined to be the edge homomorphisms of this rank spectral sequence.\end{definition}

\subsection{Constructing the rank spectral sequence} The rank spectral sequence will be that associated to a filtered object in $\DQ$, constructed in several steps.

Suppose we are given an $E_\infty^\rm{u}$-algebra $\bf{R}^+$ in spaces with $\pi_0(\bf{R}^+) \cong \bb{N}$, the commutative monoid of (nonnegative) natural numbers under addition. Sending each path component to a point induces a map of $E_\infty^\rm{u}$-algebras $\bf{R}^+ \to \bf{N}$, where $\bf{N}$ is simply the commutative monoid $\bb{N}$ considered as an $E_\infty^\rm{u}$-algebra in spaces. Letting $t \colon \bb{N} \to \ast$ be the unique symmetric monoidal functor from the discrete symmetric monoidal category $\bb{N}$ to the terminal one, we get induced functors $t_! \colon \Fun(\bb{N},\Spc) \to \Spc$ and $t_!^\alg \colon \Alg_{E_\infty^\rm{u}}(\Fun(\bb{N},\Spc)) \to \Alg_{E_\infty^\rm{u}}(\Spc)$. The latter takes the terminal algebras $\ul{\ast} \in \Alg_{E_\infty^\rm{u}}(\Fun(\bb{N},\Spc))$ to $\bf{N} \in \Alg_{E_\infty^\rm{u}}(\Spc)$, and the induced functor
\[t_! \colon \Alg_{E_\infty^\rm{u}}(\Fun(\bb{N},\Spc)) \lra \Alg_{E_\infty^\rm{u}}(\Spc)_{/\bf{N}}\]
is an equivalence. Through this, we obtain the input in (I) below:

\begin{enumerate}[\noindent (I)]
\item \emph{Providing input.} Our input will be an $E_\infty^\rm{u}$-algebra $\bf{R}^+ \in \Fun(\bb{N},\Spc)$, which comes with a unique map $\epsilon \colon \bf{R}^+ \to \bf{N}$.
\item \emph{Rationalising.} Upon rationalisation and writing $\bb{Q}[t] \coloneq \bf{N}_\bb{Q} \simeq \free_{E_\infty^\rm{u}}(1_! \bb{Q})$, we obtain a map  of $E_\infty^\rm{u}$-algebras in $\Fun(\bb{N},\DQ)$
\[\epsilon \colon \bf{R}^+_\ds{Q} \lra \bb{Q}[t].\]
\item \emph{Rank filtering.} Let $\iota \colon \bb{N} \to \bb{N}_\leq$ by the inclusion of the discrete category of nonnegative natural numbers into the poset of nonnegative natural numbers with their usual order. This is symmetric monoidal so induces a symmetric monoidal functor $\iota_! \colon \Fun(\bb{N},\DQ) \to \Fun(\bb{N}_\leq,\DQ)$ which in turn induces a functor $\smash{\iota_!^\alg}$ on categories of $E_\infty^\rm{u}$-algebras. Applying this and writing $\smash{\bb{Q}^\fil[t] \coloneq \iota_!^\alg \bb{Q}[t] \simeq \free_{E_\infty^\rm{u}}(1_! \bb{Q})}$ (where now the free $E_\infty^\rm{u}$-algebra is taken in $\Fun(\bb{N}_\leq,\DQ)$), we get
\[\iota_!^\alg\epsilon \colon \iota_!^\alg \bf{R}^+_\ds{Q} \lra \bb{Q}^\fil[t].\]
\item \emph{Filtered augmentation.} There is a \emph{filtered group completion augmentation}
\[\epsilon^\fil_\gc \colon \bb{Q}^\fil[t] \lra 1_{\Fun(\bb{N}_\leq,\DQ)}\]
adjoint to the map $1_! \bb{Q} \to \fgt_{E_\infty^\rm{u}}(1_{\Fun(\bb{N}_{\leq},\DQ})) \simeq 0_! \bb{Q}$ that is in turn adjoint to $\id \colon \bb{Q} \to 1^* 0_! \bb{Q} = \bb{Q}$. We use this to consider $\smash{\iota_!^\alg \bb{R}^+_\bb{Q}}$ as an \emph{augmented} $E_\infty^\rm{u}$-algebra in $\Fun(\bb{N}_{\leq},\DQ)$.
\item \emph{Indecomposables}. We now pass to the augmentation ideal and apply the functor $\cot_{E_\infty^\rm{nu}}$ of $E_\infty^\rm{nu}$-indecomposables to get a filtered object that we will denote
\[\cot_{E_\infty^\rm{nu}}(\iota_!^\alg \bf{R}^\gc_\bb{Q}) \in \Fun(\bb{N}_\leq,\DQ).\]
\end{enumerate}

Since this filtered object is concentrated in nonnegative filtration degrees, we have an associated half-plane spectral sequence with exiting differentials and hence this spectral sequence is always strongly convergent \cite[Theorem 6.1]{BoardmanSS}. The rank spectral sequence of \cref{thm:rank-ss-omnibus} is the one associated to this filtered object in the case $\bf{R}^+ = \BGLb(F)^+$, which indeed has $\pi_0(\BGLb(F)^+) \cong \bb{N}$. It remains to verify its properties.

\subsection{The abutment, $E^1$-page, and bidegrees of differentials} This spectral sequence can be interpreted as an instance of one constructed by the Galatius, Kupers, and Randal-Williams in \cite{GKRW18}. To see this, note that since rationalisation is a symmetric monoidal left adjoint, it commutes with the constructions in steps (III)--(V). In particular, we could instead have passed to the rank filtration and taken a filtered group completion augmentation in spaces. At this point, we may compute the $E_\infty^\rm{u}$-indecomposables of the augmentation ideal by instead taking the iterated bar spectrum and taking the quotient by a copy of the sphere spectrum \cite[Lemma 13.26]{GKRW18}. Thus we see that the rank spectral sequence of \cref{thm:rank-ss-omnibus} agrees with the rationalisation of the \emph{group completion spectral sequence} of \cite[Remark 13.30]{GKRW18}, given by
\[E^1_{n,d} = H^{E_\infty}_{n,d}(\bf{R}) \Longrightarrow H_n(\bf{R}^\rm{sp}),\]
where $\bf{R}^\rm{sp}$ is the connective spectrum whose infinite loop space $\Omega^\infty \bf{R}^\rm{sp}$ is the group completion $(\bf{R}^+)^\gc$ of $\bf{R}^+$ and whose differentials in this indexing have bidegree $(-r,-1)$. 

In the case $\bf{R}^+ = \BGLb(F)^+$, this is one of the constructions of the $K$-theory spectrum $K(F)$, so we have essentially by definition that
\[(\BGLb(F)^+)^\gc \simeq \Omega^\infty K(F),\]
and the identification of the abutment follows by recalling that the Hurewicz map from rational homotopy to rational homology is an isomorphism for spectra.

\subsection{The $d^1$-differential} For the computation of the $d^1$-differential we will use Koszul duality. Recall that if the augmentation ideal $\bf{R}_\bb{Q} \in \Alg_{E_\infty^\rm{nu}}(\Fun(\bb{N},\DQ))$ is \emph{reduced}, i.e.~$\bf{R}_\bb{Q}(0) \simeq 0$, then Koszul duality yields an equivalence
\[\bf{R}_\bb{Q} \overset{\simeq}\lra \prim_{s\,\coLie} (\indec_{E_\infty^\rm{nu}}(\bf{R}_\bb{Q})).\]
We may thus assume that $\bf{R}_\bb{Q}$ is of the form $\prim_{s\,\coLie}(\bf{L})$ for some $\bf{L} \in \Alg_{s\,\coLie}^\red(\Fun(\bb{N},\DQ))$.
Here we suppress the superscript $(-)^\dpw$ from (shifted) Lie coalgebras because we are working rationally, and we replace the superscript $(-)^\nil$ by $(-)^\red$ because we are working with reduced (shifted) Lie coalgebras which are in particular conilpotent. In this situation $\indec_{E^\rm{nu}_\infty}$ is an adjoint inverse to $\prim_{s\,\coLie}$ so the counit induces an isomorphism
\[H^{E_\infty}_{n,d}(\bf{R}_\bb{Q}) \overset{\cong}\lra H_{n,d}(\bf{L}),\]
allowing us to rewrite the $E^1$-page of the rank spectral sequence in terms of $\bf{L}$. We will for simplicity assume that $\smash{H^{E_\infty}_{2,1}(\bf{R}_\bb{Q})} = 0$, as is the case for $\bf{R} = \BGLb(F)$ (this avoids a discussion of what it means to extract a $\sigma$-component out of the cobracket $H_{2,1}(\bf{L}) \to \Lambda^2 H_{1,0}(\bf{L})$).

Since our $E_\infty^\rm{u}$-algebra $\bf{R}^+_\bb{Q}$ arose by rationalising $\bf{R}^+ \in \Alg_{E_\infty^\rm{u}}(\Fun(\bb{N},\Spc))$ that is path-connected in each rank, picking a point in rank 1 gives the left map in
\[\free_{E_\infty^\rm{u}}(1_! \ast) \overset{\sigma}\lra \bf{R}^+ \overset{\epsilon}\lra \bb{N}.\]
Rationalising the composition becomes an equivalence and we get a factorisation
\[\bb{Q}[t] \overset{\sigma}\lra \bf{R}^+_\bb{Q} \overset{\epsilon}\lra \bb{Q}[t]\]
of $\id_{\bb{Q}[t]}$. Moreover, as $\smash{\bf{R}^+_\bb{Q}}$ is connective and $\sigma$ as well as $\epsilon$ induce an isomorphism on $H_{*,0}(-)$, taking indecomposables we obtain a factorisation $\bb{Q}\{\sigma\} \to \bf{L} \to \bb{Q}\{\sigma\}$ inducing an identification $\bb{Q}\{\sigma\} \cong H_{1,0}(\bf{L})$. Thus in particular the (shifted) cobracket on $H_{*,*}(\bf{L})$ has a \emph{$\sigma$-component}
\[\delta_\sigma \colon H_{n,d}(\bf{L}) \lra H_{n-1,d-1}(\bf{L})\]
obtained from (shifted) cobracket $\delta$ by projection to the term involving $H_{1,0}(\bf{L})$.

\begin{lemma}Let $\bf{L}$ be as above. Then the $d^1$-differential
\[d^1 \colon E^1_{n,d} \lra E^1_{n-1,d-1}\]
of the rank spectral sequence for $\bf{R}^+_\bb{Q} = \prim_{s\,\coLie}(\bf{L})^+$ is given by $\delta_\sigma$ under the identification $E^1_{n,d} \cong H_{n,d}(\bf{L})$.
\end{lemma}

\begin{proof}We will permit ourselves to use the rectification results from \cref{sec:rect-dg} and use explicit chain complexes for the constructions: that is, we work in $\Fun(\bb{N},\rm{Ch}_\bb{Q})$, use the equivalence $E_\infty^\rm{u} \simeq_\bb{Q} \rm{Com}$, model $\prim_{s\,\coLie}$ by the shifted variant $\Omega^{s\,\coLie} \coloneq \Sigma^{-1} \circ \Omega^\rm{coLie} \circ \Sigma$ of the cobar construction from \cref{def:cobar-colie}, and model $\indec_{E_\infty^\rm{u}}$ by the bar construction $B^\rm{Com}$ from \cref{def:bar-comm}.

\medskip

The input to the rank spectral sequence may thus be taken to be a dg-commutative algebra of the form $\Omega^\rm{scoLie}(\bf{L}) \in \Alg_{\rm{Com}^\rm{nu}}(\Fun(\bb{N},\rm{Ch}_\bb{Q}))$. We may further assume, e.g.~using a variant of CW approximation as in \cite[Section 11]{GKRW18}, that $\bf{L} \in \coAlg_{s\,\coLie}(\Fun(\bb{N},\rm{Ch}_\bb{Q}))$ splits additively (though not necessarily as shifted dg-Lie coalgebras) as $\bb{Q}\{\sigma\} \oplus \bf{L}'$ where $|\sigma| = (1,0)$ and $\bf{L}'$ is rankwise connected. The unitalisation of shifted variant of the cobar construction is then given by
\[\Omega^{s\coLie}(\bf{L})^+ \coloneq (\rm{Com} \circ \bf{L},d_\bf{L}+d_\Omega) \cong (S^*(\bb{Q}\{\sigma\}\oplus \bf{L}'),d_\bf{L}+d_\Omega)\]
where $d_\bf{L}$ denotes the internal differential of $\bf{L}$ and $d_\Omega$ denotes the cobar differential from \cref{def:cobar-colie}. Applying $\smash{\iota_!^\alg}$ makes this graded dg-commutative algebra filtered by declaring that elements in rank $n$ lie in filtration degree $n$; for brevity we will replace $\smash{\iota^\alg_!}$ by an underline in our notation, so that 
\[\iota^\alg_! \Omega^{s\,\coLie}(\bf{L})^+ \cong (S^*(\bb{Q}\{\ul{\sigma}\} \oplus \ul{\bf{L}}'),d_\bf{L}+d_\Omega).\] 
The filtered group completion augmentation $\epsilon^\fil_\gc$ on this is then determined uniquely by sending $\ul{\sigma}$ to $1$ and $\ul{\bf{L}}'$ to $0$. 

We next pass to the augmentation ideal $I^\rm{gc}$ and apply the bar construction $B^\rm{Com}$ to this nonunital commutative algebra to obtain an explicit model for the filtered object that gives rise to the rank spectral sequence 
\[B^\rm{Com}I^\gc(\Omega^{s\coLie}(\bb{Q}\{\ul{\sigma}\} \oplus \ul{\bf{L}}')) \coloneq \Sigma^{-1}\Big(\coLie \circ \Sigma I^\gc(S^*(\bb{Q}\{\ul{\sigma}\} \oplus \ul{\bf{L}}')),d_\bf{L}+d_\Omega+d_B\Big)\]
where $d_B$ denotes the bar differential from \cref{def:bar-comm}.

\medskip 

The $d^1$-differential is extracted as follows: writing $F_r \coloneq F_r B^\rm{Com}I^\gc(\Omega^{s\,\coLie}(\bb{Q}\{\ul{\sigma}\} \oplus \ul{\bf{L}}'))$ for the $r$th filtration step, it is given by the connecting homomorphism in the short exact sequence
\[0 \lra F_{r-1}/F_{r-2} \lra F_r/F_{r-2} \lra F_r/F_{r-1} \lra 0.\]
We will now make this more explicit. When we pass to the associated graded, identifying the image of $\ul{\sigma-1}$ with $\sigma$ and the image of $\ul{\bf{L}}'$ with $\bf{L}'$, the filtered group completion augmentation yields the canonical augmentation, which cancels against unitalisation, and Koszul duality provides an equivalence
\[\eta \colon \bb{Q}\{\sigma\} \oplus \bf{L}' \overset{\simeq}\lra B^\rm{Com} \Omega^{s\,\coLie}(\bb{Q}\{\sigma\} \oplus \bf{L}') \cong \gr\,B^\rm{Com}I^\gc(\Omega^{s\,\coLie}(\bb{Q}\{\ul{\sigma}\} \oplus \ul{\bf{L}}')).\]
Let us first describe $\eta$ and its inverse up to homotopy $\varpi$. The latter is easier: in terms of
\[B^\rm{Com}\Omega^{s\,\coLie}(\bb{Q}\{\sigma\} \oplus \bf{L}') \coloneq \Sigma^{-1}\Big(\coLie \circ \Sigma (S^{*>0}(\bb{Q}\{\sigma\} \oplus \bf{L}')),d_\bf{L}+d_\Omega+d_B\Big)\]
the map $\varpi$ is the projection onto the generators $\Sigma^{-1} \rm{coLie}(1) \otimes_{\fr{S}_1} \Sigma (\bb{Q}\{\sigma\} \oplus \bf{L}') \cong \bb{Q}\{\sigma\} \oplus \bf{L}'$. The map $\eta$ is given by the formula in \cite[Lemma 2.21]{Souderes}, up to an unfortunate difference in sign conventions and normalisation of the cobracket: before the outer desuspension $\Sigma^{-1}$, $x \in \bb{Q}\{\sigma\} \oplus \bf{L}'$ is mapped to an element of the form
\[(\delta_1(x),\delta_2(x),\ldots) \qquad \text{with} \qquad \delta_k(x) \in  \rm{coLie}(k) \otimes_{\fr{S}_k} (\Sigma(\bb{Q}\{\sigma\} \oplus L'))^{\otimes k},\]
where we set $\delta_k = p_{\rm{I\!I\!I}} \circ \widetilde{\delta}_k$, using that $p_{\rm{I\!I\!I}}$ as on \cite[p.~105]{Souderes} (called $\gamma$ in \cite{HainIndec}) serves to project $ (\Sigma(\bb{Q}\{\sigma\} \oplus \bf{L}'))^{\otimes k}$ into $\rm{coLie}(k) \otimes_{\fr{S}_k} (\Sigma(\bb{Q}\{\sigma\} \oplus \bf{L}'))^{\otimes k}$, and think of the latter as a subspace of $\rm{coLie}(k) \otimes_{\fr{S}_k} (\Sigma(S^{*>0} (\bb{Q}\{\sigma\} \oplus \bf{L}')))^{\otimes k}$. Explicitly, we have that
\[\delta_1(x) = \ol{x} \quad \text{and} \quad \delta_2(x) = -\ol{\delta(x)}\]
where $\ol{a} \coloneq \Sigma a$, and the signs and scalars differ from \cite{Souderes} due to the aforementioned difference in conventions.

We next explain a lift from the associated graded to the filtered object: $I^\gc(S^*(\bb{Q}\{\ul{\sigma}\} \oplus \ul{\bf{L}}'))$ is given by the nonunital subalgebra $S^{*>0}(\bb{Q}\{\ul{\sigma-1}\} \oplus \ul{\bf{L}}')$ where $\ul{\sigma-1}$ lies in filtration degree $1$. The projection to the associated graded then admits a section induced by mapping $\sigma$ to $\ul{\sigma-1}$ and $x \in \bf{L}'$ to $\ul{x} \in \ul{\bf{L}}'$.

\medskip

Given this, we may compute the $d^1$-differential on a class $[x] \in H_*(F_r/F_{r-1})$ for $r \geq 3$. Using the additive splitting we can write 
\[\delta_k(x) = \sum_{0 \leq j \leq k} c_j \otimes_{\fr{S}_k} \left(\ol{\sigma}^{\otimes j} \otimes \delta^{(j)}_k(x)\right) \in \rm{coLie}(k) \otimes_{\fr{S}_k} (\Sigma(S^{*>0} (\bb{Q}\{\sigma\} \oplus \bf{L}')))^{\otimes k}\]
for $c_j \in \rm{coLie}(k)$ and $\delta^{(j)}_j(x) \in (\Sigma(S^{*>0} (\bf{L}')))^{\otimes k-j}$, and following the prescription for lifting $\eta(x)$ to $F_r$, each of the terms $\delta_k(x)$ is lifted to
\[\ul{\delta_k(x)} \coloneq \sum_{0 \leq j \leq k} c_j \otimes_{\fr{S}_k} \left( \ul{\ol{\sigma}}^{\otimes j} \otimes \ul{\delta^{(j)}_k(x)}\right)-c_j \otimes_{\fr{S}_j} \left(\ul{\ol{1}}^{\otimes j} \otimes \ul{\delta^{(j)}_k(x)}\right).\]
Taking the quotient by $F_{r-2}$, only the terms $j=0,1$ remain. In particular, the first two terms of the lift to $F_r/F_{r-2}$ are given by
\[(\ul{\ol{x}},- \ul{\ol{\delta(x)}}+\ul{\ol{1}} \wedge \ul{\ol{\delta_\sigma(x)}},\ldots),\] 
where $\delta_\sigma(x)$ is the $\sigma$-component, and we remind the reader that the underline denotes that we lift elements in $\bf{L}'$ in rank $r$ to filtration degree $r$ and the overline denotes a suspension that is part of the bar construction. Taking the differential, there are three contributions: 
\begin{enumerate}
\item the internal differential $d_\bf{L}$ vanishes since $x$ was a cycle, 
\item the cobar differential $d_\Omega$ maps into higher tensor powers of $S^*$ and since we will project these away momentarily when we apply $\varpi$ we may ignore them, and 
\item the bar differential $d_B$ maps terms with index $k$ into terms with index $k-1$ and since we will project those with $k \geq 2$ away momentarily when we apply $\varpi$ we may ignore all but from $k=1,2$.
\end{enumerate}
Keeping this in mind and using that $\eta(x)$ is a cycle so all terms not involving $\ul{\ol{1}}$ cancel, we compute that
\[\varpi(d(\text{lift of $\eta(x)$})) = \delta_\sigma(x) \in \bb{Q}\{\sigma\} \oplus \bf{L}'\]
proving the result.
\end{proof}

\begin{remark}
This should be compared to the description in \cite{Sun} of the $d^1$-differential in the Quillen rank spectral sequence from \cite{QuillenFiniteGeneration} in terms of the coproduct on $H_*(\GL_n;\St_n)$. We believe this formula can be computed by similar methods.
\end{remark}

\subsection{Actions and splittings} \label{sec:rank-ss-action} Let us recall from \cref{sec:action-by-scaling} that the units in $F$ act ``by scaling'' on the symmetric monoidal groupoid $\Vect$ of finite-dimensional vector spaces over $F$: for $\lambda \in F^\times$ this action is by a symmetric monoidal natural isomorphism $\varphi_\lambda$ of $\id_\Vect$, with components $\lambda \cdot \id_V \colon V \to V$. These cover the identity on $\bb{N}$, so taking classifying spaces we obtain a lift of $\BGLb(F)^+$ to a functor
\[\BGLb(F)^+ \colon \rm{B^2}F^\times \lra \Alg_{E_\infty^\rm{nu}}(\Fun(\bb{N},\Spc)).\]
The naturality of the construction of the rank spectral sequence yields a functor 
\[\cot_{E_\infty^\rm{nu}}(\iota_!^\alg \bf{R}^\gc_\bb{Q}) \colon \rm{B}^2F^\times \to \Fun(\bb{N}_{\le},\DQ)\]
with adjoint action map in $\Fun(\bb{N}_{\leq},\DQ)$ given by (recall $0_!$ denotes we place an object in filtration $0$)
\[0_! \rm{B} F_\bb{Q} \otimes \cot_{E_\infty^\rm{nu}}(\iota_!^\alg \bf{R}^\gc_\bb{Q}) \lra \cot_{E_\infty^\rm{nu}}(\iota_!^\alg \bf{R}^\gc_\bb{Q}).\]
This map of filtered objects induces an $\Lambda^* F^\times$-action on the rank spectral sequence. 

\subsubsection{Action on abutment} On the abutment this is induced by the action of $\rm{B}F^\times$ on the spectrum $K(F)$ arising from the functoriality of group completion, with adjoint action map
\[\Sigma^\infty_+ \rm{B}F^\times \otimes K(F) \lra K(F)\]
in the category of spectra. Interpreting the action as tensoring with the $1$-dimensional vector space $F$ and its automorphisms, combining \cite[Theorem IV.1.10, Theorem IV.4.6, Corollary IV.4.6.1]{Weibel} we see that this is a part of the ring spectrum structure on $K(F)$ induced by tensor products of vector spaces. In particular, restricting to rational homotopy groups this is the map $\Lambda^* F^\times_\bb{Q} \otimes K_*(F)_\bb{Q} \to K_*(F)_\bb{Q}$ induced by iterated products with elements of $K_1(F)_\bb{Q} \cong F^\times_\bb{Q}$, and thus factors over the Milnor $K$-theory in the first entry to yield the multiplication map
\[K^M_*(F)_\bb{Q} \otimes K_*(F)_\bb{Q} \lra K_*(F)_\bb{Q}.\]

\subsubsection{Action on $E^1$-page}

On the $E^1$-page this is by construction the action induced on the $E_\infty$-homology of $\BGLb(F)^+$. As this action is trivial on the level of infinite Steinberg modules since the $\GL_n$-action on $\StL_n$ factors over $\PGL_n$, with respect to the splitting 
\[H^{E_\infty}_{n,d}(\BGLb(F)_\bb{Q}) \cong H_{d-2n+2}(\GL_n;\StL_n) \underset{\cong}{\overset{\varpi_n}{\lra}} H_{d-2n+2}(\PGL_n;\StL_n) \otimes \Lambda^* F^\times,\]
it simply acts on the second term, by $x \star (a \otimes y) = n^{|x|} a \otimes (x \cdot y)$. 

Since the $d^1$-differential is compatible with the $\Lambda^* F^\times$-action, this gives that up to a nonzero scalar the $d^1$-differential is also compatible with the splittings of the $E^1$-page.

\subsection{The duality involution}\label{sec:rank-ss-involution}

\subsubsection{The duality involution on rational algebraic $K$-theory} The group completion functor
\[\Alg_{E_\infty^\rm{u}}(\Fun(\N,\Spc)) \lra \Sp_{\ge 0}\]
applied to this $\BGLb(F)^+$ endows the spectrum $K(F)$ with a $C_2$-action. By construction this is the standard duality involution, which agrees on homotopy with the Adams operation $\psi^{-1}$ \cite[Section 2.3.1]{FengGalatiusVenkatesh}. It includes a direct sum decomposition into $\pm 1$-eigenspaces
\[K_*(F)_\ds{Q} \cong K^+_*(F)_\ds{Q} \oplus K^-_*(F)_\ds{Q}.\]

To understand its effect rational algebraic $K$-theory, recall that for any $d \ge 1$ there are natural splittings \cite[Remark 5.10.1, Theorem IV.5.11]{Weibel}
\[K_d(F)_\bb{Q} \cong \bigoplus_{i=1}^d K_d^{(i)}(F)_\bb{Q}\]
where $\smash{K_d^{(i)}}(F)_\bb{Q}$ is the \emph{weight $i$ part}, determined by the property that the $k$th Adams operation $\psi^k$ acts by multiplication by $k^i$ for any $k \in \bb{Z}$. It is known that in degree $d$ the the weight $d$ part agrees with $K^M_d(F)_\bb{Q}$, and the weight $1$ part vanishes as long as $d \ge 2$ \cite[Corollaire 1]{Soule}. We conclude that:

\begin{lemma} \label{lem:involution-vs-weight}
    The splitting $K_d(F)_\bb{Q}= K_d^+(F)_\bb{Q} \oplus K_d^-(F)_\bb{Q}$ induced by the duality involution is given by
    \[K_d^+(F)_\bb{Q}= \bigoplus_{i \; \text{even}} K^{(i)}_d(F)_\bb{Q} \qquad \text{and} \qquad K_d^-(F)_\bb{Q}= \bigoplus_{i \; \text{odd}} K^{(i)}_d(F)_\bb{Q}.\]
\end{lemma}

\subsubsection{The duality involution on the $E^1$-page}
On the $E^1$-page, the duality involution is given by that $\cot_{E_\infty}(\BGLb(F))$ and we investigated this in \cref{section duality involution}. There we found the following: it acts on the entries $E^1_{n,2n-1} \cong \scr{G}_n(F)$ by $(-1)^n$ by \cref{theorem: duality}, and on the entries $\smash{E^1_{1,d}} \cong \smash{\Lambda^d F^\times}$ for $d \ge 1$ by $(-1)^d$. We add this to the following result, which helps us understand the action on the second column of the spectral sequence:

\begin{lemma}
    The duality involution acts by $+1$ on $H_*(\PGL_2(F),\StL_2(F))$. 
\end{lemma}

\begin{proof}This is a consequence of the duality involution being inner. On the one hand, this involution is given by
\[\begin{bmatrix}
        a & b \\ c & d
    \end{bmatrix} \longmapsto \frac{1}{ad-bc}\begin{bmatrix} d & -c \\ -b & a
    \end{bmatrix}.\]
On the other hand, we compute
    \[\begin{bmatrix}
        0 & 1 \\ -1 & 0
    \end{bmatrix}^{-1} \begin{bmatrix}
        a & b \\ c & d
    \end{bmatrix} \begin{bmatrix}
        0 & 1 \\ -1 & 0
    \end{bmatrix} = \begin{bmatrix} d & -c \\ -b & a
    \end{bmatrix},\]
    which means that on $\PGL_2(F)$ involution acts by conjugation with the matrix 
    $\begin{bsmallmatrix}
        0 & 1 \\ -1 & 0
    \end{bsmallmatrix}$.\end{proof}

\section{The Goncharov conjectures in weight 3}
In this section we prove \cref{theorem weight 3}. Our main tool will be the Rognes rank spectral sequence of \cref{thm:rank-ss-omnibus} and the information about it obtained there; this approach to rational algebraic $K$-theory of fields was suggested by Rognes in \cite{RognesMotivic}. 

\subsection{The $E^1$-page} 
\begin{figure}[ht]
	\begin{tikzpicture}[xscale=1.3]
	\begin{scope}
	\clip (-1,-1) rectangle ({2.5*3+1},6.5);
	\draw (0,0)--(10.5,0);
	\draw (0,0) -- (0,6.5);
	\foreach \s in {0,...,6}
	{
		\draw [dotted] (-0.5,\s)--(10.5,\s);
		\node [fill=white] at (-0.5,\s) [left] {\tiny $\s$};
	}
	\foreach \s in {0,...,4}
	{
		\draw [dotted] ({2.5*\s},-0.5)--({2.5*\s},6.5);
		\node [fill=white] at ({2.5*\s},-.5) {\tiny $\s$};
	}
	
	\node [fill=white] at (2.5,0) {\textcolor{Periwinkle!50!black}{$\ds{Q}$}};
	\node [fill=white] at (2.5,1) {\mahoganydashed{$F^\times$}};
	\node [fill=white] at (2.5,2) {\textcolor{Periwinkle!50!black}{$\Lambda^2 F^\times$}};
	\node [fill=white] at (2.5,3) {\mahoganydashed{$\Lambda^3 F^\times$}};
	\node [fill=white] at (2.5,4) {\textcolor{Periwinkle!50!black}{$\Lambda^4 F^\times$}};
	\node [fill=white] at (2.5,5) {\mahoganydashed{$\Lambda^5 F^\times$}};
	\node [fill=white] at (2.5,6) {\textcolor{Periwinkle!50!black}{$\Lambda^6 F^\times$}};

	\node [fill=white] at (5,3) {\textcolor{Periwinkle!50!black}{$B_2(F)$}};
	\node [fill=white] at (5,4) {\textcolor{Periwinkle!50!black}{$(E^1_{2,4})^+$} $\oplus$ \mahoganydashed{$B_2(F) \otimes F^\times$}};
	\node [fill=white] at (5,5) {\textcolor{Periwinkle!50!black}{$(E^1_{2,5})^+$} $\oplus$ \mahoganydashed{$(E^1_{2,4})^+ \otimes F^\times$}};
	\node [fill=white] at (5,6) {?};

	\node [fill=white] at (7.5,5) {\mahoganydashed{$B_3(F)$}};
	\node [fill=white] at (7.5,6) {?};
	
	\node at (-.5,-.5) {$\nicefrac{q}{p}$};
	\end{scope}
	\end{tikzpicture}
	\caption{The $E^1$-page $E^1_{p,q}$ of the Rognes rank spectral sequence, converging to $K_q(F)$. The colours and dashed lines denote the \textcolor{Periwinkle!50!black}{$+1$-eigenspaces} and \mahoganydashed{-1-eigenspaces}. The $d^1$-differential has bidegree $(-1,-1)$ and necessarily sends $\pm 1$-eigenspaces to $\pm 1$-eigenspaces.}
	\label{fig:e1page}
\end{figure}
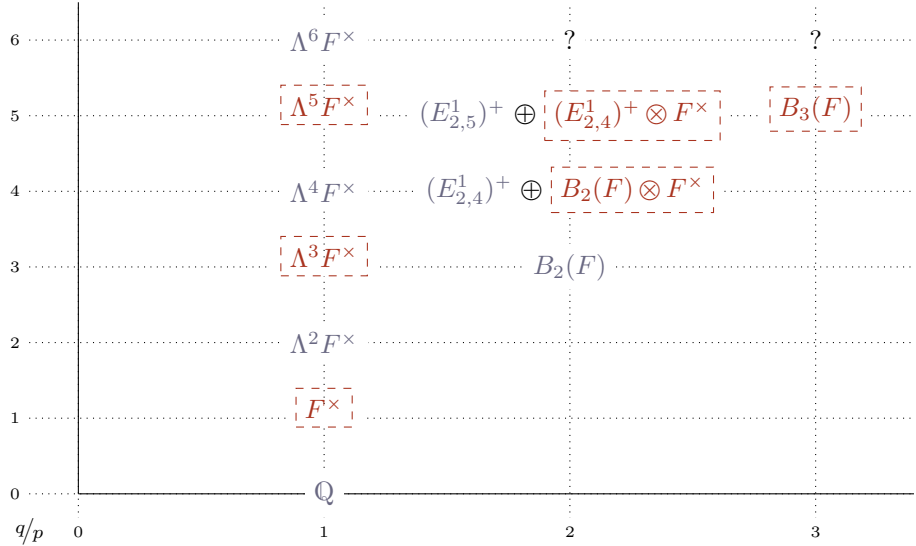

We start with a discussion of the $E^1$-page of the Rognes rank spectral sequence. \cref{fig:e1page} records the following data about its $E^1$-page, obtained in the previous section:
\begin{enumerate}[(i)]
    \item The first column is given by 
    \[E^1_{1,q} \cong \Lambda^q F^\times,\]
    and the involution acts on this by $(-1)^q$.
    \item The critical line is given by
    \[E^1_{n,2n-1} \cong \scr{G}_n(F),\]
    and the involution acts on this by $(-1)^n$. Moreover, we have identifications 
    \[F^\times \overset{\cong}\lra \scr{G}_1(F), \quad B_2(F) \overset{\cong}\lra \scr{G}_2(F), \quad \text{and} \quad B_3(F) \overset{\cong}\lra \scr{G}_3(F),\]
    with subgroups generated by ``classical'' polylogarithms.
    \item The second column is given by 
    \[E^1_{2,*} \cong H_{*+2}(\PGL_2(F);\StL_2(F)) \otimes \Lambda^* F^\times\]
    and the involution on the term $H_{*+2}(\PGL_2(F);\StL_2(F))$ is by $+1$ and on the term $\Lambda^* F^\times$ is by $(-1)^*$. 
\end{enumerate}
In particular, in order to compute rational algebraic $K$-theory of $F$ up to degree $5$, we only need to understand the $d^1$-differentials and potentially a $d^2$-differential from bidegree $(3,6)$ to $(1,5)$. 

\subsection{Computing the $d^1$ differentials in terms of the cobracket} We shall explain how to obtain most of the $d^1$-differentials. Recall that the $d^1$-differential is given by $\sigma$-component, which for $n=2$ agrees up to a sign with the cobracket and for $n \geq 3$ has the property that its projection of $\sigma$-component onto the first summand
    \[\pr_1 \delta_\sigma \colon \scr{G}_n(F) \lra \scr{G}_{n-1}(F) \otimes F^\times \subseteq \scr{G}_{n-1}(F) \otimes F^\times \oplus H_2(\PGL_{n-1}(F),\StL_{n-1}(F))\]
agrees up to a sign with the component $\delta_{n-1,1}$ of the cobracket. We have also seen that the $d^1$-differential is, up to a nonzero scalar, compatible with the splitting.

\begin{corollary} \label{cor computation of d^1}\,
    \begin{enumerate}[(i)]
        \item $d^1 \colon B_2(F) \to \Lambda^2 F^\times_\bb{Q}$ is given by $\LiG_2(x) \mapsto -(x) \wedge (1-x)$.
        \item $d^1 \colon B_2(F) \otimes F^\times_\bb{Q} \to \Lambda^3 F^\times_\bb{Q}$ is given by $\LiG_2(x) \otimes (y) \mapsto -(x) \wedge (1-x) \wedge (y)$.
        \item $d^1 \colon B_3(F) \to B_2(F) \otimes F^\times_\bb{Q}$ is given by $\LiG_3(x) \mapsto -\LiG_2(x) \otimes (x)$.
    \end{enumerate}
\end{corollary}

\begin{remark}Let us comment on the notation: when we write $-(x) \wedge (1-x)$ the minus sign thinks of $\Lambda^2 F^\times_\bb{Q}$ as being additive, so this is the additive inverse of the element $(x) \wedge (1-x)$ and is \emph{not} equal to $(-x) \wedge (1-x)$ (which in fact is equal to $(x) \wedge (1-x)$ since we work rationally).\end{remark}

As mentioned in the introduction, the complexes
\[B_2(F) \lra \Lambda^2 F^\times \quad \text{and} \quad B_3(F) \lra B_2(F) \otimes F^\times \lra \Lambda^3 F^\times\]
that appear on the $E^1$-page are isomorphic to the polylogarithmic complexes $\Gamma_2(F)$ and $\Gamma_3(F)$ studied previously by Goncharov \cite{Gon95b} (by convention, we drop $\bb{Q}$ from the notation). By convention, $B_n(F)$ is in degree $1$ and the differential increases degree.

\begin{remark}The relationship between polylogarithmic complexes of Goncharov and the Rognes rank spectral sequence is not entirely clear. The latter should be the $E^1$-page of an analogue of a Lyndon--Hochschild--Serre spectral sequence for the Lie coalgebra homology of cofibre sequences of Lie coalgebras, applied to $F^\times \to \scr{G}_n(F) \to \ol{\scr{G}}_n(F)$, together with the conjecture that $\ol{\scr{G}}_n(F)$ is cofree with cogenerators $B_n(F)$. The first term in which these complexes differ from the Lie coalgebra homology of $\scr{G}_n(F)$ is in weight $4$, in whether a term $\Lambda^2 \scr{G}_2(F)$ appears or not. Based on computations for number fields, it seems more plausible to the authors that one should relate the total differential in the Rognes rank spectral sequence (rather than just the $d^1$-differential) to the cobracket on $\scr{G}(F)$.
\end{remark}

\subsection{The first column and Milnor $K$-theory}

We now study the first column through the maps
\begin{equation}\label{eqn:first-column}\Lambda^d F^\times \cong E^1_{1,d} \lra E^2_{1,d} \overset{\pr}\lra E^\infty_{1,d} \lra K_d(F).\end{equation}
Here the left and middle map are the surjections arising from the spectral sequence, and the right map is an edge homomorphism. A similar result appears as \cite[Theorem 7.2.1]{RognesMotivic}.

\begin{lemma} \label{lem:first-column} \, \begin{enumerate}[(i)]
    \item \label{enum:first-column-i} The composition \eqref{eqn:first-column} is given by $\Lambda^* F^\times \twoheadrightarrow K_*^M(F) \hookrightarrow K_*(F)$.
    \item \label{enum:first-column-ii} The projection $\pr \colon E^2_{1,d} \to E^\infty_{1,d}$ in \eqref{eqn:first-column} is an isomorphism, or equivalently, there are no nonzero $d^r$-differentials for $r>1$ into the first column.
\end{enumerate} 
\end{lemma}

\begin{proof}By compatibility of the rank spectral sequence with the scaling action, the maps in \eqref{eqn:first-column} assemble to a map of $\Lambda^* F$-modules. For degree reasons, the maps $E^1_{1,0} \to K_0(F)$ is an isomorphism. Since $K_0(F)$ generates under the $\Lambda^* F$-action  the Milnor $K$-theory $K^M_*(F) \subseteq K_*(F)$, the composition \eqref{eqn:first-column} must be the quotient map onto the summand $K^M_d(F)_\bb{Q}$, proving the part \eqref{enum:first-column-i}.

For part \eqref{enum:first-column-ii}, we recall that the $d^1$-differential
\[E^1_{2,3} \cong B_2(F) \lra E^1_{1,2} \cong \Lambda^2 F^\times\]
is given by $\LiG_2(x) \mapsto -(x) \wedge (1-x)$, and since the $d^1$-differential is compatible with the $\Lambda^* F^\times$-action up to nonzero scalars, we see that $d^1(E^1_{2,d+3}) \subseteq E^1_{1,d+2}$ contains the span of $(x) \wedge (1-x) \wedge (y_1) \wedge \cdots \wedge (y_d)$ so that there is a factorisation $E^1_{1,d} \twoheadrightarrow K^M_d(F) \twoheadrightarrow E^2_{1,d}$. Using part \eqref{enum:first-column-i} the second map must be isomorphism, implying the part \eqref{enum:first-column-ii}.
\end{proof}

\subsection{$K$-theory groups in low degrees} 
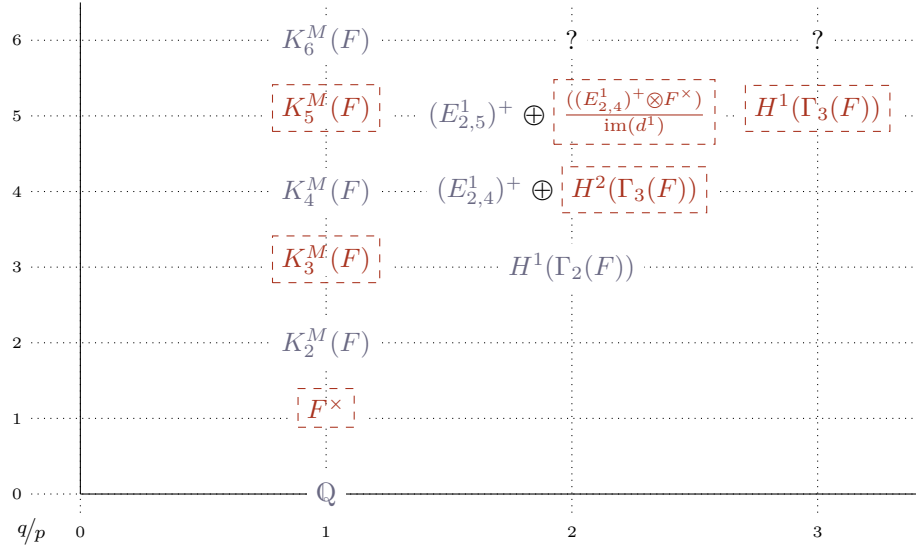
\begin{figure}[ht]
	\begin{tikzpicture}[xscale=1.3]
	\begin{scope}
	\clip (-1,-1) rectangle ({2.5*3+1},6.5);
	\draw (0,0)--(10.5,0);
	\draw (0,0) -- (0,6.5);
	\foreach \s in {0,...,6}
	{
		\draw [dotted] (-0.5,\s)--(10.5,\s);
		\node [fill=white] at (-0.5,\s) [left] {\tiny $\s$};
	}
	\foreach \s in {0,...,4}
	{
		\draw [dotted] ({2.5*\s},-0.5)--({2.5*\s},6.5);
		\node [fill=white] at ({2.5*\s},-.5) {\tiny $\s$};
	}
	
	\node [fill=white] at (2.5,0) {\textcolor{Periwinkle!50!black}{$\ds{Q}$}};
	\node [fill=white] at (2.5,1) {\mahoganydashed{$F^\times$}};
	\node [fill=white] at (2.5,2) {\textcolor{Periwinkle!50!black}{$K^M_2(F)$}};
	\node [fill=white] at (2.5,3) {\mahoganydashed{$K^M_3(F)$}};
	\node [fill=white] at (2.5,4) {\textcolor{Periwinkle!50!black}{$K^M_4(F)$}};
	\node [fill=white] at (2.5,5) {\mahoganydashed{$K^M_5(F)$}};
	\node [fill=white] at (2.5,6) {\textcolor{Periwinkle!50!black}{$K^M_6(F)$}};

	\node [fill=white] at (5,3) {\textcolor{Periwinkle!50!black}{$H^1(\Gamma_2(F))$}};
	\node [fill=white] at (5,4) {\textcolor{Periwinkle!50!black}{$(E^1_{2,4})^+$} $\oplus$ \mahoganydashed{$H^2(\Gamma_3(F))$}};
	\node [fill=white] at (5,5) {\textcolor{Periwinkle!50!black}{$(E^1_{2,5})^+$} $\oplus$ \mahoganydashed{$\frac{((E^1_{2,4})^+ \otimes F^\times)}{\rm{im}(d^1)}$}};
	\node [fill=white] at (5,6) {?};

	\node [fill=white] at (7.5,5) {\mahoganydashed{$H^1(\Gamma_3(F))$}};
	\node [fill=white] at (7.5,6) {?};
	
	\node at (-.5,-.5) {$\nicefrac{q}{p}$};
	\end{scope}
	\end{tikzpicture}
	\caption{The $E^2$-page $E^2_{p,q}$ of the rank spectral sequence, converging to $K_q(F)$. As before, the colours denote the \textcolor{Periwinkle!50!black}{$+1$-eigenspaces} and \mahoganydashed{$-1$-eigenspaces}. The $d^r$-differential has bidegree $(-r,-1)$ and necessarily sends $\pm 1$-eigenspaces to $\pm 1$-eigenspaces, so there can be only a few nonzero differentials in this range.}
	\label{fig:e2page}
\end{figure}

\cref{fig:e2page} records the results about the $E^2$-page that result from the discussion in the previous subsections. From the $E^2$-page onwards, there are no possibly nonzero differentials affecting bidegrees $q \le 5$.
Thus the three columns displayed describe the associated graded of a filtration $K_q(F)$ for $q \leq 5$, where we already know that the Milnor $K$-theory is a summand.

From the rows $q \leq 4$, we obtain the following by inspection (by convention in this section all groups are rationalised):

\begin{proposition}\label{prop:Kleq4}There are isomorphisms
\begin{enumerate}[(i)]
\item $F^\times \overset{\cong}\lra K_1(F)$. 
\item $K_2^M(F) \overset{\cong}\lra K_2(F)$. 
\item $K_3^{(2)}(F) \overset{\cong}\lra H^1(\Gamma_2(F))$. 
\item $K_4^{(3)}(F) \overset{\cong}\lra H^2(\Gamma_3(F))$. 
\item $K_4^{(2)}(F) \overset{\cong}\lra (E^1_{2,4})^+ \cong H_2(\PGL_2(F),\StL_2(F))$.
\end{enumerate}
\end{proposition}

From the rows $q=5$, we obtain the following (by convention in this section all groups are rationalised), recalling that we have
\[H^1(\Gamma_3(F))  = \ker\big[\delta \colon B_3(F) \to B_2(F) \otimes F^\times\big].\]

\begin{theorem}\label{thm:K5}
    There is an exact sequence 
    \[K_4^{(2)}(F) \otimes F^\times \overset{m_F}\lra K_5^{(3)}(F) \overset{p_F}\lra H^1(\Gamma_3(F)) \lra 0,\]
    where the first map is induced by the multiplication in $K$-theory and the second one is induced by the edge homomorphism.
\end{theorem}
 
\begin{proof}
    We know from \cref{lem:involution-vs-weight} that 
    \[K_5^-(F) = K_5^M(F) \oplus K_5^{(3)}(F)\] 
    and hence the $-1$-eigenspaces in the second and third column give a filtration for $\smash{K_5^{(3)}(F)}$. Thus, using the identification $(E^1_{2,4})^+ = K_4^{(2)}(F)$ from \cref{prop:Kleq4}, there is a short exact sequence 
    \[0 \lra (K_4^{(2)}(F) \otimes F^\times)/\rm{im}(d^1) \lra K_5^{(3)}(F) \lra H^1(\Gamma_3(F)) \lra 0,\]
    where the maps are as described.
\end{proof}

\begin{corollary}\label{cor:weight-3-beilinson-soule} If $K_4^{(2)}(F) = 0$, then the edge homomorphism induces an isomorphism 
\[K_5^{(3)}(F) \overset{\cong}\lra H^1(\Gamma_3(F)).\]
\end{corollary}

\begin{example}\label{ex:beilinson-soule-vanishing} The hypothesis for \cref{cor:weight-3-beilinson-soule} holds for number fields and fields of transcendence degree $1$ over a finite field \cite{BorelStable,Harder}. It is also closed under iterated transcendental extensions. To see this, we use the localisation theorem for algebraic $K$-theory \cite[V.6.7.1]{Weibel}, which gives split short exact sequences
\[0 \lra K_n(F) \overset{i}\lra K_n(F(t)) \overset{\partial}\lra \bigoplus_{\fr{p}} K_{n-1}(F[t]/\fr{p}) \lra 0\]
where $i$ is induced by the inclusion $F \to F(t)$ and $\partial$ is a connecting homomorphism. Each of these decompose in short exact sequences for fixed weights as \cite[Th\'eor\'eme 4]{Soule}
\[0 \lra K^{(r)}_n(F) \overset{i}\lra K^{(r)}_n(F(t)) \overset{\partial}\lra \bigoplus_{\fr{p}} K^{(r-1)}_{n-1}(F[t]/\fr{p}) \lra 0.\]
Taking $n=4$ and using that $\smash{K_3^{(1)}(E)}$ vanishes for all fields $E$ \cite[Corollaire 1]{Soule}, the inclusion  $F \to F(t)$ induces an isomorphism $\smash{K^{(2)}_4(F) \cong K^{(2)}_4(F(t))}$.
\end{example}

\subsection{The Beilinson--Soul\'e  vanishing conjecture and the Goncharov conjectures in low weight} In light of \cref{cor:weight-3-beilinson-soule}, one might ask whether validity of the Goncharov conjectures in weight 3 is equivalent to the validity of the Beilinson--Soul\'e vanishing conjecture in weight 2. In order to answer this, we first analyse how the exact sequence of \cref{thm:K5}
\[0 \lra \ker(m_F) \lra K_4^{(2)}(F) \otimes F^\times \overset{m_F}\lra K_5^{(3)}(F) \overset{p_F}\lra H^1(\Gamma_3(F)) \lra 0\]
behaves under pure transcendental extensions.

\begin{theorem} \label{thm:gamma-53}\,
\begin{enumerate}[(i)]
\item There is a short exact sequence
    \[0 \lra H^1(\Gamma_3(F)) \lra H^1(\Gamma_3(F(t))) \lra \bigoplus_\fr{p}{K_4^{(2)}(F[t]/\fr{p})/K_4^{(2)}(F)} \lra 0\]
where the sum is over all nonzero prime ideals $\fr{p}$ of $F[t]$. 
\item The map $\ker(m_F) \to \ker(m_{F(t)})$ is always an isomorphism.
\end{enumerate}
\end{theorem}

\begin{proof} We claim there is a map of short exact sequences
\[\begin{tikzcd} 0 \rar & K_4^{(2)}(F) \otimes F^\times \rar \dar{m_F} & K_4^{(2)}(F) \otimes F(t)^\times \rar \dar{m_{F(t)}} & \bigoplus_\fr{p} K_4^{(2)}(F) \rar \dar & 0 \\
0 \rar & K_5^{(3)}(F) \rar{i} & K_5^{(3)}(F(t)) \rar{\partial} &  \bigoplus_\fr{p} K_4^{(2)}(F[t]/\fr{p}) \rar & 0\end{tikzcd}\]
where the direct sums run over all nonzero prime ideals $\fr{p}$ of $F[t]$.
The bottom short exact sequence is obtained from localisation exact sequence as in \cref{ex:beilinson-soule-vanishing} for $F(t)$. To obtain the top short exact sequence, we tensor the short exact sequence 
\[0 \lra F^\times \lra F(t)^\times \lra {\textstyle \bigoplus_\fr{p}} \bb{Q} \lra 0\]
with $K_4^{(2)}(F)$. The left vertical map is the multiplication in algebraic $K$-theory and the middle vertical map is the isomorphism 
\[K_4^{(2)}(F) \otimes F(t)^\times \overset{\cong}\lra K_4^{(2)}(F(t)) \otimes F(t)^\times\] followed by multiplication; that the left square commutes is then the naturality of the multiplicative structure on algebraic $K$-theory groups. The right vertical map is given on the term $\fr{p}$ by the map induced by the inclusion $F \to F[t]/\fr{p}$. To see the right square commutes, note that the vector space $\smash{K_4^{(2)}(F(t)) \otimes F(t)^\times}$ is spanned by
\[\alpha \otimes p(t) \qquad \text{where $\alpha \in K_4^{(2)}(F)$ and $p(t) \in F[t]$ is irreducible}.\] 
The top-right composition sends this to the element that is the image of $\alpha$ in the term $\fr{p}$ where $\fr{p}$ is the ideal generated by $p(t)$. The left-bottom composition sends this to (denoting by $\cdot$ the multiplication in algebraic $K$-theory)
\[\partial(\alpha \cdot p(t)) = \alpha \cdot \partial(p(t)),\]
using the localisation sequence is a sequence of $K_*(F)$-modules  \cite[V.6.1]{Weibel}, so that in particular $\partial$ is a $K_*(F)$-module map. It is standard that in the localisation sequence $\partial(p(t))$ is given by $1$ in the term $\fr{p}$ where $\fr{p}$ is the ideal generated by $p(t)$.

Now we use that the right vertical map is injective by the existence of transfer maps \cite[V.3.3.2]{Weibel} and the identification of the left and middle vertical maps as those appearing in \cref{thm:K5}, using that $\smash{K_4^{(2)}(F)} \xrightarrow{\cong} \smash{K_4^{(2)}(F(t))}$, and apply the snake lemma to get
an exact sequence
\[\begin{tikzcd} 0 \rar & \ker(m_F) \rar & \ker(m_{F(t)}) \rar & 0 \arrow[dll,snake left] \\
& H^1(\Gamma_3(F)) \rar & H^1(\Gamma_3(F(t)) \rar & {\textstyle \bigoplus_\fr{p}} K_4^{(2)}(F[t]/\fr{p})/K_4^{(2)}(F) \rar & 0\end{tikzcd}\]
which gives the results.\end{proof}

\begin{remark}We can obtain slightly stronger results under the hypothesis of the existence of a specialisation-at-zero map: the precise claim is that
\begin{align*}\scr{G}_3(F(t)) &\lra \scr{G}_3(F) \\
\LiG_3(f(t)) &\longmapsto \begin{cases} \LiG_3(f(0)) & \text{if $f(0) \neq 0,\infty$,} \\
0 & \text{else,}\end{cases}\end{align*}
is well-defined. The existence of such specialisation maps is known on Milnor $K$-theory \cite[Theorem III.7.3]{Weibel} as well as the Bloch groups \cite[Section 3.2]{GoncharovEuler}, and amounts to verifying that the 22-term relation is sent to zero by this map; in principle a sufficiently determined reader could do so (possibly computer-aided). We strongly believe such maps exist and rather than addressing this in ad-hoc fashion here, we intend to do so in a systemic fashion in future work. In that case $\Gamma_5^3(F) \to \Gamma_5^3(F(t))$ naturally split by the specialisation-at-zero map.
\end{remark}

We obtain from this some conditions under which $H^1(\Gamma_3(F))$ is unchanged by transcendental extensions, that is, is homotopy-invariant.

\begin{corollary}\label{cor:gamma-53-iso} The map $H^1(\Gamma_3(F)) \to H^1(\Gamma_3(F(t)))$ is an isomorphism if one of the following conditions is satisfied:
\begin{enumerate}[(i)]
    \item \label{enum:gamma-53-iso-i} $F$ is algebraically closed, or
    \item \label{enum:gamma-53-iso-ii} $K_4^{(2)}(E) = 0$ for every finite extension $E$ of $F$.
\end{enumerate}
\end{corollary}

\begin{proof} We apply the short exact sequence of \cref{thm:gamma-53}. In case \eqref{enum:gamma-53-iso-i}, since $F$ is algebraically closed then all the summands of the last term are trivial. In case \eqref{enum:gamma-53-iso-ii}, since $F[t]/\fr{p}$ is a finite extension of $F$ for all $\fr{p}$, each summand vanishes.
\end{proof}

\begin{example}Part \eqref{enum:gamma-53-iso-ii} applies when $F$ is a number field or has transcendence degree $1$ over a finite field by \cref{ex:beilinson-soule-vanishing}.
\end{example}

In fact, in some sense the converse of \cref{cor:gamma-53-iso} \eqref{enum:gamma-53-iso-ii} is also true: Goncharov's homotopy-invariance conjecture in weight 3 is equivalent to the Beilinson--Soul\'e vanishing conjecture in weight 2. 

\begin{corollary}The map $H^1(\Gamma_3(F)) \to H^1(\Gamma_3(F(t)))$ is an isomorphism for every field $F$ if and only if $\smash{K_4^{(2)}(E) = 0}$ for every field $E$.
\end{corollary}

\begin{proof} The direction $\Leftarrow$ is a direct consequence of \cref{cor:gamma-53-iso} \eqref{enum:gamma-53-iso-ii}. For the direction $\Rightarrow$, we use that we have already seen that the vanishing of $\smash{K_4^{(2)}(F)}$ is closed under pure transcendental extensions in \cref{ex:beilinson-soule-vanishing}. It is also true for $F = \bb{Q}$ or $\bb{F}_q$, and so it is true for all fields if we can also show it is closed under finite extensions. It suffices to do so for simple extensions, and if $H^1(\Gamma_3(F)) \to H^1(\Gamma_3(F(t)))$ is an isomorphism for every field $F$ then from \cref{thm:gamma-53} we deduce that
\[K_4^{(2)}(F) \lra K_4^{(2)}(F[t]/\fr{p})\]
is surjective for every field $F$ and nonzero prime ideal $\fr{p}$.
\end{proof}

We also obtain the following corollary saying that the Goncharov conjecture in weight 3 imply the Beilinson--Soul\'e vanishing conjecture in weight $2$.  

\begin{corollary} \label{cor: weight 3 goncharov iff beilinson soule}
    The map $p_F \colon K_5^{(3)}(F) \to H^1(\Gamma_3(F))$ is an isomorphism for every field $F$ if and only if $K_4^{(2)}(E)=0$ for every field $E$.
\end{corollary}

\begin{proof}
    The direction $\Leftarrow$ is already done in \cref{cor:weight-3-beilinson-soule}. For the direction $\Rightarrow$ suppose that $F$ is a field such that $K_4^{(2)}(F) \neq 0$, then 
    \[K_4^{(2)}(F) \otimes F^\times \subsetneq K_4^{(2)}(F) \otimes F(t)^\times \cong K_4^{(2)}(F(t)) \otimes F(t)^\times.\]
    By \cref{thm:gamma-53} we know that $\ker(m_F) \cong \ker(m_{F(t)})$ and thus $\ker(m_{F(t)}) \subsetneq K_4^{(2)}(F(t)) \otimes F(t)^\times$ so $m_{F(t)} \neq 0$ otherwise we would not have a strict inclusion. Hence $\ker(p_{F(t)}) \neq 0$ and $p_{F(t)}$ cannot be an isomorphism.
    \end{proof}

We end with a result relating the Beilinson--Soul\'e vanishing conjecture in weight $2$ to the multiplicative structure in algebraic $K$-theory; this does not rely on our results above, but is of course related.

\begin{proposition}\label{prop:beilinson-soule-multiplication}
If there is a field $F$ so that $K_4^{(2)}(F) \neq 0$, then there exists a field $E$ of the same characteristic such that the following product map is not injective and not zero
\[m_E \colon K_4^{(2)}(E) \otimes E^\times \lra K_5^{(3)}(E).\]
\end{proposition}

\begin{proof}We first prove we can find $E$ so that it is not injective. If $K_4^{(2)}(F) \neq 0$ for some $F$ then, as algebraic $K$-theory commutes with filtered colimits, we can in fact find such an $F$ which is generated by finitely many elements over $k$, for $k=\mathbb{F}_p$ if $\rm{char}(F)=p$ and $k=\Q$ if $\rm{char}(F)=0$. In particular, we may assume $F$ has finite transcendence degree over $k$. 

Let us pick $n \in \bb{N}$ smallest so that there is a field $E_0$ with $K_4^{(2)}(E_0) \neq 0$ and $\rm{trdeg}_k(E_0)=n$. By the work of Borel, $K_4(F)$ vanishes for any number field $F$, and by the work of Quillen it vanishes (rationally) for any finite field too, so we must have $n \ge 1$. Let $F_n \coloneq k(t_1,\dots,t_n)$ be the field extension of $k$ obtained by adding $n$ pure transcendental variables, so that we can view $E_0$ as a finite degree field extension of $F_n$.

In the positive characteristic $p$ case we factor the finite extension $F_n \subset E_0$ as $F_n \subset E_1 \subset E_0$ where $F_n \subset E_1$ is finite and separable and $E_1 \subset E_0$ is finite and purely inseparable. We claim that $\smash{K_4^{(2)}(E_1)} \neq 0$ and then we can without loss of generality replace $E_0$ by $E_1$ and assume that $F_n \subset E_0$ is separable. To prove the claim, if we let $p^N = [E_0: E_1]$ be the degree then $\rm{Frob}^N(E_0) \subset E_1$ (as the extension is purely inseparable) we can consider the composition
\[K_4^{(2)}(E_0) \xrightarrow{\rm{Frob}^N_*} K_4^{(2)}(E_1) \lra K_4^{(2)}(E_0),\]
where the second map is induced by the inclusion $E_1 \subset E_0$. This composition agrees with the action of an $N$-fold composition of the Frobenius, which acts by multiplication by $p^2$ in weight 2, and hence the composition is multiplication by $p^{2N}$, showing that the second map is surjective, and hence the claim.  

Let $E$ denote the Galois closure of $E_0$ over $F_n$, which is defined by separability, then $E_0 \subseteq E$ is a finite degree extension and hence the map induced by the inclusion
\[K_4^{(2)}(E_0) \lra K_4^{(2)}(E)\]
is injective by the existence of transfer maps. In particular, we have $K_4^{(2)}(E) \neq 0$. For the rest of the proof, $G \coloneq \rm{Gal}(E/F_n)$, a finite group. We next consider the following commutative diagram 
    \[\begin{tikzcd}  & K_5^{(3)}(F_n) \dar{i_1} \rar{p_F}[swap]{\simeq} & H^1(\Gamma_3(F_n))\dar{i_2}\\[-5pt]  K_4^{(2)}(E) \otimes E^\times \rar{m_E} & K_5^{(3)}(E) \rar{p_E} & H^1(\Gamma_3(E)) \rar & 0\end{tikzcd}\]
    where the bottom short exact sequence is as in \cref{thm:gamma-53}, the top map is an isomorphism by \cref{ex:beilinson-soule-vanishing}, and the vertical maps are induced by the inclusion. Moreover, we claim that $i_2$ is injective. This is trivial in positive characteristic since $\smash{K_5^{(3)}(F_n)} \cong \smash{K_5^{(3)}(\mathbb{F}_p)}=0$ by minimality of $n$ and the localization sequence, so we just need to verify the characteristic zero case. To do so, we fix an embedding $E \subset \bb{C}$ and observe that the Borel regulator $\smash{K_5^{(3)}(\bb{C})} \to \bb{C}$ factors over $\Gamma_5^3(\bb{C})$ because the Borel regulators agree with Beilinson regulators up to a nonzero constant and the latter satisfy a product formula as  theory of Chern classes \cite[V.11]{Weibel}. Thus the composition
    \[\bb{Q} \cong K_5^{(3)}(\Q) \overset{\cong}\lra \Gamma_5^3(\Q) \overset{\cong}\lra H^1(\Gamma_3(F_n)) \overset{i_2}\lra H^1(\Gamma_3(E)) \lra H^1(\Gamma_3(\bb{C})) \lra \bb{C}\]
    with right-most map the Borel regulator, must be injective.

We use this to deduce that 
\[(K_4^{(2)}(E) \otimes E^\times)^G \subseteq \ker(m_E).\]
Given $\alpha \in (K_4^{(2)}(E) \otimes E^\times)^G$, since the multiplication on algebraic $K$-theory groups is compatible with the Galois action (because tensor product of vector spaces is), we deduce that $m_E(\alpha) \in \smash{K_5^{(3)}(E)^G}$. By the existence of transfer maps, there is a unique $\beta \in \smash{K_5^{(3)}(F_n)}$ such that $i_1(\beta)=m_E(\alpha)$.  Exactness of the bottom row implies that $p_E(i_1(\alpha))=0$ so commutativity of the diagram says that $i_2(p_F(\beta))=0$. Since $i_2$ is injective and $p_F$ is an isomorphism, $\beta=0$ so $m_E(\alpha)=i_1(\beta)=0$ and thus $\alpha \in \ker(m_E)$ as required.

Finally, we prove that $(K_4^{(2)}(E) \otimes E^\times)^G \neq 0$. Since $G$ is a finite group and $K_4^{(2)}(E) \neq 0$, the latter contains a finite irreducible $G$-representation. Since its dual occurs in the regular representation $\bb{Q}[G]$, it suffices to prove that $E^\times$ contains the regular representation; this is a consequence of \cite[Proposition 19]{BDEPS} using that $n \geq 1$.

\medskip

To get field $E$ so that the product map is not injective and not zero, take $E'$ so that $m_{E'}$ is not injective and let $E=E'(t)$. By definition the map $m_{E'}$ is not injective and hence $\ker(m_{E'}) \neq 0$ and $\smash{K_4^{(2)}(E)}\neq 0$. From \cref{thm:gamma-53}, we know $\ker(m_{E'}) \cong \ker(m_E) \neq 0$ so $m_{E}$ is not injective. But also, 
\[\ker(m_{E}) \subseteq K_4^{(2)}(E') \otimes (E')^\times \subsetneq K_4^{(2)}(E) \otimes E^\times,\]
where the second inequality uses that $K_4^{(2)}(E) \neq 0$. 
Therefore the map $m_{E}$ is not zero.
\end{proof}

\appendix

\section{Categorical foundations} To keep this appendix---and the next two---of manageable length, we opt only to include results that either are not well-known or of crucial importance to our paper. In this appendix we discuss three foundations of categeorical nature: duoidal categories (and their variants), Day convolution, and filtered/graded objects.

\begin{convention}We refer to $\infty$-categories as \emph{categories} and ordinary categories as \emph{1-categories}. We denote the former by a calligraphic font $\scr{C}$ and the latter by a roman font $\rm{C}$.
\end{convention}

\subsection{Duoidal categories and their variants} Firstly, we will have use for categories with two compatible tensor products, one monoidal and the other symmetric monoidal and oplax monoidal for the other. 

\subsubsection{Duoidal categories} \label{sec:duoidal} Let $\Mon^\rm{oplax}(\Cat)$ denote the category of monoidal categories and oplax monoidal functors. This has finite products, allowing us to make the next definition, following the work of Torii \cite{ToriiDuoidal,ToriiHigher,ToriiMult}:

\begin{definition}\label{def:duoidal} \,
\begin{enumerate}[(i)]
    \item An \emph{($E_1,E_1)$-duoidal category} is a monoid object $(\scr{C},\levi,\para)$ in $\Mon^\rm{oplax}(\Cat)$.
    \item An \emph{($E_\infty,E_1)$-duoidal category} is a commutative monoid object $(\scr{C},\levi,\para)$ in $\Mon^\rm{oplax}(\Cat)$.
\end{enumerate}
\end{definition}

Unwinding the definition, the structure of an $(E_1,E_1)$-duoidal category on $\scr{C}$ gives (a) a monoidal category structure with tensor product $\para$ and unit $1_{\para}$, and (b) a monoidal category structure with tensor product $\levi$ and unit $1_\levi$ so that all its structure is oplax monoidal for $\para$. In particular, the oplax monoidalities of the functors $1_\levi \colon \ast \to \scr{C}$ and $\levi \colon \scr{C} \times \scr{C} \to \scr{C}$ yield maps
\[1_\levi \lra 1_{\para} \qquad \text{and} \qquad 1_{\para} \lra 1_{\para} \levi 1_{\para},\]
and we say a duoidal category is \emph{normal} if these are equivalences, i.e.~if those two functors are unital. If so, we identify $1_\levi$ and $1_{\para}$ and write $1$ for them. 

More generally, the oplax monoidality on the functor $\levi$ yields an ``interchange'' natural transformation
\begin{equation}\label{eqn:duoidal-zeta}\zeta \colon (A \para B) \levi (C \para D) \lra (A \levi C) \para (B \levi D),\end{equation}
and from this we can extract a pair of natural transformations
\begin{equation}\label{eqn:duoidal-ol-ul-zeta} \begin{aligned}\bar{\zeta} &\colon A \levi B \overset{\simeq}\longleftarrow (A \para 1) \levi (1 \para B) \overset{\zeta}\lra (A \levi 1) \para (1 \levi B) \overset{\simeq}\lra A \para B, \\
\ul{\zeta} &\colon A \levi B \overset{\simeq}\longleftarrow (1 \para A) \levi (B \para 1) \overset{\zeta}\lra (1 \levi B) \para (A \levi 1) \overset{\simeq}\lra B \para A.\end{aligned}\end{equation}

For an $(E_\infty,E_1)$-duoidal category, the second tensor product $\levi$ is symmetric monoidal; that its symmetry is oplax monoidal means that the following commutes
\[\begin{tikzcd}
    (A \para B) \levi (C \para D) \rar{\zeta} \dar{\sigma} & (A \levi C) \para (B \levi D) \dar{\sigma \para \sigma} \\[-5pt]
    (C \para D) \levi (A \para B) \rar{\zeta} & (C \para A) \levi (D \para B).\end{tikzcd}\] 
This implies that $\bar{\zeta} \circ \sigma = \ul{\zeta}$.
    
\begin{remark}An $(E_1,E_1)$-duoidal structure on a 1-category is the same as a duoidal category in the sense of \cite[Sections 2.1, 2.2]{GarnerLopezFranco}, also known as 2-monoidal category \cite[Chapter 6]{AguiarMahajan}. See \cite{BataninMarkl} for a detailed descriptions of the axioms. An $(E_\infty,E_1)$-duoidal structure on a 1-category is the same as a $\levi$-symmetric one.\end{remark}

\subsubsection{A duoidal Eckmann--Hilton argument} We will give a proof of an Eckmann--Hilton argument in $(E_1,E_1)$-duoidal 1-categories. This is a slightly stronger variant of the dual of \cite[Proposition 30]{GarnerLopezFranco}:

\begin{proposition}\label{prop:duoidal-eh} Let $(\rm{C},\levi,\para)$ be a normal $(E_1,E_1)$-duoidal 1-category. Suppose that the object $C \in \rm{C}$ carries two counital operations $\delta \colon C \to C \levi C$ and $\Delta \colon C \to C \para C$ so that the following commutes
    \[\begin{tikzcd} C \rar{\delta} \arrow{dd}{\Delta} & C \levi C \dar{\Delta \levi \Delta} \\[-5pt]
    & (C \para C) \levi (C \para C) \dar{\zeta} \\[-5pt]
    C \para C \rar{\delta \para \delta} & (C \levi C) \para (C \levi C).\end{tikzcd}\]
Then $\overline{\zeta} \circ \delta = \Delta = \ul{\zeta} \circ \delta$.\end{proposition}

\begin{proof}We first verify that the counits $\epsilon \colon C \to 1$ and $E \colon C \to 1$ for $\delta$ and $\Delta$ agree. Recall that counitality implies that for any morphism $f \colon C \to D$ the following compositions
\[C \xrightarrow{\Delta} C \para C \xrightarrow{f \para E} D \para 1 \xrightarrow{\cong} D \qquad C \xrightarrow{\Delta} C \para C \xrightarrow{E \para f} 1 \para D \xrightarrow{\cong} D\]
\[C \xrightarrow{\delta} C \levi C \xrightarrow{f \levi \epsilon} D \levi 1 \xrightarrow{\cong} D \qquad C \xrightarrow{\delta} C \levi C \xrightarrow{\epsilon \levi f} 1 \levi D \xrightarrow{\cong} D\]
are equal to $f$. We claim that $\epsilon = E$ and to do so, consider next the diagram
    \[\begin{tikzcd} C \arrow{dd}{\Delta} \rar{\delta} & C \levi C \dar{\Delta \levi \Delta} \rar{\epsilon \levi \epsilon} &[30pt] 1 \levi 1 \\[-5pt]
    & (C \levi C) \para (C \levi C) \dar{\zeta} \rar{(\epsilon \levi E) \para (E \levi \epsilon)} & (1 \para 1) \levi (1 \para 1) \arrow[bend left=25]{ldd}{\cong}[swap]{\zeta} \uar{\cong} \\[-5pt]
     C \para C \dar{E \para E} \rar{\delta \para \delta} & (C \levi C) \para (C \levi C) \dar{(\epsilon \levi E) \para (E \levi \epsilon)}\\[-5pt] 1 \para 1 & (1 \levi 1) \para (1 \levi 1) \lar{\cong} \end{tikzcd}\]
where the top-left commutes by hypothesis, the bottom-left by counitality of $\delta$, the top-right by counitality of $\Delta$, and bottom-right by naturality of $\zeta$. Then by counitality of $\delta$ the top composition agrees with $\epsilon$, and by counitality of $\Delta$ the left composition agrees with $E$.

We use this to prove that $\overline{\zeta} \circ \delta = \Delta$:
    \[\begin{tikzcd} C \rar{\delta} \arrow{dd}{\Delta} & C \levi C \dar{\Delta \levi \Delta} &[30pt] \\[-5pt] 
    & (C \para C) \levi (C \para C) \dar{\zeta} \rar{(\rm{id} \para E) \levi (E \para \rm{id})} & (C \para 1) \levi (1 \para C) \arrow[bend right=10]{lu}[swap]{\cong} \dar{\zeta} \\[-5pt]
    C \para C \rar{\delta \para \delta} \arrow[bend right=10]{rrd}[swap]{\rm{id}} & (C \levi C) \para (C \levi C) \rar{(\rm{id} \levi \epsilon) \para (\epsilon \levi \rm{id})} & (C \levi 1) \para (1 \levi C) \dar{\cong} \\[-5pt]
    & & C \para C
    \end{tikzcd}\]
where the left commutes by hypothesis, the top triangle by counitality of $\Delta$, the bottom-right square by naturality of $\zeta$ and the fact that $\epsilon = E$, and the bottom triangle by counitality of $\delta$. The top-right composition is $\overline{\zeta} \circ \delta$ and the left-bottom composition is $\Delta$.

Similarly, that $\ul{\zeta} \circ \delta = \Delta$ is proven by
    \[\begin{tikzcd} C \rar{\delta} \arrow{dd}{\Delta} & C \levi C \dar{\Delta \levi \Delta} &[30pt] \\[-5pt]
    & (C \para C) \levi (C \para C) \dar{\zeta} \rar{(E \para \rm{id}) \levi (\rm{id} \para E)} & (1 \para C) \levi (C \para 1) \arrow[bend right=10]{lu}[swap]{\cong} \dar{\zeta} \\[-5pt]
    C \para C \rar{\delta \para \delta} \arrow[bend right=10]{rrd}[swap]{\rm{id}} & (C \levi C) \para (C \levi C) \rar{(\epsilon \levi \rm{id}) \para (\rm{id} \levi \epsilon)} & (1 \levi C) \para (C \levi 1) \dar{\cong} \\[-5pt]
    & & C \para C
    \end{tikzcd}\]
\end{proof}

\begin{corollary}\label{cor:duoidal-eh-sym} Suppose that $(\rm{C},\levi,\para)$ is an $(E_\infty,E_1)$-duoidal 1-category, and $C,\delta,\Delta$ are as in \cref{prop:duoidal-eh}. Then $\overline{\zeta} \circ \delta = \overline{\zeta} \circ \sigma \circ \delta$.\end{corollary}

\subsection{Day convolution and its variants}\label{sec:day} We will secondly need Day convolution monoidal structures on functor categories.

\subsubsection{Day convolution} We start recalling some well-known results, true in both the monoidal and symmetric settings, and in the unital and nonunital settings. If $\scr{C}$ has small colimits, then the construction of the functor category $\Fun(\scr{A},\scr{C})$ lifts to a functor 
\[\Fun(-,\scr{C}) \colon \Cat \lra \Cat\]
sending a functor $f \colon \scr{A} \to \scr{B}$ to its left Kan extension $f_! \colon \Fun(\scr{A},\scr{C}) \to \Fun(\scr{B},\scr{C})$, the left adjoint to the restriction functor $f^* \colon \Fun(\scr{B},\scr{C}) \to \Fun(\scr{A},\scr{C})$.

If $\scr{C}$ is moreover (symmetric) monoidal with tensor product preserving small colimits in each entry, for any (nonunital) (symmetric) monoidal category $\scr{A}$ there is a Day convolution (nonunital) (symmetric) monoidal structure on the functor category $\Fun(\scr{A},\scr{C})$ (cf.~\cite[Theorem 3.1]{BenMosheSchlank}), with tensor product and monoidal unit given by the left Kan extensions
\[\begin{tikzcd} \scr{A} \times \scr{A} \rar{F \times G} \dar[swap]{\otimes} & \scr{C} \\
 \scr{A} \arrow[dashed]{ru}[swap]{F \otimes G} & \end{tikzcd} \qquad \text{and} \qquad \begin{tikzcd} \ast \rar{1_\scr{C}} \dar[swap]{1_\scr{A}} & \scr{C} \\
 \scr{A}. \arrow[dashed]{ru}[swap]{1_{\Fun(\scr{A},\scr{C})}} & \end{tikzcd}\]
It is natural in the domain and target. For the domain we have (cf.~\cite[Proposition 3.3]{BenMosheSchlank}):

\begin{lemma}\label{lem:day-naturality-source} Let $\scr{A}$, $\scr{B}$ be (symmetric) monoidal categories and $\scr{C}$ be a symmetric monoidal category with small colimits such that the tensor product preserve small colimits in each entry.
	\begin{enumerate}
		\item If $f \colon \scr{A} \to \scr{B}$ is (nonunital) lax monoidal, then $f_! \colon \Fun(\scr{A},\scr{C}) \to \Fun(\scr{B},\scr{C})$ is (nonunital) oplax monoidal and its right adjoint $f^* \colon \Fun(\scr{B},\scr{C}) \to \Fun(\scr{A},\scr{C})$ is (nonunital) lax monoidal.
		\item If $f \colon \scr{A} \to \scr{B}$ is (nonunital) strong monoidal, then so is $f_!$.
	\end{enumerate}
\end{lemma}

For the target we have \cite[Remark 3.9 and Proposition 3.6]{BenMosheSchlank}.
 
\begin{lemma}\label{lem:day-naturality-target}Let $\scr{A}$ be a symmetric monoidal category and $\scr{C}$, $\scr{D}$ be (symmetric) monoidal categories with small colimits such that the tensor products preserve small colimits in each entry.
\begin{enumerate}
 		\item If $f \colon \scr{C} \to \scr{D}$ is lax monoidal, then so is $f_* \colon \Fun(\scr{A},\scr{C}) \to  \Fun(\scr{A},\scr{D})$.
 		\item If $f \colon \scr{C} \to \scr{D}$ is strong monoidal and preserves colimits, then $f_* \colon \Fun(\scr{A},\scr{C}) \to  \Fun(\scr{A},\scr{D})$ is strong monoidal.
 	\end{enumerate}
\end{lemma}

\subsubsection{Promonoidal categories and Day convolution}\label{sec:promonoidal-day}  It is well-known that  the construction of a Day convolution (symmetric) monoidal structure requires weaker input is than a (symmetric) monoidal category \cite[Section 1]{BGSII}. We will now focus on the promonoidal case, as it is more relevant to this paper and the symmetric promonoidal case is entirely analogous:

\begin{definition}
    A \emph{nonsymmetric promonoidal category} is an operad $\scr{C}^\otimes \to \Delta^\op$ whose restriction $\scr{C}^\otimes_\rm{act} \to \Delta^\rm{act}$ to active morphisms is flat in the sense of \cite[B.3.8]{LurieHA}.
\end{definition}

\begin{remark}It may be helpful to note that a promonoidal structure on $\scr{C}$ is the same as a monoidal structure on $\Fun(\scr{C},\Spc)$ whose tensor product preserves colimits in each entry \cite[Theorem 3.37]{LinskensNardinPol}.
\end{remark}

We let $\rm{ProMon}^\rm{lax}(\Cat) \subset \rm{Op}^\rm{ns}$ denote the full subcategory on the promonoidal categories. Promonoidal categories are closed under products as a consequence of \cite[B.3.12, B.3.16]{LurieHA}, so there is a functor
\begin{equation}\label{eqn:promon-prod}\begin{aligned}\rm{ProMon}^\rm{lax}(\Cat) \times \rm{ProMon}^\rm{lax}(\Cat) &\lra \rm{ProMon}^\rm{lax}(\Cat) \\
(\scr{A}^\otimes,\scr{C}^\otimes) &\longmapsto \scr{A}^\otimes \times_{\Delta^\op} \scr{C}^\otimes.\end{aligned}\end{equation} By \cite[2.8.3]{Hinich}, if $\scr{A}$ is a promonoidal category then $\scr{A}^\otimes \times_{\Delta^\op}(-) \colon \Op^\rm{ns} \to \Op^\rm{ns}$ admits a right adjoint $\Fun(\scr{A},-)^\otimes \colon \Op^\rm{ns} \to \Op^\rm{ns}$, and when evaluated on a monoidal category which has small colimits and whose tensor product preserves small colimits in each entry, the result is again a monoidal category by \cite[Theorem 3.2.6]{NardinShah} (by Proposition 3.1.7 loc.cit.\,their Day convolution agrees with the one used by Hinich, as it satisfies the same universal property). Let $\Mon^\rm{lax,colim}(\Cat) \subset \Mon^\rm{lax}(\Cat)$ denote the full subcategory on those monoidal categories which have small colimits and whose tensor product preserves small colimits in each entry. Extracting these right adjoints yields a functor
\begin{equation}\label{eqn:promon-day}\begin{aligned}\rm{ProMon}^\rm{lax}(\Cat)^\op \times \Mon^\rm{lax,colim}(\Cat) &\lra \rm{Mon}^\rm{lax}(\Cat) \\
(\scr{A}^\otimes,\scr{C}^\otimes) &\longmapsto \Fun(\scr{A},\scr{C})^\otimes.\end{aligned}\end{equation}
Fixing $\scr{C}^\otimes \in \Mon^\rm{lax,colim}(\Cat)$, this yields a functor
\[\Fun(-,\scr{C})^{\otimes} \colon \rm{ProMon}^\rm{lax}(\Cat)^\op \lra \Mon^\rm{lax}(\Cat)\]
whose naturality is given by restriction. Since $\scr{D}$ has small colimits, each restriction functor admits a left adjoint and we can compose with the mate correspondence $\Mon^\rm{R,lax}(\Cat)^\op \simeq \Mon^\rm{L,oplax}(\Cat)$ of \cite{HHLN} to extract a functor 
\[\Fun(-,\scr{C})^{\otimes} \colon \rm{ProMon}^\rm{lax}(\Cat) \to \Mon^\rm{oplax}(\Cat)\]
whose naturality is given by left Kan extension. 

\begin{remark}
To deal with oplax morphisms we need to work in the anti-operadic setting. Recall, e.g.~as a variant of \cite[Definition 1.2]{BGSII}, the category $\rm{aOp}^\rm{ns}$ of nonsymmetric monoidal antioperads, whose objects are maps $\scr{C}_\otimes \to \Delta$ whose opposite is an operad and whose morphisms are maps over $\Delta$ whose opposite is a map of operads. Then $\Mon^\rm{oplax}(\Cat) \subset \rm{aOp}^\rm{ns}$ is the full subcategory of cartesian fibrations, and taking opposites gives equivalences
\[(\Mon^\rm{oplax}(\Cat) \simeq \Mon^\rm{lax}(\Cat)) \subset (\rm{aOp}^\rm{ns} \simeq \rm{Op}^\rm{ns}).\]
We define $\rm{ProMon}^\rm{oplax}(\Cat) \subset \rm{aOp}^\rm{ns}$ as the full subcategory of those antioperads whose restriction to the active morphisms is flat, and taking opposites induces an equivalence $\rm{ProMon}^\rm{oplax}(\Cat) \simeq \rm{ProMon}^\rm{lax}(\Cat)$ generalising that for monoidal categories.\end{remark}

\subsubsection{Duoidal categories and Day convolution}\label{sec:duoidal-day}  The functor
\begin{align*} \Mon^\rm{lax}(\Cat)^\op \times \Mon^\rm{lax,colim}(\Cat) &\lra \Mon^\rm{lax}(\Cat) \\
(\scr{A}^\otimes,\scr{C}^\otimes) &\longmapsto \Fun(\scr{A},\scr{C})^\otimes \end{align*}
is lax symmetric monoidal with respect to the cartesian monoidal structures on domain and target. We will justify this in more generality momentarily, but see also \cite[Section 4.2]{PortaTeyssier} or, for $\scr{C} = \Spc$, \cite[4.8.1.10]{LurieHA}. After fixing a symmetric monoidal category $\scr{C}$ with small colimits such that the tensor products preserve small colimits in each entry, applying the mate correspondence, and evaluating on monoid objects, we obtain a functor
\[\Fun(-,\scr{C})^\otimes \colon \Mon(\Mon^\rm{lax}(\Cat)) \lra \Mon(\Mon^\rm{oplax}(\Cat)).\]
This is a lift of Day convolution to the duoidal setting.

\subsubsection{Produoidal categories and Day convolution} Recalling that $\rm{ProMon}^\rm{oplax}(\Cat)$ has finite products, we can define a variant of a duoidal category where the second monoidal structure is only promonoidal as a monoid object in $\rm{ProMon}^\rm{oplax}(\Cat)$.

We have phrased  \cref{sec:promonoidal-day,sec:duoidal-day} to apply essentially verbatim for this variant, as soon as we explain why
\[\Fun(-,\scr{C})^\otimes \colon \Mon(\rm{ProMon}^\rm{lax}(\Cat)) \lra \Mon(\Mon^\rm{oplax}(\Cat))\]
is lax symmetric monoidal with respect to the cartesian symmetric monoidal structures on domain and target. The starting point is that \eqref{eqn:promon-prod} is canonically symmetric monoidal with respect to this, so \eqref{eqn:promon-day} is lax symmetric monoidal as a consequence of \cite[Corollary A.5.1]{HinichRectification}. As a symmetric monoidal category $\scr{C}$ is a commutative monoid, it yields a monoid object in monoid objects by the additivity theorem, and hence fixing it we obtain that 
\[\Fun(-,\scr{C})^\otimes \colon \Mon(\rm{ProMon}^\rm{lax}(\Cat))^\op \lra \Mon(\Mon^\rm{lax}(\Cat))\]
with naturality in restriction, is lax monoidal. Finally, we apply the equivalence from the mate correspondence \cite{HHLN}.

\subsection{Filtered and graded objects} \label{sec:filtered-and-graded} We finally spell out how we think of graded and filtered objects. Let $\bb{Z}$ be the category whose objects are the integers and whose only morphisms are identities. Let $\bb{Z}_{\leq}$ be the category whose objects are the integers and there is a morphism $n \to m$ when $n \leq m$. That is, the category is given by $\cdots \to -1 \to 0 \to 1 \to 2 \to \cdots$.

\begin{definition}Let $\scr{C}$ be a category.
	\begin{itemize}
		\item The category of \emph{graded objects} in $\scr{C}$ is $\Fun(\bb{Z},\scr{C})$.
		\item The category of \emph{filtered objects} in $\scr{C}$ is $\Fun(\bb{Z}_{\leq},\scr{C})$.
	\end{itemize}
\end{definition}

Addition makes $\bb{Z}$ and $\bb{Z}_{\leq}$ into symmetric monoidal categories, inducing Day convolution symmetric monoidal structures on the categories of graded and filtered objects. Letting $\bb{N} \subset \bb{Z}$ denote the nonnegative integers, there are variants $\Fun(\bb{N},\scr{C})$ and $\Fun(\bb{N}_{\leq},\scr{C})$, which can be considered as full subcategories of $\Fun(\bb{Z},\scr{C})$ and $\Fun(\bb{Z}_{\leq},\scr{C})$, by left Kan extension along the inclusion. Because $\bb{N} \to \bb{Z}$ and $\bb{N}_\leq \to \bb{Z}_\leq$ are symmetric monoidal, so are these inclusions. In the following we focus on the $\bb{Z}$- and $\bb{Z}_{\leq}$-indexed graded and filtered objects, commenting only on when the $\bb{N}$- and $\bb{N}_{\leq}$-indexed graded and filtered objects behave differently.

\medskip

Every $n \in \bb{Z}$ gives rise to functors $n \colon \ast \to \bb{Z}$ and $n \colon \ast \to \bb{Z}_{\leq}$. Pulling back along these gives functors $n^* \colon \Fun(\bb{Z},\scr{C}) \to \scr{C}$ and $n^* \colon \Fun(\bb{Z}_{\leq},\scr{C}) \to \scr{C}$ which admit left and right adjoints denoted $n_!$ and $n_*$ as long as $\scr{C}$ has an initial object $\rm{i}$ and terminal object $\rm{t}$. These are given by
\begin{align*}n_!(X) &\simeq (\cdots \to \rm{i} \to \rm{i} \to X \xrightarrow{\id} X \xrightarrow{\id} \cdots ) \\
n_*(X) &\simeq (\cdots \xrightarrow{\id} X \xrightarrow{\id} X \to \rm{t} \to \rm{t} \to \cdots)\end{align*}
where the first, or last, $X$ appears in filtration degree $n$. If $\scr{C}$ is pointed the functor $n_!$ admits a further left adjoint $n^\dagger$ given by
\[n^\dagger(X) \simeq \rm{cofib}(X(n-1) \to \colim X).\]

\begin{lemma}\quad \label{lem:eval-lkan-sym-mon}
	\begin{enumerate}
		\item $n^* \colon \Fun(\bb{Z}_{\leq},\scr{C}) \to \scr{C}$ is a left and right adjoint, lax symmetric monoidal if $n=0$ and nonunital lax symmetric monoidal if $n < 0$.
		\item $n_! \colon \scr{C} \to \Fun(\bb{Z}_{\leq },\scr{C})$ is a left adjoint, symmetric monoidal if $n=0$ and nonunital oplax monoidal if $n < 0$.
	\end{enumerate}
\end{lemma}

\begin{proof}The functor $n$ is symmetric monoidal if and only if $n=0$, and the first parts of (1) and (2) follow from \cref{lem:day-naturality-source}. It is nonunital oplax monoidal (there is a morphism $n \leq n+n$ but no morphism $n \leq 0$) if $n>0$ and nonunital lax monoidal (there is a morphism $n+n \leq n$ but no morphisms $0 \leq n$) if $n<0$.
\end{proof}

\section{Operadic foundations} In this appendix we discuss (co)operads and (co)algebras, as well as Koszul duality, in the general setting.

\subsection{Operads, cooperads, algebras, and coalgebras}

We start with a recollection of operads and algebras, and the dual notions of cooperads and coalgebras. Details can be found in \cite{BCN,Wu,HeineMM,HeutsKoszul}. 

\subsubsection{Operads and algebras}\label{sec:operads-algebras} We fix a presentable symmetric monoidal category $\scr{C}$. Letting $\Fin^\simeq$ denote the groupoid of nonempty finite sets and bijections, the category of \emph{symmetric sequences} is defined as the functor category
\[\rm{SSeq}(\scr{C}) \coloneq \Fun(\Fin^\simeq,\scr{C}).\]
This admits two tensor products:
\begin{enumerate} 
\item A \emph{Day convolution} symmetric monoidal structure induced by disjoint union of finite sets, with tensor product denoted $\otimes$. The underlying object of $X \otimes Y$ is given by 
\[(X \otimes Y)(\ul{n}) \simeq \bigsqcup_{n = n_1+n_2} \fr{S}_{n} \times_{\fr{S}_{n_1} \times \fr{S}_{n_2}} X(\ul{n}_1) \otimes_\scr{C} Y(\ul{n}_2)\]
and monoidal unit given by $1_\scr{D}$ concentrated in arity $0$.
\item A \emph{composition} monoidal structure, with tensor product denoted $\circ$. The underlying object is given in terms of the Day convolution tensor product by
	\[X \circ Y \simeq \bigsqcup_{r \geq 0} X(r) \otimes_{\fr{S}_r} Y^{\otimes r},\]
and the monoidal unit is given by $1_\scr{D}$ concentrated in arity $1$. See \cite[p.~45]{Wu} for a formal definition. 
\end{enumerate}

\begin{definition}The category $\rm{Op}(\scr{C})$ of \emph{operads} is the category $\rm{Alg}(\rm{SSeq}(\scr{C}))$ of unital associative algebras in symmetric sequences in $\scr{C}$ under the composition product.
\end{definition}

The left Kan extension $i_{0,!}$ along the inclusion of the empty set into $\Fin^\simeq$ gives an identification of $\scr{C}$ with the full subcategory of $\rm{SSeq}(\scr{C})$ consisting of those symmetric sequences concentrated in arity zero. This subcategory is preserved by $X \circ -$, which induces a left action 
\[\rm{SSeq}(\scr{C}) \times \scr{C} \lra \scr{C},\]
or in other words, $\scr{C}$ is a left module over $\rm{SSeq}(\scr{C})$.

\begin{definition}For an operad $\scr{O} \in \rm{Op}(\scr{C})$, the category $\Alg_\scr{O}(\scr{C})$ of \emph{$\scr{O}$-algebras} is the category of left modules over $\scr{O}$ in $\scr{C}$ under the above left action.
\end{definition}

This construction is natural in $\scr{O}$: any map of operads $f \colon \scr{O} \to \scr{O}'$ induces a restriction functor $f^* \colon \Alg_{\scr{O}'}(\scr{C}) \to \Alg_\scr{O}(\scr{C})$ which is the identity on underlying objects. This admits both a left and right adjoint $f_!,f_* \colon \Alg_{\scr{O}}(\scr{C}) \to \Alg_{\scr{O}'}(\scr{C})$ \cite[4.2.3.8]{LurieHA}.

By definition an operad has a unit map $\eta \colon 1_\scr{C} \to \scr{O}$, with domain the monoidal unit of the composition tensor product. The forgetful functor $\fgt_\scr{O} \coloneq \eta^* \colon \Alg_\scr{O}(\scr{C}) \to \scr{C}$ admits a left adjoint $\free_\scr{O} \coloneq \eta_!$: using that $\fgt_{1_\scr{C}} \colon \Alg_{1_\scr{C}}(\scr{C}) \to \scr{C}$ is an equivalence we get
\[\begin{tikzcd} \scr{C} \rar[shift left=.5ex]{\free_\scr{O} \simeq \eta_!} &[20pt] \Alg_\scr{O}(\scr{C}). \lar[shift left=.5ex]{\fgt_\scr{O} \simeq \eta^*}\end{tikzcd}\]
This adjunction exhibits $\Alg_\scr{O}(\scr{C})$ as the category of algebras in $\scr{C}$ over the monad $\rm{Sym}_\scr{O} \coloneq \fgt_\scr{O} \circ \free_\scr{O} \colon \scr{C} \to \scr{C}$. It is given on underlying objects by
\[\rm{Sym}_\scr{O}(X) \simeq \bigsqcup_{r \geq 0} \scr{O}(r) \otimes_{\fr{S}_r} X^{\otimes r}.\]

This completes the discussion of the free-forgetful adjunction, but for the cotangent-trivial adjunction we need additional data. Note that the monoidal unit $1_\scr{C}$ of $\rm{SSeq}(\scr{C})$ for the composition tensor product is canonically an operad, so we can define:

\begin{definition}The category of \emph{augmented operads} is
\[\rm{Op}^\rm{aug}(\scr{C}) \coloneq \rm{Op}(\scr{C})_{/1_\scr{C}}.\]
\end{definition}

Thus an augmentation of an operad $\scr{O}$ is a map of operads $\epsilon \colon \scr{O} \to 1_\scr{C}$, which necessarily satisfies $\epsilon \circ \eta \simeq \rm{id}_{1_\scr{C}}$. We obtain a pair of adjunctions
\[\begin{tikzcd} \scr{C} \rar[shift left=.5ex]{\free_\scr{O} \simeq \eta_!} &[20pt] \Alg_\scr{O}(\scr{C}) \rar[shift left=.5ex]{\cot_\scr{O} \coloneq \epsilon_!} \lar[shift left=.5ex]{\fgt_\scr{O} \simeq \eta^*} &[20pt] \scr{C} \lar[shift left=.5ex]{\triv_\scr{O} \coloneq \epsilon^*} \end{tikzcd}\]
which satisfy 
\[\cot_\scr{O} \circ \free_\scr{O} \simeq \rm{id}_\scr{C} \quad \text{and} \quad \fgt_\scr{O} \circ \triv_\scr{O} \simeq \rm{id}_\scr{C}.\]

\subsubsection{Cooperads and coalgebras}\label{sec:cooperads-coalgebras} There is a dual story for coalgebras. It uses that the opposite of a monoidal category is itself a monoidal category by the monoidal analogue of \cite[2.4.2.7]{LurieHA}.

\begin{definition}The category $\rm{Coop}(\scr{C})$ of \emph{cooperads} is the category $\rm{Alg}(\rm{SSeq}(\scr{C})^\op)^\op$ of counital associative coalgebras in symmetric sequences in $\scr{C}$ under composition product.
\end{definition}

\begin{definition}For a cooperad $\scr{Q} \in \rm{Op}(\scr{C})$, the category $\coAlg^{{\dpw,\nil}}_\scr{O}(\scr{C})$ of \emph{conilpotent $\scr{Q}$-coalgebras with divided powers} is the category of left comodules over $\scr{Q}$ in $\scr{C}$ under the above left action.
\end{definition}

The latter construction is natural by restriction in maps of cooperads, which admit both a left and right adjoint. Using the counit map $\epsilon \colon \scr{Q} \to 1_\scr{C}$ of a cooperad, we obtain a forgetful-cofree adjunction
\[\begin{tikzcd} \coAlg^{\dpw,\nil}_\scr{Q}(\scr{C}) \rar[shift left=.5ex]{\fgt_\scr{Q} \coloneq \epsilon^*} &[20pt] \scr{C} \lar[shift left=.5ex]{\rm{cofree}_\scr{Q} \coloneq \epsilon_*}\end{tikzcd}\]
and this exhibits $\coAlg^{{\dpw,\nil}}_\scr{O}(\scr{C})$ as the category of coalgebras for the comonad $\fgt_\scr{Q} \circ \rm{cofree}_\scr{Q} \simeq \rm{Sym}_\scr{Q}$.

\begin{definition}The category of \emph{augmented cooperads} is
\[\rm{Coop}^\rm{aug}(\scr{C}) \coloneq \rm{Coop}(\scr{C})_{1_\scr{C}/}.\]
\end{definition}

Thus an augmentation of a cooperad $\scr{Q}$ is a map of cooperads $\eta \colon 1_\scr{C} \to \scr{Q}$, from which we obtain a pair of adjunctions
\[\begin{tikzcd} \scr{C} \rar[shift left=.5ex]{\rm{cotriv}_\scr{Q} \coloneq \eta^*}  &[20pt] \coAlg^{\dpw,\nil}_\scr{Q}(\scr{C}) \lar[shift left=.5ex]{\eta_*} \rar[shift left=.5ex]{\fgt_\scr{Q} \simeq \epsilon^*} &[20pt] \scr{C}. \lar[shift left=.5ex]{\rm{cofree}_\scr{Q} \simeq \epsilon_*} \end{tikzcd}\]
which satisfy 
\[\eta_* \circ \cofree_\scr{Q} \simeq \rm{id}_\scr{C} \quad \text{and} \quad \fgt_\scr{Q} \circ \rm{cotriv}_\scr{Q} \simeq \rm{id}_\scr{C}.\]

\subsubsection{Unitalisation and augmentation ideals} \label{sec:augmentation-ideals} We start with the following definition, terminology for a (co)operad concentrated in strictly positive arities:

\begin{definition}\,
    \begin{enumerate}[(i)]
    \item The category $\rm{Op}^\rm{nu}(\scr{C}) \subset \rm{Op}(\scr{C})$ of \emph{nonunitary operads} is the full subcategory of those operads $\scr{O}$ so that $\scr{O}(0)$ is initial.
    \item The category $\rm{Coop}^\rm{nu}(\scr{C}) \subset \rm{Coop}(\scr{C})$ of \emph{nonunitary cooperads} is the full subcategory of those operads $\scr{Q}$ so that $\scr{Q}(0)$ is initial.
    \end{enumerate}
\end{definition}

For the remainder of this subsection we will work with operads but there is a dual story for cooperads. 

We can alternatively define nonunitary operads in terms of nonunitary symmetric sequences. To do so, let $\rm{Fin}^{\simeq}_{\neq \varnothing} \subset \rm{Fin}$ denote the full subcategory of nonempty finite sets, and let $\rm{SSeq}^\rm{nu}(\scr{C}) \coloneq \Fun(\rm{Fin}^{\simeq}_{\neq \varnothing},\scr{C})$. The inclusion $\iota \colon \rm{Fin}^{\simeq}_{\neq \varnothing} \hookrightarrow \rm{Fin}^{\simeq}$ induces a functor $\iota_!$ identifying $\rm{SSeq}^\rm{nu}(\scr{C})$ with the full subcategory of nonunitary operads. The condition that $X(0)$ is initial is preserved by Day convolution and convolution tensor products, so these restrict to (symmetric) monoidal structures on $\iota_!$ making $\iota_!$ (symmetric) monoidal and its right adjoint $\iota^*$ lax (symmetric) monoidal. We get an induced functor $\iota_! \colon \Op^\rm{nu}(\scr{C}) \to \Op(\scr{C})$ with right adjoint $\iota^*$.

\begin{definition}Given an operad $\scr{O}$, we define $\scr{O}^\rm{nu} \coloneq \iota_! \iota^* \scr{O}$.
\end{definition}

Note that the counit of the adjunction provides with a canonical map $\upsilon \colon \scr{O}^\rm{nu} \to \scr{O}$ which is an equivalence if and only if $\scr{O}$ is nonunitary. This induces an adjunction 
\[\begin{tikzcd} \Alg_{\scr{O}^\rm{nu}}(\scr{C}) \rar[shift left=.5ex]{(-)^+ \coloneq \upsilon_!} &[20pt] \Alg_\scr{O}(\scr{C}) \lar[shift left=.5ex]{\upsilon^*}.\end{tikzcd}\]
If we suppose that $\scr{O}$ is \emph{unital}, i.e.~$\scr{O}(0) \simeq 1_\scr{C}$, then for $\bf{A} \in \Alg_{\scr{O}^\rm{nu}}(\scr{C})$ the \emph{unitalisation} $\bf{A}^+$ has underlying object $1_\scr{C} \sqcup \bf{A}$ and its $\scr{O}$-algebra structure is informally given by the $\scr{O}^\rm{nu}$-algebra structure and 0-ary operation acting by the inclusion $1_\scr{C} \to \bf{A}^+ \simeq 1_\scr{C} \sqcup \bf{A}$.

We can do better when working in a stable setting and adding augmentations. Note that if $\scr{O}$ is an augmented operad then the augmentation endows $1_\scr{D}$ with an $\scr{O}$-algebra structure, and we can make the following definition:

\begin{definition}For an augmented operad $\scr{O}$, the category of \emph{augmented $\scr{O}$-algebras} is
\[\Alg^\rm{aug}_\scr{O}(\scr{C}) \coloneq \Alg_\scr{O}(\scr{C})_{/1_{\scr{C}}}.\]\end{definition}

The following is \cite[5.4.4.10]{LurieHA}:

\begin{proposition}\label{prop:augmentation-ideal-and-unitalisation} If $\scr{C}$ is stable and $\scr{O}$ is unital, then taking the fibres of the augmentation induces an equivalence
\[\Alg^\aug_{\scr{O}}(\scr{C}) \overset{\simeq}\lra \Alg_{\scr{O}^\rm{nu}}(\scr{C})\]
with inverse given by the unitalisation $(-)^+$.
\end{proposition}

\subsubsection{Suspension and linear duals} \label{sec:operadic-suspension-linear-duals} We now recall two constructions that can be performed on (co)operads in a presentable stable symmetric monoidal category $\scr{C}$, both using the levelwise tensor product of symmetric sequences \cite[Section 3.2]{BCN}.

The first is \emph{operadic suspension}, given by taking the levelwise tensor product with the endomorphism operad of the desuspended monoidal unit $1_\scr{C}$. An $s\scr{O}$-algebra structure on $X$ is the same as an $\scr{O}$-algebra structure on $\Sigma X$, and similarly for cooperads. That is, there are commutative squares \cite[Section 3.1]{HeutsLandFormality}
\[\begin{tikzcd} \Alg_{\scr{O}}(\scr{C}) \rar{\simeq} \dar[swap]{\fgt_\scr{O}} & \Alg_{s^k\scr{O}}(\scr{C}) \dar{\fgt_{s^k\scr{O}}} \\[-5pt]
\scr{C} \rar{\Sigma^k}[swap]{\simeq} & \scr{C},\end{tikzcd} \qquad \text{and} \qquad \begin{tikzcd} \coAlg_{\scr{Q}}(\scr{C}) \rar{\simeq} \dar[swap]{\fgt_\scr{Q}} & \coAlg_{s^k\scr{Q}}(\scr{C}) \dar{\fgt_{s^k\scr{Q}}} \\[-5pt]
\scr{C} \rar{\Sigma^{-k}}[swap]{\simeq} & \scr{C}.\end{tikzcd}\]

The second is \emph{taking duals}, which provides a way to construct cooperads from operads by taking duals, and vice versa. Taking adjoints to the levelwise tensor product of symmetric sequences \cite[Proposition 3.9, Proposition 3.47]{BCN} we can define a lax monoidal linear duality functor $\rm{SSeq}(\scr{C})^\op \to \rm{SSeq}(\scr{C})$ which is symmetric monoidal on levelwise dualisable objects. Thus from a cooperad $\scr{Q}$ we obtain a dual operad $D\scr{Q}$, and for any operad $\scr{O}$ with $\scr{O}(n)$ dualisable for $n \geq 0$ we obtain a dual cooperad $D\scr{O}$.

\subsubsection{Naturality in category} Any functor induces by postcomposition a functor
\[f^\rm{SSeq} \colon \rm{SSeq}(\scr{C}) \lra \rm{SSeq}(\scr{D}),\]
which lifts to a lax monoidal functor if $f$ is a lax symmetric monoidal functor. If $f$ is colimit-preserving, then so is $f^\rm{SSeq}$ and this further lifts to a strong or oplax symmetric monoidal functor if $f$ is a strong or oplax symmetric monoidal functor.

Passing to category of associative (co)algebras, these in turn induce functors between categories of (co)operads
\begin{align*}f^\rm{Op} \colon \rm{Op}(\scr{C}) &\lra \Op(\scr{D}), && \text{if $f$ is lax symmetric monoidal,}\\ 
f^\rm{Coop} \colon \rm{Coop}(\scr{C}) &\lra \rm{Coop}(\scr{D}), && \text{\parbox[t]{5cm}{\centering if $f$ is oplax symmetric monoidal and colimit-preserving.}}
\end{align*}

Passing to category of (co)modules, we similarly get functors between categories of (co)algebras
\begin{align*}f^\rm{Alg} \colon \rm{Alg}_{\scr{O}}(\scr{C}) &\lra \Alg_{f^\rm{Op}(\scr{O})}(\scr{D}), && \text{if $f$ is lax symmetric monoidal,}\\
f^\rm{coAlg} \colon \rm{coAlg}^{\dpw,\nil}_\scr{Q}(\scr{C}) &\lra \rm{coAlg}^{\dpw,\nil}_{f^\rm{Coop}(\scr{Q})}(\scr{D}), && \text{\parbox[t]{5cm}{\centering if $f$ is oplax symmetric monoidal and colimit-preserving.}}\end{align*}
In the nonunitary setting these constructions only require that $f$ is nonunital strong, nonunital lax, or nonunital oplax symmetric monoidal.

\subsubsection{Naturality in category, continued} \label{sec:alg-nat-left-adjoint} Of particular interest is the case of an adjunction $L \dashv R$ with $R \colon \scr{D} \to \scr{C}$ lax symmetric monoidal and $L \colon \scr{C} \to \scr{D}$ given the corresponding oplax symmetric monoidality from the mate correspondence \cite{HHLN}. Then $R^\rm{SSeq}$ preserves small limits and is accessible, since limits and filtered colimits are computed in the underlying category. Using \cite[4.2.3.7]{LurieHA} and similar reasoning, the induced functor
\[R^\rm{Alg} \colon \Alg_\scr{O}(\scr{C}) \to \Alg_{R^\rm{Op}(\scr{O})}(\scr{D})\]
is a functor between presentable categories that preserves limits and is accessible. Using the adjoint functor theorem \cite[5.5.2.9]{LurieHTT}, it hence admits a left adjoint
\[L^\rm{Alg} \colon  \Alg_{R^\rm{Op}(\scr{O})}(\scr{D}) \lra \Alg_{\scr{O}}(\scr{C}).\]

If $R$ is unital, i.e.~$1_{\scr{D}} \to R(1_\scr{C})$ is an equivalence, then by using naturality in $\scr{O}$ and the monadic resolution, we see this is uniquely characterised by preserving sifted colimits and $L^\rm{Alg} \free_{R^\rm{Op}(\scr{O})} \simeq \free_{\scr{O}} L$. Indeed, one computes
\[L^\rm{Alg}(\bf{A}) \simeq |[p] \mapsto \free_\scr{O} L(\rm{Sym}_{R^\rm{Op}(\scr{O})}^p(\bf{A}))|,\]
which also makes clear that $\fgt_\scr{O}L^\rm{Alg}(\bf{A})$ rarely agrees with $L(\fgt_{R^\rm{Op}(\scr{O})}A)$. We shall give some more details on how to construct the above simplicial object. Our starting point is the natural equivalence
\[\Map_{\Alg_{R^\Op(\scr{O})}(\scr{D})}(\free_{R^\Op(\scr{O})}(X),R^\Alg \bf{B}) \simeq \Map_{\Alg_\scr{O}(\scr{C})}(\free_\scr{O}(LX),\bf{B}).\]
By Yoneda, we obtain a functor $\Alg_{R^\Op(\scr{O})}(\scr{D})^\free \to \Alg_{\scr{O}}(\scr{C})$, where the domain is the full subcategory on the free $R^\Op(\scr{O})$-algebras, which sends the object $\free_{R^\Op(\scr{O})}(X)$ to $\free_\scr{O}L(X)$. Now apply this to the simplicial object $[p] \mapsto \free_{R^\Op(\scr{O})}(\rm{Sym}^p_{R^\Op(\scr{O})}(\bf{A}))$.

\begin{remark}This can be connected to \cite[Section 3.2.2]{GKRW18}. The following should yield alternative construction of $L^\rm{Alg}$, corresponding to the construction there after passing from model categories to $\infty$-categories: there is an equivalence 
\[\rm{RMod}_T(\Fun(\scr{C},\scr{D})) \overset{\simeq}\lra \Fun^\rm{sifted}(\Alg_T(\scr{C}),\scr{D})\]
for a monad $T$ preserving sifted colimits (combine \cite[Corollary 5.29]{Heine} with the fact that for such monads the restriction $\Fun^\rm{sifted}(\Alg_T(\scr{C}),\scr{D}) \to \Fun(\Alg'_T(\scr{C}),\scr{D})$ to the essential image $\Alg'_T(\scr{C}) \subset \Alg_T(\scr{C})$ of the free $T$-algebra functor is an equivalence), which sends a right $T$-module functor $F$ to the sifted colimit
\[\bf{A} \longmapsto F^\rm{Alg}(\bf{A}) \coloneq |[p] \mapsto F(T^p(\bf{A}))|.\]
To get $L^\rm{Alg}$, we apply this to $F = \free_{\scr{O}} \circ L$ with right $\rm{Sym}_{R(\scr{O})}$-module structure arising from the natural transformation $L \rm{Sym}_{R(\scr{O})} \to \rm{Sym}_{\scr{O}} L$ induced by the oplax monoidality of $L$ and the counit of the adjunction $L \dashv R$, and the canonical right $\rm{Sym}_\scr{O}$-module structure of $\free_\scr{O}$ (which corresponds to the identity functor under the above equivalence).\end{remark}

Dually, there is an induced functor
\[L^\rm{coAlg} \colon \coAlg^{\dpw,\nil}_\scr{P}(\scr{D}) \to \coAlg^{\dpw,\nil}_{L^\rm{coOp}(\scr{Q})}(\scr{C})\]
which admits a right adjoint 
\[R^\rm{coAlg} \colon \coAlg^{\dpw,\nil}_{L^\rm{coOp}(\scr{Q})}(\scr{C}) \lra \coAlg^{\dpw,\nil}_\scr{P}(\scr{D}).\]

\subsection{Bar-cobar duality} We first recall bar-cobar duality as formulated by Lurie \cite[5.2]{LurieHA} and then apply this to (co)operads and (co)algebras. Further references include \cite[Section 3.4]{BCN}, \cite[Sections 4.1, 4.2]{BlansBlom}, and \cite{PRY} (comparing Lurie's bar-cobar duality to the more classical one in e.g.~\cite{LodayVallette}).

\subsubsection{Bar-cobar duality in general} Let $\scr{C}$ be a presentable monoidal category so that $1_\scr{C}$ is both initial and terminal; this may be arranged by slicing over and under the monoidal unit, which amounts to passing to (co)augmented (co)algebras. Bar-cobar duality is then the existence of a commutative diagram of adjunctions \cite[Theorem 3.26]{BCN}
\begin{equation}\label{eqn:bar-cobar-general} \begin{tikzcd} \rm{LMod}(\scr{C})  \arrow[shift left=.5ex]{r}{\rm{Bar}}  \dar &[5ex] \arrow[shift left=.5ex]{l}{\rm{Cobar}} \dar \rm{LComod}(\scr{C}) \eqqcolon \rm{LMod}(\scr{C}^\op)^\op \\[-5pt]
	 \arrow[shift left=.5ex]{r}{\rm{Bar}} \Alg(\scr{C}) &[5ex] \rm{coAlg}(\scr{C}) \eqqcolon \Alg(\scr{C}^\op)^\op \arrow[shift left=.5ex]{l}{\rm{Cobar}}
	 \end{tikzcd}\end{equation}
where the vertical arrows are cocartesian fibrations encoding categories of left (co)modules over (co)associative (co)algebras, the bottom horizontal arrows are (co)bar constructions with respect to the tensor product, and the top horizontal arrows are relative (co)bar constructions.

\subsubsection{Naturality in category} \label{sec:bar-cobar-naturality-general} If $R \colon \scr{C} \to \scr{D}$ is a lax monoidal functor between categories satisfying the above hypotheses and admitting a necessarily oplax monoidal left adjoint $L$, then there are functors $R^\rm{Alg} \colon \Alg(\scr{C}) \to \Alg(\scr{D})$ and $L^\rm{coAlg} \colon \Alg(\scr{D}) \to \Alg(\scr{C})$, and similarly for (co)modules. These are related to bar-cobar duality through the existence of a natural transformation of functors $\Alg(\scr{C}) \to \coAlg(\scr{C})$ which will appear in \cite{BlansBlomKupers}:
\[L^\rm{coAlg}\, \Bar_\scr{D}\, R^\rm{Alg} \lra \Bar_\scr{C}.\]
The weaker case that $R$ is symmetric monoidal and preserves geometric realisations, does already appear in the literature, e.g.~\cite[Proposition 4.1.18]{BlansBlom}.

\subsubsection{Bar-cobar duality for (co)operads and (co)algebras}\label{sec:bar-cobar-operads-algebras}
Bar-cobar duality can specialised to symmetric sequences, and there gives a commutative diagram
\begin{equation}\label{eqn:bar-cobar-operads-algebras} \begin{tikzcd} \rm{LMod}^\rm{aug}(\scr{C})  \arrow[shift left=.5ex]{r}{\rm{Bar}}  \dar &[5ex] \arrow[shift left=.5ex]{l}{\rm{Cobar}} \dar \rm{LComod}^\rm{aug}(\scr{C}) \\[-5pt]
	 \arrow[shift left=.5ex]{r}{\rm{Bar}} \rm{Op}^\rm{aug}(\scr{C}) &[5ex] \rm{Coop}^\rm{aug}(\scr{C}) \arrow[shift left=.5ex]{l}{\rm{Cobar}}
	 \end{tikzcd}\end{equation}
where the vertical arrows are cocartesian fibrations encoding categories of left (co-)modules over augmented (co-)operads, the horizontal arrows are (co-)bar constructions. Let $\rm{Op}^\rm{nu}(\scr{C}) \subset \rm{Op}(\scr{C})$ and $\rm{Coop}^\rm{nu}(\scr{C}) \subset \rm{Coop}(\scr{C})$ be the full subcategories of \emph{nonunitary} (co-)operads; alternatively, one can work with nonunitary symmetric sequences defined using nonempty finite sets in place of finite sets. The following is \cite[Theorem 3.4]{HeutsKoszul}:

\begin{theorem}Bar-cobar duality restricts to an adjoint equivalence
\[\begin{tikzcd} \rm{Op}^\rm{nu,aug}(\scr{C}) \arrow[shift left=.5ex]{r}{\rm{Bar}} &[5ex] \rm{Coop}^\rm{nu,aug}(\scr{C}) \arrow[shift left=.5ex]{l}{\rm{Cobar}}. \end{tikzcd}\]
\end{theorem}

\begin{notation}We will abbreviate $\rm{Bar}(\scr{O})$ to $B\scr{O}$ and $\rm{Cobar}(\scr{Q})$ to $\Omega \scr{Q}$.\end{notation}

The adjunction on the vertical fibres of \eqref{eqn:bar-cobar-operads-algebras} over an augmented operad $\scr{O}$ and its image $B\scr{O}$, is related to the cotangent complex as in the following proposition \cite[Corollary 3.30]{BCN}, where one uses \cite{HaugsengMonads} to extract a map of monads from a commutative square whose horizontal morphisms are left adjoints:

\begin{proposition}There is a commutative diagram of left adjoints
	\[\begin{tikzcd} \Alg_\scr{O}(\scr{C}) \arrow{rr}{\rm{Bar}} \arrow{rd}[swap]{\cot_\scr{O}} & & \coAlg^{{\dpw,\nil}}_{B\scr{O}}(\scr{C}) \arrow{ld}{\fgt_{B\scr{O}}} \\[-5pt]
		& \scr{C} & \end{tikzcd}\]
inducing an equivalence of comonads $\cot_\scr{O} \circ \triv_\scr{O} \simeq \rm{Sym}_{B\scr{O}}$.
\end{proposition}

In particular, this exhibits the adjunction $\Bar \dashv \Cobar$ as a lift of the adjunction $\cot_\scr{O} \dashv \triv_\scr{O}$, and we will use the more suggestive notation
\[\begin{tikzcd} \Alg_\scr{O}(\scr{C}) \rar[shift left=.5ex]{\indec^\nil_\scr{O}} &[20pt] \coAlg_{B\scr{O}}^{{\dpw,\nil}}(\scr{C}) \lar[shift left=.5ex]{\prim^\nil_{B\scr{O}}}.\end{tikzcd}\]

\begin{remark}
In the context of homological stability, Randal-Williams first introduced the idea that additional coalgebraic structures on indecomposables should play an important role \cite{RWchromatic}.
\end{remark}

The construction of bar-cobar duality using (co)endomorphism objects gives a different description of $\indec^\nil_\scr{O}$ \cite[Proposition 3.34]{BCN}:

\begin{lemma}\label{lem:bar-via-endomorphism-objects} There is a \emph{Koszul complex}
\[K(\scr{O}) \in \rm{RMod}_\scr{O}(\rm{LComod}_{B\scr{O}}(\rm{SSeq}(\scr{C})) \simeq  \rm{LComod}_{B\scr{O}}(\rm{RMod}_\scr{O}(\rm{SSeq}(\scr{C}))\]
which exhibits $B\scr{O}$ as a coendomorphism object of $1_\scr{C} \in \rm{RMod}_\scr{O}(\rm{SSeq}(\scr{C}))$, and we have an equivalence of functors where $\circ$ denotes the composition product
\[\indec^\nil_{B\scr{O}} \simeq K(\scr{O}) \circ_\scr{O} (-) \colon \Alg_\scr{O}(\scr{C}) \lra \coAlg_{B\scr{O}}^{\dpw,\nil}(\scr{C}).\]
\end{lemma}

\subsection{Koszul duality} Koszul duality concerns the question when the relative bar construction induces an equivalence. We discuss this now, assuming that $\scr{C}$ is a presentable stable symmetric monoidal category. The general case is addressed in \cite[Theorem 2.1]{HeutsKoszul}, where Heuts constructs a category $\coAlg^{\dpw}_\scr{Q}(\scr{C})$ of coalgebras with divided powers \cite[Appendix A]{HeutsKoszul} and a commutative diagram of left adjoints
	\begin{equation}\label{eqn:kd-diag} \begin{tikzcd} \Alg_\scr{O}(\scr{C}) \arrow{rr}{\indec_\scr{O}} \arrow{rd}[swap]{\indec^\nil_\scr{O}} & & \coAlg^{{\dpw}}_{B\scr{O}}(\scr{C}) \arrow{ld} \\[-5pt]
	& \coAlg^{{\dpw,\nil}}_{B\scr{O}}(\scr{C}). & \end{tikzcd}\end{equation}
Moreover, the top adjunction $\indec_\scr{O} \dashv \prim_{B\scr{O}}$ yields an adjoint equivalence when restricted to nilcomplete algebras on the left and conilcomplete coalgebras on the right. 

We will mainly work in the ``connected'' setting, however, where the situation simplifies significantly. Suppose that $\scr{C}$ comes with a left-compatible $t$-structure ($\otimes$ maps $\scr{C}_{\geq 0} \otimes \scr{C}_{\geq 0}$ into $\scr{C}_{\geq 0}$), then we say a nonunitary operad $\scr{O}$ is \emph{connective} if $\scr{O}(n) \in \scr{C}_{\geq 0}$ for all $n \geq 1$ and a (co)algebra is \emph{connected} if its underlying object lies in $\scr{C}_{\geq 1}$. Then \cite[Theorem 14.1, 14.6]{HeutsKoszul} says that:

\begin{theorem}[Heuts] \label{thm:kd-connected} If $\scr{O}$ is a connective nonunitary operad, then all adjunctions in \eqref{eqn:kd-diag} yield adjoint equivalences when we restrict to connected (co)algebras.\end{theorem}

We will need a comparison result between coalgebras with divided powers and the definition of coalgebras used in \cite{LurieHA}, to be proven in \cite{HeutsLand}:

\begin{proposition}[Heuts--Land] \label{prop:lurie-dp-comparison}If $\scr{O}$ is a nonunitary operad in spaces so that each $\scr{O}(n)$ a finite $\fr{S}_n$-space, then there is an equivalence 
\[\Alg_{\scr{O}}(\scr{C}^\op)^\op \simeq \coAlg^{\dpw}_{D\scr{O}}(\scr{C})\]
with target as in \eqref{eqn:kd-diag}, which is the identity on underlying objects.
\end{proposition}

\subsection{Rectification in the dg-setting} \label{sec:rect-dg} We will need to perform several explicit computations and for this it is helpful to be able to use strict (co)algebras over (co)operads and explicit methods to compute the indecomposables. We restrict our attention to 1-category $\rm{Ch}_\bb{Q}$ of (unbounded) rational chain complexes with tensor product. By \cite[1.3.5.15]{LurieHA} inverting the class $W_\rm{qiso}$ of quasi-isomorphisms yield the (unbounded) derived category \cite[1.3.5.8]{LurieHA}
\[\DQ \coloneq \rm{Ch}_\bb{Q}[W_{\rm{qiso}}^{-1}],\]
which is also equivalent to the category $\rm{Mod}_{H\bb{Q}}(\Sp)$ of $H\bb{Q}$-module spectra \cite[7.1.1.16]{LurieHA}.

We will later also need the case of functor categories $\rm{Fun}(\rm{C},\rm{Ch}_\bb{Q})$ for a symmetric monoidal 1-category $\rm{C}$, or monoidal/promonoidal/duoidal/produoidal variants, with Day convolution tensor product; there are no issues in doing so, and we will not comment on this further for the sake of brevity.

The projective model structure on the 1-category $\rm{SSeq}(\rm{Ch}_\bb{Q})$ presents $\rm{SSeq}(\DQ)$, and the composition product on the former yields the composition product on the latter \cite[Corollary 4.23, Corollary 4.32]{BCN}. We can transfer this to model structures on the 1-category $\rm{Op}(\rm{Ch}_\bb{Q})$ of 1-operads presenting $\rm{Op}(\DQ)$ \cite[Corollary 4.11]{HaugsengSS} and on the 1-category $\rm{Alg}_{\rm{O}}(\rm{Ch}_\bb{Q})$ of 1-algebras over a 1-operad presenting $\Alg_\scr{O}(\DQ)$ (where $\scr{O}$ is the operad associated to the 1-operad $\rm{O}$) \cite[Theorem 4.4]{GetzlerJones} \cite[Theorem 4.10]{HaugsengSS}. That is, we have equivalences
\[\rm{Op}(\rm{Ch}_\bb{Q})[W^{-1}] \overset{\simeq}\lra \rm{Op}(\DQ) \qquad \text{and} \qquad \Alg_\rm{O}(\rm{Ch}_\bb{Q})[W^{-1}] \overset{\simeq}\lra \Alg_\scr{O}(\DQ),\]
where we have used that in the setting of rational chain complexes any operad is $\Sigma$-cofibrant. There is a similar model structure on $\rm{Coop}(\rm{Ch}_\bb{Q})$ \cite[Theorem 2.4.1]{AubryChataur} (take $i=1$) and the chain-level operadic bar construction \cite[Section 2.1]{GetzlerJones} \cite[Section 6.5]{LodayVallette} is the left adjoint in a Quillen equivalence \cite[p.~3]{AubryChataur}
\[B \colon \rm{Op}(\rm{Ch}_\bb{Q}) \lra \rm{Coop}(\rm{Ch}_\bb{Q}).\]
Finally, while the 1-category $\rm{coAlg}^{\dpw,\nil}_\rm{Q}(\DQ)$ of 1-coalgebras over a 1-cooperad $\rm{Q}$ has a class of ``weak equivalences'' given by those maps that are quasi-isomorphisms on underlying objects; these are known to be the weak equivalence of a model structure only when we restrict to 1-coalgebras over a \emph{connected} 1-cooperad \cite[Theorem 4.7]{GetzlerJones} \cite[Theorem 3.2.3]{AubryChataur}; see \cite{PRY} for more general results without connectivity hypotheses. We highlight one feature of the rational setting: the norm maps $(X^{\otimes n})_{\mathfrak{S}_n} \to (X^{\otimes n})^{\mathfrak{S}_n}$ are equivalences, so a divided power structure on a coalgebra is no additional data:

\begin{notation}When working in a rational setting, we drop the superscript $\dpw$, unless there is a chance for confusion.
\end{notation}

The chain-level operadic bar construction for operads not only yields the aforementioned cooperad $B\rm{O} \in \rm{Coop}(\rm{Ch}_\bb{Q})$  but also a $B\rm{O}$-coalgebra $B^\rm{O} \bf{A}$ for each $\rm{O}$-algebra $\bf{A}$ \cite[Section 2.3]{GetzlerJones} \cite[Section 11.2]{LodayVallette}. The proof of \cite[Theorem 4.42]{BCN}, in the easier standard setting rather than pro-coherent one, yields a commutative diagram
\[\begin{tikzcd} \rm{Alg}_{\rm{O}}(\rm{Ch}_\bb{Q})[W^{-1}] \rar{B^\rm{O}} \dar{\simeq} & \rm{coAlg}^{\nil}_{B\rm{O}}(\rm{Ch}_\bb{Q})[W^{-1}] \dar \\[-5pt]
\Alg_\scr{O}(\DQ) \rar{\indec^\nil_\scr{O}} & \coAlg_{B\scr{O}}^{\nil}(\DQ) \end{tikzcd}\]
by using $B^\rm{O} \bb{Q}$ to exhibit $B\rm{O}$ as the coendomorphism object of the monoidal unit $\bb{Q}$ as in \cref{lem:bar-via-endomorphism-objects}; this implies that $B\scr{O}$ is the operad associated to $B\rm{O}$ and that the chain-level operadic bar construction models nil-indecomposables. Restricting to connected (co)algebras, we obtain:

\begin{proposition}There is a commutative square of equivalences
	\[\begin{tikzcd} \rm{Alg}_{\rm{O}}(\rm{Ch}_\bb{Q})[W^{-1}]_{\geq 1} \rar{B^\rm{O}}[swap]{\simeq} \dar{\simeq} & \rm{coAlg}^{\nil}_{B\rm{O}}(\rm{Ch}_\bb{Q})[W^{-1}]_{\geq 1} \dar{\simeq} \\[-5pt]
		\Alg_\scr{O}(\DQ)_{\geq 1} \rar{\indec^\nil_\scr{O}}[swap]{\simeq} & \coAlg_{B\scr{O}}^{\nil}(\DQ)_{\geq 1} \end{tikzcd}\]
where the subscripts $(-)_{\geq 1}$ indicate we restrict to the full subcategory where the underlying objects are connected.
\end{proposition}

Since the operadic cobar complex $\Omega^{B\rm{O}}$ is the right Quillen adjoint to $B^\rm{O}$, it follows that on connected coalgebras $\prim^\nil_{B\scr{O}}$ can be computed using $\Omega^{B\rm{O}}$.

\section{$E_k$-algebraic foundations} In this last appendix we specialise the theory of (co)operads, (co)algebras, and Koszul duality to the case of the $E_k$-operads. This has some special features, most importantly using the Dunn--Lurie additivity theorem.

\subsection{The $E_k$-operads} There are two variants of the $E_k$-operads for $k = 1, 2,\ldots,\infty$, cf.~\cite[Definition 12.1, 12.2]{GKRW18}. Their definition uses the notion of a rectilinear embedding between cubes, which is a map $I^k \to I^k$ of the form
\[(x_1,\ldots,x_k) \longmapsto ((b_1-a_1)x_1+a_1,\ldots,(b_k-a_k)x_k+a_k)\]
for $0 \leq a_k < b_k \leq 1$.

\begin{definition}\label{def:ek} Let $1 \leq k < \infty$.
    \begin{itemize}
        \item The \emph{unitary $E_k$-operad} $E_k^\rm{u}$ is obtained as the operadic nerve of the operad whose space of $r$-ary operations is given by the space 
        \[E_k^\rm{u}(r) \coloneq \rm{Emb}^\rm{rect}(\sqcup_r I^k,I^k)\]
        of $r$-tuples of rectilinear embeddings whose images have disjoint interior and where operadic composition is induced by composition of rectilinear embeddings.
        \item The \emph{nonunitary $E_k$-operad} $E_k^\rm{nu}$ is obtained by replacing the $0$-ary operations in $E_k^\rm{u}$, given by a single point, with the empty set.
    \end{itemize}
\end{definition}

Note that $(-) \times \rm{id}_I$ induces maps of operads $E_k^\rm{u} \to E_{k+1}^\rm{u}$ and $E_k^\rm{nu} \to E_{k+1}^\rm{nu}$.

\begin{definition}\label{def:einfty} We define the \emph{unitary $E_\infty$-operad} and \emph{nonunitary $E_\infty$-operad} as
\[E_\infty^\rm{u} \coloneq \underset{k \to \infty}{\rm{colim}}\,E_k^\rm{u} \qquad \text{and} \qquad E_\infty^\rm{nu} \coloneq \underset{k \to \infty}{\rm{colim}}\,E_k^\rm{nu}.\]
\end{definition}

\subsection{The case $k<\infty$} Passing to the category spectra by implicitly taking suspension spectra, there are for $1 \leq k < \infty$ equivalences 
\[\gamma_k \colon BE^\rm{nu}_k \overset{\simeq}\lra s^k DE^\rm{nu}_k,\]
in $\rm{Coop}(\Sp)$, where $s(-)$ is (co)operadic suspension and $D(-)$ is objectwise Spanier--Whitehead dual, see \cref{sec:operadic-suspension-linear-duals}.

\subsubsection{Koszul duality and bar-cobar duality} \label{subsection koszul bar cobar}
Combining additivity with bar-cobar duality for associative (co)algebras, Lurie established bar-cobar duality for $E_k$-(co)algebras. Suppose that $\scr{C}$ is a symmetric monoidal category and define 
\[\Alg^\aug_{E_k^u}(\scr{C}) \coloneq \Alg_{E_k^u}(\scr{C})_{/1_\scr{C}} \quad \text{and} \quad \coAlg^\aug_{E_k^u}(\scr{C}) \coloneq \left(\Alg_{E_k^u}(\scr{C}^\op)_{/1_\scr{C}}\right)^\op.\]
Assuming $\scr{C}$ has geometric realisations and totalisations, iterating bar-cobar duality as in \cite[5.2.3]{LurieHA} yields an adjunction
	\[\begin{tikzcd} \Alg^\aug_{E_k^u}(\scr{C}) \arrow[shift left=2pt]{r}{\Bar^k} &[20pt] \coAlg^\aug_{E_k^u}(\scr{C}) \arrow[shift left=2pt]{l}{\Cobar^k}.\end{tikzcd}\]
We remark that the existence of the functor $\Bar^k$ only requires $\scr{C}$ has geometric realisations, and the existence of the functor $\Cobar^k$ only requires it has totalisations.

\begin{example}\label{exam:e1-coproduct} For $k=1$, for an augmented $E^\rm{u}_1$-algebra $\epsilon \colon \bf{A} \to 1 = 1_\scr{C}$, we have $\Bar(\bf{A}) \simeq 1 \otimes_\bf{A} 1$ \cite[5.2.2.3]{LurieHA}. By \cite[p.~826]{LurieHA} the coproduct is given by
\[1 \otimes_\bf{A} 1 \simeq 1 \otimes_\bf{A}  \bf{A} \otimes_\bf{A} 1 \lra 1 \otimes_\bf{A} 1 \otimes_\bf{A} 1 \simeq (1 \otimes_\bf{A} 1) \otimes (1 \otimes_\bf{A} 1)\]
with map induced by $\epsilon$ on the middle term. Let us explain why this is the case using \cref{lem:bar-via-endomorphism-objects}. By \cite[Proposition 3.34]{BCN} the bar construction is the coendomorphism object of $1$, considered as a right $\bf{A}$-module through $\epsilon$: Recalling that $\epsilon^* \colon \scr{C} \to \rm{RMod}_\bf{A}(\scr{C})$ has a left adjoint $\epsilon_! \simeq (-) \otimes_\bf{A} 1$, the unit map $\eta \colon \rm{id} \to \epsilon^* \epsilon_!$ induces a map $\rho \coloneq \eta \epsilon^* \colon \epsilon^* 1 \to \epsilon^*\epsilon_! \epsilon^* 1$ of right $\bf{A}$-modules and this has the property that the map
\[\Map_\scr{C}(\Bar(\bf{A}),Y) \lra \Map_{\rm{RMod}_\bf{A}(\scr{C})}(\epsilon^* \Bar(\bf{A}),\epsilon^* Y) \lra \Map_{\rm{RMod}_\bf{A}(\scr{C})}(\epsilon^* 1,\epsilon^* Y)\]
given by applying $\epsilon^*$ and precomposing with $\rho$, is an equivalence. To get the coproduct from this universal property, consider $Y = \epsilon_! \epsilon^* \epsilon_! \epsilon^* 1$ with the map of right $\bf{A}$-modules
\[\epsilon^* 1 \xrightarrow{\rho} \epsilon^* \epsilon_! \epsilon^* 1 \xrightarrow{\epsilon^* \epsilon_! \rho} \epsilon^* \epsilon_! \epsilon^* \epsilon_! \epsilon^* 1\]
and recognise it arises under the above equivalence from the map $\epsilon_! \eta \epsilon^* \colon \epsilon_! \epsilon^* 1 \to \epsilon_! \epsilon^* \epsilon_! \epsilon^* 1$. Unwinding the definitions, this gives the desired coproduct.
\end{example}

The relationship between these constructions and Koszul duality is due to Heuts--Land \cite{HeutsLand}. Firstly, if $\scr{C}$ is a stable presentable symmetric monoidal category $\scr{C}$, \cref{prop:lurie-dp-comparison} provides an identification 
\[\coAlg^\aug_{E_k^u}(\scr{C}) \simeq \coAlg^{\dpw,\rm{aug}}_{E_k^\rm{u}}(\scr{C}).\] 
In terms of this identification, the following improves on \cite[Theorem 13.7]{GKRW18}, which shows the diagram commutes after composing with the forgetful functor $\coAlg_{s^{k} DE^\rm{nu}_k}(\scr{C}) \to \scr{C}$:

\begin{theorem}[Heuts--Land] \label{thm:indec-is-bar} If $\scr{C}$ is a stable presentable symmetric monoidal category, then there is a commutative square
\[\begin{tikzcd}
	\Alg^\aug_{E^\rm{u}_k}(\scr{C}) \rar{\indec_{E_k}(I(-))} \dar[swap]{\mathrm{Bar}^k(-)} &[20pt] \coAlg^{\dpw}_{BE^\rm{nu}_k}(\scr{C}) \dar{\gamma^\vee_!}[swap]{\simeq} \\[-5pt]
	\coAlg^\aug_{E^\rm{u}_k}(\scr{C}) \rar{\Sigma^{-k}I(-)}[swap]{\simeq} & \coAlg^\dpw_{s^{k} DE^\rm{nu}_k}(\scr{C})\end{tikzcd}\]
where $I$ denotes the augmentatio ideal and $s$ the (co)operadic suspension.
\end{theorem}

\subsubsection{Iterating bar constructions} \label{sec:iterated-bar-constructions} We will want to iterate bar constructions but to do so, we will need to perform these constructions in categories of (co)algebras.

\begin{lemma}\label{lem:bar-underlying} Let $\scr{C}$ be a symmetric monoidal category whose tensor product preserves sifted colimits in each entry and $\scr{O}$ be an operad in spaces.
    \begin{enumerate}[(i)]
        \item $\Alg_\scr{O}(\scr{C})$ admits the structure of symmetric monoidal category so that the forgetful functor $\fgt_\scr{O} \colon \Alg_\scr{O}(\scr{C}) \to \scr{C}$ is symmetric monoidal and creates sifted colimits.
        \item $\coAlg_\scr{O}(\scr{C})$ admits the structure of symmetric monoidal category so that the forgetful functor $\fgt_\scr{O} \colon \coAlg_\scr{O}(\scr{C}) \to \scr{C}$ is symmetric monoidal and creates colimits.
    \end{enumerate}
\end{lemma}

\begin{proof} The first part uses \cite[3.2.4.4]{LurieHA} and \cite[3.2.3.2]{LurieHA}, using that sifted colimits are created by the forgetful functor. For the second part, write $\coAlg_\scr{O}(\scr{C})^\op = \Alg_{\scr{O}}(\scr{C}^\op)$, and use \cite[3.2.4.4]{LurieHA} and \cite[3.2.2.5]{LurieHA}.
\end{proof}

\begin{remark}The underlying category of $\Alg_\scr{O}(\scr{C})$ is in fact presentable if $\scr{C}$ is presentable symmetric monoidal \cite[3.2.3.5]{LurieHA}.\end{remark}

We want to combine bar-cobar duality with the additivity theorem 
\[E^u_{k+k'} \simeq E^u_k \otimes_{\rm{BV}} E^u_{k'}\]
writing the left side as a tensor product of operads \cite[5.1.2.2]{LurieHA}. Using the adjunction $-\otimes E^u_k \dashv \Alg_{E^u_k}(-)$ that is the defining property of this tensor product, we get equivalences 
\begin{equation}\label{eqn:additivity-algebras}\Alg_{E^u_{k+k'}}(\scr{C}) \simeq \Alg_{E^u_k \otimes_\rm{BV} E^u_{k'}}(\scr{C}) \simeq \Alg_{E^u_k}(\Alg_{E^u_{k'}}(\scr{C})).\end{equation} 

\begin{proposition}\label{prop:hopf-algebra-ek-indecomposables} Suppose $\scr{C}$ is a  presentable stable symmetric monoidal category. If $\bf{R}^+$ is an augmented $E^\rm{u}_{k+k'}$-algebra in $\scr{C}$ with augmentation ideal $\bf{R}$, then
    \begin{enumerate}[(i)]
        \item $(\Sigma^k \indec_{E^\rm{nu}_k}(\bf{R}))^+$ lifts to an augmented $E^\rm{u}_k$-coalgebra in augmented $E^\rm{u}_{k'}$-algebras, and
        \item the same is true for $\Sigma^k \indec^\nil_{E^\rm{nu}_k}(\bf{R})^+$ as long as $\bf{R}$ is connected.
    \end{enumerate}
\end{proposition}

\begin{proof}For now it suffices that $\scr{C}$ is a presentable symmetric monoidal category. Slicing over and under $1_\scr{C}$ in \eqref{eqn:additivity-algebras} gives 
\[\Alg^\aug_{E^u_{k+k'}}(\scr{C}) \simeq \Alg^\aug_{E^u_k}(\Alg^\aug_{E^u_{k'}}(\scr{C})),\]
where the outer $(-)^\aug$ on the right side is tautological since the monoidal unit is terminal. Now we apply iterated bar construction to get the top map in the commutative diagram \cite[5.2.3.12]{LurieHA}
\[\begin{tikzcd}  \Alg^\aug_{E^u_k}(\Alg^\aug_{E^u_{k'}}(\scr{D})) \dar[swap]{(\fgt_{E^u_{k'}})_!} \rar{\mathrm{Bar}^k} &  \coAlg^\aug_{E^u_k}(\Alg^\aug_{E^u_{k'}}(\scr{D})) \dar{(\fgt_{E^u_{k'}})_!} \\[-5pt]
\Alg^\aug_{E^u_k}(\scr{D}) \rar{\mathrm{Bar}^k} &  \coAlg^\aug_{E^u_k}(\scr{D}) \end{tikzcd}\]
where the vertical maps are induced by the forgetful functor, which is symmetric monoidal and creates sifted colimits. 

Specialising a stable presentable symmetric monoidal category $\scr{C}$ then we apply \cref{thm:indec-is-bar} to identify the $E_k$-coalgebra $\Sigma^k \indec_{E_k} \bf{R}$ with $\rm{Bar}^k \bf{R}$, and the fact that it is in the image of the right vertical map provides the desired lift.
\end{proof}

\begin{lemma}\label{lem:bialgebra-structure} Suppose that $\bf{R}^+$ is an augmented $E^\rm{u}_{k+k'}$-algebra in $\Fun(\scr{C},\DQ)$ for $k,k' \geq 1$ with augmentation ideal $\bf{R}$. Then $H_{*,*}((\Sigma^k \indec_{E^\rm{nu}_k} \bf{R})^+)$ admits the structure of a bigraded bialgebra which is:
\begin{enumerate}[(i)]
    \item cocommutative if $k \geq 2$,
    \item commutative if $k' \geq 2$,
    \item has connected augmentation ideal if $\bf{R}$ is connected.
\end{enumerate}
\end{lemma}

\begin{proof}Applying \cref{prop:hopf-algebra-ek-indecomposables} to $\scr{C} = \Fun(\bb{N},\DQ)$ we get a lift of $(\Sigma^k \indec_{E_k^\rm{nu}}(\bf{R}))^+$ to an object of 
\[\coAlg^\aug_{E_k^\rm{u}}(\Alg^\aug_{\rm{E}_{k'}^\rm{u}}(\Fun(\bb{N},\DQ))).\] 
First supposing $k=1=k'$, we use that taking homology with rational coefficients is symmetric monoidal and $H_*(E_1^\rm{u}) \cong \rm{As}^\rm{u}$ to obtain a functor 
\[\coAlg^\aug_{E_k^\rm{u}}(\Alg^\aug_{\rm{E}_{k'}^\rm{u}}(\Fun(\bb{N},\DQ))) \lra \coAlg^\aug_{\rm{coAssoc}}(\Alg^\aug_{\rm{Assoc}}(\Fun(\bb{N},\rm{GrMod}_\bb{Q}))),\]
and the target is the category of bigraded bialgebras. To prove (i) we rather restrict along the map $\rm{Com} \to H_*(E_k^\rm{u}) \cong \rm{Pois}_{k-1}$, and similarly for (ii). Part (iii) follows by iterated bar spectral sequences.
\end{proof}

\subsection{The case $k=1$, rationally} \label{sec:ass-algebras} There is an equivalence of cooperads
\[\gamma_1 \colon BE_1^\rm{nu} \overset{\simeq}\lra s DE_1^\rm{nu}.\]
Upon passing to rational chains, we get that the $E_1$-operad is equivalent to the classical associative operad: there are equivalences
\[\rm{As}^\rm{u} \overset{\simeq}\lra E_1^\rm{u} \qquad \text{and} \qquad \rm{As}^\rm{nu} \overset{\simeq}\lra E_1^\rm{nu}.\]
Moreover, $\gamma_1$ has an inverse given by the Koszul duality equivalence 
\[s\,\rm{coAs}^\rm{nu} \overset{\simeq}\lra B\rm{As}^\rm{nu} \simeq BE_1^\rm{nu}\]
of \cite[Section 9.3]{LodayVallette}, where $\rm{coAs}^\rm{nu}$ is the linear dual of $\rm{As}^\rm{nu}$. 

We will now explain how, through the rectification results of \cref{sec:rect-dg}, we can perform Koszul duality computations for associative algebras using explicit chain complexes. More precisely, we work in the model category $\rm{C} = \rm{Fun}(\rm{A},\rm{Ch}_\bb{Q})$ of functors from a symmetric monoidal 1-category $\rm{A}$ to rational chain complexes, with associated $\infty$-category $\scr{C}$ equivalent to $\Fun(\rm{A},\DQ)$. 

We start with an explicit implementation of the indecomposables functor via bar complexes when we restrict to connected algebras: \cref{sec:rect-dg} then provides a commutative diagram of equivalences (recall our convention to drop the superscript $\dpw$ in rational settings)
	\[\begin{tikzcd} \rm{Alg}_{\rm{As}^\rm{nu}}(\rm{C})[W^{-1}]_{\geq 1} \rar{B^\rm{As}}[swap]{\simeq} \dar{\simeq} & \rm{coAlg}^{\nil}_{B\rm{As}^\rm{nu}}(\rm{C})[W^{-1}]_{\geq 1} \dar{\simeq} & \rm{coAlg}^{\nil}_{s\rm{coAs}^\rm{nu}}(\rm{C})[W^{-1}]_{\geq 1} \lar[swap]{\simeq} \dar{\simeq}\\
		\Alg_{E^\rm{nu}_1}(\scr{C})_{\geq 1} \rar{\indec^\nil_{E^\rm{nu}_1}}[swap]{\simeq} & \coAlg_{BE^\rm{nu}_1}^{\nil}(\scr{C})_{\geq 1} \rar & \lar[swap]{\simeq} \coAlg_{sDE^\rm{nu}_1}^{\nil}(\scr{C})_{\geq 1}. \end{tikzcd}\]
By the identification of the operadic bar construction for a Koszul operad in terms of a twisting morphism \cite[Section 11.2]{LodayVallette}, under the top-right equivalence the map $B^\rm{As}$ corresponds to the following classical bar construction \cite[Section 2.2.1]{LodayVallette}, which we will  denote the same:

\begin{definition}\label{def:bar-as} For a nonunital dg-algebra $\bf{A} \in \Alg_\rm{As^{nu}}(\Fun(\rm{A},\rm{Ch}_\bb{Q}))$ we define a shifted dg-coalgebra
\[B^\rm{As}(\bf{A}) \coloneq \Sigma^{-1}(\rm{coAs}^\rm{nu} \circ \Sigma,d_\bf{A}+d_B) \in \coAlg^\nil_{s\,\rm{coAs}^\rm{nu}}\Fun(\rm{A},\rm{Ch}_\bb{Q})\]
where $d_\bf{A}$ is the unique graded coderivation induced by the differential $d \colon \bf{A} \to \Sigma \bf{A}$ of $\bf{A}$ and $d_B$ is the unique graded coderivation induced by the shifted multiplication map
\[\Sigma \bf{A} \otimes \Sigma \bf{A} \cong \Sigma^2(\bf{A}^{\otimes 2}) \xrightarrow{\Sigma^2 m} \Sigma^2 \bf{A}.\]
\end{definition}
We refer to $d_\bf{A}$ as the \emph{internal differential} and $d_B$ as the \emph{bar differential}. Observe that $B^\rm{As}(\bf{A})$ is, up to a suspension, the graded conilpotent tensor coalgebra $T^c(\Sigma \bf{A}) \cong \bigoplus_{n \geq 1} (\Sigma \bf{A})^{\otimes n}$ with differentials for $\ol{a}_0,\ldots,\ol{a}_p \in \Sigma \bf{A}$ (where we abbreviate $\ol{a} \coloneq sa$) given by
\[d_\bf{A}(\ol{a}_0 \otimes \cdots \otimes \ol{a}_p) = \sum_{i=0}^p (-1)^{|\ol{a}_0|+\cdots+|\ol{a}_{i-1}|} \ol{a}_0 \otimes \cdots \otimes \ol{d(a_i)} \otimes \cdots \ol{a}_p,\]
\[d_B(\ol{a}_0 \otimes \cdots \otimes \ol{a}_p) = \sum_{i=0}^{p-1} (-1)^{|\ol{a}_0|+\cdots+|\ol{a}_{i-1}|}\ol{a}_0 \otimes \cdots \otimes \ol{a_ia_{i+1}} \otimes \cdots \otimes \ol{a}_p\]
and coproduct given by \emph{deconcatenation coproduct}, whose reduced version is
\[\overline{\Delta}(\ol{a}_0 \otimes \cdots \otimes \ol{a}_p) = \sum_{i=1}^{p} (\ol{a}_0 \otimes \cdots \otimes \ol{a}_{i-1}) \otimes (\ol{a}_i \otimes \cdots \otimes \ol{a}_p).\]

If $\bf{A}$ is commutative there is a graded-commutative product
\[B^\rm{As}(\bf{A}) \otimes B^\rm{As}(\bf{A}) \xrightarrow{\rm{EZ}} B^\rm{As}(\bf{A} \otimes \bf{A}) \lra B^\rm{As}(\bf{A})\]
where the left map is the Eilenberg--Zilber equivalence and the right map is induced by the multiplication $m \colon \bf{A} \otimes \bf{A} \to \bf{A}$, which is a map of associative algebras \cite[4.2.6]{LodayCyclic}. This yields the \emph{shuffle product} formula
\begin{align*}&\mu\big((\ol{a}_1 \otimes \cdots \otimes \ol{a}_p) \otimes (\ol{a}_{p+1} \otimes \cdots \otimes \ol{a}_{p+q}))\big) \\
&\qquad = \sum_{\sigma \in \rm{sh}_{p,q}} (-1)^{\sigma(a)} sa_{\sigma(1)} \otimes \cdots \otimes \ol{a}_{\sigma(p)} \otimes \ol{a}_{\sigma(p+1)} \otimes \cdots \otimes \ol{a}_{\sigma(p+q)}\end{align*}
where the sign $(-1)^{\sigma(a)}$ the product of $(-1)^{|\ol{a}_{i}||\ol{a}_{j}|}$ for all $1 \leq i \leq p$ and $p+1 \leq j \leq p+q$ so that $\sigma(q)<\sigma(p)$. 

If $\bf{A}$ is connected in addition to being commutative, then by \cref{prop:hopf-algebra-ek-indecomposables} the counital $E_1$-coalgebra 
\[B^\rm{As}(\bf{A})^+ \simeq (\Sigma\, \indec^\nil_{E^\rm{nu}_1}(\bf{A}))^+ \simeq \Bar(\bf{A})\] lifts to $\smash{\coAlg_{E_1^\rm{u}}(\Alg_{E_\infty^\rm{u}}(\scr{C}))}$. By the construction of the pairing underlying the additivity theorem in \cite[5.1.2.1]{LurieHA}, the underlying multiplication on the bar construction is given by
\[\Bar(\bf{A}) \otimes \Bar(\bf{A}) \overset{\simeq}\lra \Bar(\bf{A} \otimes \bf{A}) \lra \Bar(\bf{A})\]
where the left map uses that geometric realisations are sifted colimits and that the tensor product commutes with colimits in each entry and the right map is induced by the multiplication map $m \colon \bf{A} \otimes \bf{A} \to \bf{A}$. To see it agrees with the shuffle product, recall that the Eilenberg--Zilber map intertwines the diagonal of the tensor product of simplicial vector spaces and the tensor product of chain complexes.

The right adjoint of the bar construction is $B^\rm{As}$ is the cobar construction $\Omega^\rm{coAs}$ \cite[Section 2.2.5]{LodayVallette}.

\begin{definition}\label{def:cobar-as} For a nonunital coassociative dg-coalgebra $\bf{C} \in \Alg_{\rm{coAs}^\rm{nu}}(\Fun(\rm{A},\rm{Ch}_\bb{Q}))$ we define the \emph{coassociative cobar construction} as the shifted nonunital associative dg-algebra
\[\Omega^\rm{coAs}(\bf{C}) \coloneq \Sigma(\rm{As}^\rm{nu} \circ \Sigma^{-1} \bf{C},d_\bf{C}+d_\Omega)\]
where $d_\bf{C}$ is the unique graded derivation induced by the differential $\bf{C} \to \Sigma \bf{C}$ of $L$ and $d_\Omega$ is the unique graded derivation induced by the shifted coproduct
\[\Sigma^{-1} \bf{C} \xrightarrow{\Sigma^{-1} \Delta} \Sigma^{-1} (\bf{C} \otimes \bf{C}) \cong \Sigma(\rm{As}(2) \otimes_{\mathfrak{S}_2} (\Sigma^{-1} \bf{C})^{\otimes 2}).\]
\end{definition}

We refer to $d_\bf{C}$ as the \emph{internal differential} and $d_\Omega$ as the \emph{cobar differential}.

\begin{remark}It is occasionally more convenient to use a variant that is an augmented unital dg-associative algebra, by using $\rm{As}^{u}$ instead of $\rm{As}^\rm{nu}$ and with augmentation induced by the one of $\rm{As}^\rm{u}$.
\end{remark}

That is, $\Omega^\rm{coAs}(\bf{C})$ is, up to a desuspension, given by the graded tensor product $T(\bf{C}) = \bigoplus_{n \geq 1} (\Sigma^{-1} \bf{C})^{\otimes n}$ with differentials for $\ul{c}_0,\ldots,\ul{c}_p \in \Sigma^{-1} \bf{C}$ (where we abbreviate $\ul{c} = s^{-1}c$)
\[d_\bf{C}(\ul{c}_0 \otimes \cdots \otimes \ul{c}_p) = \sum_{i=0}^p (-1)^{|\ul{c}_0|+\cdots+|\ul{c}_{i-1}|} \ul{c}_0 \otimes \cdots \otimes \ul{d(c_i)} \otimes \cdots \otimes \ul{c}_p,\]
\[d_\Omega(\ul{c}_0 \otimes \cdots \otimes \ul{c}_p) = \sum_{i=0}^p (-1)^{|\ul{c}_0|+\cdots+|\ul{c}_{i-1}|} \ul{c}_0 \otimes \cdots \otimes \ul{\Delta(c_i)} \otimes \cdots \otimes \ul{c}_p.\]

Koszul duality takes the following concrete form, in terms of the unit and counit of the classical adjunctions between the bar and cobar constructions \cite[Corollary 2.3.4]{LodayVallette}: given a connected dg-associative algebra $\bf{A}$ and connected dg-associative coalgebras there are quasi-isomorphisms
\[\epsilon \colon \Omega^\rm{coAs} B^\rm{As} \bf{A} \overset{\simeq}\lra \bf{A} \qquad \text{and} \qquad \eta \colon \bf{C} \overset{\simeq}\lra B^\rm{As} \Omega^\rm{coAs} \bf{C}.\]
These have inverses on the level of chain complexes, \emph{not} compatible with (co)algebra structures, given by the inclusion of or projection onto certain terms.

\subsection{The case $k=\infty$, rationally} \label{sec:comm-colie-algebras} 
For $k=\infty$, we rather have
\[\gamma_\infty \colon BE^\rm{nu}_\infty \overset{\simeq}\lra D\rm{Lie}_\rm{Sp} \coloneq \colim_{k \to \infty} s^k DE^\rm{nu}_k,\]
where the map $s^k DE^\rm{nu}_k \to s^{k+1} DE^\rm{nu}_{k+1}$ arises a priori through Koszul duality, but is explicitly described in \cite[Section 7]{ChingSalvatore}. Its target is the Spanier--Whitehead dual of the spectral Lie operad, by definition of the latter. This is named such because, upon passing to rational chains, $\gamma_\infty$ has an inverse given by the classical Koszul duality equivalence $s\, \rm{coLie} \smash{\overset{\simeq}\lra} B \rm{Com}^\rm{nu}$ of \cite[Section 13.1.5]{LodayVallette}.

We will now explain how, as in the associative case, the rectification results of \cref{sec:rect-dg} allow us to perform Koszul duality computations using explicit chain complexes. Once more, we work in the model category $\rm{C} = \Fun(A,\rm{Ch}_\bb{Q})$ of functors from a symmetric monoidal $1$-category $\rm{A}$ to rational chain complexes, with associated $\infty$-category $\scr{C}$. Then there are explicit implementations of the indecomposables functor $\indec^\nil_{E^\rm{nu}_\infty}$ via (co)bar complexes when we restrict to connected (co)algebras: we get a commutative diagram of equivalences
	\[\begin{tikzcd} \rm{Alg}_{\rm{Com}^\rm{nu}}(\rm{C})[W^{-1}]_{\geq 1} \rar{B^\rm{Com}}[swap]{\simeq} \dar{\simeq} & \rm{coAlg}^{\nil}_{B\rm{Com}^\rm{nu}}(\rm{C})[W^{-1}]_{\geq 1} \dar{\simeq} & \rm{coAlg}^{\nil}_{s\rm{coLie}}(\rm{C})[W^{-1}]_{\geq 1} \lar[swap]{\simeq} \dar{\simeq}\\
		\Alg_{E_\infty^\rm{nu}}(\scr{C})_{\geq 1} \rar{\indec^\nil_\scr{O}}[swap]{\simeq} & \coAlg_{BE^\rm{nu}_\infty}^{\nil}(\scr{C})_{\geq 1} \rar & \lar[swap]{\simeq} \coAlg_{sDE^\rm{nu}_\infty}^{\nil}(\scr{C})_{\geq 1}.\end{tikzcd}\]
The top-right corner is given by objects that are, up to a suspension, given by conilpotent dg-Lie coalgebras $\bf{L}$: explicitly, this is an object $\bf{L} \in \rm{C}$ with a map
\[\bf{L} \lra \coLie \circ \bf{L} = \bigoplus_{n \geq 1} (\coLie(n) \otimes \bf{L}^{\otimes n})_{\mathfrak{S}_n},\]
satisfying counitality and coassociativity axioms. 

\begin{remark}A conilpotent dg-Lie coalgebra is an instance of a more classical notion of dg-Lie coalgebra used in \cite{Michaelis,LodayVallette,CharltonRadchenkoRudenko}. In these references, it is given rather by an object $\bf{L}' \in \rm{C}$ with a map
\[\delta \colon \bf{L}' \lra \bf{L}' \otimes \bf{L}'\]
satisfying $\tau \circ \delta = - \delta$ and $(1+\eta+\eta^2) \circ (1 \otimes \delta) \circ \delta$ where $\tau \colon \bf{L}' \otimes \bf{L}' \to \bf{L}' \otimes \bf{L}'$ is induced by the transposition $(1\,2)$ and $\eta \colon \bf{L}' \otimes \bf{L}' \otimes \bf{L}' \to \bf{L} \otimes \bf{L}' \otimes \bf{L}'$ is induced by the 3-cycle $(2\,3\,1)$. Using norm map isomorphism $(\coLie(n) \otimes (\bf{L}')^{\otimes n})_{\mathfrak{S}_n} \to (\coLie(n) \otimes (\bf{L}')^{\otimes n})^{\mathfrak{S}_n}$ identifying quotients with subobjects, conilpotent dg-Lie coalgebras can be identified the subcategory of those classical dg-Lie coalgebra so that $\delta$ is conilpotent.\end{remark}

By the identification of the operadic bar construction for a Koszul operad in terms of a twisting morphism, under the top-right equivalence the map $B^\rm{Com}$ corresponds up to a shift to the classical \emph{Harrison homology complex} \cite[Section 13.1.10]{LodayVallette}:

\begin{definition}\label{def:bar-comm} For a nonunital dg-commutative algebra we define the \emph{commutative bar construction} as the shifted conilpotent dg-Lie coalgebra
\[B^\rm{Com}(\bf{A}) \coloneq \Sigma^{-1}(\coLie \circ \Sigma \bf{A},d_\bf{A}+d_B)\]
where $d_\bf{A}$ is the unique graded coderivation induced by the differential $d \colon \bf{A} \to \Sigma \bf{A}$ of $\bf{A}$ and $d_B$ is the unique graded coderivation induced by the shifted product 
\[\coLie(2) \otimes_{S_2} (\Sigma \bf{A}) \cong \Sigma^2 S^2(\bf{A}) \xrightarrow{\Sigma^2 m} \Sigma^2 \bf{A}.\]
\end{definition}

We refer to $d_\bf{A}$ as the \emph{internal differential} and $d_B$ as the \emph{bar differential}. Its right adjoint $\Omega^\coLie$ models on connected Lie coalgebras the functor $\prim_{sB\rm{Com}^\rm{nu}}^\nil$ and corresponds to the classical \emph{Chevalley--Eilenberg complex}, dual to the one for Lie algebras \cite[13.2.8]{LodayVallette}:

\begin{definition}\label{def:cobar-colie} For a dg-Lie coalgebra $\bf{L}$ we define the \emph{coLie cobar construction} as the shifted nonunital dg-commutative algebra
\[\Omega^\coLie(\bf{L}) \coloneq \Sigma(\rm{Com}^\rm{nu} \circ \Sigma^{-1} \bf{L},d_\bf{L}+d_\Omega)\]
where $d_\bf{L}$ is the unique graded derivation induced by the differential $d \colon \bf{L} \to \Sigma \bf{L}$ of $\bf{L}$, and $d_\Omega$ is the unique graded derivation induced by the shifted cobracket 
\[\Sigma^{-1}\bf{L} \xrightarrow{\Sigma^{-1}\delta} \Sigma^{-1} \Lambda^2 \bf{L} \cong \Sigma \rm{Com}(2) \otimes_{S_2} (\Sigma^{-1} \bf{L})^{\otimes 2}.\]
\end{definition}

We refer to $d_\bf{L}$ as the \emph{internal differential} and $d_\Omega$ as the \emph{cobar differential}. Explicitly, we have for $\ul{x}_0,\ldots,\ul{x}_p \in \Sigma^{-1} \bf{L}$ (where we abbreviate $\ul{x} = s^{-1} x$) that
\[d_\bf{L}(\ul{x}_0 \wedge \ldots \wedge \ul{x}_p) = \sum_{i=0}^p (-1)^{|\ul{x}_0|+\cdots+|\ul{x}_{i-1}|} \ul{x}_0 \wedge \cdots \wedge \ul{d(x_i)} \wedge \cdots \wedge \ul{x}_p,\]
\[d_\Omega(\ul{x}_0 \wedge \ldots \wedge \ul{x}_p) = \sum_{i=0}^p (-1)^{|\ul{x}_0|+\cdots+|\ul{x}_{i-1}|} \ul{x}_0 \wedge \cdots \wedge \ul{\delta(x_i)} \wedge \cdots \wedge \ul{x}_p.\]

\begin{remark}Again it is occasionally more convenient to use a variant that is an augmented unital dg-commutative algebra by using $\rm{Com}^{u}$ instead of $\rm{Com}^\rm{nu}$ and with augmentation induced by the one of $\rm{Com}^\rm{u}$.
\end{remark}

Koszul duality once more takes a concrete form, in terms of the unit and counit of the classical adjunctions between the bar and cobar constructions \cite[Theorems 11.3.6, 11.3.7]{LodayVallette}: given a connected dg-commutative algebra $\bf{A}$ and connected dg-Lie coalgebra $\bf{L}$ there are quasi-isomorphisms
\[\epsilon \colon \Omega^\rm{coLie} B^\rm{Com} \bf{A} \overset{\simeq}\lra \bf{A} \qquad \text{and} \qquad \eta \colon \bf{L} \overset{\simeq}\lra B^\rm{Com} \Omega^\rm{coLie} \bf{L},\]
and we will use for an explicit formula for $\eta$ from \cite[Section 2.4]{Souderes}. These have inverses on the level of chain complexes, \emph{not} compatible with (co)algebra structures, given by the inclusion of or projection onto certain terms.

\subsubsection{Barr's splitting} One particularly useful feature of these concrete models is that they yield a splitting result, originally due to Barr \cite{Barr1968HarrisonHH}. For an $E_\infty^\rm{nu}$-algebra $\bf{A}$, we let 
\begin{align*}p_1 \colon  \fgt_{E^\rm{nu}_1}(\bf{A}) &\lra \cot_{E^\rm{nu}_1}(\bf{A}) \\
p_\infty \colon \fgt_{E^\rm{nu}_\infty}(\bf{A}) &\lra \cot_{E^\rm{nu}_\infty}(\bf{A}) \\
\pi_\bf{A} \colon \cot_{E^\rm{nu}_1}(\bf{A}) &\lra \cot_{E^\rm{nu}_\infty}(\bf{A})\end{align*}
denote the canonical maps, induced by the maps of operads $E^\rm{nu}_0 \to E^\rm{nu}_1$, $E^\rm{nu}_0 \to E^\rm{nu}_\infty$, and $E^\rm{nu}_1 \to E^\rm{nu}_\infty$. We will once more assume that $\scr{C}$ is obtained from the model category $\rm{C} = \Fun(\rm{A},\rm{Ch}_\bb{Q})$ by inverting weak equivalences.

\begin{proposition}[Barr] \label{prop:barr-splitting}
    Let $\bf{A} \in \Alg_{E^\rm{nu}_\infty}(\scr{C})_{\geq 1}$, then the map $\pi_\bf{A} \colon \cot_{E^\rm{nu}_1}(\bf{A}) \to \cot_{E^\rm{nu}_\infty}(\bf{A})$ admits a natural splitting $s_\bf{A} \colon \cot_{E^\rm{nu}_\infty}(\bf{A}) \to \rm{cot}_{E^\rm{nu}_1}(\bf{A})$ fitting in a commutative diagram
    \[\begin{tikzcd} \fgt_{E_\infty}(\bf{A}) \dar[swap]{p_\infty} \arrow{rd}{p_1} & \\
    \cot_{E^\rm{nu}_\infty}(\bf{A}) \rar[swap]{s_\bf{A}} & \cot_{E^\rm{nu}_1}(\bf{A}).\end{tikzcd}\]
\end{proposition}

\begin{proof}By the rectification results of \cref{sec:rect-dg}, we can assume that $\bf{A} \in \Alg_\rm{Com}(\scr{C})$ is a graded dg-commutative algebra in $\rm{C}$ whose underlying object is connected. By the above, $\cot_{E^\rm{nu}_1}(\bf{A})$ and $\cot_{E^\rm{nu}_\infty}(\bf{A})$ may be computed as
    \begin{align*}B^\rm{As}(\bf{A}) &\coloneq \Sigma^{-1} (\rm{coAs}^\rm{nu} \circ \Sigma, d_\bf{A}+ d_B) \\ 
    B^\rm{Com}(\bf{A}) &\coloneq \Sigma^{-1} (\rm{coLie} \circ \Sigma, d_\bf{A}+ d'_B).\end{align*}
Since the map $\pi_\bf{A}$ is induced by the map of operads $\rm{Lie} \to \rm{Ass}$, it is given in this model by the map $B^\rm{As}(\bf{A}) \to B^\rm{Com}(\bf{A})$ induces by the Koszul dual map $\rm{coAs}^\rm{nu} \to \rm{coLie}$. Similarly, we have $\fgt_{\rm{Com}^\rm{nu}}\bf{A} \simeq \Sigma^{-1}(E^\rm{nu}_0 \circ \Sigma (\bf{A}), d_\bf{A}+0)$ so that the canonical map $p_1$ is induced by the canonical map $DE^\rm{nu}_0 \to \rm{coAs}^\rm{nu}$ and similarly $p_\infty$ is induced by $(DE^\rm{nu}_0)^\vee \to \rm{coLie}$. Thus, it suffices to produce a natural splitting $j_\bf{A} \colon B^\rm{Com}(\bf{A}) \to B^\rm{As}(\bf{A})$ under $\fgt_{E^\rm{nu}_\infty}(\bf{A})$. The existence of such a splitting follows from the more general result \cite[Theorem 1.1]{Barr1968HarrisonHH} by taking coefficients to be $\bb{Q}$ viewed as a trivial $\bf{A},\bf{A}$-bimodule.
\end{proof}

\bibliographystyle{amsalpha}
\bibliography{./refs}

\newcommand{\etalchar}[1]{$^{#1}$}
\providecommand{\bysame}{\leavevmode\hbox to3em{\hrulefill}\thinspace}
\providecommand{\MR}{\relax\ifhmode\unskip\space\fi MR }
\providecommand{\MRhref}[2]{%
  \href{http://www.ams.org/mathscinet-getitem?mr=#1}{#2}
}
\providecommand{\href}[2]{#2}
\begin{thebibliography}{GKRW25b}

\bibitem[AC03]{AubryChataur}
M.~Aubry and D.~Chataur, \emph{Cooperads and coalgebras as closed model
  categories}, J. Pure Appl. Algebra \textbf{180} (2003), no.~1-2, 1--23.
  \MR{1966520}

\bibitem[AKKN18]{AKKN}
A.~Alekseev, N.~Kawazumi, Y.~Kuno, and F.~Naef, \emph{The {G}oldman-{T}uraev
  {L}ie bialgebra in genus zero and the {K}ashiwara-{V}ergne problem}, Adv.
  Math. \textbf{326} (2018), 1--53. \MR{3758425}

\bibitem[AM10]{AguiarMahajan}
M.~Aguiar and S.~Mahajan, \emph{Monoidal functors, species and {H}opf
  algebras}, CRM Monograph Series, vol.~29, American Mathematical Society,
  Providence, RI, 2010, With forewords by Kenneth Brown and Stephen Chase and
  Andr\'{e} Joyal. \MR{2724388}

\bibitem[AMP24]{AMP}
A.~Ash, J.~Miller, and P.~Patzt, \emph{Hopf algebras, {S}teinberg modules, and
  the unstable cohomology of ${SL}_n(\mathbb {Z})$ and ${GL}_n(\mathbb {Z})$},
  2024, arXiv:2404.13776.

\bibitem[Arn69]{Arnold}
V.~I. Arnold, \emph{The cohomology ring of the group of dyed braids}, Mat.
  Zametki \textbf{5} (1969), 227--231. \MR{242196}

\bibitem[AT09]{AlekseevTorossianNote}
A.~Alekseev and C.~Torossian, \emph{Flat connections and trivalent graphs},
  \url{https://web.archive.org/web/20211028202354/https://webusers.imj-prg.fr/~charles.torossian/publication/trivalent.pdf},
  2009.

\bibitem[AT12]{AlekseevTorossian}
\bysame, \emph{The {K}ashiwara-{V}ergne conjecture and {D}rinfeld's
  associators}, Ann. of Math. (2) \textbf{175} (2012), no.~2, 415--463.
  \MR{2877064}

\bibitem[Bar68]{Barr1968HarrisonHH}
M.~Barr, \emph{Harrison homology, hochschild homology and triples}, Journal of
  Algebra \textbf{8} (1968), 314--323.

\bibitem[BB24]{BlansBlom}
M.~Blans and T.~Blom, \emph{On the chain rule in {G}oodwillie calculus}, 2024,
  arXiv:2410.20504.

\bibitem[BBK]{BlansBlomKupers}
M.~Blans, T.~Blom, and A.~Kupers, \emph{Naturality of $\infty$-categorical
  bar-cobar duality}, in preparation.

\bibitem[BCGP24]{BCGP}
F.~Brown, M.~Chan, S.~Galatius, and S.~Payne, \emph{Hopf algebras in the
  cohomology of $\mathcal{A}_g$, $\mathrm{GL}_n(\mathbb{Z})$, and
  $\mathrm{SL}_n(\mathbb{Z})$}, 2024, arXiv:2405.11528.

\bibitem[BCN25]{BCN}
D.~L.~B. Brantner, R.~Campos, and J.~Nuiten, \emph{P{D} operads and explicit
  partition {L}ie algebras}, Mem. Amer. Math. Soc. \textbf{315} (2025),
  no.~1597, v+125. \MR{5003477}

\bibitem[BD94]{BD94}
A.~Beilinson and P.~Deligne, \emph{Interpr\'{e}tation motivique de la
  conjecture de {Z}agier reliant polylogarithmes et r\'{e}gulateurs}, Motives
  ({S}eattle, {WA}, 1991), Proc. Sympos. Pure Math., vol.~55, Amer. Math. Soc.,
  Providence, RI, 1994, pp.~97--121. \MR{1265552}

\bibitem[BDE{\etalchar{+}}04]{BDEPS}
N.~Berry, A.~Dubickas, N.~D. Elkies, B.~Poonen, and C.~Smyth, \emph{The
  conjugate dimension of algebraic numbers}, Q. J. Math. \textbf{55} (2004),
  no.~3, 237--252. \MR{2082091}

\bibitem[Ber14]{BerglundKoszul}
A.~Berglund, \emph{Koszul spaces}, Trans. Amer. Math. Soc. \textbf{366} (2014),
  no.~9, 4551--4569. \MR{3217692}

\bibitem[Bez94]{Bezrukavnikov}
R.~Bezrukavnikov, \emph{Koszul {DG}-algebras arising from configuration
  spaces}, Geom. Funct. Anal. \textbf{4} (1994), no.~2, 119--135. \MR{1262702}

\bibitem[BGS20]{BGSII}
C.~Barwick, S.~Glasman, and J.~Shah, \emph{Spectral {M}ackey functors and
  equivariant algebraic {$K$}-theory, {II}}, Tunis. J. Math. \textbf{2} (2020),
  no.~1, 97--146. \MR{3933393}

\bibitem[BHM93]{BHM}
M.~B\"{o}kstedt, W.~C. Hsiang, and I.~Madsen, \emph{The cyclotomic trace and
  algebraic {$K$}-theory of spaces}, Invent. Math. \textbf{111} (1993), no.~3,
  465--539. \MR{1202133}

\bibitem[BK94]{BlochKriz}
S.~Bloch and I.~K\v{r}\'{\i}\v{z}, \emph{Mixed {T}ate motives}, Ann. of Math.
  (2) \textbf{140} (1994), no.~3, 557--605. \MR{1307897}

\bibitem[Blo00]{Bloch}
S.~J. Bloch, \emph{Higher regulators, algebraic {$K$}-theory, and zeta
  functions of elliptic curves}, CRM Monograph Series, vol.~11, American
  Mathematical Society, Providence, RI, 2000. \MR{1760901}

\bibitem[BM12]{BataninMarkl}
M.~Batanin and M.~Markl, \emph{Centers and homotopy centers in enriched
  monoidal categories}, Adv. Math. \textbf{230} (2012), no.~4-6, 1811--1858.
  \MR{2927355}

\bibitem[BMS24]{BenMosheSchlank}
S.~Ben-Moshe and T.~M. Schlank, \emph{Higher semiadditive algebraic {K}-theory
  and redshift}, Compos. Math. \textbf{160} (2024), no.~2, 237--287.
  \MR{4679205}

\bibitem[Boa99]{BoardmanSS}
J.~M. Boardman, \emph{Conditionally convergent spectral sequences}, Homotopy
  invariant algebraic structures ({B}altimore, {MD}, 1998), Contemp. Math.,
  vol. 239, Amer. Math. Soc., Providence, RI, 1999, pp.~49--84. \MR{1718076}

\bibitem[Bol24]{Bol24}
V.~Bolbachan, \emph{On the {G}oncharov's conjecture in degree $m{-}1$ and
  weight $m$}, 2024, arXiv:2404.06271.

\bibitem[Bor74]{BorelStable}
A.~Borel, \emph{Stable real cohomology of arithmetic groups}, Ann. Sci.
  \'{E}cole Norm. Sup. (4) \textbf{7} (1974), 235--272 (1975). \MR{387496}

\bibitem[BPW24]{BPW}
B.~Br\"{u}ck, K.~I. Piterman, and V.~Welker, \emph{The common basis complex and
  the partial decomposition poset}, Int. Math. Res. Not. IMRN (2024), no.~18,
  12746--12760. \MR{4798648}

\bibitem[Bro94]{Brown}
K.~S. Brown, \emph{Cohomology of groups}, Graduate Texts in Mathematics,
  vol.~87, Springer-Verlag, New York, 1994, Corrected reprint of the 1982
  original. \MR{1324339}

\bibitem[Bro12]{Bro12}
F.~Brown, \emph{Mixed {T}ate motives over {$\Bbb Z$}}, Ann. of Math. (2)
  \textbf{175} (2012), no.~2, 949--976. \MR{2993755}

\bibitem[BVGS90]{BGSV90}
A.~Beilinson, A.~Varchenko, A.~Goncharov, and V.~Shekhtman, \emph{Projective
  geometry and {$K$}-theory}, Algebra i Analiz \textbf{2} (1990), no.~3,
  78--130. \MR{1073210}

\bibitem[Car07]{Cartier}
P.~Cartier, \emph{A primer of {H}opf algebras}, Frontiers in number theory,
  physics, and geometry. {II}, Springer, Berlin, 2007, pp.~537--615.
  \MR{2290769}

\bibitem[Cat93]{Cathelineau}
J.-L. Cathelineau, \emph{Homologie du groupe lin\'{e}aire et polylogarithmes
  (d'apr\`es {A}. {B}. {G}oncharov et d'autres)}, Ast\'{e}risque (1993),
  no.~216, Exp. No. 772, 5, 311--341, S\'{e}minaire Bourbaki, Vol. 1992/93.
  \MR{1246402}

\bibitem[CMRR24]{CMRR}
S.~Charlton, A.~Matveiakin, D.~Radchenko, and D.~Rudenko, \emph{The {H}opf
  algebra of formal multiple polylogarithms}, 2024.

\bibitem[CMRR26]{CMRR24}
\bysame, \emph{The {H}opf algebra of formal multiple polylogarithms}, Int.
  Math. Res. Not. IMRN (2026), no.~2, Paper No. rnaf361, 29. \MR{5013247}

\bibitem[Coh95]{Cohen}
F.~R. Cohen, \emph{On configuration spaces, their homology, and {L}ie
  algebras}, J. Pure Appl. Algebra \textbf{100} (1995), no.~1-3, 19--42.
  \MR{1344842}

\bibitem[CRR25]{CharltonRadchenkoRudenko}
S.~Charlton, D.~Radchenko, and D.~Rudenko, \emph{Multiple polylogarithms and
  the {S}teinberg module}, 2025.

\bibitem[CS22]{ChingSalvatore}
M.~Ching and P.~Salvatore, \emph{Koszul duality for topological
  {$E_n$}-operads}, Proc. Lond. Math. Soc. (3) \textbf{125} (2022), no.~1,
  1--60. \MR{4456966}

\bibitem[CV03]{ConantVogtmann}
J.~Conant and K.~Vogtmann, \emph{On a theorem of {K}ontsevich}, Algebr. Geom.
  Topol. \textbf{3} (2003), 1167--1224. \MR{2026331}

\bibitem[CZ24]{CampbellZakharevich}
J.~A. Campbell and I.~Zakharevich, \emph{Hilbert's third problem and a
  conjecture of {G}oncharov}, Adv. Math. \textbf{451} (2024), Paper No. 109757,
  57. \MR{4759410}

\bibitem[Day70]{Day}
B.~Day, \emph{On closed categories of functors}, Reports of the {M}idwest
  {C}ategory {S}eminar, {IV}, Lecture Notes in Math., Vol. 137, Springer,
  Berlin-New York, 1970, pp.~1--38. \MR{272852}

\bibitem[Del71]{Del71b}
P.~Deligne, \emph{Th\'{e}orie de {H}odge. {II}}, Inst. Hautes \'{E}tudes Sci.
  Publ. Math. (1971), no.~40, 5--57. \MR{498551}

\bibitem[DG05]{DG05}
P.~Deligne and A.~B. Goncharov, \emph{Groupes fondamentaux motiviques de {T}ate
  mixte}, Ann. Sci. \'{E}cole Norm. Sup. (4) \textbf{38} (2005), no.~1, 1--56.
  \MR{2136480}

\bibitem[DJ02]{deJeu}
R.~De~Jeu, \emph{A remark on the rank conjecture}, $K$-Theory \textbf{25}
  (2002), no.~3, 215--231. \MR{1909867}

\bibitem[{Dri}90]{Drinfeld}
V.~G. {Drinfeld}, \emph{On quasitriangular quasi-{H}opf algebras and on a group
  that is closely connected with {${\rm Gal}(\overline{\bf Q}/{\bf Q})$}},
  Algebra i Analiz \textbf{2} (1990), no.~4, 149--181, translated in Leningrad
  Math. J. 2 (1991), no. 4, 829--860. \MR{1080203}

\bibitem[Dup21]{Dup20}
C.~Dupont, \emph{Progr{\`e}s rec{\'e}nts sur la conjecture de {Z}agier et le
  programme de {G}oncharov [d'apr{\`e}s {G}oncharov, {R}udenko, {G}angl, ...]},
  S{\'e}minaire Bourbaki \textbf{73{\`e}me ann{\'e}e} (2021), no.~1176.

\bibitem[Fel18]{Felder}
M.~Felder, \emph{Internally connected graphs and the {K}ashiwara-{V}ergne {L}ie
  algebra}, Lett. Math. Phys. \textbf{108} (2018), no.~6, 1407--1441.
  \MR{3797752}

\bibitem[FGV22]{FengGalatiusVenkatesh}
T.~Feng, S.~Galatius, and A.~Venkatesh, \emph{The {G}alois action on symplectic
  {K}-theory}, Invent. Math. \textbf{230} (2022), no.~1, 225--319. \MR{4480148}

\bibitem[GHN17]{GepnerHaugsengNikolaus}
D.~Gepner, R.~Haugseng, and T.~Nikolaus, \emph{Lax colimits and free fibrations
  in {$\infty$}-categories}, Doc. Math. \textbf{22} (2017), 1225--1266.
  \MR{3690268}

\bibitem[GJ94]{GetzlerJones}
E.~Getzler and J.~D.~S. Jones, \emph{Operads, homotopy algebra and iterated
  integrals for double loop spaces}, 1994, arXiv:hep-th/9403055.

\bibitem[GK95]{GetzlerKapranov}
E.~Getzler and M.~M. Kapranov, \emph{Cyclic operads and cyclic homology},
  Geometry, topology, \& physics, Conf. Proc. Lecture Notes Geom. Topology, IV,
  Int. Press, Cambridge, MA, 1995, pp.~167--201. \MR{1358617}

\bibitem[GKRW24]{GKRW19}
S.~Galatius, A.~Kupers, and O.~Randal-Williams, \emph{{$E_\infty$}-cells and
  general linear groups of finite fields}, Ann. Sci. \'{E}c. Norm. Sup\'{e}r.
  (4) \textbf{57} (2024), no.~6, 1845--1882. \MR{4862505}

\bibitem[GKRW25a]{GKRW18}
\bysame, \emph{Cellular {$E_k$}-algebras}, Ast\'{e}risque (2025), no.~460,
  x+299. \MR{4987221}

\bibitem[GKRW25b]{GKRW20}
\bysame, \emph{{$E_\infty$}-cells and general linear groups of infinite
  fields}, Duke Math. J. \textbf{174} (2025), no.~14, 2927--3046. \MR{4974473}

\bibitem[GLF16]{GarnerLopezFranco}
R.~Garner and I.~L\'{o}pez~Franco, \emph{Commutativity}, J. Pure Appl. Algebra
  \textbf{220} (2016), no.~5, 1707--1751. \MR{3437265}

\bibitem[Gon95a]{Gon95b}
A.~Goncharov, \emph{Geometry of configurations, polylogarithms, and motivic
  cohomology}, Adv. Math. \textbf{114} (1995), no.~2, 197--318. \MR{1348706}

\bibitem[Gon95b]{Gon95}
\bysame, \emph{Polylogarithms in arithmetic and geometry}, Proceedings of the
  {I}nternational {C}ongress of {M}athematicians, {V}ol. 1, 2 ({Z}\"{u}rich,
  1994), Birkh\"{a}user, Basel, 1995, pp.~374--387. \MR{1403938}

\bibitem[Gon99]{Gon99}
Alexander Goncharov, \emph{Volumes of hyperbolic manifolds and mixed {T}ate
  motives}, J. Amer. Math. Soc. \textbf{12} (1999), no.~2, 569--618.
  \MR{1649192}

\bibitem[Gon01a]{Gon01}
A.~Goncharov, \emph{Multiple polylogarithms and mixed {T}ate motives},
  arXiv:math/0103059 [math.AG], 2001.

\bibitem[Gon01b]{Goncharov01}
A.~B. Goncharov, \emph{The dihedral {L}ie algebras and {G}alois symmetries of
  {$\pi_1^{(l)}(\Bbb P^1-(\{0,\infty\}\cup\mu_N))$}}, Duke Math. J.
  \textbf{110} (2001), no.~3, 397--487. \MR{1869113}

\bibitem[Gon08]{GoncharovEuler}
\bysame, \emph{Euler complexes and geometry of modular varieties}, Geom. Funct.
  Anal. \textbf{17} (2008), no.~6, 1872--1914. \MR{2399086}

\bibitem[Gon19]{Gon19}
A.~Goncharov, \emph{Hodge correlators}, J. Reine Angew. Math. \textbf{748}
  (2019), 1--138. \MR{3918430}

\bibitem[GR14]{GrinbergReiner}
D.~Grinberg and V.~Reiner, \emph{Hopf algebras in combinatorics}, 2014,
  arXiv:1409.8356.

\bibitem[Hai86]{HainIndec}
R.~M. Hain, \emph{On the indecomposable elements of the bar construction},
  Proc. Amer. Math. Soc. \textbf{98} (1986), no.~2, 312--316. \MR{854039}

\bibitem[Hai94]{HainClassical}
\bysame, \emph{Classical polylogarithms}, Motives ({S}eattle, {WA}, 1991),
  Proc. Sympos. Pure Math., vol.~55, Amer. Math. Soc., Providence, RI, 1994,
  pp.~3--42. \MR{1265550}

\bibitem[Har77]{Harder}
G.~Harder, \emph{Die {K}ohomologie {$S$}-arithmetischer {G}ruppen \"{u}ber
  {F}unktionenk\"{o}rpern}, Invent. Math. \textbf{42} (1977), 135--175.
  \MR{473102}

\bibitem[Hau21]{HaugsengMonads}
R.~Haugseng, \emph{On lax transformations, adjunctions, and monads in
  {$(\infty,2)$}-categories}, High. Struct. \textbf{5} (2021), no.~1, 244--281.
  \MR{4367222}

\bibitem[Hau22]{HaugsengSS}
\bysame, \emph{{$\infty$}-operads via symmetric sequences}, Math. Z.
  \textbf{301} (2022), no.~1, 115--171. \MR{4405646}

\bibitem[Hei17]{Heine}
H.~Heine, \emph{A monadicity theorem for higher algebraic structures}, 2017,
  arXiv:1712.00555.

\bibitem[Hei26]{HeineMM}
\bysame, \emph{A derived {M}ilnor-{M}oore theorem}, Selecta Math. (N.S.)
  \textbf{32} (2026), no.~2, Paper No. 39. \MR{5060174}

\bibitem[Heu24]{HeutsKoszul}
G.~Heuts, \emph{Koszul duality and a conjecture of {F}rancis--{G}aitsgory},
  2024, arXiv:2408.06173.

\bibitem[HHLN23]{HHLN}
R.~Haugseng, F.~Hebestreit, S.~Linskens, and J.~Nuiten, \emph{Lax monoidal
  adjunctions, two-variable fibrations and the calculus of mates}, Proc. Lond.
  Math. Soc. (3) \textbf{127} (2023), no.~4, 889--957. \MR{4655344}

\bibitem[HHS25]{HHS}
P.~Hanlon, P.~Hersh, and J.~Shareshian, \emph{A ${G}l_n(q)$ analogue of the
  partition lattice}, arXiv:2505.02202, 2025.

\bibitem[Hin15]{HinichRectification}
V.~Hinich, \emph{Rectification of algebras and modules}, Doc. Math. \textbf{20}
  (2015), 879--926. \MR{3404213}

\bibitem[Hin20]{Hinich}
\bysame, \emph{Yoneda lemma for enriched {$\infty$}-categories}, Adv. Math.
  \textbf{367} (2020), 107129, 119. \MR{4080581}

\bibitem[HL]{HeutsLand}
G.~S. K.~S. Heuts and M.~Land, \emph{Koszul duality of ${E}_n$-algebras and
  ${E}_n$-operads}, in preparation.

\bibitem[HL24]{HeutsLandFormality}
\bysame, \emph{Formality of $\mathbb{E}_n$-algebras and cochains on spheres},
  2024, arXiv:2407.00790.

\bibitem[KKM{\etalchar{+}}25]{KKMMW}
I.~Klang, J.~Kuijper, C.~Malkiewich, D.~Mehrle, and T.~Wittich, \emph{Higher
  spherical scissors congruence {I}: {H}opf algebra}, arxiv:2509.18009.

\bibitem[Kon93]{KontsevichFormal}
M.~Kontsevich, \emph{Formal (non)commutative symplectic geometry}, The
  {G}elfand {M}athematical {S}eminars, 1990--1992, Birkh\"{a}user Boston,
  Boston, MA, 1993, pp.~173--187. \MR{1247289}

\bibitem[KRS26]{KRS2}
A.~Kupers, D.~Rudenko, and I.~Sierra, \emph{Mixed tate motives over number
  fields}, 2026, in preparation.

\bibitem[KS14]{KahnSun}
B.~Kahn and F.~Sun, \emph{On universal modular symbols}, 2014, arXiv:1407.0475.

\bibitem[Lev93]{Levine}
M.~Levine, \emph{Tate motives and the vanishing conjectures for algebraic
  {$K$}-theory}, Algebraic {$K$}-theory and algebraic topology ({L}ake
  {L}ouise, {AB}, 1991), NATO Adv. Sci. Inst. Ser. C: Math. Phys. Sci., vol.
  407, Kluwer Acad. Publ., Dordrecht, 1993, pp.~167--188. \MR{1367296}

\bibitem[Lev94]{LevineBloch}
\bysame, \emph{Bloch's higher {C}how groups revisited}, no. 226, 1994,
  $K$-theory (Strasbourg, 1992), pp.~10, 235--320. \MR{1317122}

\bibitem[LNP25]{LinskensNardinPol}
S.~Linskens, D.~Nardin, and L.~Pol, \emph{Global homotopy theory via partially
  lax limits}, Geom. Topol. \textbf{29} (2025), no.~3, 1345--1440. \MR{4918109}

\bibitem[Lod98]{LodayCyclic}
J.-L. Loday, \emph{Cyclic homology}, second ed., Grundlehren der mathematischen
  Wissenschaften [Fundamental Principles of Mathematical Sciences], vol. 301,
  Springer-Verlag, Berlin, 1998, Appendix E by Mar\'{\i}a O. Ronco, Chapter 13
  by the author in collaboration with Teimuraz Pirashvili. \MR{1600246}

\bibitem[Lur09]{LurieHTT}
J.~Lurie, \emph{Higher topos theory}, Annals of Mathematics Studies, vol. 170,
  Princeton University Press, 2009. \MR{2522659}

\bibitem[Lur17]{LurieHA}
\bysame, \emph{{Higher algebra}}, September 2017 version (2017).

\bibitem[LV12]{LodayVallette}
J.-L. Loday and B.~Vallette, \emph{Algebraic operads}, Grundlehren der
  mathematischen Wissenschaften [Fundamental Principles of Mathematical
  Sciences], vol. 346, Springer, Heidelberg, 2012. \MR{2954392}

\bibitem[Mal20]{Mal20}
N.~Malkin, \emph{Shuffle relations for {H}odge and motivic correlators},
  arXiv:2003.06521 [math.AG], 2020.

\bibitem[Mic80]{Michaelis}
W.~Michaelis, \emph{Lie coalgebras}, Adv. in Math. \textbf{38} (1980), no.~1,
  1--54. \MR{594993}

\bibitem[Mir11]{Mirzaii}
B.~Mirzaii, \emph{Bloch-{W}igner theorem over rings with many units}, Math. Z.
  \textbf{268} (2011), no.~1-2, 329--346. \MR{2805438}

\bibitem[Mir13]{MirzaiiErratum}
\bysame, \emph{Erratum to: {B}loch-{W}igner theorem over rings with many units
  [mr2805438]}, Math. Z. \textbf{275} (2013), no.~1-2, 653--655. \MR{3101826}

\bibitem[MNP20]{MNP}
J.~Miller, R.~Nagpal, and P.~Patzt, \emph{Stability in the high-dimensional
  cohomology of congruence subgroups}, Compos. Math. \textbf{156} (2020),
  no.~4, 822--861. \MR{4079629}

\bibitem[MPP25]{MPPII}
J.~Miller, P.~Patzt, and A.~Putman, \emph{Homological vanishing for the
  {S}teinberg representation {II}: reductive groups and integral conjectures},
  2025, arXiv:2509.01559.

\bibitem[MPW23]{MPW23}
J.~Miller, P.~Patzt, and J.~C.~H. Wilson, \emph{On rank filtrations of
  algebraic {K}-theory and {S}teinberg modules}, 2023, arXiv:2303.00245.

\bibitem[MVW06]{MVW}
C.~Mazza, V.~Voevodsky, and C.~Weibel, \emph{Lecture notes on motivic
  cohomology}, Clay Mathematics Monographs, vol.~2, American Mathematical
  Society, Providence, RI; Clay Mathematics Institute, Cambridge, MA, 2006.
  \MR{2242284}

\bibitem[NS89]{NesterenkoSuslin}
Yu.~P. Nesterenko and A.~A. Suslin, \emph{Homology of the general linear group
  over a local ring, and {M}ilnor's {$K$}-theory}, Izv. Akad. Nauk SSSR Ser.
  Mat. \textbf{53} (1989), no.~1, 121--146. \MR{992981}

\bibitem[NS22]{NardinShah}
D.~Nardin and J.~Shah, \emph{Parametrized and equivariant higher algebra},
  2022, arXiv:2203.00072.

\bibitem[Pri70]{Priddy}
S.~B. Priddy, \emph{Koszul resolutions}, Trans. Amer. Math. Soc. \textbf{152}
  (1970), 39--60. \MR{265437}

\bibitem[PRiLY26]{PRY}
D.~Petersen, V.~Roca~i Lucio, and S.~Yalin, \emph{Unifying {K}oszul dualities
  via point-set models}, 2026, arXiv:2603.29910.

\bibitem[PT]{PortaTeyssier}
M.~Porta and J.-B. Teyssier, \emph{Day's convolution for
  pro-$\infty$-categories},
  \url{http://jbteyssier.com/papers/jbteyssier_day.pdf}.

\bibitem[Qui73]{QuillenFiniteGeneration}
D.~Quillen, \emph{Finite generation of the groups {$K\sb{i}$} of rings of
  algebraic integers}, Algebraic {$K$}-theory, {I}: {H}igher {$K$}-theories
  ({P}roc. {C}onf., {B}attelle {M}emorial {I}nst., {S}eattle, {W}ash., 1972),
  Lecture Notes in Math., Vol. 341, Springer, Berlin-New York, 1973,
  pp.~179--198. \MR{349812}

\bibitem[Ram82]{Ramakrishnan}
D.~Ramakrishnan, \emph{On the monodromy of higher logarithms}, Proc. Amer.
  Math. Soc. \textbf{85} (1982), no.~4, 596--599. \MR{660611}

\bibitem[Rob04]{Robinson}
A.~Robinson, \emph{Partition complexes, duality and integral tree
  representations}, Algebr. Geom. Topol. \textbf{4} (2004), 943--960.
  \MR{2100687}

\bibitem[Rog92]{Rognes}
J.~Rognes, \emph{A spectrum level rank filtration in algebraic {$K$}-theory},
  Topology \textbf{31} (1992), no.~4, 813--845. \MR{1191383}

\bibitem[Rog10]{RognesMotivic}
\bysame, \emph{Motivic complexes from the stable rank filtration}, 2010,
  \url{https://www.mn.uio.no/math/personer/vit/rognes/papers/bergen10.pdf}.

\bibitem[Rog21]{RognesWeight}
\bysame, \emph{The weight and rank filtrations}, 2021, arXiv:2110.12264.

\bibitem[RW25]{RWchromatic}
O.~Randal-Williams, \emph{A chromatic approach to homological stability}, 2025,
  arXiv:2508.20629.

\bibitem[Seg73]{SegalConfiguration}
G.~Segal, \emph{Configuration-spaces and iterated loop-spaces}, Invent. Math.
  \textbf{21} (1973), 213--221. \MR{331377}

\bibitem[Sin13]{Sinha}
D.~P. Sinha, \emph{The (non-equivariant) homology of the little disks operad},
  O{PERADS} 2009, S\'{e}min. Congr., vol.~26, Soc. Math. France, Paris, 2013,
  pp.~253--279. \MR{3203375}

\bibitem[Sou85]{Soule}
Christophe Soul\'e, \emph{Op\'erations en {$K$}-th\'eorie alg\'ebrique}, Canad.
  J. Math. \textbf{37} (1985), no.~3, 488--550. \MR{787114}

\bibitem[Sou16]{Souderes}
I.~Soud\`eres, \emph{A relative basis for mixed {T}ate motives over the
  projective line minus three points}, Commun. Number Theory Phys. \textbf{10}
  (2016), no.~1, 87--131. \MR{3521910}

\bibitem[Sun16]{Sun}
F.~Sun, \emph{Algebraic {K}-theory and modular symbols}, 2016,
  arXiv:1604.04700.

\bibitem[Sus84]{Sus84}
A.~Suslin, \emph{Homology of {${\rm GL}_{n}$}, characteristic classes and
  {M}ilnor {$K$}-theory}, vol. 165, 1984, Algebraic geometry and its
  applications, pp.~188--204. \MR{752941}

\bibitem[Sus90]{Sus90}
\bysame, \emph{{$K_3$} of a field, and the {B}loch group}, vol. 183, 1990,
  Translated in Proc. Steklov Inst. Math. {1991}, no. 4, 217--239, Galois
  theory, rings, algebraic groups and their applications (Russian),
  pp.~180--199, 229. \MR{1092031}

\bibitem[{\v{S}}W11]{SeveraWillwacher}
P.~{\v{S}}evera and T.~Willwacher, \emph{Equivalence of formalities of the
  little discs operad}, Duke Math. J. \textbf{160} (2011), no.~1, 175--206.
  \MR{2838354}

\bibitem[Tor21]{ToriiHigher}
T.~Torii, \emph{On higher monoidal $\infty$-categories}, 2021,
  arXiv:2111.00158.

\bibitem[Tor25a]{ToriiMult}
\bysame, \emph{Multiplicative structures on comodules in higher categories},
  2025.

\bibitem[Tor25b]{ToriiDuoidal}
\bysame, \emph{On duoidal {$\infty$}-categories}, J. Homotopy Relat. Struct.
  \textbf{20} (2025), no.~1, 125--162. \MR{4868047}

\bibitem[Wei13]{Weibel}
C.A. Weibel, \emph{The {$K$}-book}, Graduate Studies in Mathematics, vol. 145,
  American Mathematical Society, Providence, RI, 2013, An introduction to
  algebraic $K$-theory. \MR{3076731}

\bibitem[Whi01]{Whitehouse}
S.~Whitehouse, \emph{The integral tree representation of the symmetric group},
  J. Algebraic Combin. \textbf{13} (2001), no.~3, 317--326. \MR{1836907}

\bibitem[Woj91]{Wojtkowiak}
Z.~Wojtkowiak, \emph{The basic structure of polylogarithmic functional
  equations}, Structural properties of polylogarithms, Math. Surveys Monogr.,
  vol.~37, Amer. Math. Soc., Providence, RI, 1991, pp.~205--231. \MR{1148381}

\bibitem[Wu23]{Wu}
H.~Wu, \emph{A {H}opf algebra model for {D}wyer's tame spaces}, Ph.D. thesis,
  EFPL, 2023,
  \url{https://infoscience.epfl.ch/entities/publication/89700c3d-c155-4cd1-8ac6-0973ca4ee813}.

\bibitem[Zag07]{Zagier}
D.~Zagier, \emph{The dilogarithm function}, Frontiers in number theory,
  physics, and geometry. {II}, Springer, Berlin, 2007, pp.~3--65. \MR{2290758}

\end{thebibliography}

\bigskip

\end{document}